\documentclass[11pt,reqno]{article}
\usepackage{amsmath, amsthm, amssymb}
\usepackage{enumitem}

\newcommand{\makeheading}[2]%
        {\hspace*{-\marginparsep minus \marginparwidth}%
         \begin{minipage}[t]{\textwidth}%
                {\large \bfseries #1\hfill#2}\\[-0.15\baselineskip]%
                 \rule{\columnwidth}{1pt}%
         \end{minipage}\bigskip
}

\usepackage{fullpage}
\usepackage[many]{tcolorbox}
\usepackage{xcolor}
\usepackage{wasysym}
\usepackage{mathtools}

\usepackage{tikz}
\usetikzlibrary{cd}
\usetikzlibrary{calc}
\usetikzlibrary{decorations,decorations.pathreplacing,decorations.markings,decorations.pathmorphing}
\newcommand{\tikzmath}[2][]
{\vcenter{\hbox{\begin{tikzpicture}[#1]#2\end{tikzpicture}}}
}
\newcommand{\roundNbox}[6]{
	\draw[rounded corners=5pt, very thick, #1] ($#2+(-#3,-#3)+(-#4,0)$) rectangle ($#2+(#3,#3)+(#5,0)$);
	\coordinate (ZZa) at ($#2+(-#4,0)$);
	\coordinate (ZZb) at ($#2+(#5,0)$);
	\node at ($1/2*(ZZa)+1/2*(ZZb)$) {#6};
}
\tikzset{super thick/.style={line width=3pt}}
\tikzstyle{far>}=[decoration={markings, mark=at position 0.75 with {\arrow{>}}}, postaction={decorate}]
\tikzstyle{mid>}=[decoration={markings, mark=at position 0.55 with {\arrow{>}}}, postaction={decorate}]
\tikzstyle{mid<}=[decoration={markings, mark=at position 0.55 with {\arrow{<}}}, postaction={decorate}]
\tikzset{super thick/.style={line width=3pt}}
\tikzstyle{far>}=[decoration={markings, mark=at position 0.75 with {\arrow{>}}}, postaction={decorate}]
\tikzstyle{mid>}=[decoration={markings, mark=at position 0.5 with {\arrow{>}}}, postaction={decorate}]
\tikzstyle{mid<}=[decoration={markings, mark=at position 0.5 with {\arrow{<}}}, postaction={decorate}]
\tikzstyle{knot}=[preaction={super thick, white, draw}]
\tikzstyle{coupon}=[draw, very thick, rectangle, rounded corners=5pt]
\tikzset{Rightarrow/.style={double equal sign distance,>={Implies},->},
triplecd/.style={-,preaction={draw,Rightarrow}},
quadruplecd/.style={preaction={draw,Rightarrow,
shorten >=0pt
},
shorten >=1pt,
-,double,double
distance=0.2pt}}
\tikzset{
    tripleline/.style args={[#1] in [#2] in [#3]}{
        #1,preaction={preaction={draw,#3},draw,#2}
    }
}
\tikzstyle{triple}=[tripleline={[line width=.15mm,black] in
      [line width=.7mm,white] in
      [line width=1mm,black]}] 
\tikzset{
    quadrupleline/.style args={[#1] in [#2] in [#3] in [#4]}{
        #1,preaction={preaction={preaction={draw,#4},draw,#3}, draw,#2}
    }
}
\tikzstyle{quadruple}=[quadrupleline={[line width=.3mm,white] in
      [line width=.6mm,black] in
      [line width=1.2mm,white] in
      [line width=1.5mm,black]}]

\definecolor{violet}{RGB}{148,0,211}
\definecolor{DarkGreen}{RGB}{0,150,0}

\usepackage[pdftex,plainpages=false,hypertexnames=false,pdfpagelabels]{hyperref}
\definecolor{medium-blue}{rgb}{0,0,.8}
\hypersetup{colorlinks, linkcolor={purple}, citecolor={medium-blue}, urlcolor={medium-blue}}
\newcommand{\arxiv}[1]{\href{http://arxiv.org/abs/#1}{\tt arXiv:\nolinkurl{#1}}}
\newcommand{\arXiv}[1]{\href{http://arxiv.org/abs/#1}{\tt arXiv:\nolinkurl{#1}}}

\newcommand{\doi}[1]{\href{http://dx.doi.org/#1}{{\tt DOI:#1}}}

\newcommand{\TI}[1]{\tikzmath{
\draw[mid>] (-.3,-.6) -- (0,0);
\draw[mid>] (.3,-.6) -- (0,0);
\draw[mid<] (0,.6) -- (0,0);
\filldraw (0,0) node[left]{$\scriptstyle #1$} circle (.05cm);
}}

\newcommand{\Ta}[1]{
\tikzmath{
\draw[thick, blue, mid>] (0,0) -- (-.6,.6);
\draw[mid>] (-.6,-.6) -- (0,0);
\draw[mid>] (.6,-.6) -- (0,0);
\draw[mid>] (0,0) -- (.6,.6);
\filldraw (0,0) node[left]{$\scriptstyle #1$} circle (.05cm);
}}
\DeclareMathOperator{\Ad}{Ad}
\DeclareMathOperator{\Aut}{Aut}

\DeclareMathOperator{\End}{End}

\DeclareMathOperator{\FPdim}{FPdim}

\DeclareMathOperator{\Hom}{Hom}
\DeclareMathOperator{\id}{id}

\DeclareMathOperator{\Inv}{Inv}
\DeclareMathOperator{\Irr}{Irr}

\DeclareMathOperator{\Tr}{Tr}
\DeclareMathOperator{\tr}{tr}

\renewcommand{\Im}{\operatorname{Im}}

\newcommand{\Fib}{\mathsf{Fib}}

\newcommand{\Vect}{\mathsf{Vect}}
\newcommand{\Hilb}{\mathsf{Hilb}}
\newcommand{\Rep}{\mathsf{Rep}}

\def\semicolon{;}
\def\applytolist#1{
    \expandafter\def\csname multi#1\endcsname##1{
        \def\multiack{##1}\ifx\multiack\semicolon
            \def\next{\relax}
        \else
            \csname #1\endcsname{##1}
            \def\next{\csname multi#1\endcsname}
        \fi
        \next}
    \csname multi#1\endcsname}

\def\calc#1{\expandafter\def\csname c#1\endcsname{{\mathcal #1}}}
\applytolist{calc}QWERTYUIOPLKJHGFDSAZXCVBNM;
\def\bbc#1{\expandafter\def\csname bb#1\endcsname{{\mathbb #1}}}
\applytolist{bbc}QWERTYUIOPLKJHGFDSAZXCVBNM;
\def\bfc#1{\expandafter\def\csname bf#1\endcsname{{\mathbf #1}}}
\applytolist{bfc}QWERTYUIOPLKJHGFDSAZXCVBNM;
\def\sfc#1{\expandafter\def\csname s#1\endcsname{{\sf #1}}}
\applytolist{sfc}QWERTYUIOPLKJHGFDSAZXCVBNM;
\def\fc#1{\expandafter\def\csname f#1\endcsname{{\mathfrak #1}}}
\applytolist{fc}QWERTYUIOPLKJHGFDSAZXCVBNM;
\def\rmc#1{\expandafter\def\csname rm#1\endcsname{{\mathrm #1}}}
\applytolist{rmc}QWERTYUIOPLKJHGFDSAZXCVBNM;

\numberwithin{equation}{section}

\theoremstyle{plain}
\newtheorem{thm}[equation]{Theorem}
\newtheorem*{thm*}{Theorem}
\newtheorem{cor}[equation]{Corollary}
\newtheorem{lem}[equation]{Lemma}
\newtheorem{prop}[equation]{Proposition}

\newtheorem*{claim*}{Claim}

\newtheorem{thmalpha}{Theorem}

\newtheorem{coralpha}[thmalpha]{Corollary}
\theoremstyle{definition}
\newtheorem{defn}[equation]{Definition}

\newtheorem*{trick*}{Trick}

\newtheorem{rem}[equation]{Remark}

\makeatletter
\def\thmhead@plain#1#2#3{%
  \thmname{#1}\thmnumber{\@ifnotempty{#1}{ }\@upn{#2}}%
  \thmnote{ {\the\thm@notefont#3}}}
\let\thmhead\thmhead@plain
\makeatother

\newcommand{\bigraph}[1]{{\hspace{-3pt}\begin{array}{c}%
  \raisebox{-2.5pt}{\includegraphics[height=6mm]{#1}}% 
\end{array}\hspace{-3pt}}}

\title{Classification of $\mathbb{Z}/2\mathbb{Z}$-quadratic unitary fusion categories}
\date{\today}
\begin{document}
\author{Cain Edie-Michell, Masaki Izumi, and David Penneys
\\
\\
{\tiny Appendix by: 
Ryan Johnson, Siu-Hung Ng, David Penneys, Jolie Roat, Matthew Titsworth, and Henry Tucker}
}
\maketitle
\begin{abstract}
A unitary fusion category is called $\mathbb{Z}/2\mathbb{Z}$-quadratic if it has a $\mathbb{Z}/2\mathbb{Z}$ group of invertible objects and one other orbit of simple objects under the action of this group.
We give a complete classification of $\mathbb{Z}/2\mathbb{Z}$-quadratic unitary fusion categories.
The main tools for this classification are
skein theory,
a generalization of Ostrik's results on formal codegrees to analyze the induction of the group elements to the center, 
and a computation similar to Larson's rank-finiteness bound for $\mathbb{Z}/3\mathbb{Z}$-near group pseudounitary fusion categories.
This last computation is contained in an appendix coauthored with attendees from the 2014 AMS MRC on Mathematics of Quantum Phases of Matter and Quantum Information.
%
%This is the submitted version of \arxiv{???}
\end{abstract}

%%%%%%%%%%%%%%%%%%%%%%%%%%%%%%%%%%%%%%%%%%%%%%%%%%%%%%
%%%%%%%%%%%%%%%%%%%%%%%%%%%%%%%%%%%%%%%%%%%%%%%%%%%%%%
%%%%%%%%%%%%%%%%%%%%%%%%%%%%%%%%%%%%%%%%%%%%%%%%%%%%%%
\section{Introduction}

In the past several decades, unitary fusion categories (UFCs) have seen broad applications to many areas of mathematics and physics, including representation theory, operator algebras, topological quantum field theory (TQFT), theoretical condensed matter, and quantum information.
Given the complete list of 6j-symbols for a UFC, one can build 
unitary TQFTs 
which compute quantum invariants of links and 3-manifolds \cite{MR1357878,MR1642584},
together with
physical lattice models 
which realize these TQFTs 
\cite{PhysRevB.71.045110,2012.14424}.
These computations are increasingly difficult in the presence of \emph{multiplicity}, i.e., where there is a fusion channel with a dimension greater than 1, a.k.a.~a fusion coefficient which is larger than 1.

While many classification techniques work well for multiplicity free fusion categories, more techniques are required to achieve classification in the multiplicity setting.
We note that at the time of writing, multiplicity free UFCs have been completely classified up to rank 6 \cite{2010.10264}, while arbitrary UFCs have only been classified up to rank 3 \cite{MR1981895,1309.4822}.
For rank 4 fusion categories with a dual pair of simple objects, there is a classification of possible fusion rings in the pseudounitary setting \cite{MR3229513}; our Theorem \ref{thm:Main} (and Corollary \ref{cor:FinishLarson}) below completes the classification of rank 4 UFCs with a dual pair of simples.
The case of rank 4 with 4 self-dual objects still seems out of reach at this time.

Surprisingly, all currently known fusion categories fit into four families:
(1) those related to groups,
(2) those related to quantum groups at roots of unity \cite{MR936086,MR1090432,MR1328736},
(3) the Haagerup-Izumi quadratic categories \cite{MR1317352,MR1659954,MR1686551,MR1832764,MR2837122,MR3167494,MR3635673,MR3827808}, and
(4) the Extended Haagerup fusion categories \cite{MR2979509,1810.06076}.
Given a finite group $G$, a \emph{$G$-quadratic fusion category} is a fusion category $\cC$ with a finite group $G$ of simple objects and one other $G$-orbit $\{g\rho\}_{g\in G}$ of simple objects.
(The collection of all $G$-quadratic fusion categories over all finite groups $G$ is exactly the family (3) above.)
The term `quadratic' comes from the existence of a quadratic relation for the self-fusion of an object $\rho$ which generates the other $G$-orbit.
Surprisingly, for a fixed group $G$ beyond the trivial group, rank-finiteness is not known for $G$-quadratic fusion categories (for the trivial group, see \cite{MR1981895}).
The case $G=\bbZ/2\bbZ$ is classified in the pivotal setting in \cite{1309.4822}, and rank-finiteness for $G=\bbZ/3\bbZ$ is achieved in the pseudounitary setting in \cite{MR3229513}.

In this article, we give a complete classification of $\bbZ/2\bbZ$-quadratic unitary fusion categories.
While we do not find any new fusion categories in this article, we provide important techniques for finding 6j-symbols for fusion categories with multiplicity.
Our main theorem is as follows.

\begin{thmalpha}
\label{thm:Main}
The complete list of $\bbZ/2\bbZ$ quadratic UFCs is as follows.
\\
3 object categories:
\begin{itemize}
\item
the Ising/Tambara-Yamagami categories of the form $\cT\cY(\bbZ/2\bbZ,\chi,\pm)$ \cite{MR1659954} with $A_3$ fusion rules, of which there are exactly 2.
The case $+$ is Temperley-Lieb-Jones at $q=\exp(2\pi i/8)$, and the case $-$ is $SU(2)_2$.
\item
the three UFCs with $\Rep(S_3)$ fusion rules \cite[Theorem 5.1]{MR3635673}.
\item
the two complex cojugate UFCs with $\Ad(E_6)$ fusion rules \cite{MR1193933,MR1832764}.
These are exactly the adjoint subcategories of the exceptional quantum subgroups of Temperley-Lieb-Jones at $q=\exp(2\pi i/24)$ and $SU(2)_{10}$ \cite{MR1907188,MR1936496}.
\end{itemize}
4 object categories:
\begin{itemize}
\item
the pointed categories
$\Hilb(\bbZ/4\bbZ,\omega)$ where $\omega \in H^3(\bbZ/4\bbZ, U(1))$
and
$\Hilb(\bbZ/2\bbZ\times \bbZ/2\bbZ,\omega)$ where $\omega \in H^3(\bbZ/2\bbZ\times\bbZ/2\bbZ, U(1))/\Aut(\bbZ/2\bbZ\times\bbZ/2\bbZ ) \cong (\bbZ/2\bbZ)^3/\Aut(\bbZ/2\bbZ\times\bbZ/2\bbZ )$ \cite[Remark 4.10.4]{MR3242743}.
\item
the Deligne products $\Fib\boxtimes \Hilb(\bbZ/2\bbZ,\omega)$ for $\omega \in H^3(\bbZ/2\bbZ,U(1))$, which have $A_4$ fusion rules.
These two categories are also Temperley-Lieb-Jones at $q=\exp(2\pi i/10)$ and $SU(2)_3$.
\item
$\Ad(SU(2)_7)$, which is also equivalent to the adjoint subcategory of the $A_7$ Temperley-Lieb-Jones category with $q=\exp(2\pi i/16)$.
\item
the even parts of the two complex conjugate subfactor planar algebras with principal graphs $\cS'=\bigraph{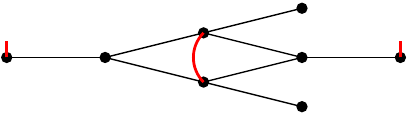}$ from \cite{MR3306607,MR3827808}.
These categories are also de-equivariantizations of $2^{\bbZ/4\bbZ}1$ near group fusion categories \cite[Ex.~9.5]{MR3635673} \cite[Ex.~2.2]{MR4079744}.
\item
the even part of the 2D2 subfactor planar algebra with principal graph
$\bigraph{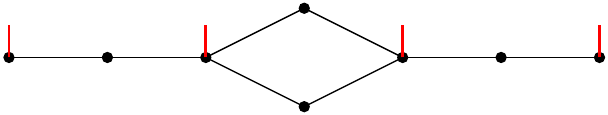}$ from \cite{MR3394622} \cite[Cor.~9.3]{MR3827808}.
This category is also a de-equivariantization of the even part of the $3^{\bbZ/4\bbZ}$ subfactor from \cite{MR3402358,MR3827808}.
\end{itemize}
All these UFCs are related to quantum groups at roots of unity or near group fusion categories \cite{MR3167494,MR3635673}.
\end{thmalpha}

The result \cite[Thm.~1.1]{MR3229513} gave a finite list of possible fusion rings for rank 4 pseudounitary fusion categories with a dual pair of simple objects, but included one fusion ring not previously known to be categorifiable (the case $c=2$ from \cite[Thm.~1.1(6)]{MR3229513}), and left open the classification of those fusion rings from \cite[Thm.~1.1]{MR3229513} which were previously known to be categorifiable.

\begin{coralpha}
\label{cor:FinishLarson}
We have a complete classification of rank 4 unitary fusion categories with a dual pair of simple objects.
In particular, there is no UFC with $c=2$ from \cite[Thm.~1.1(6)]{MR3229513}.
\end{coralpha}

One tool to prove our classification is an adaptation of Larson's rank-finiteness bound for $\mathbb{Z}/3\mathbb{Z}$-near group pseudounitary fusion categories \cite[\S4]{MR3229513}.
This adaptation appears in Appendix \ref{appendix:PseudounitaryBound} below, coauthored with attendees from the 2014 AMS MRC program on the Mathematics of Quantum Phases of Matter and Quantum Information.

Our main new technical tool to achieve Theorem \ref{thm:Main} is a generalization of Ostrik's results on formal codegrees of a spherical fusion category \cite{MR2576705,1309.4822}.
We use results of \cite[\S5]{MR1966525}, but we use the conventions of \cite{MR1782145}.
Suppose $\cC$ is a spherical fusion category, and denote by $\cF: Z(\cC) \to \cC$ the forgetful functor and let $\cI: \cC \to Z(\cC)$ be its adjoint.
Let $A$ be the tube algebra of $\cC$, and let $A_{X\leftarrow X}$ be the corner of $A$ corresponding to $X\in \Irr(\cC)$.
We pick a non-degenerate trace $\Tr_X$ on $A_{X\leftarrow X}$ given by
$$
\Tr_X\left(
\tikzmath{
    \filldraw[very thick, fill=gray!30] (0,0) circle (.2cm);
    \roundNbox{fill=white}{(-.2,.8)}{.3}{.2}{.2}{$f$}
    \draw (0,.2) -- node[right]{$\scriptstyle X$} (0,.5);
    \draw (0,1.1) arc (180:0:.3cm) -- node[right]{$\scriptstyle \overline{W}$} (.6,0) arc (0:-180:.5cm) -- node[left]{$\scriptstyle W$} (-.4,.5);
    \draw (-.4,1.1) -- node[left]{$\scriptstyle X$} (-.4,1.4);
}
\right)
:=
\delta_{W=\mathbf{1}}\dim(X)\tr_\cC(f).
$$
Given an irreducible representation $(V,\pi_V)$ of $A_{X\leftarrow X}$, its \emph{formal codegree} \cite{MR1974442,MR2576705} with respect to $\Tr_X$ is given by
$$
f_V:=\sum_b \Tr_V(\pi(b))\pi(b^\vee)
$$
where $\{b\}$ is a basis of $A_{X\leftarrow X}$ and $\{b^\vee\}$ is the dual basis with respect to the non-degenerate pairing $(a,b):= \Tr_X(ab)$ on $A_{X\leftarrow X}$.
Observe that $f_V$ is independent of the choice of basis $\{b\}$, but depends on the choice of $\Tr_X$.
We refer the reader to \S\ref{sec:FormalCodegrees} for more details.

\begin{thmalpha}\label{thm:tubecorner}
There is a bijective correspondence between irreducible representations $(V,\pi_V)$ of $A_{X\leftarrow X}$ and simple subobjects $\Gamma_V \subset \cI(X) \in Z(\cC)$. 
The formal codegree $f_V$ of $(V,\pi)$ with respect to $\Tr_X$ is a scalar, and the categorical dimension of $\Gamma_V$ is given by $\frac{\dim(\cC)}{f_V \dim(X)}$.
Moreover, if $Y\in \Irr(\cC)$ and ${}_X\pi_Y$ is the action of $A_{X\leftarrow X}$ on $A_{X\leftarrow Y}$, 
then
\[ 
\dim\Hom_{\cC}(\cF(\Gamma_V)\to Y) 
= 
\dim\Hom( \pi_V \to  {}_X\pi_Y  ).
\]
\end{thmalpha}

In the case $X=\mathbf{1}\in \Irr(\cC)$, this theorem recovers \cite[Thm.~2.13]{1309.4822}, which allowed the computation of the simple decomposition of $\cI(\mathbf{1})$ in terms of representations of the fusion algebra of $\cC$. 
Our theorem generalises this result in several ways. 
The main improvement is that this result allows us to determine the simple decomposition of $\cI(X)$ by studying the representations of the corner of the tube algebra $A_{X\leftarrow X}$. 
When $X = \mathbf{1}$, this algebra is isomorphic to the fusion algebra of $\cC$. 
However, when $X$ is non-trivial, this corner depends on certain 6j-symbols of the category involving $X$. 
One immediate application of this theorem comes from the fact that the dimensions of simple objects in $Z(\cC)$ are highly restricted, which implies the representations of $A_{X\to X}$ (which depend on the 6j-symbols) are also restricted. 
Hence we obtain obstructions based on 6j-symbols.
We make use of this application in this article to determine several non-trivial 6j-symbols involving the invertible object of a $\bbZ/2\bbZ$-quadratic UFC.

The other improvement Theorem \ref{thm:tubecorner} offers is that for each simple $\Gamma\subset \cI(X)$, we can determine $\cF(\Gamma) \in \cC$. 
This information is encoded in the action of $A_{X\leftarrow X}$ acting on the entire tube algebra $A$. 
As these algebras are semi-simple, it is easy to decompose this action into irreducibles, and hence apply Theorem \ref{thm:tubecorner}. 
A surprising application of this theorem comes from the fact that if we know the action of $A_{X\leftarrow X}$ acting on the entire tube algebra $A$ up to isomorphism, we can often determine the action on the nose.
As this action is determined by the 6j-symbols of $\cC$, this allows us to find many linear equations involving the 6j-symbols. 
We use this application later in this article to get our hands on many 6j-symbols. 
In the general setting, this result allows the combinatorial data of the forgetful functor $Z(\cC)\to \cC$ to be leveraged into the categorical data of the 6j-symbols of $\cC$. 
As the forgetful functor can often be easily determined from the fusion ring of $\cC$ \cite{MR3611056}, we expect this application to have many exciting future uses. 

%%%%%%%%%%%%%%%%%%%%%%%%%%%%%%%%%%%%%%%%%%%%%%%%%%%%%%
\subsection*{Acknowledgements}
Cain Edie-Michell was supported by NSF grant DMS 2055105 and a AMS-Simons Travel Grant.
Masaki Izumi was supported by JSPS KAKENHI Grant Number JP20H01805.
David Penneys was supported by NSF grants DMS 1654159, 1927098, and 2051170.

The authors of the appendix would like to thank the NSF-funded AMS MRC program and the organizers of the 2014 workshop on Quantum Phases of Matter and Quantum Information.

%%%%%%%%%%%%%%%%%%%%%%%%%%%%%%%%%%%%%%%%%%%%%%%%%%%%%%
%%%%%%%%%%%%%%%%%%%%%%%%%%%%%%%%%%%%%%%%%%%%%%%%%%%%%%
%%%%%%%%%%%%%%%%%%%%%%%%%%%%%%%%%%%%%%%%%%%%%%%%%%%%%%
\section{Preliminaries}
We refer the reader to \cite{MR3242743} for the basics of fusion categories.

%%%%%%%%%%%%%%%%%%%%%%%%%%%%%%%%%%%%%%%%%%%%%%%%%%%%%%
\subsection{The tube algebra}
One of the key tools in this paper is \textit{Ocneanu's tube algebra} (or equivalently the annular category) of a fusion category. 
This algebra was first introduced by \cite{MR996454} and \cite{MR1316301} in the context of subfactors, and by \cite{MR1782145,MR1832764} and \cite{MR1966525} in the context of fusion categories.

\begin{defn}
Let $\cC$ be a spherical fusion category whose spherical trace is denoted $\tr_\cC$. 
The tube algebra $A$ of $\cC$ is the finite dimensional semisimple algebra 
\[  
\bigoplus_{X,Y\in \Irr(\cC)} A_{Y \leftarrow X}    
\qquad\qquad\text{where}\qquad\qquad
A_{Y \leftarrow X} := \bigoplus_{W\in \Irr(\cC)}
\cC( W\otimes X \to Y\otimes W).\]
We graphically represent an element of $A$ as
\[
\tikzmath{
    \filldraw[very thick, fill=gray!30] (0,0) circle (.2cm);
    \roundNbox{fill=white}{(-.2,.8)}{.3}{.2}{.2}{$f$}
    \draw (0,.2) -- node[right]{$\scriptstyle X$} (0,.5);
    \draw (0,1.1) arc (180:0:.3cm) -- node[right]{$\scriptstyle \overline{W}$} (.6,0) arc (0:-180:.5cm) -- node[left]{$\scriptstyle W$} (-.4,.5);
    \draw (-.4,1.1) -- node[left]{$\scriptstyle Y$} (-.4,1.4);
}
\qquad\qquad
f\in \cC(W\otimes X \to Y\otimes W).
\]
The multiplication on $A$ is defined by composition of the tubes and applying the fusion relation obtained from semisimplicity to the strands around the outside.
In more detail, we pick a basis $\{\alpha\}\subset \cC(U\otimes V \to W)$ for all $U,V,W\in \Irr(\cC)$,
and let $\{\alpha^\vee\}\subset \cC(W\to U\otimes V)$ be the dual basis with respect to the non-degenerate pairing $(\,\cdot\,,\,\cdot\,) : \cC(U\otimes V \to W) \times \cC(W\to U\otimes V)\to \bbC$
determined by the formula $( h, k) \id_W = h\circ k \in \End_\cC(W)$.
We have the fusion relation
$$
\tikzmath{
\draw (-.2,-1.1) node[below]{$\scriptstyle U$} -- (-.2,1.1) node[above]{$\scriptstyle U$};
\draw (.2,-1.1) node[below]{$\scriptstyle V$} -- (.2,1.1) node[above]{$\scriptstyle V$};
}
=
\sum_{\substack{
W\in \Irr(\cC)
\\
\alpha
}}
\tikzmath{
\draw (-.2,-1.1) node[below]{$\scriptstyle U$} -- (-.2,-.8);
\draw (.2,-1.1) node[below]{$\scriptstyle V$} -- (.2,-.8);
\draw (-.2,1.1) node[above]{$\scriptstyle U$} -- (-.2,.8);
\draw (.2,1.1) node[above]{$\scriptstyle V$} -- (.2,.8);
\draw (0,-.3) --node[right]{$\scriptstyle W$} (0,.3);
\roundNbox{fill=white}{(0,-.5)}{.3}{.1}{.1}{$ \alpha $}
\roundNbox{fill=white}{(0,.5)}{.3}{.1}{.1}{$ \alpha^\vee$}
}
$$
which gives the following formula for composition in the tube algebra, which is independent of the choice of basis $\{\alpha\}$:
$$
\tikzmath{
    \filldraw[very thick, fill=gray!30] (0,0) circle (.2cm);
    \roundNbox{fill=white}{(-.2,.8)}{.3}{.1}{.1}{$f$}
    \draw (0,.2) -- node[right]{$\scriptstyle Y$} (0,.5);
    \draw (0,1.1) arc (180:0:.3cm) -- node[right]{$\scriptstyle \overline{U}$} (.6,0) arc (0:-180:.5cm) -- node[left]{$\scriptstyle U$} (-.4,.5);
    \draw (-.4,1.1) -- node[left]{$\scriptstyle Z$} (-.4,1.4);
}
\cdot
\tikzmath{
    \filldraw[very thick, fill=gray!30] (0,0) circle (.2cm);
    \roundNbox{fill=white}{(-.2,.8)}{.3}{.1}{.1}{$g$}
    \draw (0,.2) -- node[right]{$\scriptstyle X$} (0,.5);
    \draw (0,1.1) arc (180:0:.3cm) -- node[right]{$\scriptstyle \overline{V}$} (.6,0) arc (0:-180:.5cm) -- node[left]{$\scriptstyle V$} (-.4,.5);
    \draw (-.4,1.1) -- node[left]{$\scriptstyle Y$} (-.4,1.4);
}
:=
\tikzmath{
    \filldraw[very thick, fill=gray!30] (0,0) circle (.2cm);
    \roundNbox{fill=white}{(-.2,.8)}{.3}{.1}{.1}{$g$}
    \roundNbox{fill=white}{(-.6,1.7)}{.3}{.1}{.1}{$f$}
    \draw (0,.2) -- node[right,xshift=-.1cm]{$\scriptstyle X$} (0,.5);
    \draw (0,1.1) arc (180:0:.3cm) -- node[right,xshift=-.1cm]{$\scriptstyle \overline{V}$} (.6,0) arc (0:-180:.5cm) -- node[left,xshift=.1cm]{$\scriptstyle V$} (-.4,.5);
    \draw (-.4,1.1) -- node[left,xshift=.1cm]{$\scriptstyle Y$} (-.4,1.4);
    \draw (-.4,2) .. controls ++(90:.5cm) and ++(90:.5cm) .. (1,2) -- node[right,xshift=-.1cm]{$\scriptstyle \overline{U}$} (1,.2) arc (0:-180:.9cm) -- node[left,xshift=.1cm]{$\scriptstyle U$} (-.8,1.4);
    \draw (-.8,2) -- node[left,xshift=.1cm]{$\scriptstyle Z$} (-.8,2.5);
}
=
\sum_{\substack{
W\in \Irr(\cC)
\\
\alpha
}}
\tikzmath{
    \filldraw[very thick, fill=gray!30] (0,0) circle (.2cm);
    \roundNbox{fill=white}{(-.2,2.6)}{.25}{0}{0}{$\scriptstyle \alpha$}
    \roundNbox{fill=white}{(-.6,-.1)}{.25}{0}{0}{$\scriptstyle \alpha^\vee$}
    \roundNbox{fill=white}{(-.2,.8)}{.3}{.1}{.1}{$g$}
    \roundNbox{fill=white}{(-.6,1.7)}{.3}{.1}{.1}{$f$}
    \draw (0,.2) -- node[right,xshift=-.1cm]{$\scriptstyle X$} (0,.5);
    \draw (-.2,2.85) arc (180:0:.3cm) -- node[right,xshift=-.1cm]{$\scriptstyle \overline{W}$} (.4,-.35) arc (0:-180:.5cm) node[left, yshift=-.2cm,xshift=.1cm]{$\scriptstyle W$};
    \draw (-.4,1.1) -- node[left,xshift=.1cm]{$\scriptstyle Y$} (-.4,1.4);
    \draw (-.4,2) --node[left,xshift=.1cm]{$\scriptstyle U$} (-.3,2.35);
    \draw (0,1.1) --node[right,xshift=-.1cm]{$\scriptstyle V$} (-.1,2.35);
    \draw (-.4,.5) --node[right,xshift=-.1cm]{$\scriptstyle V$} (-.5,.15);
    \draw (-.8,1.4) --node[left,xshift=.1cm]{$\scriptstyle U$} (-.7,.15);
    \draw (-.8,2) -- node[left,xshift=.1cm]{$\scriptstyle Z$} (-.8,3.3);
}\,.
$$

There is a non-degenerate linear functional $\phi$ on $A$ given by 
\[
\tikzmath{
    \filldraw[very thick, fill=gray!30] (0,0) circle (.2cm);
    \roundNbox{fill=white}{(-.2,.8)}{.3}{.2}{.2}{$f$}
    \draw (0,.2) -- node[right]{$\scriptstyle X$} (0,.5);
    \draw (0,1.1) arc (180:0:.3cm) -- node[right]{$\scriptstyle \overline{W}$} (.6,0) arc (0:-180:.5cm) -- node[left]{$\scriptstyle W$} (-.4,.5);
    \draw (-.4,1.1) -- node[left]{$\scriptstyle Y$} (-.4,1.4);
}
\longmapsto 
\delta_{X=Y}\delta_{W=\mathbf{1}}\dim(X)\tr_\cC(f).
\]
Its restriction to $A_{X\leftarrow X}$ is tracial for all $X\in \Irr(\cC)$, and we denote it by $\Tr_X$.
\end{defn}

Note that each of the spaces $A_{X \leftarrow X}$ is the corner $1_XA1_X$ of $A$, where we cut down by the idempotent 
$$
1_X
:=
\tikzmath{
    \filldraw[very thick, fill=gray!30] (0,0) circle (.2cm);
    \draw (0,.2) -- node[right]{$\scriptstyle X$} (0,.6);
},
$$
and $A_{X \leftarrow X}$ acts on the spaces $A_{X \leftarrow Y}$ by multiplication.

The tube algebra of $\cC$ is intimately related to the Drinfeld centre $Z(\cC)$ of $\cC$. 
From the data of $Z(\cC)$, we obtain a basis of matrix units for the spaces $A_{X\leftarrow Y}$ given by
\[ 
e(\Gamma)_{(X,i),(Y,j)} 
:= 
\frac{\dim(\Gamma)}{\dim(\cC)\sqrt{\dim(X)\dim(Y)}}
\sum_{W\in \Irr{\cC}}
\dim(W) 
\tikzmath{
    \filldraw[very thick, fill=gray!30] (0,0) circle (.2cm);
    \draw (0,.2) -- node[right]{$\scriptstyle X$} (0,.5);
    \draw (0,.9) -- node[right]{$\scriptstyle \Gamma$} (0,1.2);
    \draw (0,1.8) arc (180:0:.3cm) -- node[right]{$\scriptstyle \overline{W}$} (.6,0) arc (0:-180:.5cm) -- node[left]{$\scriptstyle W$} (-.4,1.2);
    \draw (-.4,1.8) -- node[left]{$\scriptstyle \Gamma$} (-.4,2.1);
    \draw (-.4,2.5) -- node[left]{$\scriptstyle Y$} (-.4,2.8);
    \roundNbox{fill=white}{(-.4,2.3)}{.2}{.05}{.05}{$\scriptstyle j'$}
    \roundNbox{fill=white}{(0,.7)}{.2}{0}{0}{$\scriptstyle i$}
    \roundNbox{fill=white}{(-.2,1.5)}{.3}{.2}{.2}{$\beta_{W,\Gamma}$}
}
\]
where $(\Gamma, \beta_\Gamma) \in \Irr(Z(\cC))$, 
$\{i\}$ is a basis for $\cC(X\to \cF(\Gamma))$, 
and  $\{j\}$ is a basis for $\cC(Y\to \cF(\Gamma))$, 
where $\cF: Z(\cC) \to \cC$ is the forgetful functor. 
Here, $\{j'\}\subset \cC(\cF(\Gamma) \to Y)$ denotes the dual basis of $\{j\}$ with respect to the pairing $k'\circ j = \delta_{j=k}\id_{Y}$.
With respect to our functional $\phi$ on $A$, we have that
\[
\phi( e(\Gamma)_{(X,i),(Y,j)} ) = \delta_{X,Y} \delta_{i,j}   \frac{\dim(X)\dim(\Gamma)}{\dim(\cC)},    
\]
and so the dual basis with respect to the $\phi$-pairing is given by 
\[ e(\Gamma)_{(X,i),(Y,j)}^\vee = \frac{\dim(\cC)}{\dim(X)\dim(\Gamma)}  e(\Gamma)_{(Y,j),(X,i)}.
\]

The construction above shows us that $Z(\cC)$ entirely determines the structure of the tube algebra of $\cC$. 
The converse is also true. The tube algebra of $\cC$ entirely determines the Drinfeld centre of $\cC$. 
The following theorem gives a bijective correspondence between representations of the tube algebra and objects of $Z(\cC)$.
\begin{thm}[{\cite{MR1782145} and \cite[\S5]{MR1966525}}]
\label{thm:tubecentre}
There is a bijective correspondence between equivalence classes of irreducible representations of the tube algebra of $\cC$ and isomorphism classes of simple objects in $Z(\cC)$. 
This bijection sends
\[
(V,\pi) \mapsto \Gamma_V:= \bigoplus_{X\in \Irr(\cC)}  V|_{A_{X\leftarrow X}}\otimes X.      
\]
Further, we have that the minimal central projection $z_V\in A$ corresponding to $(V,\pi)$ is given by 
\[ 
z_V
=
\sum_{\substack{
X\in \Irr(\cC),
\\ 
\{i\}\subset \cC( X\to \cF( \Gamma_V ))
}} 
e(\Gamma_V)_{(X,i),(X,i)}.
\]
\end{thm}

%%%%%%%%%%%%%%%%%%%%%%%%%%%%%%%%%%%%%%%%%%%%%%%%%%%%%%%%%%%%%%%
\subsection{A new result on formal codegrees}
\label{sec:FormalCodegrees}

If one knows the full tube algebra of $\cC$, then Theorem~\ref{thm:tubecentre} essentially gives you the full data of $Z(\cC)$. 
However in many situations, such as in this article, we only know information about some sub-algebra of the tube algebra, and we wish to leverage this information into partial information about $Z(\cC)$. 
Towards this goal, we introduce the \textit{formal codegree} of a representation.
\begin{defn}[\cite{MR2576705,MR1974442}]
Let $B$ be a finite dimensional semisimple algebra equipped with a non-degenerate trace $\Tr_B$,
and let $(V,\pi)$ be a finite dimensional representation of $B$.
We define the \emph{formal codegree} of $(V,\pi)$ as 
\[
f_V := \sum_{b} \Tr_V(\pi(b))\pi(b^\vee) \in \pi(B)\subset \End(V)
\]
where we sum over a basis $\{b\}\subset B$, and $\{b^\vee\}$ denotes the dual basis with respect to the $\Tr_B$-pairing.
Observe that $f_V$ is independent of the choice of basis,
but depends on the choice of trace $\Tr_B$.
\end{defn}

The following theorem allows us to determine the simple summands of $\cI(X)\in Z(\cC)$ by classifying the representations of the subalgebra $A_{X\leftarrow X}$.
Here, $\cI: \cC \to Z(\cC)$ is the induction functor which is adjoint to the forgetful functor $\cF: Z(\cC) \to \cC$.
Moreover, we can compute categorical dimensions in terms of formal codegrees of $A_{X\leftarrow X}$ with respect to $\Tr_X$.

\begin{thm*}[\ref{thm:tubecorner}]
Let $\cC$ be a spherical fusion category, and let $A$ be the tube algebra of $\cC$. 
Fix $X\in \Irr(\cC)$. 
There is a bijective correspondence between equivalence classes of irreducible representations $(V,\pi)$ of $A_{X\leftarrow X}$ and isomorphism classes of simple subobjects $\Gamma_V \subset \cI(X) \in Z(\cC)$. 
The formal codegree $f_V$ of $(V,\pi)$ with respect to $\Tr_X$ is a scalar, and the categorical dimension of $\Gamma_V$ is given by 
$$
\dim(\Gamma_V) = \frac{\dim(\cC)}{f_V \dim(X)}.
$$
Moreover, if $Y\in \Irr(\cC)$ and ${}_X\pi_Y$ is the action of $A_{X\leftarrow X}$ on $A_{X\leftarrow Y}$, 
then
\[ 
\dim(\cC(Y\to \cF(\Gamma_V))) = \dim(\Hom(\pi_V \to  {}_X\pi_Y  )).  
\]
\end{thm*}
\begin{proof}
Let $(V,\pi)$ be an irreducible representation of $A_{X\leftarrow X}$. 
Since $A_{X\leftarrow X}$ is semisimple, $(V,\pi)$ corresponds to a simple summand of $A_{X\leftarrow X}$.
As $A_{X\leftarrow X}$ is a corner of $A$,
each simple summand of $A_{X\leftarrow X}$ is of the form $A_{X\leftarrow X}z_\Gamma$
for a simple object $(\Gamma,\beta)\in \Irr(Z(\cC))$.
Hence there is a simple $(\Gamma_V,\beta_{\Gamma_V})$ corresponding to $(V,\pi)$,
and by Theorem~\ref{thm:tubecentre}, $z_V1_X = \sum_i e(\Gamma_V)_{(X,i),(X,i)}$.
Moreover, for any other simple object $\Lambda\in \Irr(Z(\cC))$, we have that $\pi( e(\Lambda)_{(X,i),(X,i)}) = 0$ unless $\Lambda \cong \Gamma_V$. 
In particular, 
$\Tr_V(\pi( e(\Lambda)_{(X,i),(X,j)}))=0$ 
unless 
$\Lambda = \Gamma_V$ and $i=j$. 
We now compute
\begin{align*} 
f_V 
&= 
\sum_{\substack{
\Lambda \subseteq \cI(X)
\\ 
i,j
% \in \ONB(X \to \cF(\Lambda))
}}    
\Tr_V(  \pi( e(\Lambda)_{(X,i),(X,j)})) \pi(e(\Lambda)_{(X,i),(X,j)}^\vee) \\
&= 
\sum_{i
%\in \ONB(X \to \cF(\Gamma_V))
} 
\frac{\dim(\cC)}{\dim(X)\dim(\Gamma_V)}
\pi(e(\Gamma_V)_{(X,i),(X,i)})
\\
&= \frac{\dim(\cC)}{\dim(X)\dim(\Gamma_V)} \pi(z_V1_X).  
\end{align*}
Thus the formal codegree of $(V,\pi)$ is given by $f_V =\frac{\dim(\cC)}{\dim(X)\dim(\Gamma_V)}$,
and
$\dim(\Gamma_V) = \frac{\dim(\cC)}{f_V \dim(X)}$.
%Now since
Finally, we observe
\begin{align*}
\dim(\Hom(V \to {}_X\pi_Y))
&=
\dim(\Hom(V \to 1_XA1_Y))
\\&=
\dim(\Hom(V \to z_V1_XA1_Y))
\\&=
\dim(\Hom(V \to z_V1_XA1_Yz_V))
\\&=
\sum_{
j 
%\in \ONB(Y\to \cF(\Gamma))
}
\underbrace{\dim(\Hom(V \to  z_V1_X A e(\Gamma_V)_{(Y,j),(Y,j)}))}_{=1}
\\&=
\dim(\cC(Y\to \cF(\Gamma_V))).
\qedhere
\end{align*}
\end{proof}

Note that if we just consider the subalgebra $A_{\mathbf{1}\leftarrow \mathbf{1}} \cong K_0(\cC)$, the fusion algebra of $\cC$,
then the above theorem recovers \cite[Theorem 2.13]{1309.4822},
which shows that irreducible representations of $K_0(\cC)$ are in bijective correspondence with simple summands of $\cI(\mathbf{1})$. 
Thus our theorem generalises Ostrik's in two ways:
(1) it gives us the simple summands of $\cI(X)$ where $X$ is any simple object of $\cC$, 
and 
(2) it tells us the image under the forgetful functor of each of these summands.

%%%%%%%%%%%%%%%%%%%%%%%%%%%%%%%%%%%%%%%%%%%%%%%%%%%%%%%%%%%%%%%
\subsection{ \texorpdfstring{$\bbZ/2\bbZ$}{Z/2Z}-quadratic fusion categories}

A $\bbZ/2\bbZ$-quadratic fusion category is a fusion category $\cC$
whose invertible objects form the group $\bbZ/2\bbZ$, i.e., $\Inv(\cC) = \{1,\alpha\}$ with $\alpha^2\cong \mathbf{1}$,
with one other orbit of simple objects under the $\bbZ/2\bbZ$-action.
A simple associativity argument shows we have three cases:
\begin{enumerate}[label=(Q\arabic*)]
\item 
\label{Q:3Objects}
simple objects: $\mathbf{1},\alpha,\rho$;
fusion rules determined by:
$\rho^2 \cong 1\oplus m \rho \oplus \alpha$.

\item
\label{Q:4ObjectsNonSelfDual}
simple objects: $\mathbf{1},\alpha,\rho, \alpha\rho$, $\rho$ not self-dual;
fusion rules determined by:
$\rho^2 \cong m \rho\oplus n \alpha \rho\oplus \alpha$.
\item
\label{Q:4ObjectsSelfDual}
simple objects: $\mathbf{1},\alpha,\rho, \alpha\rho$, $\rho$ self-dual;
fusion rules determined by:
$\rho^2 \cong 1\oplus m \rho\oplus n \alpha \rho$.
\end{enumerate}

Note that in all three cases we have $\alpha^2\cong \mathbf{1}$ and $\alpha \rho \cong \rho \alpha$.
\subsubsection{Multiplicity bounds and categorifiability}

The case \ref{Q:3Objects} was classified in the pivotal setting in \cite[Thm.~4.1]{1309.4822}, where it was shown $m\leq 2$.
The complete classification of such unitary fusion categories was known prior to this article:
\begin{enumerate}[label=(${m=}$ \arabic*)]
\setcounter{enumi}{-1}
\item 
such a fusion category is a Tambara-Yamagami category of the form $\cT\cY(\bbZ/2\bbZ,\chi,\pm)$ \cite{MR1659954}, of which there are exactly 2.
Both such categories are unitarizable.
\item
such a fusion category has the fusion rules of $\Rep(S_3)$. There are exactly three such unitary fusion categories \cite[Theorem 5.1]{MR3635673}.
\item
such a fusion category has the fusion rules of $\Ad(E_6)$, and there are exactly 4 such fusion categories \cite{MR2559711}, all within the same Galois orbit, and each admits a spherical structure.
Two of these are unitary and complex conjugate to each other  \cite{MR1832764}.
\end{enumerate}

The case \ref{Q:4ObjectsNonSelfDual} was studied in the pseudounitary setting ($\dim(\cC) = \FPdim(\cC)$) in \cite{MR3229513}, where it was shown that $m=n\leq 2$. 
The classification of such fusion categories prior to this article is as follows:
\begin{enumerate}[label=(${m=}$ \arabic*)]
\setcounter{enumi}{-1}
\item 
such a fusion category is necessarily pointed with $\bbZ/4\bbZ$ fusion rules.
It is thus of the form $\Vect(\bbZ/4\bbZ,\omega)$ for $\omega\in H^3(\bbZ/4\bbZ,U(1))=\bbZ/4\bbZ$, of which there are 4 categories \cite[Remark 4.10.4]{MR3242743}.
\item
this case was open.
Two such unitary fusion categories which are complex conjugate were known to exist from \cite{MR3306607,MR3827808}.
\item
this case was open.
No such examples were known to exist.
\end{enumerate}
We finish this classification for unitary fusion categories in Theorem \ref{thm:NSD} below.

In Appendix \ref{appendix:PseudounitaryBound}, we adapt the results of \cite{MR3229513} in the pseudounitary setting to case \ref{Q:4ObjectsSelfDual}, where we prove the following theorem.

\begin{thm}
\label{thm:PseudounitaryBound}
Suppose $\cC$ is a pseudounitary fusion category with the fusion rules \ref{Q:4ObjectsSelfDual}.
Then $(m,n)$ must be equal to one of
$(0,0),(0,1), (1,0),(1,1),(2,2)$.
\end{thm}
\begin{proof}
By Theorem \ref{thm:TwoParameterBound} in Appendix \ref{appendix:PseudounitaryBound}, we must have $m+n\leq 5$.
If either $m$ or $n$ is zero, then there is a fusion subcategory with 2 simple objects, so $(m,n)$ must be one of $(0,0),(0,1),(1,0)$ by \cite{MR1981895}.
If $0\neq m\neq n\neq 0$, then $m+n\geq 11$ by Remark \ref{rem:K2abBetaConditions} in Appendix \ref{appendix:PseudounitaryBound}.
The result follows.
\end{proof}

The proof that $m+n \leq 5$ that appears in Appendix \ref{appendix:PseudounitaryBound} below was written by
Ryan Johnson, Siu-Hung Ng, David Penneys, Jolie Roat, Matthew Titsworth, and Henry Tucker
at the 2014 AMS MRC on The Mathematics of Quantum Phases of Matter and Quantum Information.

By \cite{MR3345186,MR3536926} (and applying Galois conjugation),
any fusion category with fusion rules \ref{Q:4ObjectsSelfDual} with $(m,n) \in \{(0,1),(1,0)\}$ factorizes as a Deligne product of a fusion category with Fibonacci fusion rules, of which there are two, namely $\Fib$ and $\mathsf{YL}$, and a $\bbZ/2\bbZ$-pointed fusion category which must be of the form $\Vect(\bbZ/2\bbZ,\omega)$ for $\omega \in H^3(\bbZ/2\bbZ,U(1))$, of which there are two.
Thus there are exactly 4 such fusion categories, and 2 are unitarizable.

When $m=n\leq 2$, the complete classification of such unitary fusion categories as in Theorem \ref{thm:PseudounitaryBound} is given in Theorem \ref{thm:SD} below.

%%%%%%%%%%%%%%%%%%%%%%%%%%%%%%%%%%%%%%%%%%%%%%%%%%%%%%
%%%%%%%%%%%%%%%%%%%%%%%%%%%%%%%%%%%%%%%%%%%%%%%%%%%%%%
%%%%%%%%%%%%%%%%%%%%%%%%%%%%%%%%%%%%%%%%%%%%%%%%%%%%%%
\section{The self-dual case}
In this section we will focus on the unitary categorification of the fusion ring with four simple objects
$\mathbf{1},\alpha, \rho, \alpha\rho$
and fusion rules
\begin{equation}
\label{eq:R(m)}
\alpha \otimes \alpha\cong \mathbf{1}
\qquad\qquad
\rho\otimes \rho \cong \mathbf{1} \oplus m \rho \oplus m\alpha \rho.
\tag{$R(m)$}   
\end{equation}
Let us write \ref{eq:R(m)} for such a fusion ring. 
By Theorem \ref{thm:PseudounitaryBound} above, we know $m\leq 2$.
Our main result of this section is as follows.

\begin{thm}\label{thm:SD}
The complete classification of unitary fusion categories $\cC_m$ with $K_0(\cC_m) \cong$ \ref{eq:R(m)} for $m\leq 2$ is as follows:
\begin{enumerate}[label=(${m=}$ \arabic*)]
\setcounter{enumi}{-1}
\item 
$\cC_0$ is pointed and thus equivalent to one of the four monoidally distinct categories $\Vect^\omega(\bbZ/2\bbZ\times\bbZ/2\bbZ)$ where $\omega \in H^3(\mathbb{Z}/2\bbZ\times\mathbb{Z}/2\bbZ, U(1))/\Aut(\bbZ/2\bbZ\times\bbZ/2\bbZ ) \cong (\bbZ/2\bbZ)^3/\Aut(\mathbb{Z}/2\bbZ\times\mathbb{Z}/2\bbZ )$ \cite[Remark 4.10.4]{MR3242743}.
\item
$\cC_1$ is equivalent to $\cC(\mathfrak{sl}_2,7)^\text{ad}$,
which is also equivalent to the even part of the $A_7$ Temperley-Lieb-Jones category with $q=\exp(2\pi i/16)$ \cite[Example 9.1]{MR3827808}.
\item
$\cC_2$ is equivalent to the even part of the $2D2$ subfactor from \cite{MR3394622,MR3827808}.
\end{enumerate}
\end{thm}
\begin{proof}
It suffices to restrict our attention to the cases of $m=1$ and $m=2$.
%We begin by outlining our general method to prove Theorem~\ref{thm:SD}. 
The first step in our proof is to provide a set of numerical data which fully classifies a categorification of \ref{eq:R(m)}; we do this in \S\ref{sec:SDNumericalData}. 
By describing a sufficient list of local relations in our category, we are able to come up with such a set of numerical data. 
This data consists of $8m^4$ complex scalars, and a collection of small roots of unity. 
This data is precisely a subset of the $6j + 4k$ symbols of such a categorification. 
Using techniques developed in the subfactor classification program, we prove that this subset of the $6j + 4k$ symbols is sufficient to describe the entire category.

In \S\ref{SD:CenterAnalysis}, in order to get a foothold on the numerical data of a categorification of \ref{eq:R(m)}, we study the Drinfeld centre of such a category. 
By studying certain small representations of the tube algebra of the category (using basic combinatorial arguments), we are able to deduce a surprising amount of numerical data of the category. 
This centre analysis tells us nearly all of the small roots of unity in our numerical data, and even gives us highly non-trivial linear equations involving the $8m^4$ complex scalars.

To reduce the $8m^4$ complex scalars down to a more manageable number, in \S\ref{sec:SDSymmetries}, we apply the tetrahedral symmetries of the $6j + 4k$ symbols. 
These symmetries only apply in the unitary setting, and give $S_4$ symmetries of these $8m^4$ complex scalars.\footnote{\label{footnote:TetrahedralSymmetry} While writing this article, the article \cite{2106.16186} was posted to the arXiv, which describes tetrahedral symmetries for general fusion categories.
It would be interesting to use their work to extend our results to the non-unitary setting.}
This essentially reduces the complexity of the problem by a factor of 24. 
For example, in the $m=2$ case, we reduce from $128$ complex scalars to roughly $10$ (some of the $S_4$ symmetries are redundant). 
These symmetries turn an intractable amount of data into a set that can easily be dealt with by hand.

To finish off this section we explicitly solve for the remaining numerical data which describes a categorification of \ref{eq:R(m)} in \S\ref{sec:SDClassification}. 
The results of the previous subsections essentially determine everything except the remaining complex scalars. 
By evaluating diagrams in our category in multiple ways, we are able to obtain equations of these complex scalars. 
In the $m=2$ case we find a single solution, which necessarily has to correspond to the even part of the $2D2$ subfactor. 
We prove this in Theorem \ref{thm:2D2} below.
%This finishes the proof of Theorem~\ref{thm:SD}.
\end{proof}

%%%%%%%%%%%%%%%%%%%%%%%%%%%%%%%%%%%%%%%%%%%%%%%%%%%%%%%%%%%%%%%%%%%
\subsection{Numerical data}
\label{sec:SDNumericalData}

Let $\cC_m$ be a unitary fusion category with $K(\cC_m) \cong$ \ref{eq:R(m)}, $m\in\{1,2\}$. 
The goal of this subsection is to give a list of numerical data which completely describes the category $\cC_m$. 
We will show that the category $\cC_m$ can be completely described by the following data:
\begin{itemize}
    \item two choices of signs $\lambda_\rho,\lambda_\alpha \in \{-1,1\}$ which are the 2nd Frobenius Schur indicators of $\alpha$ and $\rho$ respectively,
    \item a choice of sign $\mu\in \{-1,1\}$,
    \item $2m$ choices of $\chi_{\mathbf{1}, i} \in \{-1,1\}$ and $\chi_{\alpha, i} \in \{-\sqrt{\lambda_\alpha}, \sqrt{\lambda_\alpha}\}$ for $0\leq i < m$,
    \item two involution's $(\mathbf{1}, i) \mapsto (\mathbf{1}, \tilde{i})$ and $(\alpha, i) \mapsto (\alpha, \tilde{i})$,
    \item $2m$ choices of 3rd roots of unity $\omega_{\mathbf{1},i},\omega_{\alpha,i}  \in \{1, e^{2\pi i \frac{1}{3}},e^{2\pi i \frac{2}{3}} \}$ for $0\leq i < m$, and
    \item $8m^4$ complex scalars $A^{i,j}_{k,\ell},B^{i,j}_{k,\ell},C^{i,j}_{k,\ell},D^{i,j}_{k,\ell},\widehat{A}^{i,j}_{k,\ell},\widehat{B}^{i,j}_{k,\ell},\widehat{C}^{i,j}_{k,\ell},\widehat{D}^{i,j}_{k,\ell} \in \mathbb{C}$ for $0 \leq i,j,k,\ell <m$. These complex scalars are entries of the $F$-tensors $F^{\rho, \rho, \rho}_{\rho}, F^{\rho, \rho, \rho}_{\alpha \rho},F^{\alpha\rho, \rho, \rho}_{\rho}$, and $F^{\alpha\rho, \rho, \rho}_{\alpha\rho}$.
\end{itemize}
While at first glance this may seem like an intractable amount of data to classify, the following subsections will reduce this to a far more manageable list. 

%%%%%%%%%%%%%%%%%%%%%%%%%%%%%%%%%%%%%%%%%%%%%%%%%%%%%%%%%%%%%%%%
\subsubsection{Generators, basic relations, and Frobenius maps}
\label{sec:SDBasicRelations}

To simplify notation, we define $d:=\dim(\rho)$, which is the largest solution to $d^2=1+2md$. 
If $m=1$ then $d = 1+\sqrt{2}$, and if $m=2$ then $d = 2+\sqrt{5}$.
We choose orthonormal basis for our hom spaces
$$
\tikzmath{
\draw[mid>] (-.3,-.6) node[below]{$\scriptstyle \rho$} -- (0,0);
\draw[mid>] (.3,-.6) node[below]{$\scriptstyle \rho$} -- (0,0);
\draw[mid<] (0,.6) node[above]{$\scriptstyle \rho$} -- (0,0);
\filldraw (0,0) node[left]{$\scriptstyle i$} circle (.05cm);
}
\in \cC_m(\rho\otimes \rho \to \rho)
\qquad\qquad
\tikzmath{
\draw[thick, red, mid>] (0,0) -- (-.6,.6) node[above]{$\scriptstyle \alpha$};
\draw[mid>] (-.6,-.6) node[below]{$\scriptstyle \rho$} -- (0,0);
\draw[mid>] (.6,-.6) node[below]{$\scriptstyle \rho$} -- (0,0);
\draw[mid>] (0,0) -- (.6,.6) node[above]{$\scriptstyle \rho$};
\filldraw (0,0) node[left]{$\scriptstyle i$} circle (.05cm);
}
\in \cC_m(\rho\otimes \rho \to \alpha\rho)
\qquad\qquad
0 \leq i < m.
$$
We also choose unitary isomorphisms\footnote{Using the convention of switching the orientation of the $\alpha$-strand through the crossing works better for $\bbZ/2\bbZ$-equivariantization, which is related to the $3^{\bbZ/4\bbZ}$-subfactor \cite{MR3827808}.
In the non-self-dual case in \S\ref{sec:NSD} below, we will use a more natural convention from a diagrammatic point of view which does not change the orientation of the $\alpha$-strand.
}
\[  
\tikzmath{
\draw[thick, red, mid>] (0,-.5) node[below]{$\scriptstyle \alpha$} -- (0,0);
\draw[thick, red, mid>] (0,.5) node[above]{$\scriptstyle \overline{\alpha}$} -- (0,0);
\draw[thick, red] (0,0) -- (.1,0);
}\in \cC_m(\alpha\to \overline{\alpha})
\quad 
\tikzmath{
\draw[mid>] (0,-.5) node[below]{$\scriptstyle \rho$} -- (0,0);
\draw[mid>] (0,.5) node[above]{$\scriptstyle \overline{\rho}$} -- (0,0);
\draw[thick] (0,0) -- (.1,0);
}\in \cC_m(\rho\to \overline{\rho})
\quad  
\text{ and } 
\quad 
\tikzmath{
\draw[mid>] (-.5,-.5) node[below]{$\scriptstyle \rho$} -- (0,0);
\draw[mid>] (0,0) -- (.5,.5) node[above]{$\scriptstyle \rho$};
\draw[thick, red, mid>] (.5,-.5) node[below]{$\scriptstyle \alpha$} -- (0,0);
\draw[thick, red, mid>] (-.5,.5) node[above]{$\scriptstyle \overline{\alpha}$} -- (0,0);
}
\in \cC_m(\rho \otimes\alpha\to \overline{\alpha}\otimes \rho).
\]
We can unambiguously write their inverses as
\[ 
\tikzmath{
\draw[thick, red, mid<] (0,-.5) node[below]{$\scriptstyle \overline{\alpha}$} -- (0,0);
\draw[thick, red, mid<] (0,.5) node[above]{$\scriptstyle \alpha$} -- (0,0);
\draw[thick, red] (0,0) -- (.1,0);
}\in \cC_m(\overline{\alpha}\to \alpha)
,\quad 
\tikzmath{
\draw[mid<] (0,-.5) node[below]{$\scriptstyle \overline{\rho}$} -- (0,0);
\draw[mid<] (0,.5) node[above]{$\scriptstyle \rho$} -- (0,0);
\draw[thick] (0,0) -- (.1,0);
}\in \cC_m(\overline{\rho}\to \rho)
,\quad  
\text{ and } 
\quad 
\tikzmath{
\draw[mid>] (.5,-.5) node[below]{$\scriptstyle \rho$} -- (0,0);
\draw[mid>] (0,0) -- (-.5,.5) node[above]{$\scriptstyle \rho$};
\draw[thick, red, mid<] (-.5,-.5) node[below]{$\scriptstyle \overline{\alpha}$} -- (0,0);
\draw[thick, red, mid<] (.5,.5) node[above]{$\scriptstyle \alpha$} -- (0,0);
}\in\cC_m( \overline{\alpha}\otimes \rho\to\rho \otimes\alpha).
\]
The one-click rotation of these first two isomorphisms are related to their inverses respectively by the Frobenius-Schur indicators of $\alpha$ and $\rho$, via the following equations: 
\[  
\tikzmath{
\draw[thick, red, mid>] (0,-.5) node[below]{$\scriptstyle \alpha$} -- (0,0);
\draw[thick, red, mid>] (0,.5) node[above]{$\scriptstyle \overline{\alpha}$} -- (0,0);
\draw[thick, red] (0,0) -- (.1,0);
}
=
\lambda_\alpha
\tikzmath{
\draw[thick, red, mid>] (0,-.5) node[below]{$\scriptstyle \alpha$} -- (0,0);
\draw[thick, red, mid>] (0,.5) node[above]{$\scriptstyle \overline{\alpha}$} -- (0,0);
\draw[thick, red] (0,0) -- (-.1,0);
}
%\raisebox{-.5\height}{ \includegraphics[scale = .5]{graphs/U.pdf}} = \lambda \raisebox{-.5\height}{ \includegraphics[scale = .6]{graphs/uleft.pdf}} 
\qquad \qquad 
\tikzmath{
\draw[mid>] (0,-.5) node[below]{$\scriptstyle \rho$} -- (0,0);
\draw[mid>] (0,.5) node[above]{$\scriptstyle \overline{\rho}$} -- (0,0);
\draw[thick] (0,0) -- (.1,0);
}
=
\lambda_\rho
\tikzmath{
\draw[mid>] (0,-.5) node[below]{$\scriptstyle \rho$} -- (0,0);
\draw[mid>] (0,.5) node[above]{$\scriptstyle \overline{\rho}$} -- (0,0);
\draw[thick] (0,0) -- (-.1,0);
}
\qquad\qquad
\lambda_\alpha, \lambda_\rho \in \{\pm 1\}.
%\raisebox{-.5\height}{ \includegraphics[scale = .5]{graphs/Up.pdf}} = \epsilon_\rho \raisebox{-.5\height}{ \includegraphics[scale = .6]{graphs/pright.pdf}}   
\]
We can re-scale the crossing so that 
\[  
\tikzmath{
\draw[mid>] (-.8,-.8) node[below]{$\scriptstyle \rho$} -- (0,0);
\draw[mid>] (0,0) -- (.5,.5) node[above]{$\scriptstyle \rho$};
\draw[thick, red, mid<] (.5,-.5) -- (0,0);
\draw[thick, red, mid>] (.8,-.8) node[below]{$\scriptstyle \alpha$} -- (.5,-.5);
\draw[thick, red] (.5,-.5) -- (.62,-.38);
\draw[thick, red, mid<] (-.5,.5) node[above]{$\scriptstyle \alpha$} -- (0,0);
}
=
\tikzmath{
\draw[mid>] (-.5,-.5) node[below]{$\scriptstyle \rho$} -- (0,0);
\draw[mid>] (0,0) -- (.8,.8) node[above]{$\scriptstyle \rho$};
\draw[thick, red, mid>] (-.5,.5) -- (0,0);
\draw[thick, red, mid<] (-.8,.8) node[above]{$\scriptstyle \alpha$} -- (-.5,.5);
\draw[thick, red] (-.5,.5) -- (-.38,.62);
\draw[thick, red, mid>] (.5,-.5) node[below]{$\scriptstyle \alpha$} -- (0,0);
}
\]
due to the implicit inverses on both sides.
Semisimplicity gives us the following local relations.
$$
\tikzmath{
\draw[mid>] (-.3,-.6) node[below]{$\scriptstyle \rho$} -- (-.3,1.2);
\draw[mid>] (.3,-.6) node[below]{$\scriptstyle \rho$} -- (.3,1.2);
}
= 
\frac{1}{d}\tikzmath{
\draw[mid>] (-.3,-1.2) node[below]{$\scriptstyle \rho$} -- (-.3,-.6);
\draw[mid>] (.3,-1.2) node[below]{$\scriptstyle \rho$} -- (.3,-.6) arc (0:180:.3cm);
\draw[thick] (-.3,-.6) -- (-.2,-.6);
\draw[mid<] (-.3,1.2) node[above]{$\scriptstyle \rho$} -- (-.3,.6);
\draw[mid<] (.3,1.2) node[above]{$\scriptstyle \rho$} -- (.3,.6) arc (0:-180:.3cm);
\draw[thick] (-.3,.6) -- (-.2,.6);
}
+ 
\sum_{i=0}^m
\tikzmath{
\draw[mid>] (-.3,-.6) node[below]{$\scriptstyle \rho$} -- (0,0);
\draw[mid>] (.3,-.6) node[below]{$\scriptstyle \rho$} -- (0,0);
\draw[mid<] (0,.6) -- (0,0);
\filldraw (0,0) node[left]{$\scriptstyle i$} circle (.05cm);
\filldraw (0,.6) node[left]{$\scriptstyle i$} circle (.05cm);
\draw[mid<] (.3,1.2) node[above]{$\scriptstyle \rho$} -- (0,0.6);
\draw[mid<] (-.3,1.2) node[above]{$\scriptstyle \rho$} -- (0,0.6);
}
+ 
\sum_{i=0}^m
\tikzmath{
\draw[mid>] (-.3,-.6) node[below]{$\scriptstyle \rho$} -- (0,0);
\draw[mid>] (.3,-.6) node[below]{$\scriptstyle \rho$} -- (0,0);
\draw[mid>] (0,0) to[out=45,in=-45] (0,.6);
\draw[thick, red, mid>] (0,0) to[out=135,in=-135] (0,.6);
\filldraw (0,0) node[left]{$\scriptstyle i$} circle (.05cm);
\filldraw (0,.6) node[left]{$\scriptstyle i$} circle (.05cm);
\draw[mid<] (.3,1.2) node[above]{$\scriptstyle \rho$} -- (0,0.6);
\draw[mid<] (-.3,1.2) node[above]{$\scriptstyle \rho$} -- (0,0.6);
}
$$
\begin{equation}
\label{eq:SD-Spaghetti}
\left.
\begin{aligned}
\tikzmath{
\draw[thick, red, mid>] (0,0) circle (.4cm);
}
&=
\tikzmath{
\draw[thick, red, mid<] (0,0) circle (.4cm);
}
=
1
&&&
\tikzmath{
\draw[mid>] (0,0) circle (.4cm);
}
&=
\tikzmath{
\draw[mid<] (0,0) circle (.4cm);
}
=
d
\\
\tikzmath{
\draw[mid>] (-.3,-.6) node[below]{$\scriptstyle \rho$} -- (-.3,-.3) arc (180:90:.3cm);
\filldraw (0,0) node[left,yshift=.1cm]{$\scriptstyle i$} circle (.05cm);
\draw[mid>] (0,0) -- (0,.3) arc (180:0:.45cm) -- (.9,-.3) arc (0:-180:.3cm) arc (0:90:.3cm);
}
&=
\tikzmath{
\draw[mid>] (.3,-.6) node[below]{$\scriptstyle \rho$} -- (.3,-.3) arc (0:90:.3cm);
\filldraw (0,0) node[left,yshift=.12cm]{$\scriptstyle i$} circle (.05cm);
\draw[mid>] (0,0) -- (0,.3) arc (0:180:.45cm) -- (-.9,-.3) arc (-180:0:.3cm) arc (180:90:.3cm);
}
=0\,\,
&
\tikzmath{
\draw[mid>] (-.3,-.3) arc (180:90:.3cm);
\draw[mid>] (-.3,-.3) arc (-180:90:.3cm);
\draw[thick] (-.3,-.3) -- (-.2,-.3);
\filldraw (0,0) node[left,yshift=.1cm]{$\scriptstyle i$} circle (.05cm);
\draw[mid>] (0,0) -- (0,.4) node[above]{$\scriptstyle \rho$};
}
&=
0\,\,
&
\tikzmath{
\draw[mid>] (-.4,-.4) node[below]{$\scriptstyle \rho$} -- (0,0);
\draw[red, thick, mid>] (0,0) -- (-.4,.4) node[above]{$\scriptstyle \alpha$};
\draw[mid>] (0,0) .. controls ++(45:.1cm) and ++(90:.7cm) ..(.4,0) .. controls ++(-90:.7cm) and ++(-45:.1cm) .. (0,0);
\filldraw (0,0) node[left]{$\scriptstyle i$} circle (.05cm);
}
&=
\tikzmath{
\draw[mid>] (.4,-.4) node[below]{$\scriptstyle \rho$} -- (0,0);
\draw[red, thick] (0,0) -- (-.6,.8) node[above]{$\scriptstyle \alpha$};
\draw[red, thick] (-.4,.52) -- (-.32,.6);
\draw[mid>] (0,0) .. controls ++(45:.7cm) and ++(90:.7cm) ..(-.4,0) .. controls ++(-90:.7cm) and ++(-45:.1cm) ..(0,0);
\filldraw (0,0) node[left]{$\scriptstyle i$} circle (.05cm);
}
=0\,\,
&
\tikzmath{
\draw[mid>] (-.4,-.4) to[out=90,in=-135] (0,0);
\draw[mid>] (-.4,-.4) arc (-180:0:.4cm) to[out=90,in=-45] (0,0);
\draw[thick] (-.4,-.4) -- (-.3,-.4);
\draw[mid>] (0,0) -- (.4,.4) node[above]{$\scriptstyle \rho$};
\draw[red, thick, mid>] (0,0) -- (-.4,.4) node[above]{$\scriptstyle \alpha$};
\filldraw (0,0) node[left]{$\scriptstyle i$} circle (.05cm);
}
&=
0.
\end{aligned}
\,\,
\right\}
\end{equation}

\begin{defn}
Let $\mu\in \bbC^\times$  such that 
$
\tikzmath{
\draw[thick, red, mid>] (-.3,-.5) node[below]{$\scriptstyle \alpha$} -- (-.3,.5);
\draw[mid>] (0,-.5) node[below]{$\scriptstyle \rho$} -- (0,0);
\draw[mid>] (0,.5) node[above]{$\scriptstyle \overline{\rho}$} -- (0,0);
\draw[thick] (0,0) -- (.1,0);
}
= 
\mu \,
\tikzmath{
\draw[thick, red, mid>] (-.3,-1) node[below]{$\scriptstyle \alpha$} .. controls ++(90:.2cm) and ++(-135:.2cm) .. (0,-.5);% -- (-.3,1);
\draw[thick, red, mid<] (-.3,1) node[above]{$\scriptstyle \alpha$} .. controls ++(270:.2cm) and +(135:.2cm) .. (0,.5);% -- (-.3,1);
\draw[thick, red, mid<] (0,-.5)  .. controls ++(45:.2cm) and ++(270:.2cm) .. (.3,0) .. controls ++(90:.2cm) and ++(-45:.2cm) .. (0,.5);
\draw[mid>] (0,-1) node[below]{$\scriptstyle \rho$} -- (0,-.5);
\draw[mid>] (0,1) node[above]{$\scriptstyle \overline{\rho}$} -- (0,.5);
\draw[mid>] (0,-.5) -- (0,0);
\draw[mid<] (0,0) -- (0,.5);
\draw[thick] (0,0) -- (.1,0);
}
$\,.
Clearly $\mu^2=1$.
\end{defn}

In order to choose a natural basis for the spaces $\cC_m(\rho\otimes\rho\to \rho)$ and $\cC_m(\rho\otimes\rho\to \alpha\rho)$, we introduce the following linear operators on these spaces.
We often suppress the orientation on the red $\alpha$ strands, as it may be inferred from the other orientations in the diagram.
\[ 
K^\mathbf{1}\left(\tikzmath{
\draw[mid>] (-.3,-.6) -- (0,0);
\draw[mid>] (.3,-.6) -- (0,0);
\draw[mid<] (0,.6) -- (0,0);
\filldraw (0,0) node[left]{$\scriptstyle i$} circle (.05cm);
}\right)
:= 
\tikzmath{
\draw[thick, red, mid>] (-120:.6cm) arc (-120:-60:.6cm);
\draw[thick, red, mid>] (-120:.6cm) arc (-120:-180:.6cm);
\draw[thick, red, mid>] (90:.6cm) arc (90:180:.6cm);
\draw[thick, red] (-.6,0) -- (-.5,0); 
\draw[thick, red, mid>] (90:.6cm) arc (90:-60:.6cm);
\draw[mid>] (-.4,-.8) -- (0,0);
\draw[mid>] (.4,-.8) -- (0,0);
\draw[mid<] (0,.8) -- (0,0);
\filldraw (0,0) node[left]{$\scriptstyle i$} circle (.05cm);
}
\qquad \text{and} \qquad  
K^\alpha\left(\tikzmath{
\draw[thick, red, mid>] (0,0) -- (-.6,.6);
\draw[mid>] (-.6,-.6) -- (0,0);
\draw[mid>] (.6,-.6) -- (0,0);
\draw[mid>] (0,0) -- (.6,.6);
\filldraw (0,0) node[left]{$\scriptstyle i$} circle (.05cm);
}\right)
:= 
\tikzmath{
\draw[thick, red, mid>] (0,0) .. controls ++(135:.4cm) and ++(90:.4cm) .. (-.6,0);
\draw[thick, red, mid>] (45:.6cm) arc (45:-45:.6cm);
\draw[thick, red, mid<] (-45:.6cm) arc (-45:-135:.6cm);
\draw[thick, red, mid>] (225:.6) .. controls ++(135:.2cm) and ++(270:.3cm) .. (-.6,0);
\draw[thick, red] (-.6,0) -- (-.5,0);
\draw[thick, red, mid>] (45:.6cm) .. controls ++(135:.3cm) and ++(-45:.3cm) .. (-.8,1);
\draw[mid>] (-.8,-.8) -- (-.4,-.4);
\draw[mid>] (-.4,-.4) -- (0,0);
\draw[mid>] (.8,-.8) -- (.4,-.4);
\draw[mid>] (.4,-.4) -- (0,0);
\draw[mid>] (0,0) -- (.4,.4);
\draw[mid>] (.4,.4) -- (1,1);
\filldraw (0,0) node[left]{$\scriptstyle i$} circle (.05cm);
}\,.
\]
We also define the \emph{Frobenius operators}
\begin{align*}
L^\mathbf{1}\left(
\tikzmath{
\draw[mid>] (-.3,-.6) -- (0,0);
\draw[mid>] (.3,-.6) -- (0,0);
\draw[mid<] (0,.6) -- (0,0);
\filldraw (0,0) node[left]{$\scriptstyle i$} circle (.05cm);
}
\right)
&:= 
\tikzmath{
\draw[mid>] (0,0) -- (.3,.6);
\draw[mid>] (0,-.6) -- (0,0);
\draw[mid>] (-.6,-.6) -- (-.6,0);
\draw[mid>] (0,0) arc (0:180:.3cm);
\draw[thick] (-.6,0) -- (-.5,0);
\filldraw (0,0) node[left]{$\scriptstyle i$} circle (.05cm);
}
&
L^\alpha\left(
\tikzmath{
\draw[thick, red, mid>] (0,0) -- (-.6,.6);
\draw[mid>] (-.6,-.6) -- (0,0);
\draw[mid>] (.6,-.6) -- (0,0);
\draw[mid>] (0,0) -- (.6,.6);
\filldraw (0,0) node[left]{$\scriptstyle i$} circle (.05cm);
}
\right)
&:= 
\tikzmath{
\coordinate (a) at (-.45,.25);
\draw[mid>] (0,0) -- (.3,.7);
\draw[mid>] (.2,-.6) -- (0,0);
\draw[mid>] (-.6,-.6) -- (-.6,0);
\draw[mid>] (0,0) arc (0:180:.3cm);
\draw[thick] (-.6,0) -- (-.5,0);
\filldraw (0,0) node[left]{$\scriptstyle i$} circle (.05cm);
\draw[red, thick, mid<] (0,0) .. controls ++(-120:.8cm) and ++(270:.3cm) .. (a);
\draw[red, thick, mid>] (a) -- (-.45,.7);
}
\\ 
R^\mathbf{1}\left(
\tikzmath{
\draw[mid>] (-.3,-.6) -- (0,0);
\draw[mid>] (.3,-.6) -- (0,0);
\draw[mid<] (0,.6) -- (0,0);
\filldraw (0,0) node[left]{$\scriptstyle i$} circle (.05cm);
}
\right)
&:= 
\tikzmath{
\draw[mid>] (0,0) -- (-.3,.6);
\draw[mid>] (0,-.6) -- (0,0);
\draw[mid>] (.6,-.6) -- (.6,0);
\draw[mid>] (0,0) arc (180:0:.3cm);
\draw[thick] (.6,0) -- (.5,0);
\filldraw (0,0) node[left]{$\scriptstyle i$} circle (.05cm);
}
&
R^\alpha\left(
\tikzmath{
\draw[thick, red, mid>] (0,0) -- (-.6,.6);
\draw[mid>] (-.6,-.6) -- (0,0);
\draw[mid>] (.6,-.6) -- (0,0);
\draw[mid>] (0,0) -- (.6,.6);
\filldraw (0,0) node[left]{$\scriptstyle i$} circle (.05cm);
}
\right)
&:= 
\tikzmath{
\draw[mid>] (0,0) -- (-.3,.7);
\draw[mid>] (.3,-.7) -- (0,0);
\draw[mid>] (.8,-.7) -- (.8,0);
\draw[mid>] (0,0) arc (180:0:.4cm);
\draw[thick] (.8,0) -- (.7,0);
\filldraw (0,0) node[left]{$\scriptstyle i$} circle (.05cm);
\draw[red, thick] (0,0) .. controls ++(-120:1cm) and ++(270:.5cm) .. (.4,-.1) .. controls ++(90:.4cm) and ++(270:.4cm) .. (-.6,.7);
}
\end{align*}
These operators satisfy the following relations:
\begin{align*}
K^\mathbf{1} \circ K^\mathbf{1} &= \lambda_\alpha \operatorname{Id} 
&
K^\alpha \circ K^\alpha &= \operatorname{Id}\\
L^\mathbf{1}\circ L^\mathbf{1} 
&= 
\lambda_\rho \operatorname{Id} = R^\mathbf{1}\circ R^\mathbf{1} 
&
L^\alpha\circ L^\alpha &= \lambda_\rho \mu \operatorname{Id} = R^\alpha\circ R^\alpha  
\\
K^\mathbf{1} \circ L^\mathbf{1} &= \mu (L^\mathbf{1} \circ  K^\mathbf{1})
&
K^\alpha \circ L^\alpha &= \mu \lambda_\alpha (L^\alpha \circ  K^\alpha)
\\
K^\mathbf{1} \circ R^\mathbf{1} &= \mu (R^\mathbf{1} \circ  K^\mathbf{1})
&
K^\alpha \circ R^\alpha &= \mu \lambda_\alpha (R^\alpha \circ  K^\alpha)  
\\
(  R^{\mathbf{1}}\circ  L^{\mathbf{1}})^3 &= 1 & (  R^{\alpha}\circ  L^{\alpha})^3 &  = 1
\end{align*}
Note that these relations give a projective representation of $\bbZ/2\bbZ \times D_3$ on our morphism spaces. This allows us to choose a nice basis for these spaces. 

\begin{lem}\label{lem:FrobBasis}
There exists a orthonormal basis of the spaces
$\cC_m(\rho\otimes\rho\to \rho)$ and $\cC_m(\rho\otimes\rho\to \alpha\rho)$
such that 
\begin{align*}
R^\mathbf{1}\left(
\tikzmath{
\draw[mid>] (-.3,-.6) -- (0,0);
\draw[mid>] (.3,-.6) -- (0,0);
\draw[mid<] (0,.6) -- (0,0);
\filldraw (0,0) node[left]{$\scriptstyle i$} circle (.05cm);
}
\right) 
&= 
\lambda_\rho^i 
\tikzmath{
\draw[mid>] (-.3,-.6) -- (0,0);
\draw[mid>] (.3,-.6) -- (0,0);
\draw[mid<] (0,.6) -- (0,0);
\filldraw (0,0) node[left]{$\scriptstyle \widetilde{i}$} circle (.05cm);
}
&
R^\alpha\left(
\tikzmath{
\draw[thick, red, mid>] (0,0) -- (-.6,.6);
\draw[mid>] (-.6,-.6) -- (0,0);
\draw[mid>] (.6,-.6) -- (0,0);
\draw[mid>] (0,0) -- (.6,.6);
\filldraw (0,0) node[left]{$\scriptstyle i$} circle (.05cm);
}
\right) 
&= 
(\lambda_\rho\mu)^i 
\tikzmath{
\draw[thick, red, mid>] (0,0) -- (-.6,.6);
\draw[mid>] (-.6,-.6) -- (0,0);
\draw[mid>] (.6,-.6) -- (0,0);
\draw[mid>] (0,0) -- (.6,.6);
\filldraw (0,0) node[left]{$\scriptstyle \widetilde{i}$} circle (.05cm);
}
\\
L^\mathbf{1}\left(
\tikzmath{
\draw[mid>] (-.3,-.6) -- (0,0);
\draw[mid>] (.3,-.6) -- (0,0);
\draw[mid<] (0,.6) -- (0,0);
\filldraw (0,0) node[left]{$\scriptstyle i$} circle (.05cm);
}
\right ) 
&= 
\lambda_\rho^{i+1}\omega_{\mathbf{1}, i}^{-1} 
\tikzmath{
\draw[mid>] (-.3,-.6) -- (0,0);
\draw[mid>] (.3,-.6) -- (0,0);
\draw[mid<] (0,.6) -- (0,0);
\filldraw (0,0) node[left]{$\scriptstyle \widetilde{i}$} circle (.05cm);
}
&
L^\alpha\left(
\tikzmath{
\draw[thick, red, mid>] (0,0) -- (-.6,.6);
\draw[mid>] (-.6,-.6) -- (0,0);
\draw[mid>] (.6,-.6) -- (0,0);
\draw[mid>] (0,0) -- (.6,.6);
\filldraw (0,0) node[left]{$\scriptstyle i$} circle (.05cm);
}
\right) 
&= 
(\lambda_\rho\mu)^{i+1}\omega_{\alpha, i}^{-1} 
\tikzmath{
\draw[thick, red, mid>] (0,0) -- (-.6,.6);
\draw[mid>] (-.6,-.6) -- (0,0);
\draw[mid>] (.6,-.6) -- (0,0);
\draw[mid>] (0,0) -- (.6,.6);
\filldraw (0,0) node[left]{$\scriptstyle \widetilde{i}$} circle (.05cm);
}
\end{align*}
where $i\mapsto \tilde{i}$ is an order two involution on the indexing set $0 \leq i < m$ and the $\omega_{\alpha, i}$ are 3rd roots of unity indexed by $0 \leq i < m$. 
If $\lambda_\rho = -1$, then $(\mathbf{1}, \tilde{i}) = (\mathbf{1}, 1-i)$, and if $\mu\lambda_\rho = -1$, then $(\alpha, \tilde{i}) = (\alpha, 1-i)$.

We are free to exchange our distinguished basis elements, and to rescale them by \begin{align*}
\tikzmath{
\draw[mid>] (-.3,-.6) -- (0,0);
\draw[mid>] (.3,-.6) -- (0,0);
\draw[mid<] (0,.6) -- (0,0);
\filldraw (0,0) node[left]{$\scriptstyle i$} circle (.05cm);
}&\mapsto z_{\mathbf{1},i} \tikzmath{
\draw[mid>] (-.3,-.6) -- (0,0);
\draw[mid>] (.3,-.6) -- (0,0);
\draw[mid<] (0,.6) -- (0,0);
\filldraw (0,0) node[left]{$\scriptstyle i$} circle (.05cm);
}  \qquad & \qquad    \tikzmath{
\draw[mid>] (-.3,-.6) -- (0,0);
\draw[mid>] (.3,-.6) -- (0,0);
\draw[mid<] (0,.6) -- (0,0);
\filldraw (0,0) node[left]{$\scriptstyle \tilde{i}$} circle (.05cm);
}&\mapsto \overline{ z_{\mathbf{1},i}} \tikzmath{
\draw[mid>] (-.3,-.6) -- (0,0);
\draw[mid>] (.3,-.6) -- (0,0);
\draw[mid<] (0,.6) -- (0,0);
\filldraw (0,0) node[left]{$\scriptstyle \tilde{i}$} circle (.05cm);
}  \qquad z_{\mathbf{1},i} \in U(1)\\
\tikzmath{
\draw[thick, red, mid>] (0,0) -- (-.6,.6);
\draw[mid>] (-.6,-.6) -- (0,0);
\draw[mid>] (.6,-.6) -- (0,0);
\draw[mid>] (0,0) -- (.6,.6);
\filldraw (0,0) node[left]{$\scriptstyle  i$} circle (.05cm);
} &\mapsto z_{\alpha,i} \tikzmath{
\draw[thick, red, mid>] (0,0) -- (-.6,.6);
\draw[mid>] (-.6,-.6) -- (0,0);
\draw[mid>] (.6,-.6) -- (0,0);
\draw[mid>] (0,0) -- (.6,.6);
\filldraw (0,0) node[left]{$\scriptstyle i$} circle (.05cm);
} \qquad & \qquad \tikzmath{
\draw[thick, red, mid>] (0,0) -- (-.6,.6);
\draw[mid>] (-.6,-.6) -- (0,0);
\draw[mid>] (.6,-.6) -- (0,0);
\draw[mid>] (0,0) -- (.6,.6);
\filldraw (0,0) node[left]{$\scriptstyle \widetilde{i}$} circle (.05cm);
}&\mapsto \overline{z_{\alpha,i}}  \tikzmath{
\draw[thick, red, mid>] (0,0) -- (-.6,.6);
\draw[mid>] (-.6,-.6) -- (0,0);
\draw[mid>] (.6,-.6) -- (0,0);
\draw[mid>] (0,0) -- (.6,.6);
\filldraw (0,0) node[left]{$\scriptstyle \widetilde{i}$} circle (.05cm);
}\qquad z_{\alpha,i} \in U(1).
\end{align*}
\end{lem}
\begin{proof}
As the maps $R^\mathbf{1}, L^\mathbf{1},R^\alpha$, and $L^\alpha$ give projective $D_3$ representations on the spaces $\cC_m(\rho\otimes\rho\to\rho)$ and $\cC_m(\rho\otimes\rho\to\alpha\rho)$, we can use \cite[Section VI]{MR1657800} or \cite[Lemma 7.3]{MR3635673} to get orthonormal basis of these spaces which satisfy the first part of this lemma.

By rescaling each of the basis elements as in the statement of the lemma, we preserve all of the relations.
\end{proof}

%%%%%%%%%%%%%%%%%%%%%%%%%%%%%%%%%%%%%%%%%%%%%%%%%%%%%%%%%%%%%%%%
\subsubsection{Jellyfish relations}
\label{sec:SDJellyfishRelations}

With our particular choices of generators we can determine local jellyfish relations in $\cC_m$. 
As a consequence of the nice basis we choose in Lemma~\ref{lem:FrobBasis}, we obtain simple relations for pulling an $\alpha$ strand through our vertices.
\begin{lem}[($\alpha$ Jellyfish)]
\label{lem:alphaJellyfish}
We have the local relations 
$$
\tikzmath{
\draw[thick, red, mid>] (-.6,-.6) node[below]{$\scriptstyle \alpha$} -- (-.6,.6);
\draw[mid>] (-.3,-.6) node[below]{$\scriptstyle \rho$} -- (0,0);
\draw[mid>] (.3,-.6) node[below]{$\scriptstyle \rho$} -- (0,0);
\draw[mid<] (0,.6) node[above]{$\scriptstyle \rho$} -- (0,0);
\filldraw (0,0) node[left]{$\scriptstyle i$} circle (.05cm);
}
=
\chi_{\mathbf{1}, i} 
\tikzmath{
\draw[thick, red, mid<] (-120:.6cm) .. controls ++(30:.2cm) and ++(-150:.3cm) .. (-60:.4cm);
\draw[thick, red, mid<] (-120:.6cm) .. controls ++(-150:.2cm) and ++(90:.3cm) .. (-.7,-.8);
\draw[thick, red, mid<] (90:.4cm) arc (90:-60:.4cm);
\draw[thick, red,mid<]  (90:.4cm) arc (270:180:.4cm);
\draw[thick, red] (-.4,.8) -- (-.3,.8);
\draw[thick, red, mid>] (-.4,.8) -- (-.4,1.2);
\draw[mid>] (-.4,-.8) -- (0,0);
\draw[mid>] (.4,-.8) -- (0,0);
\draw[mid<] (0,1.2) -- (0,0);
\filldraw (0,0) node[left]{$\scriptstyle i$} circle (.05cm);
}
\qquad \text{ and } \qquad  
\tikzmath{
\draw[thick, red, mid>] (-.9,-.6) node[below]{$\scriptstyle \alpha$} -- (-.9,.6);
\draw[thick, red, mid>] (0,0) -- (-.6,.6) node[above]{$\scriptstyle \alpha$};
\draw[mid>] (-.6,-.6) node[below]{$\scriptstyle \rho$} -- (0,0);
\draw[mid>] (.6,-.6) node[below]{$\scriptstyle \rho$} -- (0,0);
\draw[mid>] (0,0) -- (.6,.6) node[above]{$\scriptstyle \rho$};
\filldraw (0,0) node[left]{$\scriptstyle i$} circle (.05cm);
}   
= 
\chi_{\alpha, i} 
\tikzmath{
\draw[thick, red, mid>] (0,0) .. controls ++(135:.2cm) and ++(180:.3cm) .. (0,.6);
\draw[thick, red] (0,.6) arc (90:45:.6cm);
\draw[thick, red, mid<] (45:.6cm) arc (45:-45:.6cm);
\draw[thick, red, mid>] (-45:.6cm) arc (-45:-135:.6cm);
\draw[thick, red, mid<] (225:.6) .. controls ++(135:.3cm) and ++(90:.3cm) .. (-1.1,-.8);
\draw[thick, red, mid>] (-1.1,.8) -- (-1.1,.6) arc (-180:-90:.3cm);
\draw[thick, red] (-.8,.3) -- (-.8,.4);
\draw[thick, red, mid>] (-.5,.8) -- (-.5,.6) arc (0:-90:.3cm);
\draw[mid>] (-.8,-.8) -- (-.4,-.4);
\draw[mid>] (-.4,-.4) -- (0,0);
\draw[mid>] (.8,-.8) -- (.4,-.4);
\draw[mid>] (.4,-.4) -- (0,0);
\draw[mid>] (0,0) -- (.4,.4);
\draw[mid>] (.4,.4) -- (.8,.8);
\filldraw (0,0) node[left]{$\scriptstyle i$} circle (.05cm);
}
$$
where $\chi_{\mathbf{1}, i}^2 = \lambda_\alpha$, and $\chi_{\alpha, i}^2 = 1$. This data satisfies the relations
$\chi_{\mathbf{1}, i} = \lambda_\alpha \mu  \chi_{\mathbf{1}, \widetilde{i}}$
and
$\chi_{\alpha, i} = \lambda_\alpha \mu  \chi_{\alpha, \widetilde{i}}$.
\end{lem}
\begin{proof}
From Lemma~\ref{lem:FrobBasis} the linear operators $L^\mathbf{1}\circ R^\mathbf{1}$ and $L^\alpha\circ R^\alpha$ act diagonally. 
As the linear operators $K^\mathbf{1}$ and $K^\alpha$ commute with these operators, they must be diagonal as well. 
Hence we can choose scalars $\chi_{\mathbf{1}, i}$ and $\chi_{\alpha, i}$ which are the diagonal entries of these diagonal operators. 
As $K^\mathbf{1}\circ K^\mathbf{1} = \lambda_\alpha \operatorname{Id}$ we see that $\chi_{\mathbf{1}, i}^2 = \lambda_\alpha$, and $K^\alpha\circ K^\alpha =  \operatorname{Id}$ implies $\chi_{\alpha, i}^2 = 1$.

From the relation $K^{\mathbf{1}}\circ R^\mathbf{1} = \mu R^\mathbf{1}\circ K^{\mathbf{1}}$ we get that 
$
\chi_{\mathbf{1}, \tilde{i}} 
= 
\mu \overline{\chi_{\mathbf{1}, i} } 
= 
\mu \lambda_\alpha \chi_{\mathbf{1}, i}
$,
and the relation $K^{\alpha}\circ R^\alpha = \mu \lambda_\alpha R^\alpha\circ K^{\alpha}$ gives 
$\chi_{\alpha, \tilde{i}} =  \mu \lambda_\alpha \chi_{\alpha, i}$.
\end{proof}
As a corollary of these local relations we can determine the scalars $\lambda_\alpha$, $\lambda_\rho$, and $\mu $ in the case where $m=1$.

\begin{cor}\label{cor:SDm1}
If $m=1$, then $\lambda_\alpha=   \lambda_\rho = \mu = 1$.
\end{cor}
\begin{proof}
As $m=1$ we must have that $(\mathbf{1}, \tilde{0}) = (\mathbf{1}, 0)$, and $(\alpha, \tilde{0}) = (\alpha, 0)$, which implies that $   \lambda_\rho = \mu = 1$. The equation $\chi_{\alpha, 0}=\chi_{\alpha, \tilde{0}} =  \mu \lambda_\alpha \chi_{\alpha, 0}$ shows that $\lambda_\alpha = 1$.
\end{proof}

\begin{lem}[($\rho$ Jellyfish)]
\label{lem:SDJellyfishRelations}
There exist scalars
\[
A^{i,j}_{k,\ell},B^{i,j}_{k,\ell},C^{i,j}_{k,\ell},D^{i,j}_{k,\ell},
\widehat{A}^{i,j}_{k,\ell},\widehat{B}^{i,j}_{k,\ell},\widehat{C}^{i,j}_{k,\ell},\widehat{D}^{i,j}_{k,\ell} 
\in \mathbb{C} 
\qquad\qquad
0 \leq i,j,k,\ell <m
\]
such that the following local relations hold in $\cC_m$:
\begin{footnotesize}
\begin{align*}
\tikzmath{
\draw[mid>] (-.3,-.5) node[below]{$\scriptstyle \rho$} -- (-.3,.5);
\draw[mid>] (0,-.5) node[below]{$\scriptstyle \rho$} -- (0,0);
\draw[mid>] (0,.5) node[above]{$\scriptstyle \overline{\rho}$} -- (0,0);
\draw[thick] (0,0) -- (.1,0);
}
&=
\frac{\lambda_\rho}{\delta}
\tikzmath{
\draw[mid>] (.2,.5) arc (0:-180:.2cm);
\draw[mid>] (-.2,-.5) -- (-.2,-.1);
\draw[mid>] (.2,-.5) -- (.2,-.1) arc (0:180:.2cm);
\draw[thick] (-.2,-.1) -- (-.1,-.1);
}
+
\sum_i
\lambda_\rho^i
\tikzmath{
\draw[mid>] (-.2,-.5) -- (0,-.1);
\draw[mid>] (.2,-.5) -- (0,-.1);
\draw[mid>] (0,-.1) -- (.4,.5);
\draw[mid>] (.8,1.1) -- (.8,.5) arc (0:-180:.2cm);
\draw[mid>] (.4,.5) -- (.4,1.1);
\filldraw (0,-.1) node[left]{$\scriptstyle i$} circle (.05cm);
\filldraw (.4,.5) node[left]{$\scriptstyle \widetilde{i}$} circle (.05cm);
}
+
\mu(\lambda_\rho\mu)^i
\tikzmath{
\draw[red,thick] (.4,.5) .. controls ++(135:.6cm) and ++(90:.3cm) .. (.7,.5) .. controls ++(270:.3cm) and ++(135:.6cm) .. (0,-.1);
\draw[mid>] (-.2,-.5) -- (0,-.1);
\draw[mid>] (.2,-.5) -- (0,-.1);
\draw[mid>] (0,-.1) -- (.4,.5);
\draw[mid>] (1,1.1) -- (1,.5) arc (0:-180:.3cm);
\draw[far>] (.4,.5) -- (.6,1.1);
\filldraw (0,-.1) node[left]{$\scriptstyle i$} circle (.05cm);
\filldraw (.4,.5) node[left]{$\scriptstyle \widetilde{i}$} circle (.05cm);
}
\\
\tikzmath{
\draw[mid>] (-.6,-.6) node[below]{$\scriptstyle \rho$} -- (-.6,.6);
\draw[mid>] (-.3,-.6) node[below]{$\scriptstyle \rho$} -- (0,0);
\draw[mid>] (.3,-.6) node[below]{$\scriptstyle \rho$} -- (0,0);
\draw[mid<] (0,.6) node[above]{$\scriptstyle \rho$} -- (0,0);
\filldraw (0,0) node[left]{$\scriptstyle \ell$} circle (.05cm);
}
&=
\frac{\lambda_\rho^{\ell+1} \omega_{\mathbf{1},\ell}}{d}\,
\tikzmath{
\draw[mid>] (-.4,-.6) -- (-.4,-.2);
\draw[mid>] (0,-.6) -- (0,-.2);
\draw (-.4,-.2) arc (180:2:.2cm) -- (0,-.6);
\draw[thick] (-.4,-.2) -- (-.3,-.2);
\draw[mid>] (.4,-.6) -- (.4,0);
\draw[mid>] (.4,0) -- (0,.4);
\draw[mid>] (.4,0) -- (.8,.4);
\filldraw (.4,0) node[left, yshift=-.1cm]{$\scriptstyle \widetilde{\ell}$} circle (.05cm);
}
+
\frac{\lambda_\rho}{\omega_{\mathbf{1},\ell}}
\tikzmath{
\draw[mid>] (-.4,-.8) -- (-.4,-.4);
\draw[thick] (-.4,-.4) -- (-.3,-.4);
\draw[mid>] (0,-1.2) -- (0,-.4);
\draw (-.4,-.4) arc (180:0:.2cm);
\draw[mid>] (-.2,-1.2) -- (-.4,-.8);
\draw[mid>] (-.6,-1.2) -- (-.4,-.8);
\draw (-.4,.2) arc (-180:0:.2cm);
\draw[thick] (-.4,.2) -- (-.3,.2);
\draw[mid>] (0,.2) -- (0,.6);
\draw[mid>] (-.4,.2) -- (-.4,.6);
\filldraw (-.4,-.8) node[left]{$\scriptstyle \ell$} circle (.05cm);
}
+
\sum_{i,j,k}
A_{k,\ell}^{i,j}
\tikzmath{
\draw[mid>] (-.4,-.8) -- (-.2,-.2);
\draw[mid>] (.4,-1.2) -- (-.2,-.2);
\draw[mid>] (0,-1.2) -- (-.4,-.8);
\draw[mid>] (-.8,-1.2) -- (-.4,-.8);
\draw[mid<] (-.2,.4)-- (-.2,-.2);
\draw[mid>] (-.2,.4) -- (.2,.8);
\draw[mid>] (-.2,.4) -- (-.6,.8);
\filldraw (-.2,.4) node[left]{$\scriptstyle k$} circle (.05cm);
\filldraw (-.2,-.2) node[left]{$\scriptstyle j$} circle (.05cm);
\filldraw (-.4,-.8) node[left]{$\scriptstyle i$} circle (.05cm);
}
+
B_{k,\ell}^{i,j}
\tikzmath{
\draw[thick, red, mid<] (-.4,.4) .. controls ++(90:.3cm) and ++(90:.4cm) .. (.4,0) .. controls ++(-90:.4cm) and ++(135:.6cm) .. (-.4,-.8);
\draw[thick, red, mid>] (0,0) to[out=135,in=-90] (-.4,.4);
\draw[mid>] (-.4,-.8) -- (0,0);
\draw[mid>] (.4,-1.2) -- (0,0);
\draw[mid>] (0,-1.2) -- (-.4,-.8);
\draw[mid>] (-.8,-1.2) -- (-.4,-.8);
\draw[mid<] (0,.8) to[out=-90,in=45] (0,0);
\draw[mid>] (0,.8) -- (.4,1.2);
\draw[mid>] (0,.8) -- (-.4,1.2);
\filldraw (0,.8) node[left]{$\scriptstyle k$} circle (.05cm);
\filldraw (0,0) node[left]{$\scriptstyle j$} circle (.05cm);
\filldraw (-.4,-.8) node[left]{$\scriptstyle i$} circle (.05cm);
}
+
C_{k,\ell}^{i,j}
\tikzmath{
\draw[thick, red, mid<] (-.2,.4) to[out=-135,in=135] (-.2,-.2);
\draw[mid>] (-.4,-.8) -- (-.2,-.2);
\draw[mid>] (.4,-1.2) -- (-.2,-.2);
\draw[mid>] (0,-1.2) -- (-.4,-.8);
\draw[mid>] (-.8,-1.2) -- (-.4,-.8);
\draw[mid<] (-.2,.4) to[out=-45,in=45] (-.2,-.2);
\draw[mid>] (-.2,.4) -- (.2,.8);
\draw[mid>] (-.2,.4) -- (-.6,.8);
\filldraw (-.2,.4) node[left]{$\scriptstyle k$} circle (.05cm);
\filldraw (-.2,-.2) node[left]{$\scriptstyle j$} circle (.05cm);
\filldraw (-.4,-.8) node[left]{$\scriptstyle i$} circle (.05cm);
}
+
D_{k,\ell}^{i,j}
\tikzmath{
\draw[thick, red, mid<] (0,.5) .. controls ++(-135:.6cm) and ++(90:.4cm) .. (.4,0) .. controls ++(-90:.4cm) and ++(135:.6cm) .. (-.4,-.8);
\draw[thick, red] (-.12,.34) -- (-.02,.34);
\draw[mid>] (-.4,-.8) -- (0,0);
\draw[mid>] (.4,-1.2) -- (0,0);
\draw[mid>] (0,-1.2) -- (-.4,-.8);
\draw[mid>] (-.8,-1.2) -- (-.4,-.8);
\draw[mid<] (0,.5) to[out=-45,in=45] (0,0);
\draw[mid>] (0,.5) -- (.4,.9);
\draw[mid>] (0,.5) -- (-.4,.9);
\filldraw (0,.5) node[left]{$\scriptstyle k$} circle (.05cm);
\filldraw (0,0) node[left]{$\scriptstyle j$} circle (.05cm);
\filldraw (-.4,-.8) node[left]{$\scriptstyle i$} circle (.05cm);
}
\\
\tikzmath{
\draw[thick, red, mid>] (-.9,-.6) node[below]{$\scriptstyle \alpha$} -- (-.9,.6);
\draw[mid>] (-.6,-.6) node[below]{$\scriptstyle \rho$} -- (-.6,.6);
\draw[mid>] (-.3,-.6) node[below]{$\scriptstyle \rho$} -- (0,0);
\draw[mid>] (.3,-.6) node[below]{$\scriptstyle \rho$} -- (0,0);
\draw[thick, red, mid<] (-.3,.6) node[above]{$\scriptstyle \alpha$} -- (0,0);
\draw[mid<] (.3,.6) node[above]{$\scriptstyle \rho$} -- (0,0);
\filldraw (0,0) node[left]{$\scriptstyle \ell$} circle (.05cm);
}
&=
\frac{(\lambda_\rho \mu)^{\ell+1}\omega_{\alpha, \ell}}{d}
\tikzmath{
\draw[mid>] (-.4,-.6) -- (-.4,-.2);
\draw[mid>] (0,-.6) -- (0,-.2);
\draw (-.4,-.2) arc (180:2:.2cm);
\draw[thick] (-.4,-.2) -- (-.3,-.2);
\draw[thick, red] (.4,0) .. controls ++(-135:.6cm) and ++(90:.3cm) .. (-.6,-.6);
\draw[mid>] (.8,-.6) -- (.4,0);
\draw[mid>] (.4,0) -- (0,.4);
\draw[mid>] (.4,0) -- (.8,.4);
\draw[thick, red] (-.2,.4) arc (-180:0:.2cm);
\filldraw (.4,0) node[left, xshift=-.1cm]{$\scriptstyle \widetilde{\ell}$} circle (.05cm);
}
+
\frac{\lambda_\rho}{d \omega_{\alpha,\ell}}
\tikzmath{
\draw[thick, red] (-.4,-.8) .. controls ++(135:.4cm) and ++(90:.3cm) .. (-.2,-.8) .. controls ++(-90:.3cm) and ++(90:.3cm) .. (-.8,-1.2);
\draw[thick, red] (-.6,.6) arc (-180:0:.2cm);
\draw[mid>] (-.4,-.8) to[out=45,in=-90] (-.2,-.4);
\draw[thick] (-.2,-.4) -- (-.1,-.4);
\draw[mid>] (.2,-1.2) -- (.2,-.4);
\draw (-.2,-.4) arc (180:0:.2cm);
\draw[mid>] (-.2,-1.2) -- (-.4,-.8);
\draw[mid>] (-.6,-1.2) -- (-.4,-.8);
\draw (-.4,.2) arc (-180:0:.2cm);
\draw[thick] (-.4,.2) -- (-.3,.2);
\draw[mid>] (0,.2) -- (0,.6);
\draw[mid>] (-.4,.2) -- (-.4,.6);
\filldraw (-.4,-.8) node[left]{$\scriptstyle \ell$} circle (.05cm);
}
+
\sum_{i,j,k}
\widehat{A}_{k,\ell}^{i,j}
\tikzmath{
\draw[thick, red] (-.4,-.8) .. controls ++(135:.5cm) and ++(90:.4cm) .. (-.2,-.8) .. controls ++(-90:.3cm) and ++(90:.3cm) .. (-.8,-1.2);
\draw[thick, red, mid<] (-.2,.4) to[out=-135,in=135] (-.2,-.2);
\draw[thick, red] (-.8,.8) arc (-180:0:.2cm);
\draw[mid>] (-.4,-.8) to[out=45,in=-135] (-.2,-.2);
\draw[mid>] (.4,-1.2) -- (-.2,-.2);
\draw[mid>] (-.2,-1.2) -- (-.4,-.8);
\draw[mid>] (-.6,-1.2) -- (-.4,-.8);
\draw[mid<] (-.2,.4) to[out=-45,in=45] (-.2,-.2);
\draw[mid>] (-.2,.4) -- (.2,.8);
\draw[mid>] (-.2,.4) -- (-.6,.8);
\filldraw (-.2,.4) node[left]{$\scriptstyle k$} circle (.05cm);
\filldraw (-.2,-.2) node[left]{$\scriptstyle j$} circle (.05cm);
\filldraw (-.4,-.8) node[left]{$\scriptstyle i$} circle (.05cm);
}
+
\widehat{B}_{k,\ell}^{i,j}
\tikzmath{
\draw[thick, red] (0,.5) .. controls ++(-135:.6cm) and ++(90:.6cm) .. (.4,-.4) .. controls ++(-90:.6cm) and ++(90:.3cm) .. (-.8,-1.2);
\draw[thick, red] (-.6,.9) arc (-180:0:.2cm);
\draw[mid>] (-.4,-.8) -- (0,0);
\draw[mid>] (.4,-1.2) -- (0,0);
\draw[mid>] (-.2,-1.2) -- (-.4,-.8);
\draw[mid>] (-.6,-1.2) -- (-.4,-.8);
\draw[mid<] (0,.5) to[out=-45,in=90] (0,0);
\draw[mid>] (0,.5) -- (.4,.9);
\draw[mid>] (0,.5) -- (-.4,.9);
\filldraw (0,.5) node[left]{$\scriptstyle k$} circle (.05cm);
\filldraw (0,0) node[left]{$\scriptstyle j$} circle (.05cm);
\filldraw (-.4,-.8) node[left]{$\scriptstyle i$} circle (.05cm);
}
+
\widehat{C}_{k,\ell}^{i,j}
\tikzmath{
\draw[thick, red] (-.4,-.8) .. controls ++(135:.5cm) and ++(90:.4cm) .. (-.2,-.8) .. controls ++(-90:.3cm) and ++(90:.3cm) .. (-.8,-1.2);
\draw[thick, red] (-.8,.8) arc (-180:0:.2cm);
\draw[mid>] (-.4,-.8) to[out=45,in=-135] (-.2,-.2);
\draw[mid>] (.4,-1.2) -- (-.2,-.2);
\draw[mid>] (-.2,-1.2) -- (-.4,-.8);
\draw[mid>] (-.6,-1.2) -- (-.4,-.8);
\draw[mid<] (-.2,.4) -- (-.2,-.2);
\draw[mid>] (-.2,.4) -- (.2,.8);
\draw[mid>] (-.2,.4) -- (-.6,.8);
\filldraw (-.2,.4) node[left]{$\scriptstyle k$} circle (.05cm);
\filldraw (-.2,-.2) node[left]{$\scriptstyle j$} circle (.05cm);
\filldraw (-.4,-.8) node[left]{$\scriptstyle i$} circle (.05cm);
}
+
\widehat{D}_{k,\ell}^{i,j}
\tikzmath{
\draw[thick, red] (-.2,.-.2) .. controls ++(135:.5cm) and ++(90:.4cm) ..(.2,-.2) .. controls ++(-90:.8cm) and ++(90:.3cm) .. (-.8,-1.2) -- (-.8,-1.4);
\draw[thick, red] (-.8,.8) arc (-180:0:.2cm);
\draw[thick, red] (-.8,-1.25) -- (-.7,-1.25); 
\draw[mid>] (-.4,-.8) -- (-.2,-.2);
\draw[mid>] (.4,-1.4) -- (-.2,-.2);
\draw[mid>] (-.2,-1.4) -- (-.4,-.8);
\draw[mid>] (-.6,-1.4) -- (-.4,-.8);
\draw[mid<] (-.2,.4) to[out=-90,in=45] (-.2,-.2);
\draw[mid>] (-.2,.4) -- (.2,.8);
\draw[mid>] (-.2,.4) -- (-.6,.8);
\filldraw (-.2,.4) node[left]{$\scriptstyle k$} circle (.05cm);
\filldraw (-.2,-.2) node[left]{$\scriptstyle j$} circle (.05cm);
\filldraw (-.4,-.8) node[left]{$\scriptstyle i$} circle (.05cm);
}
\end{align*}
\end{footnotesize}
\end{lem}
\begin{proof}
We provide the proof for the third relation, and the other two are left to the reader.
We compute
\begin{footnotesize}
\begin{align*}
\tikzmath{
\draw[thick, red, mid>] (-.9,-.6) node[below]{$\scriptstyle \alpha$} -- (-.9,.6);
\draw[mid>] (-.6,-.6) node[below]{$\scriptstyle \rho$} -- (-.6,.6);
\draw[mid>] (-.3,-.6) node[below]{$\scriptstyle \rho$} -- (0,0);
\draw[mid>] (.3,-.6) node[below]{$\scriptstyle \rho$} -- (0,0);
\draw[thick, red, mid<] (-.3,.6) node[above]{$\scriptstyle \alpha$} -- (0,0);
\draw[mid<] (.3,.6) node[above]{$\scriptstyle \rho$} -- (0,0);
\filldraw (0,0) node[left]{$\scriptstyle \ell$} circle (.05cm);
\draw[dotted, thick, blue, rounded corners=5pt] (-.75,-.5) rectangle (0,-.2);
}
&=
\frac{1}{d}\,
\tikzmath{
\draw[thick, red, mid>] (-.9,-1) -- (-.9,.6);
\draw[mid>] (-.6,0) -- (-.6,.6);
\draw[mid>] (-.6,0) .. controls ++(-90:.4cm) and ++(-135:.3cm) .. (0,0);
\draw[mid>] (-.6,-1) -- (-.6,-.6);
\draw[thick] (-.6,-.6) -- (-.5,-.6);
\draw[mid>] (-.2,-1) -- (-.2,-.6);
\draw (-.6,-.6) arc (180:0:.2cm);
\draw[thick] (-.6,0) -- (-.5,0);
\draw[mid>] (.3,-1) -- (0,0);
\draw[thick, red, mid<] (-.3,.6) -- (0,0);
\draw[mid<] (.3,.6) -- (0,0);
\filldraw (0,0) node[left]{$\scriptstyle \ell$} circle (.05cm);
}
+
\sum_i
\tikzmath{
\draw[thick, red, mid>] (-.9,-1) -- (-.9,.6);
\draw[mid>] (-.3,-.3) to[out=135,in=-90] (-.6,.6);
\draw[mid>] (-.3,-.7) -- (-.3,-.3);
\draw[mid>] (-.3,-.3) -- (0,0);
\draw[mid>] (-.6,-1) -- (-.3,-.7);
\draw[mid>] (0,-1) -- (-.3,-.7);
\draw[mid>] (.3,-1) -- (0,0);
\draw[thick, red, mid<] (-.3,.6) -- (0,0);
\draw[mid<] (.3,.6) -- (0,0);
\filldraw (0,0) node[left]{$\scriptstyle \ell$} circle (.05cm);
\filldraw (-.3,-.3) node[left]{$\scriptstyle i$} circle (.05cm);
\filldraw (-.3,-.7) node[left]{$\scriptstyle i$} circle (.05cm);
}
+
\tikzmath{
\draw[thick, red, mid>] (-.9,-1) -- (-.9,.6);
\draw[mid>] (-.3,-.3) to[out=135,in=-90] (-.6,.6);
\draw[mid>] (-.3,-.7) to[out=45,in=-45] (-.3,-.3);
\draw[thick, red] (-.3,-.7) to[out=135,in=-135] (-.3,-.3);
\draw[mid>] (-.3,-.3) -- (0,0);
\draw[mid>] (-.6,-1) -- (-.3,-.7);
\draw[mid>] (0,-1) -- (-.3,-.7);
\draw[mid>] (.3,-1) -- (0,0);
\draw[thick, red, mid<] (-.3,.6) -- (0,0);
\draw[mid<] (.3,.6) -- (0,0);
\filldraw (0,0) node[left]{$\scriptstyle \ell$} circle (.05cm);
\filldraw (-.3,-.3) node[left]{$\scriptstyle i$} circle (.05cm);
\filldraw (-.3,-.7) node[left]{$\scriptstyle i$} circle (.05cm);
}
\\&=
\frac{\mu}{d}\,
\tikzmath{
\draw[thick, red, mid>] (-1.1,-1) -- (-1.1,.6);
\draw[mid>] (-.6,0) -- (-.6,.6);
\draw[mid>] (-.6,0) .. controls ++(-90:.4cm) and ++(-135:.3cm) .. (0,0);
\draw[mid>] (-.6,-1) -- (-.6,-.6);
\draw[thick] (-.6,-.6) -- (-.5,-.6);
\draw[mid>] (-.2,-1) -- (-.2,-.6);
\draw (-.6,-.6) arc (180:0:.2cm);
\draw[thick] (-.6,0) -- (-.5,0);
\draw[mid>] (.3,-1) -- (0,0);
\draw[thick, red] (0,0) .. controls ++(135:.3cm) and ++(90:.2cm) .. (-.4,0)
arc (0:-180:.25cm) .. controls ++(90:.4cm) and ++(-90:.3cm) .. (-.4,.6)
;
\draw[mid<] (.3,.6) -- (0,0);
\filldraw (0,0) node[left]{$\scriptstyle \ell$} circle (.05cm);
}
+
\sum_i
\tikzmath{
\draw[thick, red] (-.8,.8) arc (-180:0:.2cm);
\draw[mid>] (-.3,-.3) to[out=135,in=-90] (-.6,.8);
\draw[mid>] (-.3,-.7) -- (-.3,-.3);
\draw[mid>] (-.3,-.3) -- (0,0);
\draw[mid>] (-.6,-1) -- (-.3,-.7);
\draw[mid>] (0,-1) -- (-.3,-.7);
\draw[mid>] (.3,-1) -- (0,0);
\draw[mid>] (0,0) -- (.3,.8);
\draw[thick, red] (0,0) .. controls ++(135:.7cm) and ++(90:.3cm) .. (-.9,-1);
\filldraw (0,0) node[left]{$\scriptstyle \ell$} circle (.05cm);
\filldraw (-.3,-.3) node[left]{$\scriptstyle i$} circle (.05cm);
\filldraw (-.3,-.7) node[left]{$\scriptstyle i$} circle (.05cm);
\draw[thick, dotted, blue, rounded corners=5pt] (-.8,.2) rectangle (.4,.5);
}
+
\tikzmath{
\draw[thick, red] (-.8,.8) arc (-180:0:.2cm);
\draw[mid>] (-.3,-.3) to[out=135,in=-90] (-.6,.8);
\draw[mid>] (-.3,-.7) to[out=45,in=-45] (-.3,-.3);
\draw[thick, red] (-.3,-.7) to[out=135,in=-135] (-.3,-.3);
\draw[mid>] (-.3,-.3) -- (0,0);
\draw[mid>] (-.6,-1) -- (-.3,-.7);
\draw[mid>] (0,-1) -- (-.3,-.7);
\draw[mid>] (.3,-1) -- (0,0);
\draw[mid>] (0,0) -- (.3,.8);
\draw[thick, red] (0,0) .. controls ++(135:.7cm) and ++(90:.3cm) .. (-.9,-1);
\filldraw (0,0) node[left]{$\scriptstyle \ell$} circle (.05cm);
\filldraw (-.3,-.3) node[left]{$\scriptstyle i$} circle (.05cm);
\filldraw (-.3,-.7) node[left]{$\scriptstyle i$} circle (.05cm);
\draw[thick, dotted, blue, rounded corners=5pt] (-.8,.2) rectangle (.4,.5);
}
\\
&=
\frac{1}{d}\,
\tikzmath{
\draw[thick, red] (-.8,.6) arc (-180:0:.2cm);
\draw[mid>] (-.6,0) -- (-.6,.6);
\draw[mid>] (-.6,0) .. controls ++(-90:.4cm) and ++(-135:.3cm) .. (0,0);
\draw[mid>] (-.6,-1) -- (-.6,-.6);
\draw[thick] (-.6,-.6) -- (-.5,-.6);
\draw[mid>] (-.2,-1) -- (-.2,-.6);
\draw (-.6,-.6) arc (180:0:.2cm);
\draw[thick] (-.6,0) -- (-.5,0);
\draw[mid>] (.3,-1) -- (0,0);
\draw[thick, red] (0,0) .. controls ++(135:.3cm) and ++(90:.2cm) .. (-.4,0) -- (-.4,-.5)
.. controls ++(-90:.3cm) and ++(90:.3cm) .. (-.9,-1)
;
\draw[mid<] (.3,.6) -- (0,0);
\filldraw (0,0) node[left]{$\scriptstyle \ell$} circle (.05cm);
\draw[thick, dotted, blue, rounded corners=5pt] (-.8,-.3) rectangle (.4,.3) node[right, yshift=-.3cm]{$\scriptstyle L^\alpha(\ell)^*$};
}
+
\sum_i
\frac{1}{d}
\underbrace{
\tikzmath{
\draw[thick, red] (-.8,1.2) arc (-180:0:.2cm);
\draw[mid>] (-.6,.8) -- (-.6,1.2);
\draw[mid>] (0,.8) -- (0,1.2);
\draw[thick] (-.6,.8) -- (-.5,.8);
\draw (-.6,.8) arc (-180:0:.3cm);
\draw[mid>] (-.3,-.7) -- (-.3,-.3);
\draw[mid>] (-.3,-.3) -- (0,0);
\draw[mid>] (-.6,-1) -- (-.3,-.7);
\draw[mid>] (0,-1) -- (-.3,-.7);
\draw[mid>] (.3,-1) -- (0,0);
\draw[mid>] (-.3,-.3) to[out=135, in=-90] (-.5,.2);
\draw[thick] (-.5,.2) -- (-.4,.2);
\draw[mid>] (0,0) .. controls ++(45:.5cm) and ++(90:.3cm) .. (-.5,.2);
\draw[thick, red] (0,0) .. controls ++(135:.7cm) and ++(90:.3cm) .. (-.9,-1);
\filldraw (0,0) node[left]{$\scriptstyle \ell$} circle (.05cm);
\filldraw (-.3,-.3) node[left]{$\scriptstyle i$} circle (.05cm);
\filldraw (-.3,-.7) node[left]{$\scriptstyle i$} circle (.05cm);
}
}_{=0}
+
\frac{1}{d}
\tikzmath{
\draw[thick, red] (-.8,1.2) arc (-180:0:.2cm);
\draw[mid>] (-.6,.8) -- (-.6,1.2);
\draw[mid>] (0,.8) -- (0,1.2);
\draw[thick] (-.6,.8) -- (-.5,.8);
\draw (-.6,.8) arc (-180:0:.3cm);
\draw[mid>] (-.3,-.7) to[out=45,in=-45] (-.3,-.3);
\draw[thick, red] (-.3,-.7) to[out=135,in=-135] (-.3,-.3);
\draw[mid>] (-.3,-.3) -- (0,0);
\draw[mid>] (-.6,-1) -- (-.3,-.7);
\draw[mid>] (0,-1) -- (-.3,-.7);
\draw[mid>] (.3,-1) -- (0,0);
\draw[mid>] (-.3,-.3) to[out=135, in=-90] (-.5,.2);
\draw[thick] (-.5,.2) -- (-.4,.2);
\draw[mid>] (0,0) .. controls ++(45:.5cm) and ++(90:.3cm) .. (-.5,.2);
\draw[thick, red] (0,0) .. controls ++(135:.7cm) and ++(90:.3cm) .. (-.9,-1);
\filldraw (0,0) node[left]{$\scriptstyle \ell$} circle (.05cm);
\filldraw (-.3,-.3) node[left]{$\scriptstyle i$} circle (.05cm);
\filldraw (-.3,-.7) node[left]{$\scriptstyle i$} circle (.05cm);
}
+
\sum_{i,k}
\tikzmath{
\draw[thick, red] (-.8,1.2) arc (-180:0:.2cm);
\draw[mid>] (-.3,.7) -- (-.6,1.2);
\draw[mid>] (-.3,.7) -- (0,1.2);
\draw[mid>] (-.3,-.7) -- (-.3,-.3);
\draw[mid>] (-.3,-.3) -- (0,0);
\draw[mid>] (-.6,-1) -- (-.3,-.7);
\draw[mid>] (0,-1) -- (-.3,-.7);
\draw[mid>] (.3,-1) -- (0,0);
\draw[mid>] (-.3,-.3) to[out=135, in=-135] (-.3,.3);
\draw[mid>] (0,0) to[out=45,in=-45] (-.3,.3);
\draw[mid>] (-.3,.3) -- (-.3,.7);
\draw[thick, red] (0,0) .. controls ++(135:.7cm) and ++(90:.3cm) .. (-.9,-1);
\filldraw (0,0) node[left]{$\scriptstyle \ell$} circle (.05cm);
\filldraw (-.3,-.3) node[left]{$\scriptstyle i$} circle (.05cm);
\filldraw (-.3,-.7) node[left]{$\scriptstyle i$} circle (.05cm);
\filldraw (-.3,.3) node[left]{$\scriptstyle k$} circle (.05cm);
\filldraw (-.3,.7) node[left]{$\scriptstyle k$} circle (.05cm);
}
+
\tikzmath{
\draw[thick, red] (-.8,1.2) arc (-180:0:.2cm);
\draw[mid>] (-.3,.7) -- (-.6,1.2);
\draw[mid>] (-.3,.7) -- (0,1.2);
\draw[mid>] (-.3,-.7) -- (-.3,-.3);
\draw[mid>] (-.3,-.3) -- (0,0);
\draw[mid>] (-.6,-1) -- (-.3,-.7);
\draw[mid>] (0,-1) -- (-.3,-.7);
\draw[mid>] (.3,-1) -- (0,0);
\draw[mid>] (-.3,-.3) to[out=135, in=-135] (-.3,.3);
\draw[mid>] (0,0) to[out=45,in=-45] (-.3,.3);
\draw[mid>] (-.3,.3) to[out=45, in=-45] (-.3,.7);
\draw[thick, red] (-.3,.3) to[out=135, in=-135] (-.3,.7);
\draw[thick, red] (0,0) .. controls ++(135:.7cm) and ++(90:.3cm) .. (-.9,-1);
\filldraw (0,0) node[left]{$\scriptstyle \ell$} circle (.05cm);
\filldraw (-.3,-.3) node[left]{$\scriptstyle i$} circle (.05cm);
\filldraw (-.3,-.7) node[left]{$\scriptstyle i$} circle (.05cm);
\filldraw (-.3,.3) node[left]{$\scriptstyle k$} circle (.05cm);
\filldraw (-.3,.7) node[left]{$\scriptstyle k$} circle (.05cm);
}
+
\tikzmath{
\draw[thick, red] (-.8,1.2) arc (-180:0:.2cm);
\draw[mid>] (-.3,.7) -- (-.6,1.2);
\draw[mid>] (-.3,.7) -- (0,1.2);
\draw[mid>] (-.3,-.7) to[out=45,in=-45] (-.3,-.3);
\draw[thick, red] (-.3,-.7) to[out=135,in=-135] (-.3,-.3);
\draw[mid>] (-.3,-.3) -- (0,0);
\draw[mid>] (-.6,-1) -- (-.3,-.7);
\draw[mid>] (0,-1) -- (-.3,-.7);
\draw[mid>] (.3,-1) -- (0,0);
\draw[mid>] (-.3,-.3) to[out=135, in=-135] (-.3,.3);
\draw[mid>] (0,0) to[out=45,in=-45] (-.3,.3);
\draw[mid>] (-.3,.3) -- (-.3,.7);
\draw[thick, red] (0,0) .. controls ++(135:.7cm) and ++(90:.3cm) .. (-.9,-1);
\filldraw (0,0) node[left]{$\scriptstyle \ell$} circle (.05cm);
\filldraw (-.3,-.3) node[left]{$\scriptstyle i$} circle (.05cm);
\filldraw (-.3,-.7) node[left]{$\scriptstyle i$} circle (.05cm);
\filldraw (-.3,.3) node[left]{$\scriptstyle k$} circle (.05cm);
\filldraw (-.3,.7) node[left]{$\scriptstyle k$} circle (.05cm);
}
+
\tikzmath{
\draw[thick, red] (-.8,1.2) arc (-180:0:.2cm);
\draw[mid>] (-.3,.7) -- (-.6,1.2);
\draw[mid>] (-.3,.7) -- (0,1.2);
\draw[mid>] (-.3,-.7) to[out=45,in=-45] (-.3,-.3);
\draw[thick, red] (-.3,-.7) to[out=135,in=-135] (-.3,-.3);
\draw[mid>] (-.3,-.3) -- (0,0);
\draw[mid>] (-.6,-1) -- (-.3,-.7);
\draw[mid>] (0,-1) -- (-.3,-.7);
\draw[mid>] (.3,-1) -- (0,0);
\draw[mid>] (-.3,-.3) to[out=135, in=-135] (-.3,.3);
\draw[mid>] (0,0) to[out=45,in=-45] (-.3,.3);
\draw[mid>] (-.3,.3) to[out=45, in=-45] (-.3,.7);
\draw[thick, red] (-.3,.3) to[out=135, in=-135] (-.3,.7);
\draw[thick, red] (0,0) .. controls ++(135:.7cm) and ++(90:.3cm) .. (-.9,-1);
\filldraw (0,0) node[left]{$\scriptstyle \ell$} circle (.05cm);
\filldraw (-.3,-.3) node[left]{$\scriptstyle i$} circle (.05cm);
\filldraw (-.3,-.7) node[left]{$\scriptstyle i$} circle (.05cm);
\filldraw (-.3,.3) node[left]{$\scriptstyle k$} circle (.05cm);
\filldraw (-.3,.7) node[left]{$\scriptstyle k$} circle (.05cm);
}
\\
&=
\frac{(\lambda_\rho \mu)^{\ell+1}\omega_{\alpha, \ell}}{d}
\tikzmath{
\draw[mid>] (-.4,-.6) -- (-.4,-.2);
\draw[mid>] (0,-.6) -- (0,-.2);
\draw (-.4,-.2) arc (180:2:.2cm);
\draw[thick] (-.4,-.2) -- (-.3,-.2);
\draw[thick, red] (.4,0) .. controls ++(-135:.6cm) and ++(90:.3cm) .. (-.6,-.6);
\draw[mid>] (.8,-.6) -- (.4,0);
\draw[mid>] (.4,0) -- (0,.4);
\draw[mid>] (.4,0) -- (.8,.4);
\draw[thick, red] (-.2,.4) arc (-180:0:.2cm);
\filldraw (.4,0) node[left, xshift=-.1cm]{$\scriptstyle \widetilde{\ell}$} circle (.05cm);
}
+
\frac{\lambda_\rho}{d}
\sum_i
\tikzmath{
\draw[thick, red] (-.8,1.2) arc (-180:0:.2cm);
\draw[mid>] (-.6,.8) -- (-.6,1.2);
\draw[mid>] (0,.8) -- (0,1.2);
\draw[thick] (-.6,.8) -- (-.5,.8);
\draw (-.6,.8) arc (-180:0:.3cm);
\draw[mid>] (-.3,-1) to[out=45,in=-45] (-.3,-.7);
\draw[thick, red] (-.3,-1) to[out=135,in=-135] (-.3,-.7);
\draw[mid>] (-.3,-.7) -- (0,0);
\draw[mid>] (-.6,-1.3) -- (-.3,-1);
\draw[mid>] (0,-1.3) -- (-.3,-1);
\draw[mid>] (.4,0) .. controls ++(-90:.4cm) and ++(-45:.2cm) .. (0,0);
\draw[thick] (.4,0) -- (.3,0);
\draw[mid>] (.4,0) arc (180:0:.2cm);
\draw[thick] (.8,0) -- (.7,0);
\draw[mid>] (.8,-1.3) -- (.8,0);
\draw[mid>] (-.3,-.7) to[out=135, in=-90] (-.5,.2);
\draw[thick] (-.5,.2) -- (-.4,.2);
\draw[mid>] (0,0) .. controls ++(45:.5cm) and ++(90:.3cm) .. (-.5,.2);
\draw[thick, red] (0,0) .. controls ++(135:.2cm) and ++(90:.2cm) .. (-.3,0) .. controls ++(-90:.2cm) and ++(90:.2cm) .. (0,-.3) .. controls ++(-90:.3cm) and ++(90:.9cm) .. (-.9,-1.3);
\filldraw (0,0) node[left]{$\scriptstyle \ell$} circle (.05cm);
\filldraw (-.3,-.7) node[left]{$\scriptstyle i$} circle (.05cm);
\filldraw (-.3,-1) node[left]{$\scriptstyle i$} circle (.05cm);
\draw[thick, dotted, blue, rounded corners=5pt] (-.8,-.4) rectangle (.6,.45) node[above, xshift=-.2cm]{\scriptsize\eqref{eq:RotatedFrobeniusComposite}} ;
}
\\&\hspace{3cm}
+
\sum_{i,j,k}
\widehat{A}_{k,\ell}^{i,j}
\tikzmath{
\draw[thick, red] (-.4,-.8) .. controls ++(135:.5cm) and ++(90:.4cm) .. (-.2,-.8) .. controls ++(-90:.3cm) and ++(90:.3cm) .. (-.8,-1.2);
\draw[thick, red, mid<] (-.2,.4) to[out=-135,in=135] (-.2,-.2);
\draw[thick, red] (-.8,.8) arc (-180:0:.2cm);
\draw[mid>] (-.4,-.8) to[out=45,in=-135] (-.2,-.2);
\draw[mid>] (.4,-1.2) -- (-.2,-.2);
\draw[mid>] (-.2,-1.2) -- (-.4,-.8);
\draw[mid>] (-.6,-1.2) -- (-.4,-.8);
\draw[mid<] (-.2,.4) to[out=-45,in=45] (-.2,-.2);
\draw[mid>] (-.2,.4) -- (.2,.8);
\draw[mid>] (-.2,.4) -- (-.6,.8);
\filldraw (-.2,.4) node[left]{$\scriptstyle k$} circle (.05cm);
\filldraw (-.2,-.2) node[left]{$\scriptstyle j$} circle (.05cm);
\filldraw (-.4,-.8) node[left]{$\scriptstyle i$} circle (.05cm);
}
+
\widehat{B}_{k,\ell}^{i,j}
\tikzmath{
\draw[thick, red] (0,.5) .. controls ++(-135:.6cm) and ++(90:.6cm) .. (.4,-.4) .. controls ++(-90:.6cm) and ++(90:.3cm) .. (-.8,-1.2);
\draw[thick, red] (-.6,.9) arc (-180:0:.2cm);
\draw[mid>] (-.4,-.8) -- (0,0);
\draw[mid>] (.4,-1.2) -- (0,0);
\draw[mid>] (-.2,-1.2) -- (-.4,-.8);
\draw[mid>] (-.6,-1.2) -- (-.4,-.8);
\draw[mid<] (0,.5) to[out=-45,in=90] (0,0);
\draw[mid>] (0,.5) -- (.4,.9);
\draw[mid>] (0,.5) -- (-.4,.9);
\filldraw (0,.5) node[left]{$\scriptstyle k$} circle (.05cm);
\filldraw (0,0) node[left]{$\scriptstyle j$} circle (.05cm);
\filldraw (-.4,-.8) node[left]{$\scriptstyle i$} circle (.05cm);
}
+
\widehat{C}_{k,\ell}^{i,j}
\tikzmath{
\draw[thick, red] (-.4,-.8) .. controls ++(135:.5cm) and ++(90:.4cm) .. (-.2,-.8) .. controls ++(-90:.3cm) and ++(90:.3cm) .. (-.8,-1.2);
\draw[thick, red] (-.8,.8) arc (-180:0:.2cm);
\draw[mid>] (-.4,-.8) to[out=45,in=-135] (-.2,-.2);
\draw[mid>] (.4,-1.2) -- (-.2,-.2);
\draw[mid>] (-.2,-1.2) -- (-.4,-.8);
\draw[mid>] (-.6,-1.2) -- (-.4,-.8);
\draw[mid<] (-.2,.4) -- (-.2,-.2);
\draw[mid>] (-.2,.4) -- (.2,.8);
\draw[mid>] (-.2,.4) -- (-.6,.8);
\filldraw (-.2,.4) node[left]{$\scriptstyle k$} circle (.05cm);
\filldraw (-.2,-.2) node[left]{$\scriptstyle j$} circle (.05cm);
\filldraw (-.4,-.8) node[left]{$\scriptstyle i$} circle (.05cm);
}
+
\widehat{D}_{k,\ell}^{i,j}
\tikzmath{
\draw[thick, red] (-.2,.-.2) .. controls ++(135:.5cm) and ++(90:.4cm) ..(.2,-.2) .. controls ++(-90:.8cm) and ++(90:.3cm) .. (-.8,-1.2) -- (-.8,-1.4);
\draw[thick, red] (-.8,.8) arc (-180:0:.2cm);
\draw[thick, red] (-.8,-1.25) -- (-.7,-1.25); 
\draw[mid>] (-.4,-.8) -- (-.2,-.2);
\draw[mid>] (.4,-1.4) -- (-.2,-.2);
\draw[mid>] (-.2,-1.4) -- (-.4,-.8);
\draw[mid>] (-.6,-1.4) -- (-.4,-.8);
\draw[mid<] (-.2,.4) to[out=-90,in=45] (-.2,-.2);
\draw[mid>] (-.2,.4) -- (.2,.8);
\draw[mid>] (-.2,.4) -- (-.6,.8);
\filldraw (-.2,.4) node[left]{$\scriptstyle k$} circle (.05cm);
\filldraw (-.2,-.2) node[left]{$\scriptstyle j$} circle (.05cm);
\filldraw (-.4,-.8) node[left]{$\scriptstyle i$} circle (.05cm);
}
\end{align*}
\end{footnotesize}
For the second diagram in the last line above, we used the relation that for $f\in \cC_m(\rho^2 \to \alpha \rho)$,
\begin{equation}
\label{eq:RotatedFrobeniusComposite}
\tikzmath{
\draw[mid>] (-.8,-.9) -- (-.8,-.3);
\draw[thick] (-.8,-.1) -- (-.7,.-.1);
\draw[mid>] (.2,.3)  .. controls ++(90:.5cm) and ++(90:.5cm) .. (-.8,.3)  -- (-.8,-.3);
\draw[mid>] (.8,.3) -- (.8,-.3) arc (0:-180:.3cm);
\draw[mid>] (.8,.3) -- (.8,.9);
\draw[thick] (.8,.1) -- (.7,.1);
\draw[mid>] (-.2,-.9) -- (-.2,-.3);
\draw[thick, red] (-.2,.3) arc (0:180:.2cm) -- (-.6,-.3)  .. controls ++(-90:.3cm) and ++(90:.3cm) .. (.2,-.9);
\roundNbox{fill=white}{(0,0)}{.3}{.1}{.1}{$f$}
}
=
\tikzmath{
\draw[mid>] (.4,.3) -- (.4,.9);
\draw[mid>] (-.4,-.9) -- (-.4,-.3);
\draw[mid>] (.4,-.9) -- (.4,-.3);
\draw[thick, red] (-.4,.3) .. controls ++(90:.5cm) and ++(90:.5cm) .. (1,.3) -- (1,-.9);
\roundNbox{fill=white}{(0,0)}{.3}{.5}{.5}{$\scriptstyle (L^\alpha\circ R^\alpha)(f)$}
}
\end{equation}
which can be verified by a straightforward diagrammatic calculation.
We leave the final simplification of this second diagram using the $\alpha$ Jellyfish Relation of Lemma \ref{lem:alphaJellyfish} to the reader.
\end{proof}

\begin{rem}
\label{rem:SD6j}
Recall that the associator $F$-tensors of a unitary fusion category are determined by the formula
\[  
\tikzmath{
\draw[mid>] (.4,-.7) node[below]{$\scriptstyle Z$} -- (.2,-.1);
\draw[mid>] (0,-.7) node[below]{$\scriptstyle Y$} -- (.2,-.1);
\draw[mid>] (-.4,-.7) node[below]{$\scriptstyle X$} -- (0,.5);
\draw[mid>] (.2,-.1) -- node[right]{$\scriptstyle U$} (0,.5);
\draw[mid>] (0,.5) -- (0,1) node[above]{$\scriptstyle W$};
\filldraw (.2,-.1) node[left]{$\scriptstyle\ell$} circle (.05cm);
\filldraw (0,.5) node[left]{$\scriptstyle k$} circle (.05cm);
}
=   
\sum_{\substack{U\in \Irr(\cC)   \\ 0\leq i < \dim\Hom(X\otimes Y\to V) \\0\leq j < \dim\Hom(V\otimes Z\to W)}}  \left(F^{X,Y,Z}_W\right)^{(V; i,j)}_{(U;k,l)} 
\tikzmath{
\draw[mid>] (.4,-.7) node[below]{$\scriptstyle Z$} -- (0,.5);
\draw[mid>] (0,-.7) node[below]{$\scriptstyle Y$} -- (-.2,-.1);
\draw[mid>] (-.4,-.7) node[below]{$\scriptstyle X$} -- (-.2,-.1);
\draw[mid>] (-.2,-.1) -- node[left]{$\scriptstyle V$} (0,.5);
\draw[mid>] (0,.5) -- (0,1) node[above]{$\scriptstyle W$};
\filldraw (-.2,-.1) node[left]{$\scriptstyle i$} circle (.05cm);
\filldraw (0,.5) node[left]{$\scriptstyle j$} circle (.05cm);
}.
\]
We have the following identification between the above $8m^4$ complex scalars and certain $F$-tensors of the category $\cC_m$:
\begin{align*}
A^{i,j}_{k,\ell} &= \left(F^{\rho, \rho, \rho}_{\rho}\right)^{(\rho; i,j)}_{(\rho; k,\ell)} 
& 
B^{i,j}_{k,\ell} &= \left(F^{\rho, \rho, \rho}_{\rho}\right)^{(\alpha\rho; i,j)}_{(\rho; k,\ell)}
&
C^{i,j}_{k,\ell} &= \left(F^{\rho, \rho, \rho}_{\alpha\rho}\right)^{(\rho; i,j)}_{(\rho; k,\ell)} 
& 
D^{i,j}_{k,\ell} &= \left(F^{\rho, \rho, \rho}_{\alpha\rho}\right)^{(\alpha\rho; i,j)}_{(\rho; k,\ell)}
\\
\widehat{A}^{i,j}_{k,\ell} &= \left(F^{\alpha\rho, \rho, \rho}_{\alpha\rho}\right)^{(\rho; i,j)}_{(\alpha\rho; k,\ell)} 
&
\widehat{B}^{i,j}_{k,\ell} &= \left(F^{\alpha\rho, \rho, \rho}_{\alpha\rho}\right)^{(\alpha\rho; i,j)}_{(\alpha\rho; k,\ell)}
&
\widehat{C}^{i,j}_{k,\ell} &= \left(F^{\alpha\rho, \rho, \rho}_{\rho}\right)^{(\rho; i,j)}_{(\alpha\rho; k,\ell)} 
& 
\widehat{D}^{i,j}_{k,\ell} &= \left(F^{\alpha\rho, \rho, \rho}_{\rho}\right)^{(\alpha\rho; i,j)}_{(\alpha\rho; k,\ell)}.
\end{align*}
In the name of readability, we will not use this $F$-tensor notation in this article.
\end{rem}

%%%%%%%%%%%%%%%%%%%%%%%%%%%%%%%%%%%%%%%%%%%%%%%%%%%%%%%%%%%%%%%%
\subsubsection{Absorption relations}

Using the nomenclature from \cite{MR2979509}, a closed diagram in our generators is said to be in \emph{jellyfish form} if all trivalent and tetravalent vertices and their labels appear on the external boundary of the closed diagram.
By a slight abuse of nomenclature, we will say that a morphism in a hom space is in jellyfish form (or a \emph{train} in the nomenclature of \cite{MR3157990}) if all labels of trivalent and tetravalent vertices in the morphism meet the leftmost region of the morphism.

\begin{lem}[(Absorption)]
\label{lem:AbsorptionRelations}
Using the relations from \S\ref{sec:SDBasicRelations} and \S\ref{sec:SDJellyfishRelations}, any two trivalent/tetravalent vertices in jellyfish form connected by two of their $\rho$ strands so that the composite is still in jellyfish form
may be simplified into a diagram with no trivalent/tetravalent vertices.
\end{lem}
\begin{proof}
There are 16 words of length 2 on the symbols
$$
\left\{
\tikzmath{
\draw[mid>] (-.3,-.6) -- (0,0);
\draw[mid>] (.3,-.6) -- (0,0);
\draw[mid<] (0,.6) -- (0,0);
\filldraw (0,0) circle (.05cm);
}
\,,\,
\tikzmath{
\draw[mid>] (-.3,.6) -- (0,0);
\draw[mid>] (.3,.6) -- (0,0);
\draw[mid<] (0,-.6) -- (0,0);
\filldraw (0,0) circle (.05cm);
}
\,,\,
\tikzmath{
\draw[thick, red, mid>] (0,0) -- (-.6,.6);
\draw[mid>] (-.6,-.6) -- (0,0);
\draw[mid>] (.6,-.6) -- (0,0);
\draw[mid>] (0,0) -- (.6,.6);
\filldraw (0,0) circle (.05cm);
}
\,,\,
\tikzmath{
\draw[thick, red, mid<] (0,0) -- (-.6,-.6);
\draw[mid>] (.6,-.6) -- (0,0);
\draw[mid<] (-.6,.6) -- (0,0);
\draw[mid>] (0,0) -- (.6,.6);
\filldraw (0,0) circle (.05cm);
}
\right\},
$$
and up to adjoints, 10 are distinct.
Given any word of length 2, there is a unique composite in jellyfish form with two $\rho$ strands connected, up to labels and moving tags through crossings. 
There are thus 10 cases to consider:
\begin{align*}
&
\tikzmath{
\draw[mid>] (-.3,-.8) -- (-.3,-.3) arc (180:90:.3cm);
\draw[mid>] (.9,-.3) arc (0:-180:.3cm) arc (0:90:.3cm);
\draw[thick] (.9,-.3) -- (.8,-.3);
\draw[mid>] (0,0) -- (0,.15) arc (180:90:.45cm);
\draw[mid>] (.9,-.3) -- (.9,.15) arc (0:90:.45cm);
\draw[mid>] (.45,.6) -- (.45,1);
\filldraw (0,0) node[left,yshift=.1cm]{$\scriptstyle i$} circle (.05cm);
\filldraw (.45,.6) node[left,yshift=.1cm]{$\scriptstyle j$} circle (.05cm);
}
\qquad
\tikzmath{
\draw[mid>] (0,.3) -- (0,.7);
\draw[mid>] (0,-.7) -- (0,-.3);
\draw[mid>] (0,-.3) arc (-90:90:.3cm);
\draw[mid>] (0,-.3) arc (270:90:.3cm);
\filldraw (0,-.3) node[left,yshift=-.1cm]{$\scriptstyle i$} circle (.05cm);
\filldraw (0,.3) node[left,yshift=.1cm]{$\scriptstyle j$} circle (.05cm);
}
\qquad
\tikzmath{
\draw[mid>] (-.3,-.8) -- (-.3,-.3) to[out=90,in=-135] (0,0);
\draw[mid>] (.9,-.3) arc (0:-180:.3cm) to[out=90,in=-45] (0,0);
\draw[thick] (.9,-.3) -- (.8,-.3);
\draw[mid>] (0,0) to[out=45,in=-135] (.45,1);
\draw[mid>] (.9,-.3) -- (.9,.15) to[out=90,in=-45] (.45,1);
\draw[mid>] (.45,1) -- (.45,1.4);
\draw[red, thick] (0,0) to[out=135,in=-90] (-.2,.2) to[out=90,in=180] (1.2,.6);
\filldraw (0,0) node[left]{$\scriptstyle i$} circle (.05cm);
\filldraw (.45,1) node[left,yshift=.1cm]{$\scriptstyle j$} circle (.05cm);
}
\qquad
\tikzmath{
\draw[mid>] (0,.3) -- (0,.7);
\draw[mid>] (.4,-1.3) to[out=90,in=-45] (0,-.3);
\draw[red, thick] (.8,-1) to[out=180,in=-135] (-.3,-.6) -- (0,-.3);
\draw[mid>] (0,-.3) to[out=135,in=-90] (-.3,0) arc (180:90:.3cm); 
\draw[mid>] (0,-.3) to[out=45,in=-90] (.3,0) arc (0:90:.3cm); 
\filldraw (0,-.3) node[left]{$\scriptstyle i$} circle (.05cm);
\filldraw (0,.3) node[left,yshift=.1cm]{$\scriptstyle j$} circle (.05cm);
}
\qquad
\tikzmath[yscale=-1]{
\draw[mid<] (-.3,-.8) -- (-.3,-.3) arc (180:90:.3cm);
\draw[mid<] (.9,-.3) arc (0:-180:.3cm) arc (0:90:.3cm);
\draw[thick] (.9,-.3) -- (.8,-.3);
\draw[mid<] (0,0) -- (0,.15) arc (180:90:.45cm);
\draw[mid<] (.9,-.3) -- (.9,.15) arc (0:90:.45cm);
\draw[mid<] (.45,.6) -- (.45,1);
\filldraw (0,0) node[left,yshift=-.1cm]{$\scriptstyle j$} circle (.05cm);
\filldraw (.45,.6) node[left,yshift=-.1cm]{$\scriptstyle i$} circle (.05cm);
}
\\
&
\tikzmath{
\draw[mid>] (-.3,-.8) -- (-.3,-.3) to[out=90,in=-135] (0,0);
\draw[mid>] (.4,1) .. controls ++(45:.7cm) and ++(90:.3cm) .. (.8,.5) .. controls ++(-90:.3cm) and ++(-45:.7cm) .. (0,0);
\draw[mid>] (0,0) to[out=45,in=-90] (.4,1);
\draw[mid>] (.4,1) -- (.2,1.4);
\draw[red, thick] (0,0) to[out=135,in=-90] (-.2,.2) to[out=90,in=180] (1,.6);
\filldraw (0,0) node[left]{$\scriptstyle i$} circle (.05cm);
\filldraw (.4,1) node[left,yshift=-.1cm]{$\scriptstyle j$} circle (.05cm);
}
\qquad
\tikzmath{
\draw[mid>] (-.3,.3) arc (270:180:.3cm) -- (-.6,.9);
\draw[mid>] (.4,-1.3) to[out=90,in=-45] (0,-.3);
\draw[red, thick] (.8,-1) to[out=180,in=-135] (-.3,-.6) -- (0,-.3);
\draw[mid>] (0,-.3) to[out=135,in=-90] (-.3,.3); 
\draw[mid>] (0,-.3) to[out=45,in=-90] (.3,.3);
\draw[thick] (.3,.3) -- (.2,.3);
\draw[mid<] (.3,.3) -- (.3,.5) arc (0:180:.2cm) arc (0:-90:.2cm); 
\filldraw (0,-.3) node[left]{$\scriptstyle i$} circle (.05cm);
\filldraw (-.3,.3) node[left,yshift=-.1cm]{$\scriptstyle j$} circle (.05cm);
}
\qquad
\tikzmath{
\draw[mid<] (.4,1.3) to[out=-90,in=45] (0,.3);
\draw[red, thick] (.8,1) to[out=180,in=135] (-.3,.6) -- (0,.3);
\draw[mid>] (.4,-1.3) to[out=90,in=-45] (0,-.3);
\draw[red, thick] (.8,-1) to[out=180,in=-135] (-.3,-.6) -- (0,-.3);
\draw[mid>] (0,-.3) to[out=135,in=-90] (-.3,0) to[out=90,in=-135] (0,.3); 
\draw[mid>] (0,-.3) to[out=45,in=-90] (.3,0) to[out=90,in=-45] (0,.3); 
\filldraw (0,-.3) node[left]{$\scriptstyle i$} circle (.05cm);
\filldraw (0,.3) node[left]{$\scriptstyle j$} circle (.05cm);
}
\qquad
\tikzmath{
\draw[mid<] (.8,2) to[out=-90,in=45] (.4,1);
\draw[red, thick] (1.2,1.7) to[out=180,in=135] (.1,1.3) -- (.4,1);
\draw[mid>] (-.3,-.8) -- (-.3,-.3) to[out=90,in=-135] (0,0);
\draw[mid<] (0,0) to[out=-45, in=90] (.3,-.3) arc (-180:0:.3cm);
\draw[mid>] (.9,-.3) to[out=90,in=-45] (.4,1);
\draw[mid>] (0,0) to[out=45,in=-135] (.4,1);
\draw[thick] (.9,-.3) -- (.8,-.3);
\draw[red, thick] (0,0) to[out=135,in=-90] (-.2,.2) to[out=90,in=180] (1.2,.6);
\filldraw (0,0) node[left]{$\scriptstyle i$} circle (.05cm);
\filldraw (.4,1) node[left]{$\scriptstyle j$} circle (.05cm);
}
\qquad
\tikzmath{
\draw[mid<] (-.3,1.6) --  (0,1.2);
\draw[mid>] (-.3,-.4) -- (0,0);
\draw[mid>] (0,1.2) .. controls ++(45:.7cm) and ++(90:.3cm) .. (.8,.6) .. controls ++(-90:.3cm) and ++(-45:.7cm) .. (0,0);
\draw[mid>] (0,0) .. controls ++(45:.4cm) and ++(-45:.4cm) .. (0,1.2);
\draw[red, thick] (0,0) .. controls ++(135:.7cm) and ++(180:.3cm) .. (1.2,.4);
\draw[red, thick] (0,1.2) .. controls ++(-135:.7cm) and ++(180:.3cm) .. (1.2,.8);
\filldraw (0,0) node[left]{$\scriptstyle i$} circle (.05cm);
\filldraw (0,1.2) node[left]{$\scriptstyle j$} circle (.05cm);
}\,.
\end{align*}
We give a full proof for the last case, and the others are similar and omitted:
\[
\tikzmath{
\draw[mid<] (-.3,1.6) --  (0,1.2);
\draw[mid>] (-.3,-.4) -- (0,0);
\draw[mid>] (0,1.2) .. controls ++(45:.7cm) and ++(90:.3cm) .. (.8,.6) .. controls ++(-90:.3cm) and ++(-45:.7cm) .. (0,0);
\draw[mid>] (0,0) .. controls ++(45:.4cm) and ++(-45:.4cm) .. (0,1.2);
\draw[red, thick] (0,0) .. controls ++(135:.7cm) and ++(180:.3cm) .. (1.2,.4);
\draw[red, thick] (0,1.2) .. controls ++(-135:.7cm) and ++(180:.3cm) .. (1.2,.8);
\filldraw (0,0) node[left]{$\scriptstyle i$} circle (.05cm);
\filldraw (0,1.2) node[left]{$\scriptstyle j$} circle (.05cm);
}
=
\tikzmath{
\draw[mid<] (-.3,2) --  (0,1.2);
\draw[mid>] (-.3,-.8) -- (0,0);
\draw[mid>] (0,1.2) .. controls ++(45:.7cm) and ++(90:.3cm) .. (.8,1.2);
\draw[mid<] (0,0) .. controls ++(-45:.7cm) and ++(-90:.3cm) .. (.8,0);
\draw[mid>] (0,0) .. controls ++(45:.4cm) and ++(-45:.4cm) .. (0,1.2);
\draw[mid>] (.8,0) -- (.8,1.2);
\draw[thick] (.8,0) -- (.7,0);
\draw[thick] (.8,1.2) -- (.7,1.2);
\draw[red, thick] (0,0) .. controls ++(135:.4cm) and ++(90:.4cm) .. (.2,0) .. controls ++(-90:.3cm) and ++(90:.3cm) .. (-.5,-.4) .. controls ++(-90:.3cm) and ++(180:.3cm) .. (1.2,-.6);
\draw[red, thick] (0,1.2) .. controls ++(-135:.4cm) and ++(-90:.4cm) .. (.2,1.2) .. controls ++(90:.3cm) and ++(-90:.3cm) .. (-.5,1.6) .. controls ++(90:.3cm) and ++(180:.3cm) .. (1.2,1.8);
\filldraw (0,0) node[left]{$\scriptstyle i$} circle (.05cm);
\filldraw (0,1.2) node[left]{$\scriptstyle j$} circle (.05cm);
}
=
(\lambda_\rho\mu)^{i+j}
\tikzmath{
\draw[mid<] (.4,1.3) to[out=-90,in=45] (0,.3);
\draw[red, thick] (.8,1) to[out=180,in=135] (-.3,.6) -- (0,.3);
\draw[mid>] (.4,-1.3) to[out=90,in=-45] (0,-.3);
\draw[red, thick] (.8,-1) to[out=180,in=-135] (-.3,-.6) -- (0,-.3);
\draw[mid>] (0,-.3) to[out=135,in=-90] (-.3,0) to[out=90,in=-135] (0,.3); 
\draw[mid>] (0,-.3) to[out=45,in=-90] (.3,0) to[out=90,in=-45] (0,.3); 
\filldraw (0,-.3) node[left]{$\scriptstyle i$} circle (.05cm);
\filldraw (0,.3) node[left]{$\scriptstyle j$} circle (.05cm);
}
=
\delta_{i=j}
(\lambda_\rho\mu)^{i+j}
\,\,
\tikzmath{
\draw[mid>] (0,-.5) -- (0,.5);
\draw[thick, red, mid<] (.6,-.3) arc (270:90:.3cm);
}\,.
\qedhere
\]
\end{proof}

%%%%%%%%%%%%%%%%%%%%%%%%%%%%%%%%%%%%%%%%%%%%%%%%%%%%%%%%%%%%%%%%
\subsubsection{Evaluation algorithm}

With these local relations in hand, we can show that the numerical data we have described uniquely determines the category $\cC_m$.

\begin{prop}\label{prop:UniqueSD}
There is at most one unitary fusion category $\cC_m$ realising each tuple of data
\[(\lambda_\alpha,\lambda_\rho,\mu, \chi, \omega,\tilde{\quad}, A, B, C, D, \widehat{A},\widehat{B},\widehat{C},\widehat{D}).\]
\end{prop}
\begin{proof}
The proof is an adaptation of Bigelow's jellyfish algorithm \cite{MR2577673,MR2979509}.
Given any closed diagram in our generators, we show it can be evaluated to a scalar using our relations.
This immediately implies the stated result.

By the jellyfish relations from Lemmas \ref{lem:alphaJellyfish} and  \ref{lem:SDJellyfishRelations}, it suffices to show we can evaluate any closed diagram in \emph{jellyfish form}, in which all trivalent and tetravalent vertices and their labels appear on the external boundary of the closed diagram.
There are 3 cases for such a diagram:
\begin{enumerate}[label=\underline{Case \arabic*:}]
\item 
there are no vertices at all in the closed diagram.
Then we may use \eqref{eq:SD-Spaghetti} to evaluate the closed diagram to a scalar.
\item
there is a trivalent/tetravalent vertex connected to itself.
Then we may use \eqref{eq:SD-Spaghetti} to show that this closed diagram is equal to zero.
\item
there are two neighboring trivalent/tetravalent vertices that are connected by at least 2 of their $\rho$ strands.
Then using the absoprtion relations from Lemma \ref{lem:AbsorptionRelations}, we can express our closed diagram in jellyfish form as a linear combination of diagrams with strictly fewer vertices, which are still in jellyfish form.
\end{enumerate}
We are finished by a simple induction argument on the number of vertices in our closed diagram in jellyfish form.
\end{proof}

\begin{rem}
We wish to point out that we can also give an existence result for the categories $\cC_m$ by realising them as actions by endomorphisms on the Cuntz algebras $O_{2m+1}\rtimes \mathbb{Z}_{2}$. 
To obtain existence one needs to verify a finite list of polynomial equations that the above tuple needs to satisfy. 
As we can conclude existence of the examples in this article from the existing literature, we will not 
include the details of this existence result.
\end{rem}
%%%%%%%%%%%%%%%%%%%%%%%%%%%%%%%%%%%%%%%%%%%%%%%%%%%%%%%%%%%%%%%%%%%
\subsection{Centre analysis}
\label{SD:CenterAnalysis}

In this subsection we will analyse the Drinfeld centre of the categories $\cC_2$ in order to pin down the values of our data $\chi$, $\mu$, and $\tilde{\quad}$. We restrict our attention to the $m=2$ case here as with Corollary~\ref{cor:SDm1} in hand, the $m=1$ case is already covered by \cite[Example 9.1]{MR3827808}.

\begin{lem}\label{lem:maincentre}
If $m=2$, then
$\mu = 1$,
and 
\[\chi_{\mathbf{1}, 0} = \sqrt{\lambda_\alpha}, \quad  \chi_{\mathbf{1}, 1} = -\sqrt{\lambda_\alpha} ,\quad   \chi_{\alpha, 0} =1, \quad \text{and} \quad  \chi_{\alpha, 1}= -1.\]
There are two cases depending on $\lambda_\alpha$:
\begin{itemize}
\item
if $\lambda_\alpha = 1$, then 
$(\mathbf{1}, \tilde{i}) =(\mathbf{1}, i)$
and
$(\alpha, \tilde{i}) =(\alpha, i)$, and 
\item
if $\lambda_\alpha = -1$, then 
$(\mathbf{1}, \tilde{i}) =(\mathbf{1}, 1-i)$ and $(\alpha, \tilde{i}) =(\alpha, 1-i)$.
\end{itemize}
\end{lem}

Our main tool to prove Lemma~\ref{lem:maincentre} is Theorem~\ref{thm:tubecentre}. To remind the reader this result states that for an object $X \in \cC$, the irreducible representations $V$ of $A_{X \to X}$ are in bijective correspondence with simple summands $\Gamma_V \subset \cI(X) \in Z(\cC)$, the dimension of $\Gamma_V$ is given by $\frac{\dim(\cC)}{\dim(X) f_V}$ where $f_V$ is the formal codegree of the representation $V$, and the multiplicity of $Y \in \cF(\Gamma_V)$ is equal to the multiplicity of $V$ in the left-action of $A_{X\to X}$ on $A_{X\to Y}$. 

With this tool in mind, we aim to study the tube algebra for $\cC_2$:
\[A \quad  = \quad  
    \begin{tabular}{|c|c|c|c|}
     \hline $A_{\mathbf{1}\leftarrow\mathbf{1}}$ & $A_{\mathbf{1}\leftarrow\alpha}$ &$A_{\mathbf{1}\leftarrow\rho}$  &$A_{\mathbf{1}\leftarrow\alpha\rho}$  \\\hline
      
      \hline $A_{\alpha\leftarrow\mathbf{1}}$ & $A_{\alpha\leftarrow\alpha}$ &$A_{\alpha\leftarrow\rho}$  &$A_{\alpha\leftarrow\alpha\rho}$  \\\hline
      $A_{\rho\leftarrow\mathbf{1}}$ & $A_{\rho\leftarrow\alpha}$ &$A_{\rho\leftarrow\rho}$  &$A_{\rho\leftarrow\alpha\rho}$  \\\hline
      
      $A_{\alpha\rho\leftarrow\mathbf{1}}$ & $A_{\alpha\rho\leftarrow\alpha}$ &$A_{\alpha\rho\leftarrow\rho}$  &$A_{\alpha\rho\leftarrow\alpha\rho}$  \\\hline
    \end{tabular}
    \]
By determining the irreducible representations of the sub-algebras $A_{X\leftarrow X}$, and their  multiplicities in the left action of $A_{X\leftarrow X}$ on $A_{X\leftarrow Y}$, we can determine the simple objects of $Z(\cC_2)$, and their images under the forgetful functor.

Performing this computation over all of the tube algebra is far too computationally taxing. Instead we restrict our attention to the sub-algebra
\[  
A_{\mathbf{1}\leftarrow\mathbf{1}}
\oplus
A_{\alpha\leftarrow\alpha}
\oplus
A_{\mathbf{1}\leftarrow\rho}
\oplus
A_{\alpha\leftarrow\rho}
\oplus
A_{\mathbf{1}\leftarrow\alpha\rho}.
\]
We pick the following bases for these spaces:
\begin{align*}
    A_{\mathbf{1}\leftarrow \mathbf{1}} &=\operatorname{span} \left\{\,\,  
    \tikzmath{\filldraw[very thick, fill=gray!30] (0,0) circle (.2cm);}
    \,,\,\,
    \tikzmath{\filldraw[very thick, fill=gray!30] (0,0) circle (.2cm);
    \draw[thick, red, mid<] (0,0) circle (.4cm);}
    \,,\,\,
    \tikzmath{\filldraw[very thick, fill=gray!30] (0,0) circle (.2cm);
    \draw[mid<] (0,0) circle (.4cm);}
    \,,\,\,
    \tikzmath{\filldraw[very thick, fill=gray!30] (0,0) circle (.2cm);
    \draw[red, thick, mid<] (0,0) circle (.65cm);
    \draw[mid<] (0,0) circle (.4cm);}\,\,  
    \right\}
    \\
    A_{\mathbf{1}\leftarrow \rho} 
    &= 
    \operatorname{span}\left\{\,\,
    \tikzmath{\filldraw[very thick, fill=gray!30] (0,-.2) circle (.2cm);
    \draw[mid<] (0,0) ellipse (.4cm and .6cm);
    \draw[mid>] (0,0) .. controls ++(90:.2cm) and ++(-45:.3cm) .. (135:.46cm);
    \filldraw (135:.46cm) node[left]{$\scriptstyle 0$} circle (.05cm);
    }
    \,,\,\,
    \tikzmath{\filldraw[very thick, fill=gray!30] (0,-.2) circle (.2cm);
    \draw[mid<] (0,0) ellipse (.4cm and .6cm);
    \draw[mid>] (0,0) .. controls ++(90:.2cm) and ++(-45:.3cm) .. (135:.46cm);
    \filldraw (135:.46cm) node[left]{$\scriptstyle 1$} circle (.05cm);
    }\,,\,\,     
    \tikzmath{\filldraw[very thick, fill=gray!30] (0,-.2) circle (.2cm);
    \draw[mid<] (0,0) ellipse (.4cm and .6cm);
    \draw[thick, red, mid<] (0,0) ellipse (.7cm and .9cm);
    \draw[mid>] (0,0) .. controls ++(90:.2cm) and ++(-45:.3cm) .. (135:.46cm);
    \filldraw (135:.46cm) node[left]{$\scriptstyle 0$} circle (.05cm);
    }
    \,,\,\,
    \tikzmath{\filldraw[very thick, fill=gray!30] (0,-.2) circle (.2cm);
    \draw[mid<] (0,0) ellipse (.4cm and .6cm);
    \draw[thick, red, mid<] (0,0) ellipse (.7cm and .9cm);
    \draw[mid>] (0,0) .. controls ++(90:.2cm) and ++(-45:.3cm) .. (135:.46cm);
    \filldraw (135:.46cm) node[left]{$\scriptstyle 1$} circle (.05cm);
    }\,\,    
    \right\}\\
    A_{\mathbf{1}\leftarrow \alpha\rho} &= \operatorname{span}\left\{
    \tikzmath{
    \coordinate (a) at (115:.52cm);
    \coordinate (b) at (145:.45cm);
    \filldraw[very thick, fill=gray!30] (0,-.2) circle (.2cm);
    \draw[mid<] (0,0) ellipse (.4cm and .6cm);
    \draw[mid>] (0,0) .. controls ++(90:.2cm) and ++(-45:.3cm) .. (a);
    \draw[red, thick, mid>] ($ (135:.2cm) + (0,-.2) $) .. controls ++(135:.2cm) and ++(-45:.3cm) .. (b);
    \draw[thick, red, far>] (a) .. controls ++(135:.6cm) and ++(135:.6cm) .. (b);
    \filldraw (a) node[left]{$\scriptstyle 0$} circle (.05cm);
    }
    \,,\,\,
    \tikzmath{
    \coordinate (a) at (115:.52cm);
    \coordinate (b) at (145:.45cm);
    \filldraw[very thick, fill=gray!30] (0,-.2) circle (.2cm);
    \draw[mid<] (0,0) ellipse (.4cm and .6cm);
    \draw[mid>] (0,0) .. controls ++(90:.2cm) and ++(-45:.3cm) .. (a);
    \draw[red, thick, mid>] ($ (135:.2cm) + (0,-.2) $) .. controls ++(135:.2cm) and ++(-45:.3cm) .. (b);
    \draw[thick, red, far>] (a) .. controls ++(135:.6cm) and ++(135:.6cm) .. (b);
    \filldraw (a) node[left]{$\scriptstyle 1$} circle (.05cm);
    }
    \,,\,\, 
    \tikzmath{
    \coordinate (b) at (115:.82cm);
    \coordinate (a) at (145:.65cm);
    \filldraw[very thick, fill=gray!30] (0,-.2) circle (.2cm);
    \draw[mid<] (0,0) ellipse (.6cm and .9cm);
    \draw[mid>] (0,0) .. controls ++(90:.2cm) and ++(-45:.3cm) .. (a);
    \draw[red, thick, mid>] ($ (135:.2cm) + (0,-.2) $) .. controls ++(135:.2cm) and ++(90:.2cm) .. (-.4,-.2) arc (-180:0:.4cm) .. controls ++(90:.4cm) and ++(-45:.3cm) .. (b);
    \draw[thick, red, far>] (a) .. controls ++(135:.6cm) and ++(135:.6cm) .. (b);
    \filldraw (a) node[left]{$\scriptstyle 0$} circle (.05cm);
    }
    \,,\,\,
    \tikzmath{
    \coordinate (b) at (115:.82cm);
    \coordinate (a) at (145:.65cm);
    \filldraw[very thick, fill=gray!30] (0,-.2) circle (.2cm);
    \draw[mid<] (0,0) ellipse (.6cm and .9cm);
    \draw[mid>] (0,0) .. controls ++(90:.2cm) and ++(-45:.3cm) .. (a);
    \draw[red, thick, mid>] ($ (135:.2cm) + (0,-.2) $) .. controls ++(135:.2cm) and ++(90:.2cm) .. (-.4,-.2) arc (-180:0:.4cm) .. controls ++(90:.4cm) and ++(-45:.3cm) .. (b);
    \draw[thick, red, far>] (a) .. controls ++(135:.6cm) and ++(135:.6cm) .. (b);
    \filldraw (a) node[left]{$\scriptstyle 1$} circle (.05cm);
    }   
    \right\}\\
    A_{\alpha\leftarrow\alpha} 
    &= 
    \operatorname{span}\left\{\,\,  
    \tikzmath{\filldraw[very thick, fill=gray!30] (0,0) circle (.2cm);
    \draw[red, thick, mid>] (0,.2) -- (0,.8);}\,,\,\, 
    \tikzmath{\filldraw[very thick, fill=gray!30] (0,0) circle (.2cm);
    \draw[red, thick, mid>] (0,.2) arc (180:0:.2cm) -- (.4,0) arc (0:-180:.4cm) .. controls ++(90:.3cm) and ++(270:.3cm) .. (0,.8);}\,,\,\,
    \tikzmath{\filldraw[very thick, fill=gray!30] (0,0) circle (.2cm);
    \draw[red, thick, mid>] (0,.2) -- (0,.6);
    \draw[mid<] (0,0) circle (.6cm);
    \draw[thick, red, mid<] (0,.6) -- (0,1);
    \draw[thick, red] (0,1) -- (.1,1);
    \draw[thick, red, mid>] (0,1) -- (0,1.4);
    }\,,\,\,
    \tikzmath{\filldraw[very thick, fill=gray!30] (0,0) circle (.2cm);
    \draw[red, thick, mid>] (0,.2) arc (180:0:.2cm) -- (.4,0) arc (0:-180:.4cm) .. controls ++(90:.3cm) and ++(270:.3cm) .. (0,.8);
    \draw[thick, red] (0,.8) -- (.1,.8);
    \draw[thick, red,mid>] (0,.8) -- (0,1.2);
    \draw[mid<] (0,0) circle (.6cm);}\,\,
    \right\}
    \\
    A_{\alpha\leftarrow\rho} &= \operatorname{span}\left\{ 
    \tikzmath{
    \coordinate (a) at (135:.46cm);
    \filldraw[very thick, fill=gray!30] (0,-.2) circle (.2cm);
    \draw[mid<] (0,0) ellipse (.4cm and .6cm);
    \draw[mid>] (0,0) .. controls ++(90:.2cm) and ++(-45:.3cm) .. (a);
    \draw[thick, red, mid>] (a) .. controls ++(135:.2cm) and ++(270:.3cm) .. (-.6,.8);
    \filldraw (a) node[left, yshift=-.1cm]{$\scriptstyle 0$} circle (.05cm);
    }\,,\,\,
    \tikzmath{
    \coordinate (a) at (135:.46cm);
    \filldraw[very thick, fill=gray!30] (0,-.2) circle (.2cm);
    \draw[mid<] (0,0) ellipse (.4cm and .6cm);
    \draw[mid>] (0,0) .. controls ++(90:.2cm) and ++(-45:.3cm) .. (a);
    \draw[thick, red, mid>] (a) .. controls ++(135:.2cm) and ++(270:.3cm) .. (-.6,.8);
    \filldraw (a) node[left, yshift=-.1cm]{$\scriptstyle 1$} circle (.05cm);
    }\,,\,\,
    \tikzmath{
    \coordinate (a) at (135:.46cm);
    \filldraw[very thick, fill=gray!30] (0,-.2) circle (.2cm);
    \draw[mid<] (0,0) ellipse (.4cm and .6cm);
    \draw[mid>] (0,0) .. controls ++(90:.2cm) and ++(-45:.3cm) .. (a);
    \draw[thick, red, mid>] (a) .. controls ++(135:1cm) and ++(90:1cm) .. (.6,.4) .. controls ++(270:.4cm) and ++(0:.6cm) .. (0,-.8) .. controls ++(180:.6cm) and ++(270:.6cm) .. (-.8,1.2);
    \filldraw (a) node[left, yshift=-.1cm]{$\scriptstyle 0$} circle (.05cm);
    }\,,\,\,
    \tikzmath{
    \coordinate (a) at (135:.46cm);
    \filldraw[very thick, fill=gray!30] (0,-.2) circle (.2cm);
    \draw[mid<] (0,0) ellipse (.4cm and .6cm);
    \draw[mid>] (0,0) .. controls ++(90:.2cm) and ++(-45:.3cm) .. (a);
    \draw[thick, red, mid>] (a) .. controls ++(135:1cm) and ++(90:1cm) .. (.6,.4) .. controls ++(270:.4cm) and ++(0:.6cm) .. (0,-.8) .. controls ++(180:.6cm) and ++(270:.6cm) .. (-.8,1.2);
    \filldraw (a) node[left, yshift=-.1cm]{$\scriptstyle 1$} circle (.05cm);
    }\,\,
    \right\}.
\end{align*}
By direct computation we obtain that: 
\begin{itemize}
    \item The irreducible representations of $A_{\mathbf{1}\leftarrow\mathbf{1}}$ are: 
$$
    \begin{array}{c|c c c c}
         &
         \tikzmath{\filldraw[very thick, fill=gray!30] (0,0) circle (.1cm);}
         &
        \tikzmath{\filldraw[very thick, fill=gray!30] (0,0) circle (.1cm);
        \draw[thick, red, mid<] (0,0) circle (.2cm);}
        &
        \tikzmath{\filldraw[very thick, fill=gray!30] (0,0) circle (.1cm);
        \draw[mid<] (0,0) circle (.2cm);}
        &
        \tikzmath{\filldraw[very thick, fill=gray!30] (0,0) circle (.1cm);
        \draw[red, thick, mid<] (0,0) circle (.35cm);
        \draw[mid<] (0,0) circle (.2cm);}
         \\ \hline
         \chi_0 &1 & 1& 2 + \sqrt{5}&2 + \sqrt{5} \\
         \chi_1 &1 & 1& 2 - \sqrt{5}& 2 - \sqrt{5}\\
         \chi_2 &1 &-1 & 1&-1 \\
         \chi_3 &1 & -1&-1 & 1
    \end{array}
$$
    Hence $\cI(\mathbf{1})$ contains 4 simple objects $X_i$ with dimensions
\[ \dim(X_0) = 1 , \quad \dim(X_1) = 9 + 4\sqrt{5},\quad  \text{and} \quad \dim(X_2) = \dim(X_3) = 5 + 2\sqrt{5}.\]
\item  The irreducible representations of $A_{\alpha\leftarrow\alpha}$ are:
\[  
\begin{array}{c|c c c c}
    & 
    \tikzmath{\filldraw[very thick, fill=gray!30] (0,0) circle (.2cm);
    \draw[red, thick, mid>] (0,.2) -- (0,.8);}
    &
    \tikzmath{\filldraw[very thick, fill=gray!30] (0,0) circle (.2cm);
    \draw[red, thick, mid>] (0,.2) arc (180:0:.2cm) -- (.4,0) arc (0:-180:.4cm) .. controls ++(90:.3cm) and ++(270:.3cm) .. (0,.8);}
    &
    \tikzmath{\filldraw[very thick, fill=gray!30] (0,0) circle (.2cm);
    \draw[red, thick, mid>] (0,.2) -- (0,.6);
    \draw[mid<] (0,0) circle (.6cm);
    \draw[thick, red, mid<] (0,.6) -- (0,1);
    \draw[thick, red] (0,1) -- (.1,1);
    \draw[thick, red, mid>] (0,1) -- (0,1.4);
    \filldraw[white] (0,-.75) circle (.01cm);
    }
    &
    \tikzmath{\filldraw[very thick, fill=gray!30] (0,0) circle (.2cm);
    \draw[red, thick, mid>] (0,.2) arc (180:0:.2cm) -- (.4,0) arc (0:-180:.4cm) .. controls ++(90:.3cm) and ++(270:.3cm) .. (0,.8);
    \draw[thick, red] (0,.8) -- (.1,.8);
    \draw[thick, red,mid>] (0,.8) -- (0,1.2);
    \draw[mid<] (0,0) circle (.6cm);}
    \\\hline
         \tau_{0} &1 &   \sqrt{\lambda_\alpha}& \frac{a\sqrt{\lambda_\alpha}  +\sqrt{a^2\lambda_\alpha +4\mu}}{2}& \frac{a\lambda_\alpha  +\sqrt{a^2 +4\lambda_\alpha\mu}}{2} \\
         \tau_{1} &1 &   -\sqrt{\lambda_\alpha}& \frac{-a\sqrt{\lambda_\alpha}  +\sqrt{a^2\lambda_\alpha +4\mu}}{2}& \frac{a\lambda_\alpha  -\sqrt{a^2 +4\lambda_\alpha\mu}}{2} \\
            \tau_{2} &1 &   \sqrt{\lambda_\alpha}& \frac{a\sqrt{\lambda_\alpha}  -\sqrt{a^2\lambda_\alpha +4\mu}}{2}& \frac{a\lambda_\alpha  -\sqrt{a^2 +4\lambda_\alpha\mu}}{2}
            \\ \tau_{3} &1 &   -\sqrt{\lambda_\alpha}& \frac{-a\sqrt{\lambda_\alpha}  -\sqrt{a^2\lambda_\alpha +4\mu}}{2}& \frac{a\lambda_\alpha  +\sqrt{a^2 +4\lambda_\alpha\mu}}{2} 
\end{array}
\]
where $a := \lambda_\alpha\sqrt{\lambda_\alpha}(\chi_{\mathbf{1}, 1} + \chi_{\mathbf{1}, 2}) +\lambda_\alpha( \chi_{\alpha, 1} + \chi_{\alpha, 2})$. 
(Note that we always have $a\in \{0,\pm 2,\pm4\}$.)
Hence $\cI(\alpha)$ contains 4 simple objects $Y_i$ with dimensions
\begin{align*}
\dim(Y_0) &= \dim(Y_3) =\frac{\dim(\cC)}{  4+a^2\lambda_\alpha\mu   +  a\sqrt{a^2 +4\lambda_\alpha\mu}}  
\\
\dim(Y_1) &= \dim(Y_2) =\frac{\dim(\cC)}{ 4+a^2\lambda_\alpha\mu   -  a\sqrt{a^2 +4\lambda_\alpha\mu}}.
\end{align*}

\item 
Let $_\mathbf{1} \pi _\rho$ be the action of $A_{\mathbf{1} \leftarrow\mathbf{1}}$ on  $A_{\mathbf{1} \leftarrow\rho}$. Then 
\[
\resizebox{.9\textwidth}{!}{
$
{}_\mathbf{1}\pi_\rho
\left(\,
    \tikzmath{\filldraw[very thick, fill=gray!30] (0,0) circle (.2cm);
    \draw[thick, red, mid<] (0,0) circle (.4cm);}
\,\right)
= 
\begin{bmatrix}
0 &0 & 1 & 0\\
0 & 0 &0& 1 \\
1 & 0 & 0 &0 \\
0 & 1 & 0 &0
\end{bmatrix} 
,\,\,
{}_\mathbf{1}\pi_\rho
\left(\,
    \tikzmath{\filldraw[very thick, fill=gray!30] (0,0) circle (.2cm);
    \draw[mid<] (0,0) circle (.4cm);}
\,\right) 
= 
\begin{bmatrix}
\phi & \phi'\\
\phi' & \phi
\end{bmatrix}
, \text{ and } 
{}_\mathbf{1}\pi_\rho
\left(\,
    \tikzmath{\filldraw[very thick, fill=gray!30] (0,0) circle (.2cm);
    \draw[red, thick, mid<] (0,0) circle (.65cm);
    \draw[mid<] (0,0) circle (.4cm);}
\,\right) 
= 
\begin{bmatrix}
\phi' & \phi\\
\phi & \phi'
\end{bmatrix}
$}
\]
where $\phi$ and $\phi'$ are the operators on $\Hom(\rho\otimes\rho\leftarrow\rho)$ defined by
\begin{align*}
\phi\left(
    \tikzmath{
    \draw[mid>] (-.3,-.6) -- (0,0);
    \draw[mid>] (.3,-.6) -- (0,0);
    \draw[mid<] (0,.6) -- (0,0);
    \filldraw (0,0) node[left]{$\scriptstyle i$} circle (.05cm);
    }
\right)
&=
\sum_{j} 
\tikzmath{
\draw[mid>] (-.2,-.4) -- (.2,.8);
\draw[mid>] (-.3,-.8) -- (-.2,-.4);
\draw[mid>] (.2,.8) -- (.2,1.2);
\draw[mid>] (-.2,-.4) -- (.4,0);
\draw[mid>] (.4,0) -- (.2,.8);
\draw[mid>] (.6,-.8) -- (.4,0);
\filldraw (-.2,-.4) node[left]{$\scriptstyle j$} circle (.05cm);
\filldraw (.2,.8) node[left]{$\scriptstyle j$} circle (.05cm);
\filldraw (.4,0) node[left, yshift=.05cm]{$\scriptstyle i$} circle (.05cm);
}
= 
\sum_{j,k} A^{j,k}_{j,i}
    \tikzmath{
    \draw[mid>] (-.3,-.6) -- (0,0);
    \draw[mid>] (.3,-.6) -- (0,0);
    \draw[mid<] (0,.6) -- (0,0);
    \filldraw (0,0) node[left]{$\scriptstyle k$} circle (.05cm);
    }
\\
\phi'\left(
    \tikzmath{
    \draw[mid>] (-.3,-.6) -- (0,0);
    \draw[mid>] (.3,-.6) -- (0,0);
    \draw[mid<] (0,.6) -- (0,0);
    \filldraw (0,0) node[left]{$\scriptstyle i$} circle (.05cm);
    }
\right)
&=
\sum_{j} 
\tikzmath{
\draw[thick, red] (-.2,-.4) .. controls ++(-135:.6cm) and ++(270:.5cm) .. (.1,-.4) .. controls ++(90:.3cm) and ++(270:.3cm) .. (-.3,.2) .. controls ++(90:.3cm) and ++(270:.3cm) .. (.4,.8) .. controls ++(90:.5cm) and ++(135:.6cm) .. (.2,.8);
\draw[mid>] (-.2,-.4) -- (.2,.8);
\draw[far>] (-.05,-1) -- (-.2,-.4);
\draw[far>] (.2,.8) -- (.3,1.5);
\draw[far>] (-.2,-.4) -- (.4,0);
\draw[mid>] (.4,0) -- (.2,.8);
\draw[mid>] (.65,-1) -- (.4,0);
\filldraw (-.2,-.4) node[left]{$\scriptstyle j$} circle (.05cm);
\filldraw (.2,.8) node[left]{$\scriptstyle j$} circle (.05cm);
\filldraw (.4,0) node[left, yshift=.1cm]{$\scriptstyle i$} circle (.05cm);
}
=
\chi_{\mathbf{1}, i}\sum_{j,k} D^{j,k}_{j,i}
    \tikzmath{
    \draw[mid>] (-.3,-.6) -- (0,0);
    \draw[mid>] (.3,-.6) -- (0,0);
    \draw[mid<] (0,.6) -- (0,0);
    \filldraw (0,0) node[left]{$\scriptstyle k$} circle (.05cm);
    }
\end{align*}
which we can naturally identify as operators on the two spaces:
\[
\left\{ \,\,
\tikzmath{\filldraw[very thick, fill=gray!30] (0,-.2) circle (.2cm);
    \draw[mid<] (0,0) ellipse (.4cm and .6cm);
    \draw[mid>] (0,0) .. controls ++(90:.2cm) and ++(-45:.3cm) .. (135:.46cm);
    \filldraw (135:.46cm) node[left]{$\scriptstyle 0$} circle (.05cm);
    }
    \,,\,\,
    \tikzmath{\filldraw[very thick, fill=gray!30] (0,-.2) circle (.2cm);
    \draw[mid<] (0,0) ellipse (.4cm and .6cm);
    \draw[mid>] (0,0) .. controls ++(90:.2cm) and ++(-45:.3cm) .. (135:.46cm);
    \filldraw (135:.46cm) node[left]{$\scriptstyle 1$} circle (.05cm);
    }\,\,
\right\}
\qquad\text{and}\qquad
\left\{\,\,
    \tikzmath{\filldraw[very thick, fill=gray!30] (0,-.2) circle (.2cm);
    \draw[mid<] (0,0) ellipse (.4cm and .6cm);
    \draw[thick, red, mid<] (0,0) ellipse (.7cm and .9cm);
    \draw[mid>] (0,0) .. controls ++(90:.2cm) and ++(-45:.3cm) .. (135:.46cm);
    \filldraw (135:.46cm) node[left]{$\scriptstyle 0$} circle (.05cm);
    }
    \,,\,\,
    \tikzmath{\filldraw[very thick, fill=gray!30] (0,-.2) circle (.2cm);
    \draw[mid<] (0,0) ellipse (.4cm and .6cm);
    \draw[thick, red, mid<] (0,0) ellipse (.7cm and .9cm);
    \draw[mid>] (0,0) .. controls ++(90:.2cm) and ++(-45:.3cm) .. (135:.46cm);
    \filldraw (135:.46cm) node[left]{$\scriptstyle 1$} circle (.05cm);
    }\,\, 
\right\}
\]
by local insertion, i.e., the elements of $A_{\mathbf{1}\leftarrow \mathbf{1}}$ which involve $\phi,\phi'$ above act on $A_{\mathbf{1}\leftarrow \rho}$ by applying $\phi,\phi'$ locally on the trivalent vertices in our standard basis of $A_{\mathbf{1}\leftarrow \rho}$.

\item 
Let $_\mathbf{1} \pi _{\alpha\rho}$ be the action of $A_{\mathbf{1} \leftarrow\mathbf{1}}$ on  $A_{\mathbf{1} \leftarrow\alpha\rho}$. Then 
\[ 
\resizebox{.9\textwidth}{!}{
$
{}_\mathbf{1} 
\pi_{\alpha\rho}
\left(\,
    \tikzmath{\filldraw[very thick, fill=gray!30] (0,0) circle (.2cm);
    \draw[thick, red, mid<] (0,0) circle (.4cm);}
\,\right)
= 
\begin{bmatrix}
0 &0 & 1 & 0\\
0 & 0 &0& 1 \\
1 & 0 & 0 &0 \\
0 & 1 & 0 &0
\end{bmatrix} \qquad 
{}_\mathbf{1} \pi _{\alpha\rho}
\left(\,
    \tikzmath{\filldraw[very thick, fill=gray!30] (0,0) circle (.2cm);
    \draw[mid<] (0,0) circle (.4cm);}
\,\right)
= 
\begin{bmatrix}
\psi' & \psi\\
\psi & \psi'
\end{bmatrix}
\quad \text{and} \qquad 
{}_\mathbf{1} \pi _{\alpha\rho}
\left(\,
    \tikzmath{\filldraw[very thick, fill=gray!30] (0,0) circle (.2cm);
    \draw[red, thick, mid<] (0,0) circle (.65cm);
    \draw[mid<] (0,0) circle (.4cm);}
\,\right)
= 
\begin{bmatrix}
\psi & \psi'\\
\psi' & \psi
\end{bmatrix}  
$
}
\]
where $\psi$ and $\psi'$ are the operators on $\Hom(\rho\otimes\rho\to \alpha\rho)$ defined by
\begin{align*}
\psi\left(
    \tikzmath{
    \draw[thick, red, mid>] (0,0) -- (-.6,.6);
    \draw[mid>] (-.6,-.6) -- (0,0);
    \draw[mid>] (.6,-.6) -- (0,0);
    \draw[mid>] (0,0) -- (.6,.6);
    \filldraw (0,0) node[left]{$\scriptstyle i$} circle (.05cm);
    }
\right)
&=
\sum_{j} 
\tikzmath{
\draw[thick, red] (.4,0) .. controls ++(135:.4cm) and ++(-90:.4cm) .. (.6,1) .. controls ++(90:.4cm) and ++(270:.3cm) .. (-.2,1.8);
\draw[thick, red] (-.2,-.4) .. controls ++(-135:.6cm) and ++(270:.5cm) .. (.1,-.4) .. controls ++(90:.3cm) and ++(270:.3cm) .. (-.3,.2) .. controls ++(90:.3cm) and ++(270:.3cm) .. (.4,.8) .. controls ++(90:.5cm) and ++(135:.6cm) .. (.2,.8);
\draw[mid>] (-.2,-.4) -- (.2,.8);
\draw[far>] (-.05,-1) -- (-.2,-.4);
\draw[mid>] (.2,.8) -- (.3,1.8);
\draw[far>] (-.2,-.4) -- (.4,0);
\draw[mid>] (.4,0) .. controls ++(45:.4cm) and ++(-45:.2cm) .. (.2,.8);
\draw[mid>] (.65,-1) -- (.4,0);
\filldraw (-.2,-.4) node[left]{$\scriptstyle j$} circle (.05cm);
\filldraw (.2,.8) node[left]{$\scriptstyle j$} circle (.05cm);
\filldraw (.4,0) node[left, yshift=.1cm]{$\scriptstyle i$} circle (.05cm);
}
=
\sum_{j,k} \widehat{A}^{j,k}_{j,i}
    \tikzmath{
    \draw[thick, red, mid>] (0,0) -- (-.6,.6);
    \draw[mid>] (-.6,-.6) -- (0,0);
    \draw[mid>] (.6,-.6) -- (0,0);
    \draw[mid>] (0,0) -- (.6,.6);
    \filldraw (0,0) node[left]{$\scriptstyle k$} circle (.05cm);
    }
\\ 
\psi'\left(
    \tikzmath{
    \draw[thick, red, mid>] (0,0) -- (-.6,.6);
    \draw[mid>] (-.6,-.6) -- (0,0);
    \draw[mid>] (.6,-.6) -- (0,0);
    \draw[mid>] (0,0) -- (.6,.6);
    \filldraw (0,0) node[left]{$\scriptstyle i$} circle (.05cm);
    }
\right)
&=
\sum_{j} 
\tikzmath{
\draw[thick, red] (.4,0) .. controls ++(135:.4cm) and ++(-90:.3cm) .. (.4,.8) .. controls ++(90:.3cm) and ++(270:.3cm) .. (-.2,1.5);
\draw[mid>] (-.2,-.4) -- (.2,.8);
\draw[mid>] (-.3,-.8) -- (-.2,-.4);
\draw[mid>] (.2,.8) -- (.2,1.5);
\draw[mid>] (-.2,-.4) -- (.4,0);
\draw[mid>] (.4,0) .. controls ++(45:.4cm) and ++(-45:.2cm) .. (.2,.8);
\draw[mid>] (.6,-.8) -- (.4,0);
\filldraw (-.2,-.4) node[left]{$\scriptstyle j$} circle (.05cm);
\filldraw (.2,.8) node[left]{$\scriptstyle j$} circle (.05cm);
\filldraw (.4,0) node[left, yshift=.05cm]{$\scriptstyle i$} circle (.05cm);
}
= 
\lambda_\alpha  \chi_{\alpha,i}\sum_{j,k} \widehat{D}^{j,k}_{j,i}
    \tikzmath{
    \draw[thick, red, mid>] (0,0) -- (-.6,.6);
    \draw[mid>] (-.6,-.6) -- (0,0);
    \draw[mid>] (.6,-.6) -- (0,0);
    \draw[mid>] (0,0) -- (.6,.6);
    \filldraw (0,0) node[left]{$\scriptstyle k$} circle (.05cm);
    }.
\end{align*}
As before, we can naturally identify $\psi,\psi'$ as operators on the following two spaces by local insertion:
\[
\left\{
    \tikzmath{
    \coordinate (a) at (115:.52cm);
    \coordinate (b) at (145:.45cm);
    \filldraw[very thick, fill=gray!30] (0,-.2) circle (.2cm);
    \draw[mid<] (0,0) ellipse (.4cm and .6cm);
    \draw[mid>] (0,0) .. controls ++(90:.2cm) and ++(-45:.3cm) .. (a);
    \draw[red, thick, mid>] ($ (135:.2cm) + (0,-.2) $) .. controls ++(135:.2cm) and ++(-45:.3cm) .. (b);
    \draw[thick, red, far>] (a) .. controls ++(135:.6cm) and ++(135:.6cm) .. (b);
    \filldraw (a) node[left]{$\scriptstyle 0$} circle (.05cm);
    }
    \,,\,\,
    \tikzmath{
    \coordinate (a) at (115:.52cm);
    \coordinate (b) at (145:.45cm);
    \filldraw[very thick, fill=gray!30] (0,-.2) circle (.2cm);
    \draw[mid<] (0,0) ellipse (.4cm and .6cm);
    \draw[mid>] (0,0) .. controls ++(90:.2cm) and ++(-45:.3cm) .. (a);
    \draw[red, thick, mid>] ($ (135:.2cm) + (0,-.2) $) .. controls ++(135:.2cm) and ++(-45:.3cm) .. (b);
    \draw[thick, red, far>] (a) .. controls ++(135:.6cm) and ++(135:.6cm) .. (b);
    \filldraw (a) node[left]{$\scriptstyle 1$} circle (.05cm);
    }
\right\}
\qquad\text{and}\qquad
\left\{
    \tikzmath{
    \coordinate (b) at (115:.82cm);
    \coordinate (a) at (145:.65cm);
    \filldraw[very thick, fill=gray!30] (0,-.2) circle (.2cm);
    \draw[mid<] (0,0) ellipse (.6cm and .9cm);
    \draw[mid>] (0,0) .. controls ++(90:.2cm) and ++(-45:.3cm) .. (a);
    \draw[red, thick, mid>] ($ (135:.2cm) + (0,-.2) $) .. controls ++(135:.2cm) and ++(90:.2cm) .. (-.4,-.2) arc (-180:0:.4cm) .. controls ++(90:.4cm) and ++(-45:.3cm) .. (b);
    \draw[thick, red, far>] (a) .. controls ++(135:.6cm) and ++(135:.6cm) .. (b);
    \filldraw (a) node[left]{$\scriptstyle 0$} circle (.05cm);
    }
    \,,\,\,
    \tikzmath{
    \coordinate (b) at (115:.82cm);
    \coordinate (a) at (145:.65cm);
    \filldraw[very thick, fill=gray!30] (0,-.2) circle (.2cm);
    \draw[mid<] (0,0) ellipse (.6cm and .9cm);
    \draw[mid>] (0,0) .. controls ++(90:.2cm) and ++(-45:.3cm) .. (a);
    \draw[red, thick, mid>] ($ (135:.2cm) + (0,-.2) $) .. controls ++(135:.2cm) and ++(90:.2cm) .. (-.4,-.2) arc (-180:0:.4cm) .. controls ++(90:.4cm) and ++(-45:.3cm) .. (b);
    \draw[thick, red, far>] (a) .. controls ++(135:.6cm) and ++(135:.6cm) .. (b);
    \filldraw (a) node[left]{$\scriptstyle 1$} circle (.05cm);
    }
\right\}.
\]
\item 
Denoting by ${}_\alpha\pi _{\rho}$ the action of $A_{\alpha \to\alpha}$ on  $A_{\alpha \leftarrow\rho}$,
we have
$
{}_\alpha\pi_{\rho}\left(\,
    \tikzmath{\filldraw[very thick, fill=gray!30] (0,0) circle (.2cm);
    \draw[red, thick, mid>] (0,.2) arc (180:0:.2cm) -- (.4,0) arc (0:-180:.4cm) .. controls ++(90:.3cm) and ++(270:.3cm) .. (0,.8);}\,
    \right) 
=   
\begin{bmatrix}
0 &0 & 1 & 0\\
0 & 0 &0& 1 \\
-1 & 0 & 0 &0 \\
0 & -1 & 0 &0
\end{bmatrix}.  
$
\end{itemize}

With these computations in hand, we can now determine a large amount of structure of $\cZ(\cC_2)$.
    
\begin{lem}\label{lem:SDidInd}
There exists $b \in \{0,1,2\}$ such that 
\begin{align*}
    \cF(X_0) &= \mathbf{1}\\
    \cF(X_1) &= \mathbf{1} \oplus 2\rho \oplus 2\alpha\rho\\
    \cF(X_2) &= \mathbf{1} \oplus b\rho \oplus (2-b)\alpha \rho\\
    \cF(X_3) &= \mathbf{1} \oplus (2-b)\rho \oplus b\alpha \rho.
\end{align*}
Furthermore, if $b\in \{0,2\}$, then the operators $\phi$ and $\psi$ are both the same scalar 
\[ \phi = \psi = \frac{1+b-\sqrt{5}}{2},
\]
and if $b=1$, the the operators $\phi$ and $\psi$ have the two eigenvalues
$\frac{1-\sqrt{5}}{2}$ and $\frac{3-\sqrt{5}}{2}$.
\end{lem}
\begin{proof}
First note that as $X_0$ is the tensor unit of $Z(\cC)$, we have that $_\mathbf{1} \pi _{\rho}$ and $_\mathbf{1} \pi _{\alpha\rho}$ contain no copies of $\chi_0$.
From the above computations, we have that 
$$\Tr\left({}_\mathbf{1}\pi_\rho
\left(\,
    \tikzmath{\filldraw[very thick, fill=gray!30] (0,0) circle (.2cm);
    \draw[thick, red, mid<] (0,0) circle (.4cm);}
\,\right)
\right)= 0.
$$ 
As $_\mathbf{1}\pi_\rho$ is $4$-dimensional, and $\chi_0$ is not a sub-representation, we must have that 
\[ 
{}_\mathbf{1}\pi_\rho \cong 2\chi_1 \oplus b\chi_2 \oplus (2-b)\chi_3   
\qquad\qquad
\text{where }b\in \{0,1,2\}.
\]
Thus $\cF(X_1)$ contains two copies of $\rho$, and a dimension count shows that 
\[ 
\cF(X_1) = \mathbf{1} \oplus 2\rho \oplus 2\alpha \rho.  
\]
From this we can deduce three possibilities for the restrictions of $X_2$ and $X_3$. 
\begin{enumerate}[label=\underline{Case \arabic*:}]
    \item $\cF(X_2) = \cF(X_3) = \mathbf{1} \oplus \rho \oplus \alpha\rho$, in which case $_\mathbf{1}\pi_\rho \cong    {{}_\mathbf{1}\pi_{\alpha\rho}} \cong 2\chi_1 \oplus \chi_2 \oplus \chi_3$, and in particular
    \[ \Tr\left(  {}_\mathbf{1}\pi_\rho
    \left(\,
    \tikzmath{\filldraw[very thick, fill=gray!30] (0,0) circle (.2cm);
    \draw[mid<] (0,0) circle (.4cm);}
    \,\right) 
    \right)= 4-2\sqrt{5},  
    \]
    \item $\cF(X_2) = \mathbf{1} \oplus 2\rho$ and $\cF(X_3) = \mathbf{1} \oplus 2\alpha\rho$, in which case ${}_\mathbf{1}\pi_\rho \cong 2\chi_1 \oplus  2\chi_2$ and ${}_\mathbf{1}\pi_{\alpha\rho}=2\chi_1 \oplus  2\chi_3$, and in particular
    \[  
    \Tr\left(  {}_\mathbf{1}\pi_\rho
    \left(\,
    \tikzmath{\filldraw[very thick, fill=gray!30] (0,0) circle (.2cm);
    \draw[mid<] (0,0) circle (.4cm);}
    \,\right)
    \right)= 6-2\sqrt{5},  
    \]
    \item $\cF(X_2) = \mathbf{1} \oplus 2\alpha\rho$ and $\cF(X_3) = \mathbf{1} \oplus 2\rho$, in which case ${}_\mathbf{1}\pi_\rho=2\chi_1 + 2\chi_3$ and ${}_\mathbf{1}\pi_{\alpha\rho}=2\chi_1 + 2\chi_2$, and in particular
    \[   
    \Tr\left(  {}_\mathbf{1}\pi_\rho
    \left(\,
    \tikzmath{\filldraw[very thick, fill=gray!30] (0,0) circle (.2cm);
    \draw[mid<] (0,0) circle (.4cm);}
    \,\right)
    \right)= 2-2\sqrt{5}. 
    \]
\end{enumerate}

We now aim to deduce more information about the operator $\psi$. 
Note that 
\[z_\mathbf{1} :=   \frac{1}{\dim(\cC)}\left(\tikzmath{\filldraw[very thick, fill=gray!30] (0,0) circle (.2cm);} +\tikzmath{\filldraw[very thick, fill=gray!30] (0,0) circle (.2cm);
    \draw[thick, red, mid<] (0,0) circle (.4cm);} + (2+\sqrt{5}) \tikzmath{\filldraw[very thick, fill=gray!30] (0,0) circle (.2cm);
    \draw[mid<] (0,0) circle (.4cm);} + (2+\sqrt{5})  \tikzmath{\filldraw[very thick, fill=gray!30] (0,0) circle (.2cm);
    \draw[red, thick, mid<] (0,0) circle (.65cm);
    \draw[mid<] (0,0) circle (.4cm);}\right)  \]
is the minimal central idempotent corresponding to the representation $\chi_0$. i.e.
$z_{\mathbf{1}}\cdot x = \chi_{0}(x) \cdot z_{\mathbf{1}}$.
As $\pi_\rho$ contains no copies of $\chi_0$, we get ${}_\mathbf{1}\pi_\rho
\left(  z_{\mathbf{1}}\right) =0$, and so
\[ \begin{bmatrix}
1 & 0\\
0 & 1
\end{bmatrix} +\left(2 + \sqrt{5}\right)(\phi + \phi') = 0
\qquad\Longrightarrow\qquad
\phi' = -\left(\phi +\begin{bmatrix}
\frac{1}{2+\sqrt{5}} & 0\\
0 & \frac{1}{2+\sqrt{5}}
\end{bmatrix} \right).
\]
To solve for $\phi$, we use the fusion $\rho^2 = \mathbf{1} \oplus 2\rho \oplus 2\alpha \rho$ to get
\[ 
\phi^2 +\left(\sqrt{5}-2\right)\phi  =\begin{bmatrix}
\frac{1}{2+\sqrt{5}} & 0\\
0 & \frac{1}{2+\sqrt{5}}
\end{bmatrix}. 
\]
Together with knowing the trace of $\pi_\rho(\rho)$ in each of the above cases, we can solve to get the statement of the lemma.

To obtain the statement about $\psi$ we repeat the above analysis with $_{\mathbf{1}}\pi_{\alpha\rho}$.
\end{proof}

\begin{rem}\label{rem:phi}
In the following subsection we are able to show that the operators $\phi$ and $\psi$ are diagonal. Thus the above lemma shows that there are only three possibilities for these operators:
\[  
\begin{bmatrix}
\frac{1-\sqrt{5}}{2} & 0\\
0 &\frac{1-\sqrt{5}}{2}
\end{bmatrix} 
\,, \qquad 
\begin{bmatrix}
\frac{3-\sqrt{5}}{2} & 0\\
0 &\frac{1-\sqrt{5}}{2}
\end{bmatrix}
\,,\qquad \text{ or }\qquad 
\begin{bmatrix}
\frac{3-\sqrt{5}}{2} & 0\\
0 &\frac{3-\sqrt{5}}{2}
\end{bmatrix} . 
\]
Hence the expressions $A^{0,0}_{0,0} + A^{1,0}_{1,0}$, $A^{0,1}_{0,1} + A^{1,1}_{1,1}$, $\widehat{A}^{0,0}_{0,0} + \widehat{A}^{1,0}_{1,0}$, and $\widehat{A}^{0,1}_{0,1} + \widehat{A}^{1,1}_{1,1}$ take values in $\left\{\frac{1-\sqrt{5}}{2}, \frac{3-\sqrt{5}}{2}    \right\}.$
    \end{rem}
We also aim to understand the restriction of the objects $Y_i$. 
For reasons that will become apparent soon, we only need to consider the case where $\chi_{\mathbf{1}, 1} = -\chi_{\mathbf{1}, 2}\quad \text{and} \quad \chi_{\alpha, 1} = -\chi_{\alpha, 2}$.

\begin{lem}\label{lem:SDalphaInd}
Suppose that $\chi_{\mathbf{1}, 1} = -\chi_{\mathbf{1}, 2}$ and $\chi_{\alpha, 1} = -\chi_{\alpha, 2}$. 
Then we have that
$$
\begin{aligned}
    \cF(Y_0) & = \alpha \oplus c_0 \rho \oplus (2-c_0)\alpha \rho\\
     \cF(Y_1) & = \alpha \oplus (2-c_0) \rho \oplus (c_0)\alpha \rho\\
    \cF(Y_2) & = \alpha \oplus c_2 \rho \oplus (2-c_2)\alpha \rho\\
     \cF(Y_3) & = \alpha \oplus (2-c_2) \rho \oplus c_2\alpha \rho.
\end{aligned}
\qquad\qquad
\text{where }c_0, c_2 \in \{0,1,2\}.
$$
\end{lem}
\begin{proof}
As $\chi_{\mathbf{1}, 1} = -\chi_{\mathbf{1}, 2}$ and $\chi_{\alpha, 1} = -\chi_{\alpha, 2}$, we have that $a=0$, and so each of the objects $Y_i$ has dimension $5 + 2\sqrt{5}$.
From the earlier computations
\[\Tr\left({_\alpha\pi_{\rho}}\left(\,
    \tikzmath{\filldraw[very thick, fill=gray!30] (0,0) circle (.2cm);
    \draw[red, thick, mid>] (0,.2) arc (180:0:.2cm) -- (.4,0) arc (0:-180:.4cm) .. controls ++(90:.3cm) and ++(270:.3cm) .. (0,.8);}\,
    \right)\right) = 0.\]
Therefore
\[{_\alpha\pi_{\rho}}    \cong c_0 \tau_0 \oplus (2-c_0)\tau_1 \oplus c_2\tau_2 \oplus (2-c_2)\tau_3.\]
Counting dimensions of the objects $Y_i$ gives the statement of the lemma.
\end{proof}

With the restrictions of the objects $X_i$ and $Y_i$ now understood, we can give a fairly explicit formula for the even Frobenius-Schur indicators of $\rho$. 
This formula will come in handy at several points later in this article.
Let $Z_i$ denote the remaining objects in $Z(\cC)$ beyond the $X_i,Y_i$, i.e., the objects in $Z(\cC)$ such that applying the forgetful functor, we have 
\[
\cF(Z_i) = p_i \rho \oplus q_i \alpha\rho.
\] 

\begin{lem}\label{lem:2nfrob}
Suppose that $\chi_{\mathbf{1}, 1} = -\chi_{\mathbf{1}, 2}$ and $\chi_{\alpha, 1} = -\chi_{\alpha, 2}$. 
Then we have that
\[ \nu_{2n}(\rho) \dim(\cC)=28+ 12 \sqrt{5} +\lambda_\alpha^n( 20 + 8\sqrt{5}) +    (2+ \sqrt{5})\sum p_i(p_i + q_i) \theta_{Z_i}^{2n},   \]
where $p_i$ and $q_i$ are integers satisfying 
$\sum p_i(p_i+q_i) = 16$,
and the $\theta_i$ are roots of unity.
\end{lem}
\begin{proof}
From Lemmas~\ref{lem:SDidInd} and \ref{lem:SDalphaInd}, we know the image under the forgetful functor of each of the simple objects appearing in $\cI(\mathbf{1})$ and $\cI(\alpha)$, up to some small integers $b, c_0,c_2$. 
Then we can write
\begin{align*}
    \cI(\rho) &= 2X_1 \oplus bX_2 \oplus (2-b)X_3 \oplus c_0 Y_0 \oplus (2-c_0)Y_1 \oplus c_2 Y_2 \oplus (2-c_2)Y_3 \oplus \bigoplus p_i Z_i\\
    \cI(\alpha\rho) &= 2X_1 \oplus (2-b)X_2 \oplus bX_3 \oplus (2-c_0) Y_0 \oplus c_0Y_1 \oplus (2-c_2) Y_2 \oplus c_2Y_3 \oplus \bigoplus q_i Z_i.
\end{align*}
Using the fact that $\cF(\cI(\rho)) \cong \bigoplus_{X \in \Irr(\cC)} X \rho X^*$ we obtain
\begin{align*}
    20 &= \dim\Hom(\cI(\rho), \cI(\rho)) = 4 + b^2 + (2-b)^2+ c_0^2  + (2-c_0)^2 + c_2^2  + (2-c_2)^2 + \sum  p_i^2\\
     16 &= \dim\Hom(\cI(\alpha\rho), \cI(\rho)) = 4 + 2b(2-b)+ 2c_0(2-c_0) + 2c_2(2-c_2) + \sum  p_iq_i,
\end{align*}
so $36 = \dim\Hom(\cI(\rho\oplus \alpha\rho), \cI(\rho)) = 20 + \sum  p_i^2+ \sum  p_iq_i$, and thus $\sum p_i(p_i+q_i) = 16$.

For a simple $W\in Z(\cC)$, we write $\theta_{W}$ for its twist. 
We have $\theta_{X_i}^2 =1$ and $ \theta_{Y_i}^2 = \lambda_\alpha$. 
We can use the Ng-Shauenburg formula for the 2n-th Frobenius-Schur indicator \cite[Theorem 4.1]{MR2313527} 
to obtain
\begin{align*}
\nu_{2n}(\rho)\dim(\cC) 
&= 
\sum_{W\in Z(\cC)} \dim\cC(\cF(W)\to \rho) \dim(W) \theta_W^{2n} 
\\&=
28+ 12 \sqrt{5} +\lambda_\alpha^n( 20 + 8\sqrt{5}) +   (2+ \sqrt{5})\sum p_i(p_i + q_i) \theta_{Z_i}^{2n} .  
\qedhere
\end{align*}
\end{proof}

We now have all the required pieces to prove the main result of this subsection.
\begin{proof}[Proof of Lemma~\ref{lem:maincentre}]
\item[\underline{Case 1:}]

First suppose that $\lambda_\alpha = \mu$.

If $a = \pm 2$, then two of the objects in $\cI(\alpha)$ would have dimension $\frac{2 \sqrt{5}+5}{\sqrt{2}+2}$, which is impossible.

If $a = \pm 4$, then two of the objects in $\cI(\alpha)$ would have dimension $1$, and so $\alpha$ would have two lifts to the centre. This would imply that $\mathbf{1}$ also has two lifts to the centre by \cite[Cor.~3.7]{MR3354332}, which contradicts our knowledge of $\cI(\mathbf{1})$.

If $a=0$, then we have $\chi_{\mathbf{1}, 0} = - \chi_{\mathbf{1}, 1}$ and $\chi_{\alpha, 0} = - \chi_{\alpha, 1}$. From Lemma~\ref{lem:FrobBasis} we have
\[ \chi_{\mathbf{1}, 0} = \lambda_\alpha\mu \chi_{\mathbf{1}, \tilde{0}} = \chi_{\mathbf{1}, \tilde{0}}\quad \text{and} \quad \chi_{\alpha, 0} = \lambda_\alpha\mu \chi_{\alpha, \tilde{0}} = \chi_{\alpha, \tilde{0}}.\]
Therefore $(\mathbf{1}, \tilde{0})=(\mathbf{1}, 0)$ and $(\alpha, \tilde{0})=(\alpha, 0)$. 
If $\mu = -1$, then by exchanging $\rho$ and $\alpha\rho$ we can arrange so that $\lambda_\rho = -1$. 
However Lemma~\ref{lem:FrobBasis} then implies that $(\mathbf{1}, \tilde{0})=(\mathbf{1}, 1)$ which gives a contradiction. 
Thus $\mu= 1$ in this case.

Summarising, if $\lambda_\alpha = \mu$, we can only have $\lambda_\alpha = \mu = 1$, in which case
$\chi_{\mathbf{1}, 0} = -\chi_{\mathbf{1}, 1}$, $\chi_{\alpha, 0} = -\chi_{\alpha, 1}$, $(\mathbf{1}, \tilde{0})=(\mathbf{1}, 0)$, and $(\alpha, \tilde{0})=(\alpha, 0)$.

\item[\underline{Case 2:}]
If $\lambda_\alpha = -\mu$, then equations $\chi_{\mathbf{1}, 0} = \lambda_\alpha\mu \chi_{\mathbf{1}, \tilde{0}} =-\chi_{\mathbf{1}, \tilde{0}} $ and $\chi_{\alpha, 0} = \lambda_\alpha\mu \chi_{\alpha, \tilde{0}}=-\chi_{\alpha, \tilde{0}}$ from Lemma~\ref{lem:FrobBasis} imply that $(\mathbf{1}, \tilde{0})=(\mathbf{1}, 1)$ and $(\alpha, \tilde{0})=(\alpha, 1)$, and hence $\chi_{\mathbf{1}, 0} = -\chi_{\mathbf{1}, 1}$ and $\chi_{\alpha, 0} = -\chi_{\alpha, 1}$. 
Thus all that remains to be shown is that $\mu = 1$. 
Suppose for contradiction that $\mu =-1$, which implies that $\lambda_\alpha=1$. 
As $\mu= -1$, we can exchange $\rho$ and $\alpha\rho$ if necessary to arrange $\lambda_\rho = -1$.

We can now use Lemma~\ref{lem:2nfrob}, along with the fact that the 2nd Frobenius-Schur indicator of $\rho$ is $\lambda_\rho$ to get the equation 
\[-20 - 8\sqrt{5}= 48+ 20 \sqrt{5}+ (2+ \sqrt{5})\sum p_i(p_i + q_i) \theta_{Z_i}^2 ,  \]
where $\sum p_i(p_i+q_i) = 16. $
Thus
\[\sum p_i(p_i + q_i) \theta_{Z_i}^2 = -4 - 12\sqrt{5}.   \]
However Theorem~\ref{thm:LarsonRootsOfUnityAB} implies that it takes at least 32 roots of unity to write $ -4 - 12\sqrt{5}$, contradicting $\sum p_i(p_i+q_i) = 16$. Hence $\mu = 1$.

Summarising, if $\lambda_\alpha = -\mu$, then $\lambda_\alpha = -1$ and $\mu = 1$, in which case
$\chi_{\mathbf{1}, 1} = -\chi_{\mathbf{1}, 2}$, $\chi_{\alpha, 1} = -\chi_{\alpha, 2}$, $(\mathbf{1}, \tilde{0})=(\mathbf{1}, 1)$, and $(\alpha, \tilde{0})=(\alpha, 1)$.

In either case, we have that $\mu = 1$, and 
$\chi_{\mathbf{1}, 1} = -\chi_{\mathbf{1}, 2}$ and $\chi_{\alpha, 1} = -\chi_{\alpha, 2}$.
As $ \chi_{\mathbf{1}, i} \in \{\pm 1\}$ and  $\chi_{\alpha, i} \in \{\pm \sqrt{\lambda_\alpha}\}$ we can swap our basis elements to arrange that
\[
\chi_{\mathbf{1}, 0} = \sqrt{\lambda_\alpha}, \quad  \chi_{\mathbf{1}, 1} = -\sqrt{\lambda_\alpha} ,\quad   \chi_{\alpha, 0} =1, \quad \text{and} \quad  \chi_{\alpha, 1}= -1.
\qedhere
\]
\end{proof}

%%%%%%%%%%%%%%%%%%%%%%%%%%%%%%%%%%%%%%%%%%%%%%%%%%%%%%%%%%%%%%%%%%%
\subsection{Symmetries}
\label{sec:SDSymmetries}

With the results of the last subsection in hand, the major task in front of us is to determine the $8m^4$ complex scalars:
\[ 
A^{i,j}_{k,\ell},\quad  B^{i,j}_{k,\ell} ,\quad C^{i,j}_{k,\ell},\quad D^{i,j}_{k,\ell} ,\quad \widehat{A}^{i,j}_{k,\ell}, \quad \widehat{B}^{i,j}_{k,\ell} , \quad\widehat{C}^{i,j}_{k,\ell},\quad\text{and} \quad  \widehat{D}^{i,j}_{k,\ell}.    
\]
In theory, we could begin evaluating diagrams in our category in multiple ways in order to obtain equations of these variables. 
However, in practice this task is too complicated, especially in the $m=2$ case where we have $128$ unknowns. 
To make our task of pinning down these scalars easier, we aim to find symmetries between them, and to show that many of them must in fact vanish.
The symmetries of these scalars come from the tetrahedral symmetries of the $6j$ symbols, which were rigorously studied in \cite{MR1657800}, and have been used in previous works of the second author \cite{MR3635673,MR3827808}.
(See also Footnote \ref{footnote:TetrahedralSymmetry}.)

The main result of this subsection is as follows.
\begin{lem}\label{lem:mainSym}
The scalars 
$B^{i,j}_{k,\ell} , C^{i,j}_{k,\ell}, \widehat{B}^{i,j}_{k,\ell} , \widehat{C}^{i,j}_{k,\ell}, \widehat{D}^{i,j}_{k,\ell}$
can be expressed in terms of the $D^{i,j}_{k,\ell}$ as:
\begin{align*}
 B^{i,j}_{k,\ell} &= \lambda_\rho^{1+i+k}\lambda_\alpha \sqrt{\lambda_\alpha}(-1)^\ell \omega_{\mathbf{1}, \ell}D^{j,\tilde{k}}_{\tilde{i},\ell},\qquad &&\widehat{B}^{i,j}_{k,\ell}= \lambda_\rho^{1+j+\ell}\lambda_\alpha \sqrt{\lambda_\alpha}(-1)^i \omega_{\mathbf{1}, j}\omega_{\mathbf{1}, i}^2 \omega_{\alpha, \ell}D^{\tilde{\ell},i}_{k,\tilde{j}},\\
    C^{i,j}_{k,\ell} &= \lambda_\rho^{1+j+k}\lambda_\alpha \sqrt{\lambda_\alpha}(-1)^\ell \omega_{\mathbf{1}, \ell}^2D^{\tilde{k},i}_{\tilde{j},\ell},\qquad   &&\widehat{C}^{i,j}_{k,\ell} =\lambda_\rho^{1+j+k}\sqrt{\lambda_\alpha}(-1)^k \omega_{\mathbf{1}, k}\omega_{\mathbf{1}, j}^2 \omega_{\alpha, i}D^{i,\tilde{k}}_{\ell,\tilde{j}},\\
     & &&   \widehat{D}^{i,j}_{k,\ell} = \sqrt{\lambda_\alpha}(-1)^{k+j} \omega_{\mathbf{1}, i}\omega_{\mathbf{1}, k}^2 \omega_{\alpha, \ell}\omega_{\alpha, j}^2D^{j,i}_{\ell,k}.
\end{align*}
The scalars $A^{i,j}_{k,\ell}$ and $\widehat{A}^{i,j}_{k,\ell}$ satisfy $S_4$ symmetries generated by the order three rotation:
\begin{align*}
 A^{i,j}_{k,\ell} &=\lambda_\rho^{1+i+k} \omega_{\mathbf{1}, \ell}A^{j,\tilde{k}}_{\tilde{i},\ell}=\lambda_\rho^{1+j+k} \omega_{\mathbf{1}, \ell}^2A^{\tilde{k},i}_{\tilde{j},\ell}
 &
  \widehat{A}^{i,j}_{k,\ell} &= \lambda_\rho^{1+i+k}\omega_{\alpha, \ell}\widehat{A}^{j,\tilde{k}}_{\tilde{i},\ell}=\lambda_\rho^{1+j+k} \omega_{\alpha, \ell}^2\widehat{A}^{\tilde{k},i}_{\tilde{j},\ell}
\end{align*}
and the order two flips:
\begin{align*}
     A^{i,j}_{k,\ell} &= \omega_{\mathbf{1}, k}\omega_{\mathbf{1}, i}^2 A^{\tilde{k},\tilde{\ell}}_{\tilde{i},\tilde{j}} = \lambda_\rho^{j+\ell}\overline{A^{k,\tilde{j}}_{i,\tilde{\ell}}}
     &
     \widehat{A}^{i,j}_{k,\ell} &= \omega_{\alpha, k}\omega_{\alpha, i}^2 \widehat{A}^{\tilde{k},\tilde{\ell}}_{\tilde{i},\tilde{j}} = \lambda_\rho^{j+\ell}\overline{\widehat{A}^{k,\tilde{j}}_{i,\tilde{\ell}}}.
\end{align*}
The $D^{i,j}_{k,\ell}$ scalars satisfy the $\bbZ/2\bbZ \times \bbZ/2\bbZ$ symmetries generated by:
\[ D^{i,j}_{k,\ell}=   \lambda_\rho^{j+l}\overline{ D^{k,\tilde{j}}_{i,\tilde{\ell}}} = \lambda_\alpha(-1)^{j+\ell}\omega_{\alpha, k}\omega_{\alpha, i}^2D^{\tilde{k},\tilde{\ell}}_{\tilde{i},\tilde{j}}= \lambda_\rho^{i+k}\lambda_\alpha(-1)^{j+\ell}\omega_{\alpha, k}\omega_{\alpha, i}^2\overline{D^{\tilde{i},\ell}_{\tilde{k},j}}.\]
Finally, we have
\[ 
A^{i,j}_{k,\ell} = \widehat{A}^{i,j}_{k,\ell} = D^{i,j}_{k,\ell} = 0 \quad \text{if}\quad i+j+k+\ell \not\equiv 0 \pmod m.    
\]
\end{lem}

When $m=2$, this result reduces the number of complex scalars to solve for down to $11$ in the $\lambda_\alpha = 1$ case, and $7$ in the $\lambda_\alpha = -1$ case. 
This simplification makes it feasible to solve for these scalars in the next subsection.

\begin{lem}
We have
\[ A^{i,j}_{k,\ell} = \widehat{A}^{i,j}_{k,\ell} = D^{i,j}_{k,\ell} = 0 \quad \text{if}\quad i+j+k+\ell \not\equiv 0 \pmod 2.\]
\end{lem}
\begin{proof}
Note that if $m=1$, then the statement of the lemma is trivially true. 
Thus it suffices to restrict our attention to $m=2$.
We will prove the statement of the lemma in the case of the $A^{i,j}_{k,\ell}$ coefficients, as the remaining two cases are nearly identical.
We have
\[
A^{i,j}_{k,\ell}
\,
\tikzmath{
\draw[thick, red, mid>] (-.2,-.5) -- (-.2,.5);
\draw[mid>] (.2,-.5) -- (.2,.5);
}
=
\tikzmath{
\draw[thick, red, mid>] (-.6,-1) -- (-.6,1);
\coordinate (b) at (-.2,-.2);
\coordinate (a) at (.2,-.6);
\coordinate (d) at (-.2,.6);
\coordinate (c) at (.2,.2);
\draw[mid>] (d) -- ($ (d) + (0,.4) $);
\draw[mid<] (a) -- ($ (a) + (0,-.4) $);
\draw[mid>] (a) to (b);
\draw[mid>] (b) to (c);
\draw[mid>] (c) to (d);
\draw[mid>] (b) to[out=135,in=-135] (d);
\draw[mid>] (a) to[out=45,in=-45] (c);
\filldraw (a) node[left]{$\scriptstyle j$} circle (.05cm);
\filldraw (b) node[left]{$\scriptstyle i$} circle (.05cm);
\filldraw (c) node[left]{$\scriptstyle \ell$}circle (.05cm);
\filldraw (d) node[left]{$\scriptstyle k$} circle (.05cm);
}
=
\frac{\chi_{\mathbf{1},k}\chi_{\mathbf{1},\ell}}{\chi_{\mathbf{1},i}\chi_{\mathbf{1},j}}
\,
\tikzmath{
\draw[thick, red] (-.6,-1.2) .. controls ++(90:.4cm) and ++(-90:1.2cm) .. (.7,0) .. controls ++(90:1cm) and ++(-90:.4cm) .. (-.6,1.2);
\coordinate (b) at (-.2,-.2);
\coordinate (a) at (.2,-.6);
\coordinate (d) at (-.2,.6);
\coordinate (c) at (.2,.2);
\draw[mid>] (d) -- ($ (d) + (0,.6) $);
\draw[mid<] (a) -- ($ (a) + (0,-.6) $);
\draw[mid>] (a) to (b);
\draw[mid>] (b) to (c);
\draw[mid>] (c) to (d);
\draw[mid>] (b) to[out=135,in=-135] (d);
\draw[mid>] (a) to[out=45,in=-45] (c);
\filldraw (a) node[left]{$\scriptstyle j$} circle (.05cm);
\filldraw (b) node[left]{$\scriptstyle i$} circle (.05cm);
\filldraw (c) node[left]{$\scriptstyle \ell$}circle (.05cm);
\filldraw (d) node[left]{$\scriptstyle k$} circle (.05cm);
}
= 
\frac{\chi_{\mathbf{1},k}\chi_{\mathbf{1},\ell}}{\chi_{\mathbf{1},i}\chi_{\mathbf{1},j}}A^{i,j}_{k,\ell}
\,
\tikzmath{
\draw[thick, red, mid>] (-.2,-.5) -- (-.2,.5);
\draw[mid>] (.2,-.5) -- (.2,.5);
}
\qquad
\Longrightarrow
\qquad
A^{i,j}_{k,\ell} =\frac{\chi_{\mathbf{1},k}\chi_{\mathbf{1},\ell}}{\chi_{\mathbf{1},i}\chi_{\mathbf{1},j}}A^{i,j}_{k,\ell}.    \]
Recall from Lemma~\ref{lem:maincentre} of the previous section that $\chi_{\mathbf{1}, i} = (-1)^{ i  }\sqrt{\lambda_\alpha}$. 
Thus if $i+j+k+\ell \not\equiv 0 \pmod 2$, then $\frac{\chi_{\mathbf{1},k}\chi_{\mathbf{1},\ell}}{\chi_{\mathbf{1},i}\chi_{\mathbf{1},j}} \neq 1$, which implies $A^{i,j}_{k,\ell}=0$.
\end{proof}

Now that we know that half of our coefficients vanish, we move on to describing the symmetries between them. 
As mentioned before, these symmetries are the standard tetrahedral symmetries of the 6j-symbols. 
This completes the proof of the main result of this section.

\begin{proof}[Proof of Lemma~\ref{lem:mainSym}]
We include enough examples to illuminate the necessary techniques, all of which involve using the Frobenius maps defined in \S\ref{sec:SDBasicRelations}.
The symmetries of the $A^{i,j}_{k,\ell}$ coefficients are the easiest, as the diagrams only involve $\rho$ strands. 
We compute the following symmetries:
\begin{footnotesize}
\begin{align*}
A^{i,j}_{k,\ell} 
&=
\frac{1}{d}
\,
\tikzmath{
\coordinate (b) at (-.2,-.2);
\coordinate (a) at (.2,-.6);
\coordinate (d) at (-.2,.6);
\coordinate (c) at (.2,.2);
\draw[mid>] (d) .. controls ++(90:.7cm) and ++(90:.5cm) .. (.8,0) .. controls ++(-90:.5cm) and ++(-90:.7cm) .. (a);
\draw[mid>] (a) to (b);
\draw[mid>] (b) to (c);
\draw[mid>] (c) to (d);
\draw[mid>] (b) to[out=135,in=-135] (d);
\draw[mid>] (a) to[out=45,in=-45] (c);
\filldraw (a) node[left]{$\scriptstyle j$} circle (.05cm);
\filldraw (b) node[left]{$\scriptstyle i$} circle (.05cm);
\filldraw (c) node[left]{$\scriptstyle \ell$}circle (.05cm);
\filldraw (d) node[left]{$\scriptstyle k$} circle (.05cm);
}
=
\frac{\lambda_\rho^{\ell+1}}{d}
\tikzmath{
\coordinate (t) at (.5,.5);
\coordinate (b) at (-.2,-.2);
\coordinate (a) at (.2,-.6);
\coordinate (d) at (-.2,.6);
\coordinate (c) at (.2,.2);
\draw[mid>] (d) .. controls ++(90:1cm) and ++(90:1.5cm) .. (1.2,0) .. controls ++(-90:1cm) and ++(-90:.7cm) .. (a);
\draw[mid>] (a) to (b);
\draw[mid>] (b) to (c);
\draw[mid>] (c) to (d);
\draw[mid>] (b) to[out=135,in=-135] (d);
\draw[mid>] (a) to[out=45,in=-90] ($ (t) + (.4,0) $);
\draw[mid>] (c) to[out=45,in=-90] (t);
\draw (t) arc (180:0:.2cm);
\draw[thick] (t) -- ($ (t) + (.1,0) $);
\filldraw (a) node[left]{$\scriptstyle j$} circle (.05cm);
\filldraw (b) node[left]{$\scriptstyle i$} circle (.05cm);
\filldraw (c) node[left]{$\scriptstyle \widetilde{\ell}$}circle (.05cm);
\filldraw (d) node[left]{$\scriptstyle k$} circle (.05cm);
}
=
\frac{\lambda_\rho^{\ell}}{d}
\tikzmath{
\coordinate (t) at (.5,.9);
\coordinate (b) at (.2,-.2);
\coordinate (a) at (-.2,-.6);
\coordinate (d) at (.2,.6);
\coordinate (c) at (-.2,.2);
\draw[mid>] (d) .. controls ++(90:1.5cm) and ++(90:1.5cm) .. (1.2,0) .. controls ++(-90:1cm) and ++(-90:.7cm) .. (a);
\draw[mid>] (a) to (b);
\draw[mid>] (b) to (c);
\draw[mid>] (c) to (d);
\draw[mid>] (b) to[out=45,in=-90] ($ (t) + (.4,0)$);
\draw (t) arc (180:0:.2cm);
\draw[mid>] (a) to[out=135,in=-135] (c);
\draw[mid>] (d) to[out=45,in=-90] (t);
\draw[thick] (t) -- ($ (t) + (.1,0) $);
\filldraw (a) node[left]{$\scriptstyle i$} circle (.05cm);
\filldraw (b) node[left, xshift=-.1cm]{$\scriptstyle \widetilde{\ell}$} circle (.05cm);
\filldraw (c) node[left]{$\scriptstyle k$}circle (.05cm);
\filldraw (d) node[left]{$\scriptstyle j$} circle (.05cm);
}
=
\frac{\lambda_\rho^{\ell+j}}{d}
\tikzmath{
\coordinate (b) at (.2,-.2);
\coordinate (a) at (-.2,-.6);
\coordinate (d) at (.2,.6);
\coordinate (c) at (-.2,.2);
\draw[mid>] (d) .. controls ++(90:.7cm) and ++(90:.5cm) .. (.8,0) .. controls ++(-90:.5cm) and ++(-90:.7cm) .. (a);
\draw[mid>] (a) to (b);
\draw[mid>] (b) to (c);
\draw[mid>] (c) to (d);
\draw[mid>] (b) to[out=45,in=-45] (d);
\draw[mid>] (a) to[out=135,in=-135] (c);
\filldraw (a) node[left]{$\scriptstyle i$} circle (.05cm);
\filldraw (b) node[left, xshift=-.1cm]{$\scriptstyle \widetilde{\ell}$} circle (.05cm);
\filldraw (c) node[left]{$\scriptstyle k$}circle (.05cm);
\filldraw (d) node[left]{$\scriptstyle \widetilde{j}$} circle (.05cm);
}
=
\lambda_\rho^{\ell+j}
\overline{A^{k,\tilde{j}}_{i,\tilde{\ell}}}, 
\\   
A^{i,j}_{k,\ell}
&= 
\frac{1}{d}\,
\tikzmath{
\coordinate (b) at (-.2,-.2);
\coordinate (a) at (.2,-.6);
\coordinate (d) at (-.2,.6);
\coordinate (c) at (.2,.2);
\draw[mid>] (d) .. controls ++(90:.7cm) and ++(90:.5cm) .. (.8,0) .. controls ++(-90:.5cm) and ++(-90:.7cm) .. (a);
\draw[mid>] (a) to (b);
\draw[mid>] (b) to (c);
\draw[mid>] (c) to (d);
\draw[mid>] (b) to[out=135,in=-135] (d);
\draw[mid>] (a) to[out=45,in=-45] (c);
\filldraw (a) node[left]{$\scriptstyle j$} circle (.05cm);
\filldraw (b) node[left]{$\scriptstyle i$} circle (.05cm);
\filldraw (c) node[left]{$\scriptstyle \ell$}circle (.05cm);
\filldraw (d) node[left]{$\scriptstyle k$} circle (.05cm);
}
=
\frac{\lambda_\rho^{1+\ell}}{d} \omega_{\mathbf{1}, \ell} 
\tikzmath{
\coordinate (t) at (-.2,.5);
\coordinate (b) at (-.4,-.2);
\coordinate (a) at (.2,-.6);
\coordinate (d) at (-.2,1.2);
\coordinate (c) at (.4,.2);
\draw[mid>] (d) .. controls ++(90:1cm) and ++(90:1.5cm) .. (.8,0) .. controls ++(-90:1cm) and ++(-90:.7cm) .. (a);
\draw[mid>] (a) to (b);
\draw[mid>] (a) to (c);
\draw[mid>] (c) to[out=45,in=-45] (d);
\draw[mid>] (b) to[out=135,in=-135] (d);
\draw[mid>] (c) to[out=135,in=-90] ($ (t) + (.4,0) $);
\draw[mid>] (b) to[out=45,in=-90] (t);
\draw (t) arc (180:0:.2cm);
\draw[thick] (t) -- ($ (t) + (.1,0) $);
\filldraw (a) node[left]{$\scriptstyle j$} circle (.05cm);
\filldraw (b) node[left]{$\scriptstyle i$} circle (.05cm);
\filldraw (c) node[left, xshift=-.1cm]{$\scriptstyle \widetilde{\ell}$}circle (.05cm);
\filldraw (d) node[left]{$\scriptstyle k$} circle (.05cm);
}
=
\frac{\lambda_\rho^{1+\ell}}{d} \omega_{\mathbf{1},\ell} 
\tikzmath{
\coordinate (t) at (-.2,.5);
\coordinate (b) at (-.4,.2);
\coordinate (a) at (.2,-.6);
\coordinate (d) at (-.2,1.2);
\coordinate (c) at (.4,-.2);
\draw[mid>] (d) .. controls ++(90:1cm) and ++(90:1.5cm) .. (.8,0) .. controls ++(-90:1cm) and ++(-90:.7cm) .. (a);
\draw[mid>] (a) to (b);
\draw[mid>] (a) to (c);
\draw[mid>] (c) to[out=45,in=-45] (d);
\draw[mid>] (b) to[out=135,in=-135] (d);
\draw[mid>] (c) to[out=135,in=-90] ($ (t) + (.4,0) $);
\draw[mid>] (b) to[out=45,in=-90] (t);
\draw (t) arc (180:0:.2cm);
\draw[thick] (t) -- ($ (t) + (.1,0) $);
\filldraw (a) node[left]{$\scriptstyle j$} circle (.05cm);
\filldraw (b) node[left]{$\scriptstyle i$} circle (.05cm);
\filldraw (c) node[left, xshift=-.1cm]{$\scriptstyle \widetilde{\ell}$}circle (.05cm);
\filldraw (d) node[left]{$\scriptstyle k$} circle (.05cm);
}
=
\frac{\lambda_\rho^{1+i+\ell}}{d} \omega_{\mathbf{1}, \ell} 
\tikzmath{
\coordinate (b) at (.2,-.2);
\coordinate (a) at (-.2,-.6);
\coordinate (d) at (.2,.6);
\coordinate (c) at (-.2,.2);
\draw[mid>] (d) .. controls ++(90:.7cm) and ++(90:.5cm) .. (.8,0) .. controls ++(-90:.5cm) and ++(-90:.7cm) .. (a);
\draw[mid>] (a) to (b);
\draw[mid>] (b) to (c);
\draw[mid>] (c) to (d);
\draw[mid>] (b) to[out=45,in=-45] (d);
\draw[mid>] (a) to[out=135,in=-135] (c);
\filldraw (a) node[left]{$\scriptstyle j$} circle (.05cm);
\filldraw (b) node[left, xshift=-.1cm]{$\scriptstyle \widetilde{\ell}$} circle (.05cm);
\filldraw (c) node[left]{$\scriptstyle \widetilde{i}$}circle (.05cm);
\filldraw (d) node[left]{$\scriptstyle k$} circle (.05cm);
}
=
\lambda_\rho^{1+i+\ell}\omega_{\mathbf{1}, \ell}
\overline{A^{\tilde{i},k}_{j,\tilde{\ell}}} ,
\\
A^{i,j}_{k,\ell} 
&= 
\frac{1}{d}\,
\tikzmath{
\coordinate (b) at (-.2,-.2);
\coordinate (a) at (.2,-.6);
\coordinate (d) at (-.2,.6);
\coordinate (c) at (.2,.2);
\draw[mid>] (d) .. controls ++(90:.7cm) and ++(90:.5cm) .. (-.8,0) .. controls ++(-90:.5cm) and ++(-90:.7cm) .. (a);
\draw[mid>] (a) to (b);
\draw[mid>] (b) to (c);
\draw[mid>] (c) to (d);
\draw[mid>] (b) to[out=135,in=-135] (d);
\draw[mid>] (a) to[out=45,in=-45] (c);
\filldraw (a) node[left]{$\scriptstyle j$} circle (.05cm);
\filldraw (b) node[left]{$\scriptstyle i$} circle (.05cm);
\filldraw (c) node[left]{$\scriptstyle \ell$}circle (.05cm);
\filldraw (d) node[left]{$\scriptstyle k$} circle (.05cm);
}
=
\frac{1}{d} 
\,
\tikzmath{
\coordinate (t1) at (-.6,.6);
\coordinate (t2) at (-1,.6);
\coordinate (b) at (-.2,-.2);
\coordinate (a) at (.2,-.6);
\coordinate (d) at (-.2,.6);
\coordinate (c) at (.2,.2);
\draw[mid>] (d) .. controls ++(90:.7cm) and ++(90:1cm) .. (-1.2,.4) .. controls ++(-90:1cm) and ++(-90:.7cm) .. (a);
\draw[mid>] (a) to (b);
\draw[mid>] (b) to (c);
\draw[mid>] (c) to (d);
\draw[mid>] (a) to[out=45,in=-45] (c);
\draw[mid>] (b) to[out=135, in=-90] (t2);
\draw[thick] (t1) -- ($ (t1) + (.1,0) $);
\draw[thick] (t2) -- ($ (t2) - (.1,0) $);
\draw (t2) arc (180:0:.2cm);
\draw[mid>] (t1) arc (-180:0:.2cm);
\filldraw (a) node[left]{$\scriptstyle j$} circle (.05cm);
\filldraw (b) node[left]{$\scriptstyle i$} circle (.05cm);
\filldraw (c) node[left]{$\scriptstyle \ell$}circle (.05cm);
\filldraw (d) node[left]{$\scriptstyle k$} circle (.05cm);
}
=
\frac{\lambda_\rho^{k}}{d} \omega_{\mathbf{1}, k}
\tikzmath{
\coordinate (t) at (-.6,-.2);
\coordinate (b) at (-.2,-.2);
\coordinate (a) at (.2,-.6);
\coordinate (d) at (-.2,.6);
\coordinate (c) at (.2,.2);
\draw[mid>] (d) .. controls ++(45:.7cm) and ++(90:1cm) .. (-1.2,.4) .. controls ++(-90:1cm) and ++(-90:.7cm) .. (a);
\draw[mid>] (a) to (b);
\draw[mid>] (b) to (c);
\draw[mid>] (c) to (d);
\draw[mid>] (a) to[out=45,in=-45] (c);
\draw[thick] (t) -- ($ (t) + (.1,0) $);
\draw[mid>] (d) .. controls ++(135:.7cm) and ++(90:1cm) .. (-.9,0) .. controls ++(-90:.4cm) and ++(-90:.3cm) .. (t);
\draw[mid<] (t) arc (180:0:.2cm);
\filldraw (a) node[left]{$\scriptstyle j$} circle (.05cm);
\filldraw (b) node[left]{$\scriptstyle i$} circle (.05cm);
\filldraw (c) node[left]{$\scriptstyle \ell$}circle (.05cm);
\filldraw (d) node[left]{$\scriptstyle \widetilde{k}$} circle (.05cm);
}
=
\frac{\lambda_\rho^{i+k}}{d}\frac{\omega_{\mathbf{1}, k}}{\omega_{\mathbf{1}, i}}
\tikzmath{
\coordinate (b) at (-.2,-.2);
\coordinate (a) at (.2,-.6);
\coordinate (d) at (-.2,.6);
\coordinate (c) at (.2,.2);
\draw[mid>] (d) .. controls ++(45:.7cm) and ++(90:1cm) .. (-1.2,.4) .. controls ++(-90:1cm) and ++(-90:.7cm) .. (a);
\draw[mid>] (a) to (b);
\draw[mid>] (b) to (c);
\draw[mid>] (c) to (d);
\draw[mid>] (a) to[out=45,in=-45] (c);
\draw[mid>] (d) .. controls ++(135:.7cm) and ++(90:1cm) .. (-.9,0) .. controls ++(-90:.4cm) and ++(-135:.5cm) .. (b);
\filldraw (a) node[left]{$\scriptstyle j$} circle (.05cm);
\filldraw (b) node[left, yshift=.1cm]{$\scriptstyle \widetilde{i}$} circle (.05cm);
\filldraw (c) node[left]{$\scriptstyle \ell$}circle (.05cm);
\filldraw (d) node[left]{$\scriptstyle \widetilde{k}$} circle (.05cm);
}
=
\frac{\lambda_\rho^{i+k}}{d}\frac{\omega_{\mathbf{1}, k}}{\omega_{\mathbf{1}, i}}
\,
\tikzmath{
\coordinate (b) at (.2,-.2);
\coordinate (a) at (-.2,-.6);
\coordinate (d) at (.2,.6);
\coordinate (c) at (-.2,.2);
\draw[mid>] (d) .. controls ++(90:.7cm) and ++(90:.5cm) .. (-.8,0) .. controls ++(-90:.5cm) and ++(-90:.7cm) .. (a);
\draw[mid>] (a) to (b);
\draw[mid>] (b) to (c);
\draw[mid>] (c) to (d);
\draw[mid>] (b) to[out=45,in=-45] (d);
\draw[mid>] (a) to[out=135,in=-135] (c);
\filldraw (a) node[left]{$\scriptstyle \widetilde{k}$} circle (.05cm);
\filldraw (b) node[left]{$\scriptstyle j$} circle (.05cm);
\filldraw (c) node[left]{$\scriptstyle \widetilde{i}$}circle (.05cm);
\filldraw (d) node[left]{$\scriptstyle \ell$} circle (.05cm);
}
=
\lambda_\rho^{i+k}\frac{\omega_{\mathbf{1}, k}}{\omega_{\mathbf{1}, i}}
\overline{A^{\tilde{i},\ell}_{\tilde{k},j}} 
\end{align*}
\end{footnotesize}

Together this shows that 
\[A^{i,j}_{k,\ell} =\lambda_\rho^{i+k}\omega_{\mathbf{1}, \ell}A^{j,\tilde{k}}_{\tilde{i},\ell}=\lambda_\rho^{1+j+k}\omega_{\mathbf{1}, \ell}^2A^{\tilde{k},i}_{\tilde{j},\ell}\quad \text{and}\quad A^{i,j}_{k,\ell} = \lambda_\rho^{i+j+k+\ell}\frac{\omega_{\mathbf{1}, i}}{\omega_{\mathbf{1}, k}}A^{\tilde{k},\tilde{\ell}}_{\tilde{i},\tilde{j}} = \frac{\omega_{\mathbf{1}, i}}{\omega_{\mathbf{1}, k}}A^{\tilde{k},\tilde{\ell}}_{\tilde{i},\tilde{j}}\]
as claimed.
These three tricks work to determine all of the symmetries in the statement of the lemma. In order to show how to deal with $\alpha$ strands, we include one final example.
\begin{footnotesize}
\[ 
\widehat{A}^{i,j}_{k,\ell} 
=
\frac{1}{d}
\tikzmath{
\coordinate (b) at (-.2,-.2);
\coordinate (a) at (.2,-.6);
\coordinate (d) at (-.2,.6);
\coordinate (c) at (.2,.2);
\draw[mid>] (d) .. controls ++(45:.7cm) and ++(90:.5cm) .. (.8,0) .. controls ++(-90:.5cm) and ++(-45:.7cm) .. (a);
\draw[thick, red, mid>] (d) .. controls ++(135:.7cm) and ++(90:1.5cm) .. (1,0) .. controls ++(-90:1.5cm) and ++(-135:.7cm) .. (a);
\draw[mid>] (a) to (b);
\draw[mid>] (b) to[out=45, in=-45] (c);
\draw[mid>] (c) to[out=45, in=-45] (d);
\draw[mid>] (b) to[out=135,in=-135] (d);
\draw[mid>] (a) to[out=45,in=-45] (c);
\draw[thick, red] (b) .. controls ++(-135:.3cm) and ++(-90:.5cm) .. (0,-.2) .. controls ++(90:.3cm) and ++(-90:.3cm) .. (-.6,.2) .. controls ++(90:.3cm) and ++(135:.3cm) .. (c);
\filldraw (a) node[left]{$\scriptstyle j$} circle (.05cm);
\filldraw (b) node[left]{$\scriptstyle i$} circle (.05cm);
\filldraw (c) node[left]{$\scriptstyle \ell$}circle (.05cm);
\filldraw (d) node[left]{$\scriptstyle k$} circle (.05cm);
}
=
\frac{\lambda_\rho^{\ell+1}}{d}
\tikzmath{
\coordinate (t) at (.5,.5);
\coordinate (b) at (-.2,-.2);
\coordinate (a) at (.2,-.6);
\coordinate (d) at (-.2,.6);
\coordinate (c) at (.2,.2);
\draw[mid>] (d) .. controls ++(45:1cm) and ++(90:1.5cm) .. (1.2,0) .. controls ++(-90:1cm) and ++(-45:.7cm) .. (a);
\draw[thick, red, mid>] (d) .. controls ++(135:.7cm) and ++(90:2.5cm) .. (1.4,0) .. controls ++(-90:2cm) and ++(-135:.7cm) .. (a);
\draw[mid>] (a) to (b);
\draw[mid>] (b) to[out=45, in=-45] (c);
\draw[mid>] (c) to (d);
\draw[mid>] (b) to[out=135,in=-135] (d);
\draw[mid>] (a) to[out=45,in=-90] ($ (t) + (.4,0) $);
\draw[mid>] (c) to[out=45,in=-90] (t);
\draw (t) arc (180:0:.2cm);
\draw[thick] (t) -- ($ (t) + (.1,0) $);
\draw[thick, red] (b) .. controls ++(-135:.3cm) and ++(-90:.5cm) .. (0,-.2) .. controls ++(90:.3cm) and ++(-90:.3cm) .. (-.6,.2) .. controls ++(90:.3cm) and ++(90:.3cm) .. (.4,.2) .. controls ++(-90:.3cm) and ++(-135:.3cm) .. (c);
\filldraw (a) node[left]{$\scriptstyle j$} circle (.05cm);
\filldraw (b) node[left]{$\scriptstyle i$} circle (.05cm);
\filldraw (c) node[left]{$\scriptstyle \widetilde{\ell}$}circle (.05cm);
\filldraw (d) node[left]{$\scriptstyle k$} circle (.05cm);
}
=
\frac{\lambda_\rho^{\ell}}{d}
\tikzmath{
\coordinate (t) at (.5,1.3);
\coordinate (b) at (.2,-.2);
\coordinate (a) at (-.2,-.6);
\coordinate (d) at (.2,1);
\coordinate (c) at (-.2,.2);
\draw[mid>] (d) .. controls ++(135:.7cm) and ++(90:3cm) .. (1.2,0) .. controls ++(-90:1cm) and ++(-45:.7cm) .. (a);
\draw[thick, red] (d) .. controls ++(-135:.3cm) and ++(-90:.5cm) .. (.4,1) .. controls ++(90:.3cm) and ++(-90:.3cm) .. (-.1,1.4) .. controls ++(90:.7cm) and ++(90:.7cm) .. (.8,1.4) .. controls ++(-90:1.2cm) and ++(135:.5cm) .. (c);
\draw[mid>] (a) to[out=45,in=-45] (b);
\draw[mid>] (b) to (c);
\draw[mid>] (c) to[out=45,in=-45] (d);
\draw[mid>] (b) to[out=45,in=-90] ($ (t) + (.4,0)$);
\draw (t) arc (180:0:.2cm);
\draw[mid>] (a) to[out=135,in=-135] (c);
\draw[mid>] (d) to[out=45,in=-90] (t);
\draw[thick] (t) -- ($ (t) + (.1,0) $);
\draw[thick, red] (a) .. controls ++(-135:.3cm) and ++(-90:.5cm) .. (0,-.6) .. controls ++(90:.3cm) and ++(-90:.3cm) .. (-.6,-.2) .. controls ++(90:.3cm) and ++(90:.3cm) .. (.4,-.2) .. controls ++(-90:.3cm) and ++(-135:.3cm) .. (b);
\filldraw (a) node[left]{$\scriptstyle i$} circle (.05cm);
\filldraw (b) node[left, xshift=-.1cm]{$\scriptstyle \widetilde{\ell}$} circle (.05cm);
\filldraw (c) node[left]{$\scriptstyle k$}circle (.05cm);
\filldraw (d) node[left]{$\scriptstyle j$} circle (.05cm);
}
=
\frac{\lambda_\rho^{\ell+j}}{d}
\tikzmath{
\coordinate (b) at (.2,-.2);
\coordinate (a) at (-.2,-.6);
\coordinate (d) at (.2,.6);
\coordinate (c) at (-.2,.2);
\draw[mid>] (d) .. controls ++(45:.7cm) and ++(90:.5cm) .. (.8,0) .. controls ++(-90:.5cm) and ++(-45:.7cm) .. (a);
\draw[thick, red] (c) .. controls ++(135:.5cm) and ++(-90:.3cm) .. (.4,.6) .. controls ++(90:.3cm) and ++(135:.5cm) .. (d);
\draw[mid>] (a) to[out=45, in=-45] (b);
\draw[mid>] (b) to (c);
\draw[mid>] (c) to (d);
\draw[mid>] (b) to[out=45,in=-45] (d);
\draw[mid>] (a) to[out=135,in=-135] (c);
\draw[thick, red] (a) .. controls ++(-135:.3cm) and ++(-90:.5cm) .. (0,-.6) .. controls ++(90:.3cm) and ++(-90:.3cm) .. (-.6,-.2) .. controls ++(90:.3cm) and ++(90:.3cm) .. (.4,-.2) .. controls ++(-90:.3cm) and ++(-135:.3cm) .. (b);
\filldraw (a) node[left]{$\scriptstyle i$} circle (.05cm);
\filldraw (b) node[left, xshift=-.1cm]{$\scriptstyle \widetilde{\ell}$} circle (.05cm);
\filldraw (c) node[left]{$\scriptstyle k$}circle (.05cm);
\filldraw (d) node[left]{$\scriptstyle \widetilde{j}$} circle (.05cm);
}
=
\frac{\lambda_\rho^{\ell+j}}{d}
\tikzmath{
\coordinate (b) at (.2,-.2);
\coordinate (a) at (-.2,-.6);
\coordinate (d) at (.2,.6);
\coordinate (c) at (-.2,.2);
\draw[mid>] (d) .. controls ++(45:.7cm) and ++(90:.5cm) .. (.8,0) .. controls ++(-90:.5cm) and ++(-45:.7cm) .. (a);
\draw[thick, red, mid>] (d) .. controls ++(135:.7cm) and ++(90:1.5cm) .. (1,0) .. controls ++(-90:1.5cm) and ++(-135:.7cm) .. (a);
\draw[mid>] (a) to[out=45, in=-45] (b);
\draw[mid>] (b) to (c);
\draw[mid>] (c) to (d);
\draw[mid>] (b) to[out=45,in=-45] (d);
\draw[mid>] (a) to[out=135,in=-135] (c);
\draw[thick, red] (b) .. controls ++(-135:.5cm) and ++(-90:.3cm) .. (-.6,-.2) .. controls ++(90:.3cm) and ++(-90:.5cm) .. (0,.2) .. controls ++(90:.3cm) and ++(135:.5cm) .. (c);
\filldraw (a) node[left]{$\scriptstyle i$} circle (.05cm);
\filldraw (b) node[left, xshift=-.1cm]{$\scriptstyle \widetilde{\ell}$} circle (.05cm);
\filldraw (c) node[left]{$\scriptstyle k$}circle (.05cm);
\filldraw (d) node[left]{$\scriptstyle \widetilde{j}$} circle (.05cm);
}
=
\lambda_\rho^{\ell+j}
\overline{\widehat{A}^{k,\tilde{j}}_{i,\tilde{\ell}}}.
\]
\end{footnotesize}
We leave the verification of the remaining identities to the reader.
\end{proof}

To finish off this section, we explicitly compute the 4th Frobenius-Schur indicator of $\rho$ in terms of our free variables. This formula will be useful in the next section.

\begin{lem}\label{lem:4th}
We have that 
\[\nu_4(\rho) =\frac{1}{d} + \lambda_\rho\sum_{i,j}\omega_{\mathbf{1}, i}\omega_{\mathbf{1}, j}A^{i,j}_{i,j} + \lambda_\rho\lambda_\alpha\sum_{i,j}(-1)^{i+j}{\omega_{\alpha, i}}{\omega_{\alpha, j}}\widehat{A}^{i,j}_{i,j}. \]
\end{lem}
\begin{proof}
We pick the following orthonormal basis of $\cC_2(\rho^{\otimes 4}\to \mathbf{1})$:
\[   
\left\{
\frac{1}{d}\,
\tikzmath{
\draw[mid>] (0,0) -- (0,.4);
\draw[mid>] (.4,0) -- (.4,.4) arc (0:180:.2cm);
\draw[thick] (0,.4) -- (.1,.4);
\draw[mid>] (.6,0) -- (.6,.4);
\draw[mid>] (1,0) -- (1,.4) arc (0:180:.2cm);
\draw[thick] (.6,.4) -- (.7,.4);
}
\right \}\cup  \left\{
\frac{1}{\sqrt{d}}\,
\tikzmath{
\draw[mid>] (0,0) -- (0,.4);
\draw[mid>] (.8,0) -- (.8,.4) arc (0:180:.4cm);
\draw[thick] (0,.4) -- (.1,.4);
\draw[mid>] (-.3,-.3) arc (180:90:.3cm);
\draw[mid>] (.3,-.3) arc (0:90:.3cm);
\draw[mid>] (.5,-.3) arc (180:90:.3cm);
\draw[mid>] (1.1,-.3) arc (0:90:.3cm);
\filldraw (0,0) node[left,yshift=.1cm]{$\scriptstyle i$} circle (.05cm);
\filldraw (.8,0) node[left,yshift=.1cm]{$\scriptstyle j$} circle (.05cm);
}
\right \}_{i,j}\cup\left\{
\frac{1}{\sqrt{d}}\,
\tikzmath{
\draw[mid>] (0,0) -- (.5,.5) to[out=45,in=90] (1.2,.4);
\draw[mid>] (.8,0) to[out=45,in=-90] (1.2,.4);
\draw[thick] (1.2,.4) -- (1.1,.4);
\draw[red, thick, mid>] (.8,0) -- (.3,.5) to[out=135,in=90] (-.4,.4);
\draw[red, thick, mid>] (0,0) to[out=135,in=-90] (-.4,.4);
\draw[thick, red] (-.4,.4) -- (-.3,.4);
\draw[mid>] (-.3,-.3)  -- (0,0);
\draw[mid>] (.3,-.3) -- (0,0);
\draw[mid>] (.5,-.3) -- (.8,0);
\draw[mid>] (1.1,-.3) -- (.8,0);
\filldraw (0,0) node[left]{$\scriptstyle i$} circle (.05cm);
\filldraw (.8,0) node[left]{$\scriptstyle j$} circle (.05cm);
}
\right \}_{i,j} .
\]
With this basis we compute
\begin{align*}
\nu_4(\rho)
&=  
\frac{1}{d^2}
\tikzmath{
\draw[mid>] (1.2,0) arc (-180:0:.2cm) -- (1.6,.2) arc (0:180:1cm) -- (-.4,0) arc (-180:0:.2cm) -- (0,.4);
\draw[mid>] (.4,0) -- (.4,.4) arc (0:180:.2cm);
\draw[thick] (0,.4) -- (.1,.4);
\draw[mid>] (.4,0) arc (-180:0:.2cm) -- (.8,.4);
\draw[mid>] (1.2,0) -- (1.2,.4) arc (0:180:.2cm);
\draw[thick] (.8,.4) -- (.9,.4);
\draw[thick] (.4,0) -- (.5,0);
\draw[thick] (1.2,0) -- (1.3,0);
}
+ 
\frac{1}{d}\sum_{i,j}
\tikzmath{
\coordinate (a) at (0,0);
\coordinate (b) at (.6,0);
\coordinate (c) at (-.3,.6);
\coordinate (d) at (.3,.6);
\draw[mid>] (a) to (c);
\draw[mid>] (a) to (d);
\draw[mid>] (b) to (d);
\draw[mid<] (c) .. controls ++(-135:.3cm) and ++(-90:.3cm) .. (-.9,.6) .. controls ++(90:1.5cm) and ++(90:1.5cm) .. (.9,.6) to[out=-90,in=45] (b);
\draw[mid<] (a) to ($ (a) + (0,-.4) $);
\draw[mid<] (b) to ($ (b) + (0,-.4) $);
\draw[mid>] (c) to ($ (c) + (0,.4) $);
\draw[mid>] (d) to ($ (d) + (0,.4) $);
\draw ($ (a) + (0,-.4) $) arc (-180:0:.3cm);
\draw ($ (c) + (0,.4) $) arc (180:0:.3cm);
\draw[thick] ($ (a) + (0,-.4) $) -- ($ (a) + (.1,-.4) $);
\draw[thick] ($ (c) + (0,.4) $) -- ($ (c) + (.1,.4) $);
\filldraw (a) node[left]{$\scriptstyle i$} circle (.05cm);
\filldraw (b) node[left]{$\scriptstyle j$} circle (.05cm);
\filldraw (c) node[left]{$\scriptstyle i$}circle (.05cm);
\filldraw (d) node[left]{$\scriptstyle j$} circle (.05cm);
}
+ 
\frac{1}{d}\sum_{i,j}
\tikzmath{
\coordinate (a) at (0,0);
\coordinate (b) at (.6,0);
\coordinate (c) at (-.3,.6);
\coordinate (d) at (.3,.6);
\draw[mid>] (a) to (c);
\draw[mid>] (a) to (d);
\draw[mid>] (b) to (d);
\draw[mid<] (c) .. controls ++(-135:.3cm) and ++(-90:.3cm) .. (-.9,.6) .. controls ++(90:1.5cm) and ++(45:2.5cm) .. (b);
\draw[mid<] (a) to ($ (a) + (.3,-.3) $);
\draw[mid<] (b) to ($ (b) + (.3,-.3) $);
\draw[mid>] (c) to ($ (c) + (.3,.3) $);
\draw[mid>] (d) to ($ (d) + (.3,.3) $);
\draw ($ (c) + (.3,.3) $) .. controls ++(45:.5cm) and ++(45:.4cm) .. ($ (d) + (.3,.3) $);
\draw ($ (a) + (.3,-.3) $) .. controls ++(-45:.5cm) and ++(-45:.4cm) .. ($ (b) + (.3,-.3) $);
\draw[thick, red] (c) .. controls ++(135:.9cm) and ++(135:1cm) .. (d);
\draw[thick, red] (a) .. controls ++(-135:.9cm) and ++(-110:1cm) .. (b);
\draw[thick] ($ (a) + (.3,-.3) $) -- ($ (a) + (.38,-.22) $);
\draw[thick] ($ (d) + (.3,.3) $) -- ($ (d) + (.22,.38) $);
\filldraw (a) node[left]{$\scriptstyle i$} circle (.05cm);
\filldraw (b) node[left]{$\scriptstyle j$} circle (.05cm);
\filldraw (c) node[left]{$\scriptstyle i$}circle (.05cm);
\filldraw (d) node[left]{$\scriptstyle j$} circle (.05cm);
}
\displaybreak[1]\\&= 
\frac{1}{d} +  \frac{\lambda_\rho}{d}\sum_{i,j}
\tikzmath{
\coordinate (a) at (0,0);
\coordinate (b) at (.6,0);
\coordinate (c) at (-.3,.6);
\coordinate (d) at (.3,.6);
\draw[mid>] (a) to (c);
\draw[mid>] (a) to (d);
\draw[mid>] (b) to (d);
\draw[mid>] (c) .. controls ++(-135:.3cm) and ++(-90:.3cm) .. (-.9,.6) .. controls ++(90:1.5cm) and ++(90:1cm) .. (1.5,.4) arc (0:-180:.2cm);
\draw[mid>] (b) to[out=45,in=-90] (.7,.4);
\draw (.7,.4) arc (180:0:.2cm);
\draw[mid<] (a) to ($ (a) + (0,-.4) $);
\draw[mid<] (b) to ($ (b) + (0,-.4) $);
\draw[mid>] (c) to ($ (c) + (0,.4) $);
\draw[mid>] (d) to ($ (d) + (0,.4) $);
\draw ($ (a) + (0,-.4) $) arc (-180:0:.3cm);
\draw ($ (c) + (0,.4) $) arc (180:0:.3cm);
\draw[thick] ($ (a) + (0,-.4) $) -- ($ (a) + (.1,-.4) $);
\draw[thick] ($ (c) + (0,.4) $) -- ($ (c) + (.1,.4) $);
\draw[thick] (-.9,.6) -- (-.8,.6);
\draw[thick] (.7,.4) -- (.8,.4);
\filldraw (a) node[left]{$\scriptstyle i$} circle (.05cm);
\filldraw (b) node[left]{$\scriptstyle j$} circle (.05cm);
\filldraw (c) node[left]{$\scriptstyle i$}circle (.05cm);
\filldraw (d) node[left]{$\scriptstyle j$} circle (.05cm);
}
+ 
\frac{\lambda_\rho}{d} \sum_{i,j}
\tikzmath{
\coordinate (a) at (0,0);
\coordinate (b) at (.6,0);
\coordinate (c) at (-.3,.6);
\coordinate (d) at (.3,.6);
\draw[thick, red] (c) .. controls ++(135:.1cm) and ++(90:.2cm) .. (-.5,.6) .. controls ++(-90:.4cm) and ++(-90:.4cm) .. (-.1,.6) .. controls ++(90:.3cm) and ++(-90:.2cm) .. (-.7,1.3) .. controls ++(90:.5cm) and ++(90:.5cm) .. (.1,.9) .. controls ++(-90:1.2cm) and ++(-90:.5cm) .. (-1.1,.6) .. controls ++(90:1.5cm) and ++(90:1cm) .. (.2,1) to[out=-90,in=135] (d);
;
\draw[thick, red] (b) .. controls ++(-135:.3cm) and ++(-90:.2cm) .. (.8,0) .. controls ++(90:.4cm) and ++(90:.4cm) .. (.4,0) .. controls ++(-90:.3cm) and ++(90:.2cm) .. (.6,-.5) .. controls ++(-90:.3cm) and ++(-90:.3cm) .. (.2,-.3) .. controls ++(90:1.2cm) and ++(90:.5cm) .. (1.3,.4) .. controls ++(-90:2cm) and ++(-135:1cm) .. (a);
;
\draw[mid>] (a) to (c);
\draw[mid>] (a) to (d);
\draw[mid>] (b) to (d);
\draw[mid>] (c) .. controls ++(-135:.3cm) and ++(-90:.3cm) .. (-.9,.6) .. controls ++(90:1.5cm) and ++(90:1cm) .. (1.5,.4) arc (0:-180:.2cm);
\draw[mid>] (b) to[out=45,in=-90] (.7,.4);
\draw (.7,.4) arc (180:0:.2cm);
\draw[mid<] (a) to ($ (a) + (.3,-.3) $);
\draw[mid<] (b) to ($ (b) + (.3,-.3) $);
\draw[mid>] (c) to ($ (c) + (.3,.3) $);
\draw[mid>] (d) to ($ (d) + (.3,.3) $);
\draw ($ (c) + (.3,.3) $) .. controls ++(45:.5cm) and ++(45:.4cm) .. ($ (d) + (.3,.3) $);
\draw ($ (a) + (.3,-.3) $) .. controls ++(-45:.5cm) and ++(-45:.4cm) .. ($ (b) + (.3,-.3) $);
\draw[thick] ($ (a) + (.3,-.3) $) -- ($ (a) + (.38,-.22) $);
\draw[thick] ($ (c) + (.3,.3) $) -- ($ (c) + (.38,.22) $);
\draw[thick] (-.9,.6) -- (-.8,.6);
\draw[thick] (.7,.4) -- (.8,.4);
\filldraw (a) node[left]{$\scriptstyle i$} circle (.05cm);
\filldraw (b) node[left]{$\scriptstyle j$} circle (.05cm);
\filldraw (c) node[left]{$\scriptstyle i$}circle (.05cm);
\filldraw (d) node[left]{$\scriptstyle j$} circle (.05cm);
}
\displaybreak[1]\\&= 
\frac{1}{d} +\frac{\lambda_\rho}{d}\sum_{i,j}\omega_{\mathbf{1}, i}\omega_{\mathbf{1}, j}
\tikzmath{
\coordinate (b) at (-.2,-.2);
\coordinate (a) at (.2,-.6);
\coordinate (d) at (-.2,.6);
\coordinate (c) at (.2,.2);
\draw[mid>] (d) .. controls ++(90:.7cm) and ++(90:.5cm) .. (.8,0) .. controls ++(-90:.5cm) and ++(-90:.7cm) .. (a);
\draw[mid>] (a) to (b);
\draw[mid>] (b) to (c);
\draw[mid>] (c) to (d);
\draw[mid>] (b) to[out=135,in=-135] (d);
\draw[mid>] (a) to[out=45,in=-45] (c);
\filldraw (a) node[left]{$\scriptstyle j$} circle (.05cm);
\filldraw (b) node[left]{$\scriptstyle i$} circle (.05cm);
\filldraw (c) node[left]{$\scriptstyle j$}circle (.05cm);
\filldraw (d) node[left]{$\scriptstyle i$} circle (.05cm);
}
+ 
\frac{\lambda_\rho}{d} \sum_{i,j}{\omega_{\alpha, i}}{\omega_{\alpha, j}}
\tikzmath{
\coordinate (a) at (.2,-.6);
\coordinate (b) at (-.2,-.2);
\coordinate (c) at (.2,.2);
\coordinate (d) at (-.2,.6);
\draw[mid>] (d) .. controls ++(45:.7cm) and ++(90:.5cm) .. (.8,0) .. controls ++(-90:.5cm) and ++(-45:.7cm) .. (a);
\draw[red, thick] (d) .. controls ++(135:.3cm) and ++(90:.3cm) .. (0,.6) .. controls ++(-90:.4cm) and ++(-90:.4cm) .. (-.5,.6) .. controls ++(90:.5cm) and ++(90:.5cm) ..  (.2,.6) .. controls ++(-90:.2cm) and ++(135:.2cm) .. (c);
\draw[mid>] (a) to (b);
\draw[mid>] (b) to[out=45, in=-45] (c);
\draw[mid>] (c) to[out=45, in=-45] (d);
\draw[mid>] (b) to[out=135,in=-135] (d);
\draw[mid>] (a) to[out=45,in=-45] (c);
\draw[thick, red] (a) .. controls ++(-135:.7cm) and ++(90:.7cm) .. (.4,-.6) .. controls ++(-90:.7cm) and ++(-135:.7cm) .. (b);
\filldraw (a) node[left]{$\scriptstyle j$} circle (.05cm);
\filldraw (b) node[left]{$\scriptstyle i$} circle (.05cm);
\filldraw (c) node[left]{$\scriptstyle j$}circle (.05cm);
\filldraw (d) node[left]{$\scriptstyle i$} circle (.05cm);
}
\displaybreak[1]\\&= 
\frac{1}{d}+\frac{\lambda_\rho}{d}\sum_{i,j}\omega_{\mathbf{1}, i}\omega_{\mathbf{1}, j}
A_{i,j}^{i,j}
+ 
\frac{\lambda_\rho\lambda_\alpha}{d} \sum_{i,j}\chi_{\alpha,i}\chi_{\alpha,j}{\omega_{\alpha, i}}{\omega_{\alpha, j}}
\tikzmath{
\coordinate (b) at (-.2,-.2);
\coordinate (a) at (.2,-.6);
\coordinate (d) at (-.2,.6);
\coordinate (c) at (.2,.2);
\draw[mid>] (d) .. controls ++(45:.7cm) and ++(90:.5cm) .. (.8,0) .. controls ++(-90:.5cm) and ++(-45:.7cm) .. (a);
\draw[thick, red, mid>] (d) .. controls ++(135:.7cm) and ++(90:1.5cm) .. (1,0) .. controls ++(-90:1.5cm) and ++(-135:.7cm) .. (a);
\draw[mid>] (a) to (b);
\draw[mid>] (b) to[out=45, in=-45] (c);
\draw[mid>] (c) to[out=45, in=-45] (d);
\draw[mid>] (b) to[out=135,in=-135] (d);
\draw[mid>] (a) to[out=45,in=-45] (c);
\draw[thick, red] (b) .. controls ++(-135:.3cm) and ++(-90:.5cm) .. (0,-.2) .. controls ++(90:.3cm) and ++(135:.3cm) .. (c);
\filldraw (a) node[left]{$\scriptstyle j$} circle (.05cm);
\filldraw (b) node[left]{$\scriptstyle i$} circle (.05cm);
\filldraw (c) node[left]{$\scriptstyle j$}circle (.05cm);
\filldraw (d) node[left]{$\scriptstyle i$} circle (.05cm);
}
\displaybreak[1]\\&= 
\frac{1}{d} + \lambda_\rho\sum_{i,j}\omega_{\mathbf{1}, i}\omega_{\mathbf{1}, j}A^{i,j}_{i,j} + \lambda_\rho\lambda_\alpha\sum_{i,j}(-1)^{i+j}{\omega_{\alpha, i}}{\omega_{\alpha, j}}\widehat{A}^{i,j}_{i,j}.
\qedhere
\end{align*}
\end{proof}

%%%%%%%%%%%%%%%%%%%%%%%%%%%%%%%%%%%%%%%%%%%%%%%%%%%%%%%%%%%%%
\subsection{Classification}
\label{sec:SDClassification}

In this final subsection, we complete the classification result in the self-dual case \ref{Q:4ObjectsSelfDual}, i.e., we complete the proof of Theorem~\ref{thm:SD}, and classify all categorifications of the rings \ref{eq:R(m)}. 
We have several cases to consider
depending on $m\in \{1,2\}$ and $\lambda_\alpha = \pm1$.
Note that from Corollary~\ref{cor:SDm1}, we know that if $m=1$ then $\lambda_\alpha = 1$. 
This case has already been covered in \cite[Example 9.1]{MR3827808} where the result as in the statement of Theorem~\ref{thm:SD} is found. 
Hence it suffices to consider the case $m=2$.

%%%%%%%%%%%%%%%%%%%%%%%%%%%%%%%%%%%%%%%%%%%%%%%%%%%%%%%%%%%%%
\subsubsection{The case \texorpdfstring{$m=2$ and $\lambda_\alpha = 1$}{m=2 and la=1}}

In the case of $\lambda_\alpha=1$ we have determined that  \[\lambda_\rho = 1, \qquad \mu = 1,\qquad \tilde{i} = i, \quad \text{and} \quad 
    \chi_{\mathbf{1},i} =\chi_{\alpha,i}= (-1)^{i}.\]
Thus all that remains is to deduce the 3rd roots of unity 
$\omega_{\mathbf{1},0},  \omega_{\mathbf{1},1}, \omega_{\alpha,0}, \omega_{\alpha,1}$,
along with the free variables 
$A^{i,j}_{k,\ell},  \widehat{A}^{i,j}_{k,\ell},  D^{i,j}_{k,\ell}$. 
We express these free variables in the matrix form:
\begin{equation}
\label{eq:MatrixForm}
\begin{bmatrix}
X^{0, 0}_{0, 0} & X^{0, 0}_{0, 1} & X^{0, 0}_{1, 0} & X^{0,0}_{1, 1}	\\[.3em]
X^{0, 1}_{0, 0} & X^{0, 1}_{0, 1} & X^{0, 1}_{1, 0} & X^{0,1}_{1, 1}	\\[.3em]
X^{1, 0}_{0, 0} & X^{1, 0}_{0, 1} & X^{1, 0}_{1, 0} & X^{1,0}_{1, 1}	\\[.3em]
X^{1, 1}_{0, 0} & X^{1, 1}_{0, 1} & X^{1, 1}_{1, 0} & X^{1,1}_{1, 1}	
\end{bmatrix}
\qquad\qquad
X=A,\widehat{A},D.
\end{equation}
By applying the symmetries of Lemma~\ref{lem:mainSym}, we have that our free variables are of the form
\[
\resizebox{\hsize}{!}{
$A = 
\begin{bmatrix}
a_0 &0 & 0 & \omega_{\mathbf{1},0}a_2	\\
0 & a_2 & \omega_{\mathbf{1},0}^2a_2 &0		\\
0 &  \omega_{\mathbf{1},0}^2a_2& a_2 & 0	\\ 
 \omega_{\mathbf{1},0}a_2 &0 &0 & a_1
\end{bmatrix}
\quad 
\widehat{A} = 
\begin{bmatrix}
\widehat{a}_0 &0 & 0 & \omega_{\alpha,0}\widehat{a}_2	\\
0 & \widehat{a}_2 & \omega_{\alpha,0}^2\widehat{a}_2 &0		\\
0 &  \omega_{\alpha,0}^2\widehat{a}_2& \widehat{a}_2 & 0	\\ 
 \omega_{\alpha,0}\widehat{a}_2 &0 &0 & \widehat{a}_1
\end{bmatrix}
\quad 
D =
\begin{bmatrix}
d_0 &0 & 0 &  -\frac{\omega_{\alpha,1}}{\omega_{\alpha,0}}\overline{d_4}	\\
0 & d_2 & d_4 &0		\\
0 &  -\frac{\omega_{\alpha,0}}{\omega_{\alpha,1}}d_4& d_3 & 0	\\ 
\overline{d_4} &0 &0 & d_1
\end{bmatrix}
$}  \]
all of which are real apart from $d_4$. If these free coefficients are non-zero, then this implies conditions on our twists. We have
\begin{align*}
    a_0 \neq 0 &\implies \omega_{\mathbf{1}, 0} = 1 
    &
    a_1 \neq 0&\implies \omega_{\mathbf{1}, 1} = 1
    \\
    \widehat{a}_0 \neq 0 &\implies \omega_{\alpha, 0} = 1
    &
    \widehat{a}_1 \neq 0 &\implies \omega_{\alpha, 1} = 1
    \\
    a_2 \neq 0 &\implies \omega_{\mathbf{1}, 0} = \omega_{\mathbf{1}, 1}
    &
    \widehat{a}_2 \neq 0 &\implies \omega_{\alpha, 0} = \omega_{\alpha, 1}.
\end{align*}
In order to solve for these complex variables, we evaluate certain morphisms in our categories in two ways to obtain equations of these variables. We compute
\begin{align*}
\delta_{k,k'}\delta_{\ell,\ell'}
\,
\tikzmath{\draw[mid>](0,-.5) -- (0,.5);}
&=
\delta_{\ell,\ell'}
\,
\tikzmath{
\draw[mid>] (0,-.9) -- (0,-.4);
\draw[mid>] (0,.4) -- (0,.9);
\draw[mid>] (0,-.4) arc (240:120:.47cm);
\draw[mid>] (0,-.4) arc (-60:60:.47cm);
\filldraw (0,-.4) node[left]{$\scriptstyle k'$} circle (.05cm);
\filldraw (0,.4) node[left]{$\scriptstyle k$} circle (.05cm);
}
=
\tikzmath{
\draw[mid>] (-.2,-1.5) -- (-.2,-1);
\draw[mid>] (-.2,1) -- (-.2,1.5);
\draw[mid>] (-.2,-1) -- (0,-.4);
\draw[mid>] (0,.4) -- (-.2,1);
\draw[mid>] (0,-.4) arc (240:120:.47cm);
\draw[mid>] (0,-.4) arc (-60:60:.47cm);
\draw[mid>] (-.2,-1) arc (240:120:1.15cm);
\filldraw (0,-.4) node[left]{$\scriptstyle \ell'$} circle (.05cm);
\filldraw (0,.4) node[left]{$\scriptstyle \ell$} circle (.05cm);
\filldraw (-.2,-1) node[left]{$\scriptstyle k'$} circle (.05cm);
\filldraw (-.2,1) node[left]{$\scriptstyle k$} circle (.05cm);
\draw[rounded corners=5pt, dotted, red, thick] (-1,-.2) rectangle (.4,.75);
}
\\&=
\frac{\omega_{\mathbf{1}, \ell}}{2+\sqrt{5}}
\tikzmath{
\draw[mid>] (.2,1.5) -- (.2,2);
\draw[mid>] (.2,1) to[out=45,in=-45] (.2,1.5);
\draw[mid>] (.2,1) to[out=135,in=-135] (.2,1.5);
\draw[mid>] (0,-.5) -- (0,0);
\draw[mid>] (0,0) -- (.4,.5);
\draw[mid>] (.4,.5) to[out=45,in=-90] (.2,1);
\draw[mid>] (0,0) to[out=135, in=-90] (-.3,.5) ;
\draw[mid>] (.4,.5) .. controls ++(135:.5cm) and ++(90:.3cm) .. (-.3,.5);
\draw[thick] (-.3,.5) -- (-.2,.5);
\filldraw (.4,.5) node[left]{$\scriptstyle \ell'$} circle (.05cm);
\filldraw (0,0) node[left]{$\scriptstyle k'$} circle (.05cm);
\filldraw (.2,1) node[left]{$\scriptstyle \ell$} circle (.05cm);
\filldraw (.2,1.5) node[left]{$\scriptstyle k$} circle (.05cm);
}
+
\sum_{i,j}A^{i,j}_{k,\ell}
\tikzmath{
\draw[mid>] (.2,2) -- (.2,2.5);
\draw[mid>] (.2,1.5) to[out=45,in=-45] (.2,2);
\draw[mid>] (.2,1.5) to[out=135,in=-135] (.2,2);
\draw[mid>] (.2,1) -- (.2,1.5);
\draw[mid>] (0,-1) -- (0,-.5);
\draw[mid>] (0,.5) -- (.2,1);
\draw[mid>] (.4,0) -- (0,.5);
\draw[mid>] (0,-.5) -- (.4,0);
\draw[mid>] (0,-.5) to[out=135,in=-135] (0,.5);
\draw[mid>] (.4,0) to[out=45,in=-45] (.2,1);
\filldraw (.4,0) node[left]{$\scriptstyle \ell'$} circle (.05cm);
\filldraw (0,.5) node[left]{$\scriptstyle i$} circle (.05cm);
\filldraw (0,-.5) node[left]{$\scriptstyle k'$} circle (.05cm);
\filldraw (.2,1) node[left]{$\scriptstyle j$} circle (.05cm);
\filldraw (.2,1.5) node[left]{$\scriptstyle k$} circle (.05cm);
\filldraw (.2,2) node[left]{$\scriptstyle k$} circle (.05cm);
\draw[dotted, thick, red, rounded corners=5pt] (-.4,-.3) rectangle (.7,.3);
}
+ 
B^{i,j}_{k,\ell}
\tikzmath{
\draw[mid>] (.2,2) -- (.2,2.5);
\draw[mid>] (.2,1.5) to[out=45,in=-45] (.2,2);
\draw[mid>] (.2,1.5) to[out=135,in=-135] (.2,2);
\draw[mid>] (.2,1) to[out=45, in=-90]  (.2,1.5);
\draw[thick, red] (.2,1) .. controls ++(135:.5cm) and ++(90:.3cm) .. (.4,1) .. controls ++(-90:.3cm) and ++(135:.5cm) .. (0,.5);
\draw[mid>] (0,-1) -- (0,-.5);
\draw[mid>] (0,.5) -- (.2,1);
\draw[mid>] (.4,0) -- (0,.5);
\draw[mid>] (0,-.5) -- (.4,0);
\draw[mid>] (0,-.5) to[out=135,in=-135] (0,.5);
\draw[mid>] (.4,0) to[out=45,in=-45] (.2,1);
\filldraw (.4,0) node[left]{$\scriptstyle \ell'$} circle (.05cm);
\filldraw (0,.5) node[left]{$\scriptstyle i$} circle (.05cm);
\filldraw (0,-.5) node[left]{$\scriptstyle k'$} circle (.05cm);
\filldraw (.2,1) node[left]{$\scriptstyle j$} circle (.05cm);
\filldraw (.2,1.5) node[left]{$\scriptstyle k$} circle (.05cm);
\filldraw (.2,2) node[left]{$\scriptstyle k$} circle (.05cm);
\draw[dotted, thick, red, rounded corners=5pt] (-.4,-.3) rectangle (.7,.3);
}
\\
&= \delta_{k,\ell}\delta_{k',\ell'}\frac{\omega_{\mathbf{1},\ell}\omega_{\mathbf{1},\ell'}^2}{2+\sqrt{5}}
\,
\tikzmath{\draw[mid>](0,-.5) -- (0,.5);}
+
\sum_{i,j}A^{i,j}_{k,\ell}\overline{A^{i,j}_{k',\ell'}}
\tikzmath{
\draw[mid>] (.2,-2.5) -- (.2,-2);
\draw[mid>] (.2,-2) to[out=45,in=-45] (.2,-1.5);
\draw[mid>] (.2,-2) to[out=135,in=-135] (.2,-1.5);
\draw[mid>] (.2,-1.5) -- (.2,-1);
\draw[mid>] (.2,1) -- (.2,1.5);
\draw[mid>] (.2,-1) -- (0,-.4);
\draw[mid>] (0,.4) -- (.2,1);
\draw[mid>] (0,-.4) arc (240:120:.47cm);
\draw[mid>] (0,-.4) arc (-60:60:.47cm);
\draw[mid>] (.2,-1) arc (-60:60:1.15cm);
\filldraw (0,-.4) node[left]{$\scriptstyle i$} circle (.05cm);
\filldraw (0,.4) node[left]{$\scriptstyle i$} circle (.05cm);
\filldraw (.2,-1) node[left]{$\scriptstyle j$} circle (.05cm);
\filldraw (.2,1) node[left]{$\scriptstyle j$} circle (.05cm);
\filldraw (.2,-1.5) node[left]{$\scriptstyle k'$} circle (.05cm);
\filldraw (.2,-2) node[left]{$\scriptstyle k'$} circle (.05cm);
}
+ 
B^{i,j}_{k,\ell}\overline{B^{i,j}_{k',\ell'}}
\tikzmath{
\draw[thick, red] (.2,-1) .. controls ++(-135:.5cm) and ++(-90:.3cm) .. (.4,-1) .. controls ++(90:.3cm) and ++(-135:.5cm) .. (0,-.4);
\draw[mid>] (.2,-2.5) -- (.2,-2);
\draw[mid>] (.2,-2) to[out=45,in=-45] (.2,-1.5);
\draw[mid>] (.2,-2) to[out=135,in=-135] (.2,-1.5);
\draw[mid>] (.2,-1.5) to[out=90, in=-45]  (.2,-1);
\draw[thick, red] (.2,1) .. controls ++(135:.5cm) and ++(90:.3cm) .. (.4,1) .. controls ++(-90:.3cm) and ++(135:.5cm) .. (0,.4);
\draw[mid>] (.2,1) -- (.5,1.5);
\draw[mid>] (.2,-1) -- (0,-.4);
\draw[mid>] (0,.4) -- (.2,1);
\draw[mid>] (0,-.4) arc (240:120:.47cm);
\draw[mid>] (0,-.4) arc (-60:60:.47cm);
\draw[mid>] (.2,-1) arc (-60:60:1.15cm);
\filldraw (0,-.4) node[left]{$\scriptstyle i$} circle (.05cm);
\filldraw (0,.4) node[left]{$\scriptstyle i$} circle (.05cm);
\filldraw (.2,-1) node[left]{$\scriptstyle j$} circle (.05cm);
\filldraw (.2,1) node[left]{$\scriptstyle j$} circle (.05cm);
\filldraw (.2,-1.5) node[left]{$\scriptstyle k'$} circle (.05cm);
\filldraw (.2,-2) node[left]{$\scriptstyle k'$} circle (.05cm);
}
\\&= 
\left(\delta_{k,\ell}\delta_{k',\ell'}\frac{\omega_{\mathbf{1},\ell}\omega_{\mathbf{1},\ell'}^2}{2+\sqrt{5}} +\sum_{i',j'}A^{i,j}_{k,\ell}\overline{A^{i,j}_{k',\ell'}}+ B^{i,j}_{k,\ell}\overline{B^{i,j}_{k',\ell'}}\right)
\,
\tikzmath{\draw[mid>](0,-.5) -- (0,.5);}
\\&= 
\left(\delta_{k,\ell}\delta_{k',\ell'}\frac{\omega_{\mathbf{1},\ell}\omega_{\mathbf{1},\ell'}^2}{2+\sqrt{5}} +\sum_{i',j'}A^{i,j}_{k,\ell}\overline{A^{i,j}_{k',\ell'}}+(-1)^{\ell+\ell'}\omega_{\mathbf{1},\ell}\omega_{\mathbf{1},\ell'}^2 D^{j,k}_{i,\ell}\overline{D^{j,k'}_{i,\ell'}}\right)
\,
\tikzmath{\draw[mid>](0,-.5) -- (0,.5);}
\\&= 
\left(\delta_{k,\ell}\delta_{k',\ell'}\frac{\omega_{\mathbf{1},\ell}\omega_{\mathbf{1},\ell'}^2}{2+\sqrt{5}} +\omega_{\mathbf{1},\ell}\omega_{\mathbf{1},\ell'}^2\sum_{i,j}A^{j,k}_{i,\ell}\overline{A^{j,k'}_{i,\ell'}}+(-1)^{\ell+\ell'}\omega_{\mathbf{1},\ell}\omega_{\mathbf{1},\ell'}^2 D^{j,k}_{i,\ell}\overline{D^{j,k'}_{i,\ell'}}\right) 
\,
\tikzmath{\draw[mid>](0,-.5) -- (0,.5);}
\,.
\end{align*}
Note that if $\ell\neq \ell'$ then the left hand side vanishes, and we can cancel the $\omega_{\mathbf{1},\ell}\omega_{\mathbf{1},\ell'}^2$ terms. If $\ell= \ell'$ then $\omega_{\mathbf{1},\ell}\omega_{\mathbf{1},\ell'}^2=1$. In either case, we can remove the $\omega_{\mathbf{1},\ell}\omega_{\mathbf{1},\ell'}^2$ terms from the above equation. This leaves us with the equation
\[
  \sum_{i,j}A^{j,k}_{i,\ell}\overline{A^{j,k'}_{i,\ell'}}+(-1)^{\ell+\ell'} \sum_{i,j}D^{\ell,i}_{k,j}\overline{D^{\ell',i}_{k',j}}- +(2-\sqrt{5})\delta_{k,\ell}\delta_{k',\ell'} = \delta_{k,k'}\delta_{\ell,\ell'} .
\]
In a similar fashion, we can evaluate the diagrams
\[
\tikzmath{
\draw[mid>] (-.2,1) -- (.1,1.5);
\draw[thick, red, mid>] (-.2,1) -- (-.5,1.5);
\draw[mid>] (.1,-1.5) -- (-.2,-1);
\draw[thick, red, mid>] (-.5,-1.5) -- (-.2,-1);
\draw[mid>] (-.2,-1) .. controls ++(45:.3cm) and ++(-45:.3cm) .. (0,-.4);
\draw[mid>] (0,.4) .. controls ++(45:.3cm) and ++(-45:.3cm) .. (-.2,1);
\draw[thick, red] (0,-.4) .. controls ++(-135:.9cm) and ++(-90:.5cm) .. (-1,0) .. controls ++(90:.5cm) and ++(135:.9cm) .. (0,.4);
\draw[mid>] (0,-.4) arc (240:120:.47cm);
\draw[mid>] (0,-.4) arc (-60:60:.47cm);
\draw[mid>] (-.2,-1) arc (240:120:1.15cm);
\filldraw (0,-.4) node[left]{$\scriptstyle \ell'$} circle (.05cm);
\filldraw (0,.4) node[left]{$\scriptstyle \ell$} circle (.05cm);
\filldraw (-.2,-1) node[left]{$\scriptstyle k'$} circle (.05cm);
\filldraw (-.2,1) node[left]{$\scriptstyle k$} circle (.05cm);
}
\,\qquad 
\tikzmath{
\draw[mid>] (-.2,1) -- (.1,1.5);
\draw[thick, red, mid>] (-.2,1) -- (-.5,1.5);
\draw[mid>] (.1,-1.5) -- (-.2,-1);
\draw[thick, red, mid>] (-.5,-1.5) -- (-.2,-1);
\draw[mid>] (-.2,-1) -- (0,-.4);
\draw[mid>] (0,.4) -- (-.2,1);
\draw[mid>] (0,-.4) arc (240:120:.47cm);
\draw[mid>] (0,-.4) arc (-60:60:.47cm);
\draw[mid>] (-.2,-1) arc (240:120:1.15cm);
\filldraw (0,-.4) node[left]{$\scriptstyle \ell'$} circle (.05cm);
\filldraw (0,.4) node[left]{$\scriptstyle \ell$} circle (.05cm);
\filldraw (-.2,-1) node[left]{$\scriptstyle k'$} circle (.05cm);
\filldraw (-.2,1) node[left]{$\scriptstyle k$} circle (.05cm);
}
\,\qquad
\tikzmath{
\draw (0,-.7) arc (-90:0:.3cm);
\draw[mid>] (.3,-.4) -- (.3,-.3) arc (0:90:.3cm);
\draw (0,-.7) arc (270:180:.3cm); 
\draw[mid>] (-.3,-.4) -- (-.3,-.3) arc (180:90:.3cm);
\draw[thick] (-.3,-.4) -- (-.2,-.4);
\draw[mid>] (-.9,-1) -- (-.9,0) arc (180:90:.45cm);
\draw[mid>] (0,0) arc (0:90:.45cm);
\draw[mid>] (-.45,.45) -- (-.45,1);
\filldraw (0,0) node[left, yshift=.1cm]{$\scriptstyle \ell$} circle (.05cm);
\filldraw (-.45,.45) node[left, yshift=.1cm]{$\scriptstyle k$} circle (.05cm);
}
\quad \text{and} \quad  
\tikzmath{
\draw (0,-.7) arc (-90:0:.3cm);
\draw[mid>] (.3,-.4) -- (.3,-.3) to[out=90,in=-45] (0,0);
\draw (0,-.7) arc (270:180:.3cm); 
\draw[mid>] (-.3,-.4) -- (-.3,-.3) to[out=90,in=-135] (0,0);
\draw[thick] (-.3,-.4) -- (-.2,-.4);
\draw[mid>] (-.6,-1) -- (-.6,0) to[out=90, in=-135] (-.4,1);
\draw[mid>] (0,0) to[out=45,in=-135] (.2,.2) to[out=45, in=-45] (-.4,1);
\draw[thick, red] (0,0) .. controls ++(135:.5cm) and ++(90:.5cm) .. (.5,-.3) arc (0:-90:.6cm) to[out=180,in=-90] (-.8,0) to[out=90,in=180] (-.4,.7) arc (-90:0:.4cm) .. controls ++(90:.5cm) and ++(135:.5cm) .. (-.4,1);
\draw[mid>] (-.4,1) -- (.1,1.5);
\filldraw (0,0) node[left]{$\scriptstyle \ell$} circle (.05cm);
\filldraw (-.4,1) node[left]{$\scriptstyle k$} circle (.05cm);
}
\]
in two ways\footnote{\label{footnote:EvaluateInTwoWays}When we say we evaluate a diagram in two ways to obtain a relation, one way is trivial, and the other uses the jellyfish relations from Lemmas \ref{lem:alphaJellyfish} and \ref{lem:SDJellyfishRelations}.
For the non-self dual case in \S\ref{sec:NSD} below, we use the jellyfish relations from \eqref{eq:NSDExtraJellyfish}, \eqref{eq:NSDalphaJellyfish}, and Lemma \ref{lem:NSDJellyfishRelations} instead.}
to obtain
\begin{align*}
   \sum_{i,j}\widehat{A}^{i,k}_{j,\ell}\overline{\widehat{A}^{i,k'}_{j,\ell'}} +\sum_{i,j} D^{j,k}_{i,\ell}\overline{D^{j,k'}_{i,\ell'}} - (2-\sqrt{5})\delta_{k,\ell}\delta_{k',\ell'} 
   &= 
   \delta_{k,k'}\delta_{\ell,\ell'} 
   \\
   \sum_{i,j}\widehat{D}^{i,j}_{k,\ell}\overline{\widehat{D}^{i,j}_{k',\ell'}} +(-1)^{\ell+\ell'}\omega_{\mathbf{1}, \ell'}\omega_{\mathbf{1}, \ell}^2\sum_{i,j} D^{i,j}_{k',\ell'   }\overline{D^{i,j}_{k,\ell}}
   &= 
   \delta_{k,k'}\delta_{\ell,\ell'}  
   \\
   \sum_{i}A^{i,k}_{i,\ell} + (-1)^\ell\sum_{i}D^{i,k}_{i,\ell} - (2-\sqrt{5})\delta_{k,\ell}
   &= 
   0
   \\
    \sum_{i}\widehat{A}^{i,k}_{i,\ell} + \sum_{i}(-1)^iD^{\ell,i}_{k,i} - (2-\sqrt{5})\delta_{k,\ell}
    &=
    0.
\end{align*}
In terms of our free variables, this gives us the equations:
\begin{align}
    a_0^2 {+} a_2^2 {+} d_0^2 {+} d_3^2 = a_1^2 {+} a_2^2 {+} d_1^2 {+} d_2^2 = \widehat{a}_0^2 {+} \widehat{a}_2^2 {+} d_0^2 {+} d_2^2 =\widehat{a}_1^2 {+} \widehat{a}_2^2 {+} d_1^2 {+} d_3^2 &=3 - \sqrt{5}\label{eq:coef1}\\
    a_2^2  {+} |d_4|^2 = \widehat{a}_2^2{+} |d_4|^2 =d_0^2  {+} |d_4|^2= d_1^2  {+} |d_4|^2=d_2^2  {+} |d_4|^2=d_3^2  {+} |d_4|^2&=\frac{1}{2}\qquad \label{eq:coef2}\\
    (a_0 {+} a_1)a_2 - d_0d_2 - d_1d_3 = (\widehat{a}_0 {+} \widehat{a}_1)\widehat{a}_2 {+} d_0d_3 {+} d_1d_2 &= 2-\sqrt{5}\label{eq:coef3}\\
    a_0 {+} a_2 {+} d_0 {+} d_3 =a_1 {+} a_2 - d_1 - d_2 =  \widehat{a}_0 {+} \widehat{a}_2 {+} d_0 - d_2 =\widehat{a}_1 {+} \widehat{a}_2 - d_1 {+} d_3 &=2-\sqrt{5}\label{eq:coef4}\\
    (\omega_{\alpha,0} {+} \omega_{\alpha,0}^2)a_2^2{+} \omega_{\alpha,0}\omega_{\alpha,1}^2 d_4^2 {+} \omega_{\alpha,0}^2\omega_{\alpha,1}\overline{d_4}^2 =  (\omega_{\alpha,0} {+} \omega_{\alpha,0}^2)\widehat{a}_2^2{+} \omega_{\alpha,0}^2\omega_{\alpha,1} d_4^2 {+}     \overline{d_4}^2  &=0\label{eq:coef5}\qquad\\
   (1-\omega_{\alpha,0}\omega_{\alpha,1}^2)(d_2\overline{d_4} - d_3^2)  =d_4(d_1 {+} d_0) - \overline{d_4}(\omega_{\alpha,0}^2\omega_{\alpha,1}d_0 {+}\omega_{\alpha,0}\omega_{\alpha,1}^2d_1) &=0. \qquad\label{eq:coef6}
\end{align}
While we could begin solving these equations directly, instead we opt for a more measured approach, and use our previous centre analysis to simplify our solution.

\begin{lem}
\label{lem:SD-existsTau}
There exists a $\tau \in \{-1,1\}$ such that
\begin{align*}
a_0 = a_1 = \widehat{a}_0 = \widehat{a}_1 &= \frac{2 +3\tau- \sqrt{5}}{4}
\qquad \text{and}\qquad
a_2 = \widehat{a}_2 =  d_0 = -d_1 = -d_2 = d_3 &=  \frac{2-\tau -\sqrt{5}}{4}.
\end{align*}
In particular, as $a_0,a_1,\widehat{a}_0$, and $\widehat{a}_1$ are all non-zero, we have 
$\omega_{\mathbf{1}, 0} =\omega_{\mathbf{1}, 1}=\omega_{\alpha, 0} = \omega_{\alpha, 1}=1$. 
\end{lem}
\begin{proof}
We first observe from Equation~\eqref{eq:coef2} that
\[     a_2^2 = \widehat{a}_2^2 = d_0^2 = d_1^2 = d_2^2 = d_3^2 = \frac{1}{2} - |d_4|^2,\]
and in particular we have that $a_2, \widehat{a}_2, d_0,d_1,d_2$, and $d_3$ are real numbers which are equal up to sign. With this information in hand, we can now see from Equation~\eqref{eq:coef1} that 
\[   a_0^2 = a_1^2 = \widehat{a_0}^2 =  \widehat{a_1}^2 =3- \sqrt{5} - 3d_0^2.   \]

To make additional progress on solving these equations, we recall the operators $\phi$ and $\psi$. 
In our case, via Equation~\eqref{eq:coef4}, we have that
\begin{align*}
\phi &=\begin{bmatrix}
a_0 + a_2 & 0\\
0 &a_1+a_2
\end{bmatrix}  = \begin{bmatrix}
2-\sqrt{5}-d_0-d_3 & 0\\
0 &2-\sqrt{5}+d_1 + d_2
\end{bmatrix} \\ 
\psi &=\begin{bmatrix}
\widehat{a}_0 + \widehat{a}_2 & 0\\
0 &\widehat{a}_1+\widehat{a}_2
\end{bmatrix}  = \begin{bmatrix}
2-\sqrt{5}-d_0+d_2 & 0\\
0 &2-\sqrt{5}+d_1 -d_3
\end{bmatrix}.
\end{align*}
From Lemma~\ref{lem:SDidInd} and Remark~\ref{rem:phi}, we know that $\phi$ and $\psi$ have entries in $\{\frac{3-\sqrt{5}}{2},\frac{1-\sqrt{5}}{2}\}$, so
\[  d_0 + d_3 , -d_1 - d_2, d_0 - d_2 ,d_3-d_1\in\left\{\frac{3-\sqrt{5}}{2},\frac{1-\sqrt{5}}{2} \right\}. \]
In particular, as the values $d_0,d_1,d_2$, and $d_3$ are real numbers which are the same up to sign, we have that 
\[  d_0 = -d_1 = -d_2 = d_3 = \frac{2-\sqrt{5}-\tau}{4}, \]
for some $\tau \in \{-1,1\}$.

From Equations~\eqref{eq:coef4} we can deduce that
$a_0 = a_1$ and $\widehat{a}_0 = \widehat{a}_1$.
We know that $a_2$ and $c_0$ are the same up to sign. If we have $a_2 = -c_0$, then Equation~\eqref{eq:coef4} would imply that
\[  a_0 = 2 - \sqrt{5}-d_0. \]
Plugging this value of $a_0$ into Equation~\eqref{eq:coef1} gives a contradiction. Thus $a_2 = c_0$, and so Equation~\eqref{eq:coef4} gives
\[ a_0 = 2-\sqrt{5} - 3 d_0 =\frac{2 +3\tau- \sqrt{5}}{4}. \]
A similar argument shows that $\widehat{a_2} = c_0$, and thus $\widehat{a_0} = a_0$.
\end{proof}

To pin down the value of $\tau$, we return to our analysis of the centre of $\cC_2$. 
By computing the 4th Frobenius-Schur indicator of $\rho$ in two ways, we can show that $\tau= 1$.

\begin{lem}
\label{lem:SD-tau=1}
We have that $\tau =1$.
\end{lem}
\begin{proof}
From Lemma~\ref{lem:4th}, we have that
$\nu_4(\rho) = 3\tau$.
On the other hand, we can use Lemma~\ref{lem:2nfrob} to obtain
\[\nu_4(\rho)(20 + 8\sqrt{5}) =  48  + 20\sqrt{5} + (2 + \sqrt{5})\sum p_i(p_i + q_i)  \theta_i^4   \]
where $\sum p_i(p_i + q_i)=16$ and the $\theta_i$'s are roots of unity. Thus
\[ \sum p_i(p_i + q_i)  \theta_i^4  = 4(-1 - 2\sqrt{5} + 3\tau \sqrt{5}).  \]
If $\tau = -1$, then Theorem~\ref{thm:LarsonRootsOfUnityAB} implies that it would take at least 24 roots of unity to write $4(-1 - 2\sqrt{5} + 3\tau \sqrt{5})$, and hence $\sum p_i(p_i + q_i)\geq 24$, giving a contradiction. Thus we must have $\tau = 1$.
\end{proof}
Now that we know all of our real free variables, we can solve for $d_4$, the one complex variable.
\begin{lem}\label{lem:d4}
We have that
\[
d_4 = \eta_1 \frac{1}{2} + \eta_2\frac{1}{2}\sqrt{\frac{-1+\sqrt{5}}{2}}
\qquad\qquad
\text{where }\eta_1,\eta_2 \in  \{-1,1\}.
\]
\end{lem}
\begin{proof}
From Lemmas \ref{lem:SD-existsTau} and \ref{lem:SD-tau=1}, we have $a_2=\frac{1-\sqrt{5}}{4}$.
By Equations~\eqref{eq:coef2} and \eqref{eq:coef5}, we have 
\[
|d_4|^2 = \frac{1+\sqrt{5}}{8}
\qquad\text{and}\qquad
d_4^2 + \overline{d_4}^2 = \frac{\sqrt{5}-3}{4}.  
\]
The 4 intersection points of this hyperbola and circle yield the statement of the lemma.
\end{proof}
Now that we have pinned down all of our variables, we can prove part of our main theorem which states that there is no fusion category when $\lambda_\alpha=1$.

\begin{thm}
There is no unitary fusion category that categorifies $R(2)$ with $\lambda_\alpha = 1$.
\end{thm}
\begin{proof}
By evaluating the diagram
\[ 
\tikzmath{
\draw[mid>] (-.2,1) -- (.1,1.5);
\draw[mid>] (-.2,-1.5) -- (-.2,-1);
\draw[mid>] (-.2,-1) .. controls ++(45:.3cm) and ++(-45:.3cm) .. (0,-.4);
\draw[mid>] (0,.4) -- (-.2,1);
\draw[thick, red] (0,-.4) .. controls ++(-135:.9cm) and ++(-90:.5cm) .. (-1,-.2) .. controls ++(90:.5cm) and ++(-90:.5cm) .. (.6,.2) .. controls ++(90:1cm) and ++(135:1cm) .. (-.2,1);
\draw[mid>] (0,-.4) arc (240:120:.47cm);
\draw[mid>] (0,-.4) arc (-60:60:.47cm);
\draw[mid>] (-.2,-1) arc (240:120:1.15cm);
\filldraw (0,-.4) node[left]{$\scriptstyle \ell'$} circle (.05cm);
\filldraw (0,.4) node[left]{$\scriptstyle \ell$} circle (.05cm);
\filldraw (-.2,-1) node[left]{$\scriptstyle k'$} circle (.05cm);
\filldraw (-.2,1) node[left]{$\scriptstyle k$} circle (.05cm);
}
\]
in two ways 
(see Footnote \ref{footnote:EvaluateInTwoWays}),
we obtain the equation
\[ (-1)^{k+\ell}\sum_{i,j}\overline{D}_{i,j}^{k,\ell}\overline{D}_{\ell',k'}^{i,j} + \sum_{i,j} D_{i,j}^{\ell',k'}D^{i,j}_{k,\ell}=0.\]
Taking $k= \ell' = 0$ and $\ell = k' =1$ we see that
$\sum_{i,j}D^{0,1}_{i,j}D^{i,j}_{0,1} = d_2^2 - d_4^2  \in \bbR$.
Since $d_2\in \bbR$ by Lemma \ref{lem:SD-existsTau}, this means 
$d_4^2\in \bbR$, which contradicts Lemma~\ref{lem:d4}.
\end{proof}

%%%%%%%%%%%%%%%%%%%%%%%%%%%%%%%%%%%%%%%%%%%%%%%%%%%%%%%%%%%
\subsubsection{The case \texorpdfstring{$m=2$ and $\lambda_\alpha = -1$}{m=2 and la=-1}}

In the case of $\lambda_\alpha = -1$ we have determined that
\[ \mu = 1,\qquad \tilde{i} = 1-i,\qquad \omega_{\mathbf{1}, 0} = \omega_{\mathbf{1}, 1},   \qquad\omega_{\alpha, 0} = \omega_{\alpha, 1},\qquad \chi_{\mathbf{1}, i} = (-1)^i,\quad \text{ and } \quad \chi_{\alpha, i} =(-1)^i \mathbf{i}. \]

Thus all that remains is to deduce $\lambda_\rho$, the 3rd roots of unity
$\omega_{\mathbf{1},0}$ and $\omega_{\mathbf{1},1}$,
along with the free variables
$A^{i,j}_{k,\ell}$, $\widehat{A}^{i,j}_{k,\ell},$, and $D^{i,j}_{k,\ell}$.
By studying the 4th Frobenius-Schur indicator of $\rho$, we are able to show that $\lambda_\rho = 1$ and $\omega_{\alpha,0} = \omega_{\mathbf{1}, 0}^2$, along with the values of several of our free variables.

\begin{lem}
We have that $\lambda_\rho = 1$, and $\omega_{\alpha,0} = \omega_{\mathbf{1}, 0}^2$. Further, we have that
$$
A^{0,0}_{0,0}
= 
\frac{3-\sqrt{5}}{2(1 + \omega_{\mathbf{1}, 0})}
\qquad\text{ and }\qquad
\widehat{A}^{0,0}_{0,0}
= 
\frac{3-\sqrt{5}}{2(1 + \omega_{\mathbf{1}, 0}^2)}.
$$
\end{lem}
\begin{proof}

Recall the operators $\phi$ and $\psi$. By applying the symmetries of Lemma~\ref{lem:mainSym} we have that
\begin{align*}
   \phi &= \begin{bmatrix}
A^{0,0}_{0,0}+A^{1,0}_{1,0} & 0\\
0 &A^{0,1}_{0,1}+A^{1,1}_{1,1}
\end{bmatrix} = \begin{bmatrix}
 A^{0,0}_{0,0}(1+\lambda_\rho \omega_{\mathbf{1},0}) & 0\\
0 &A^{0,0}_{0,0}(1+\lambda_\rho \omega_{\mathbf{1},0})
\end{bmatrix} \\ 
 \psi &= \begin{bmatrix}
\widehat{A}^{0,0}_{0,0}+\widehat{A}^{1,0}_{1,0} & 0\\
0 &\widehat{A}^{0,1}_{0,1}+\widehat{A}^{1,1}_{1,1}
\end{bmatrix} = \begin{bmatrix}
 \widehat{A}^{0,0}_{0,0}(1+\lambda_\rho \omega_{\alpha,0}) & 0\\
0 &\widehat{A}^{0,0}_{0,0}(1+\lambda_\rho \omega_{\alpha,0})
\end{bmatrix}.
\end{align*}
Thus the operators $\phi$ and $\psi$ are scalars, and  Lemma~\ref{lem:SDidInd} and Remark~\ref{rem:phi} tell us that
\[
     A^{0,0}_{0,0} = \frac{2 - \sqrt{5} + \tau}{2(1 + \lambda_\rho\omega_{\mathbf{1}, 0})}\quad \text{ and } \quad
     \widehat{A}^{0,0}_{0,0} = \frac{2 - \sqrt{5} + \tau}{2(1 + \lambda_\rho\omega_{\alpha, 0})}
\]
for some $\tau \in \{-1,1\}$.

From Lemma~\ref{lem:4th}, we can write the 4th Frobenius-Schur indicator of $\rho$ as
\begin{align*}\nu_4(\rho) &=  \sqrt{5}-2 +\lambda_\rho\omega_{\mathbf{1}, 0}^2(A^{0,0}_{0,0} + A^{0,1}_{0,1} + A^{1,0}_{1,0} + A^{1,1}_{1,1}) - \lambda_\rho\omega_{\alpha,0}^2(\widehat{A}^{0,0}_{0,0}-\widehat{A}^{0,1}_{0,1} - \widehat{A}^{1,0}_{1,0} +\widehat{A}^{1,1}_{1,1})\\
&=\sqrt{5}-2 + \lambda_\rho\omega_{\mathbf{1}, 0}^2 A^{0,0}_{0,0}(2 + 2\lambda_\rho\omega_{\mathbf{1}, 0}) - \lambda_\rho\omega_{\alpha,0}^2\widehat{A}^{0,0}_{0,0}( 2 - 2\lambda_\rho\omega_{\alpha, 0}) \\
&= \sqrt{5} - 2 + \lambda_\rho\omega_{\mathbf{1}, 0}^2(2 - \sqrt{5}+\tau) - \lambda_\rho\omega_{\alpha,0}^2(2 - \sqrt{5}+\tau)\frac{1- \lambda_\rho\omega_{\alpha,0}}{1+ \lambda_\rho\omega_{\alpha,0}}\\
&= \left(-2 +\lambda_\rho\omega_{\mathbf{1},0}^2(2+\tau) - \lambda_\rho\omega_{\alpha,0}^2(2 + \tau)\frac{1- \lambda_\rho\omega_{\alpha,0}}{1+ \lambda_\rho\omega_{\alpha,0}}     \right) +\sqrt{5}\left(1 - \lambda_\rho\omega_{\mathbf{1}, 0}^2+\lambda_\rho\omega_{\alpha, 0}^2\frac{1- \lambda_\rho\omega_{\alpha,0}}{1+ \lambda_\rho\omega_{\alpha,0}}  \right).
\end{align*}
As $\nu_4(\rho) \in \mathbb{Z}[i]$,  
$\lambda_\rho$ is a second root of unity, and 
$\omega_{\mathbf{1},0}$ and $\omega_{\alpha,0}$ are third roots of unity, 
we have 
\[ 
1 - \lambda_\rho\omega_{\mathbf{1}, 0}^2+\lambda_\rho\omega_{\alpha, 0}^2\frac{1- \lambda_\rho\omega_{\alpha,0}}{1+ \lambda_\rho\omega_{\alpha,0}}
=
0.  
\]
This implies that $\lambda_\rho =1$, and that $\omega_{\alpha, 0} = \omega_{\mathbf{1}, 0}^{-1}$. By simplifying the formula for $\nu_4(\rho)$ further, we find that $\nu_4(\rho) = \tau$.

To determine $\tau$ we use Lemma~\ref{lem:2nfrob} to write
\[  \tau (20 +8\sqrt{5})  = \nu_4(\rho) = 48 + 20\sqrt{5} + (2 + \sqrt{5}) \sum p_i(p_i+q_i)\theta_i^4  \]
where $\sum p_i(p_i + q_i) = 16$, and the $\theta_i$'s are roots of unity. 
If $\tau = -1$, then we have
\[ 
\sum p_i(p_i + q_i)\theta_i^4 
= 
-4 - 12\sqrt{5}.  
\]
However, Theorem~\ref{thm:LarsonRootsOfUnityAB} implies that it takes at least 48 roots of unity to write $-4 - 12\sqrt{5}$, giving a contradiction. 
Thus $\tau = 1$, which gives
\[
A^{0,0}_{0,0}
= 
\frac{3-\sqrt{5}}{2(1 + \omega_{\mathbf{1}, 0})}
\qquad\text{and}\qquad
\widehat{A}^{0,0}_{0,0}
= 
\frac{3-\sqrt{5}}{2(1 + \omega_{\mathbf{1}, 0}^2)}.
\qedhere
\]
\end{proof}

Now that we know $\lambda_\rho=1$, the symmetries of Lemma~\ref{lem:mainSym} become much simpler. Using the same matrix notation as in the $\lambda_\alpha=1$ case from \eqref{eq:MatrixForm}, we can use these symmetries to express our free variables as
\[
\resizebox{\hsize}{!}{
$A = 
\begin{bmatrix}
\omega_{\mathbf{1},0}r &0 & 0 & \overline{a_1}	\\
0 &  \omega_{\mathbf{1},0}^2r &r&0		\\
0 &  r&  \omega_{\mathbf{1},0}^2r & 0	\\ 
 a_1 &0 &0 &  \omega_{\mathbf{1},0}r
\end{bmatrix}
\quad
\widehat{A} = 
\begin{bmatrix}
\omega_{\mathbf{1},0}^2r &0 & 0 & \overline{\widehat{a}_1}	\\
0 &  \omega_{\mathbf{1},0}r &r&0		\\
0 &  r&  \omega_{\mathbf{1},0}r & 0	\\
 \widehat{a}_1 &0 &0 &  \omega_{\mathbf{1},0}^2r
\end{bmatrix}
\quad 
D = 
\begin{bmatrix}
d_0 &0 & 0 & \overline{d_2}	\\
0 & \overline{d_0} & d_1 &0		\\
0 &  \overline{d_1}&-\overline{d_0} & 0	\\
d_2 &0 &0 & -d_0
\end{bmatrix}$
}  
\]
where $r = \frac{1}{\omega_{\mathbf{1},0}+\omega_{\mathbf{1},0}^2}\frac{3- \sqrt{5}}{2}\in \mathbb{R}$, and if either of $a_1$ or $\widehat{a}_1$ are non-zero, then we have that $\omega_{\mathbf{1}, 0} = 1$.

Now that we have reduced our free variables down to $5$ complex variables, all that remains is to solve for these variables, and to determine the 3rd root of unity $\omega_{\mathbf{1}, 0}$. 
As in the $\lambda_\alpha=1$ case, we get equations of these variables by evaluating the diagrams:
\[ 
\tikzmath{
\draw[mid>] (-.2,1) -- (-.2,1.5);
\draw[mid>] (-.2,-1.5) -- (-.2,-1);
\draw[mid>] (-.2,-1) -- (0,-.4);
\draw[mid>] (0,.4) -- (-.2,1);
\draw[mid>] (0,-.4) arc (240:120:.47cm);
\draw[mid>] (0,-.4) arc (-60:60:.47cm);
\draw[mid>] (-.2,-1) arc (240:120:1.15cm);
\filldraw (0,-.4) node[left]{$\scriptstyle \ell'$} circle (.05cm);
\filldraw (0,.4) node[left]{$\scriptstyle \ell$} circle (.05cm);
\filldraw (-.2,-1) node[left]{$\scriptstyle k'$} circle (.05cm);
\filldraw (-.2,1) node[left]{$\scriptstyle k$} circle (.05cm);
}
\,,\qquad
\tikzmath{
\draw[mid>] (-.2,1) -- (.1,1.5);
\draw[thick, red, mid>] (-.2,1) -- (-.5,1.5);
\draw[mid>] (.1,-1.5) -- (-.2,-1);
\draw[thick, red, mid>] (-.5,-1.5) -- (-.2,-1);
\draw[mid>] (-.2,-1) .. controls ++(45:.3cm) and ++(-45:.3cm) .. (0,-.4);
\draw[mid>] (0,.4) .. controls ++(45:.3cm) and ++(-45:.3cm) .. (-.2,1);
\draw[thick, red] (0,-.4) .. controls ++(-135:.9cm) and ++(-90:.5cm) .. (-1,0) .. controls ++(90:.5cm) and ++(135:.9cm) .. (0,.4);
\draw[mid>] (0,-.4) arc (240:120:.47cm);
\draw[mid>] (0,-.4) arc (-60:60:.47cm);
\draw[mid>] (-.2,-1) arc (240:120:1.15cm);
\filldraw (0,-.4) node[left]{$\scriptstyle \ell'$} circle (.05cm);
\filldraw (0,.4) node[left]{$\scriptstyle \ell$} circle (.05cm);
\filldraw (-.2,-1) node[left]{$\scriptstyle k'$} circle (.05cm);
\filldraw (-.2,1) node[left]{$\scriptstyle k$} circle (.05cm);
}
\,,\qquad
\tikzmath{
\draw[mid>] (-.2,1) -- (.1,1.5);
\draw[thick, red, mid>] (-.2,1) -- (-.5,1.5);
\draw[mid>] (.1,-1.5) -- (-.2,-1);
\draw[thick, red, mid>] (-.5,-1.5) -- (-.2,-1);
\draw[mid>] (-.2,-1) -- (0,-.4);
\draw[mid>] (0,.4) -- (-.2,1);
\draw[mid>] (0,-.4) arc (240:120:.47cm);
\draw[mid>] (0,-.4) arc (-60:60:.47cm);
\draw[mid>] (-.2,-1) arc (240:120:1.15cm);
\filldraw (0,-.4) node[left]{$\scriptstyle \ell'$} circle (.05cm);
\filldraw (0,.4) node[left]{$\scriptstyle \ell$} circle (.05cm);
\filldraw (-.2,-1) node[left]{$\scriptstyle k'$} circle (.05cm);
\filldraw (-.2,1) node[left]{$\scriptstyle k$} circle (.05cm);
}
\,,\qquad
\tikzmath{
\draw[mid>] (-.2,1) -- (-.2,1.5);
\draw[mid>] (.1,-1.5) -- (-.2,-1);
\draw[mid>] (-.2,-1) .. controls ++(45:.3cm) and ++(-45:.3cm) .. (0,-.4);
\draw[mid>] (0,.4) -- (-.2,1);
\draw[thick, red] (0,-.4) .. controls ++(-135:.9cm) and ++(-90:.5cm) .. (-1,-.2) .. controls ++(90:.5cm) and ++(90:.5cm) .. (.6,-.2) .. controls ++(-90:1cm) and ++(-135:1cm) .. (-.2,-1);
\draw[mid>] (0,-.4) arc (240:120:.47cm);
\draw[mid>] (0,-.4) arc (-60:60:.47cm);
\draw[mid>] (-.2,-1) arc (240:120:1.15cm);
\filldraw (0,-.4) node[left]{$\scriptstyle \ell'$} circle (.05cm);
\filldraw (0,.4) node[left]{$\scriptstyle \ell$} circle (.05cm);
\filldraw (-.2,-1) node[left]{$\scriptstyle k'$} circle (.05cm);
\filldraw (-.2,1) node[left]{$\scriptstyle k$} circle (.05cm);
}
\,,\qquad 
\tikzmath{
\draw (0,-.7) arc (-90:0:.3cm);
\draw[mid>] (.3,-.4) -- (.3,-.3) arc (0:90:.3cm);
\draw (0,-.7) arc (270:180:.3cm); 
\draw[mid>] (-.3,-.4) -- (-.3,-.3) arc (180:90:.3cm);
\draw[thick] (-.3,-.4) -- (-.2,-.4);
\draw[mid>] (-.9,-1) -- (-.9,0) arc (180:90:.45cm);
\draw[mid>] (0,0) arc (0:90:.45cm);
\draw[mid>] (-.45,.45) -- (-.45,1);
\filldraw (0,0) node[left, yshift=.1cm]{$\scriptstyle \ell$} circle (.05cm);
\filldraw (-.45,.45) node[left, yshift=.1cm]{$\scriptstyle k$} circle (.05cm);
}
\,,\qquad
\tikzmath{
\draw (0,-.7) arc (-90:0:.3cm);
\draw[mid>] (.3,-.4) -- (.3,-.3) to[out=90,in=-45] (0,0);
\draw (0,-.7) arc (270:180:.3cm); 
\draw[mid>] (-.3,-.4) -- (-.3,-.3) to[out=90,in=-135] (0,0);
\draw[thick] (-.3,-.4) -- (-.2,-.4);
\draw[mid>] (-.6,-1) -- (-.6,0) to[out=90, in=-135] (-.4,1);
\draw[mid>] (0,0) to[out=45,in=-135] (.2,.2) to[out=45, in=-45] (-.4,1);
\draw[thick, red] (0,0) .. controls ++(135:.5cm) and ++(90:.5cm) .. (.5,-.3) arc (0:-90:.6cm) to[out=180,in=-90] (-.8,0) to[out=90,in=180] (-.4,.7) arc (-90:0:.4cm) .. controls ++(90:.5cm) and ++(135:.5cm) .. (-.4,1);
\draw[mid>] (-.4,1) -- (.1,1.5);
\filldraw (0,0) node[left]{$\scriptstyle \ell$} circle (.05cm);
\filldraw (-.4,1) node[left]{$\scriptstyle k$} circle (.05cm);
}
\]  
in two ways (see Footnote \ref{footnote:EvaluateInTwoWays}). 
This gives us the equations:
\begin{align*}
    \sum_{i,j}A^{i,k}_{j,\ell}\overline{A^{i,k'}_{j,\ell'}} +  (-1)^{\ell+\ell'}\sum_{i,j}D^{i,k}_{j,\ell}\overline{D^{i,k'}_{j,\ell'}} - (2-\sqrt{5})\delta_{k,\ell}\delta_{k',\ell'}
    &= 
    \delta_{k,k'}\delta_{\ell,\ell'} 
    \\
     \sum_{i,j}\widehat{A}^{i,k}_{j,\ell}\overline{\widehat{A}^{i,k'}_{j,\ell'}} +  (-1)^{\ell+\ell'}\sum_{i,j}D^{1-\ell',i}_{1-k',j}\overline{D^{1-\ell,i}_{1-k,j}} - (2-\sqrt{5})\delta_{k,\ell}\delta_{k',\ell'}
     &= 
     \delta_{k,k'}\delta_{\ell,\ell'} 
     \\
      \sum_{i,j}D^{i,j}_{k,\ell}\overline{D^{i,j}_{k',\ell'}} +  (-1)^{\ell+\ell'}\sum_{i,j}D^{i,j}_{1-k',1-\ell'}\overline{D^{i,j}_{1-k,1-\ell}}
      &= 
      \delta_{k,k'}\delta_{\ell,\ell'}
 \\
     \sum_{i,j}(-1)^iA^{i,j}_{k,\ell}\overline{D^{1-\ell', i}_{k', 1-j} }+ (-1)^{\ell+1}\sum_{i,j}D^{j,1-k}_{1-i,\ell}\overline{\widehat{A}^{i,j}_{k,\ell}}
     - \omega_{\mathbf{1}, 0}^2(\sqrt{5}-2)\mathbf{i}\delta_{k,1-\ell}\delta_{k',1-\ell'}
     &= 
     0.
     \\
    \sum_{i}A^{i,k}_{i,\ell} + (-1)^{\ell+1}\mathbf{i}\sum_iD^{i,k}_{i,\ell} - (2-\sqrt{5})\delta_{k,\ell}&=0
    \\
    \sum_{i}\widehat{A}^{i,k}_{i,\ell} + \mathbf{i}\sum_i (-1)^{i+1}D^{1-\ell,i}_{1-k,i} - (2-\sqrt{5})\delta_{k,\ell}&=0
\end{align*}
In terms of our free variables, this gives us the equations:
\begin{align}
    \Im(d_0) = \frac{1-\sqrt{5}}{4}&=-\frac{\omega_{\mathbf{1},0}^2}{1 + \omega_{\mathbf{1},0}^2} \frac{-1+\sqrt{5}}{2}, \label{eq:-1eq1}\\
    r^2 + |d_0|^2 &= \frac{3-\sqrt{5}}{2},\label{eq:-1eq2}\\
    r^2 + |a_1|^2 + |d_1|^2 +|d_2|^2 =  r^2 + |\widehat{a}_1|^2 + |d_1|^2 +|d_2|^2 &= 1,\label{eq:-1eq3}\\
    |d_1|^2 = |d_2|^2& =     \frac{1}{2} - |d_0|^2,\label{eq:-1eq4}\\
    (\omega_{\mathbf{1},0} + \omega_{\mathbf{1},0}^2)r^2 - (d_0^2 + \overline{d_0}^2) &= 2-\sqrt{5},\label{eq:-1eq5}\\
    ra_1&= d_1d_2\label{eq:-1eq6}\\
    r\widehat{a}_1 &= -\overline{d}_1d_2\label{eq:-1eq7}\\
    d_2\overline{a_1}+ \widehat{a}_1\overline{d_2} = d_1\overline{a_1}+ \overline{\widehat{a}_1}\overline{d_1} &= 0\label{eq:-1eq8}
\end{align}

\begin{rem}
From Equation~\eqref{eq:-1eq1} we see that $  \frac{\omega_{\mathbf{1},0}^2}{1 + \omega_{\mathbf{1},0}^2} =1/2$, which implies that $\omega_{\mathbf{1}, 0}=1$.
\end{rem}
    
It is now straightforward to solve the above system of equations.
    
\begin{lem}
A general solution to Equations~\eqref{eq:-1eq1} -- \eqref{eq:-1eq8} is given by:
\begin{align*}
a_0 = \widehat{a}_0 & = \frac{3-\sqrt{5}}{4}
&
a_1 & = (3 + \sqrt{5})d_1d_2
&
\widehat{a}_1 & = -(3 + \sqrt{5})\overline{d_1}d_2
\\
|d_1|^2 = |d_2|^2 &= \frac{-1+\sqrt{5}}{8}
&
d_0 &= -\frac{1}{2} + \mathbf{i}\frac{1-\sqrt{5}}{4}.
\end{align*}
\end{lem}

With this lemma in hand, we can show the existence and uniqueness of the unitary fusion category with fusion ring $R(2)$.

\begin{thm}
\label{thm:2D2}
There exists a unique fusion category categorifying the ring $R(2)$ with $\lambda_\alpha = -1$. This unitary fusion category can be realised as the even part of the $2D2$ subfactor.
\end{thm}
\begin{proof}
Note that from Lemma~\ref{lem:FrobBasis}, we are free to re-scale our basis elements of $\cC_2(\rho\otimes\rho \to \rho)$ and $\cC_2(\rho\otimes\rho \to \alpha\rho)$ by
\begin{align*}
    \tikzmath{
    \draw[mid>] (-.3,-.6) -- (0,0);
    \draw[mid>] (.3,-.6) -- (0,0);
    \draw[mid<] (0,.6) -- (0,0);
    \filldraw (0,0) node[left]{$\scriptstyle 0$} circle (.05cm);
    }
    &
    \mapsto z_{\mathbf{1}} \tikzmath{
    \draw[mid>] (-.3,-.6) -- (0,0);
    \draw[mid>] (.3,-.6) -- (0,0);
    \draw[mid<] (0,.6) -- (0,0);
    \filldraw (0,0) node[left]{$\scriptstyle 0$} circle (.05cm);
    } \,,
    &
    \tikzmath{
    \draw[mid>] (-.3,-.6) -- (0,0);
    \draw[mid>] (.3,-.6) -- (0,0);
    \draw[mid<] (0,.6) -- (0,0);
    \filldraw (0,0) node[left]{$\scriptstyle 1$} circle (.05cm);
    }
    &
    \mapsto \overline{ z_{\mathbf{1}}} \tikzmath{
    \draw[mid>] (-.3,-.6) -- (0,0);
    \draw[mid>] (.3,-.6) -- (0,0);
    \draw[mid<] (0,.6) -- (0,0);
    \filldraw (0,0) node[left]{$\scriptstyle 1$} circle (.05cm);
    }  \,,
    &
    \tikzmath{
    \draw[thick, red, mid>] (0,0) -- (-.6,.6);
    \draw[mid>] (-.6,-.6) -- (0,0);
    \draw[mid>] (.6,-.6) -- (0,0);
    \draw[mid>] (0,0) -- (.6,.6);
    \filldraw (0,0) node[left]{$\scriptstyle  0$} circle (.05cm);
    } 
    &\mapsto z_{\alpha} 
    \tikzmath{
    \draw[thick, red, mid>] (0,0) -- (-.6,.6);
    \draw[mid>] (-.6,-.6) -- (0,0);
    \draw[mid>] (.6,-.6) -- (0,0);
    \draw[mid>] (0,0) -- (.6,.6);
    \filldraw (0,0) node[left]{$\scriptstyle 0$} circle (.05cm);
    } \,,
    & 
    \tikzmath{
    \draw[thick, red, mid>] (0,0) -- (-.6,.6);
    \draw[mid>] (-.6,-.6) -- (0,0);
    \draw[mid>] (.6,-.6) -- (0,0);
    \draw[mid>] (0,0) -- (.6,.6);
    \filldraw (0,0) node[left]{$\scriptstyle 1$} circle (.05cm);
    }
    &\mapsto \overline{z_{\alpha}}  
    \tikzmath{
    \draw[thick, red, mid>] (0,0) -- (-.6,.6);
    \draw[mid>] (-.6,-.6) -- (0,0);
    \draw[mid>] (.6,-.6) -- (0,0);
    \draw[mid>] (0,0) -- (.6,.6);
    \filldraw (0,0) node[left]{$\scriptstyle 1$} circle (.05cm);
    }
  \end{align*}
where $z_\mathbf{1},z_\alpha \in U(1)$. This re-scaling changes the phase of our free variables $d_1$ and $d_2$ by $z_{\mathbf{1}}^{-2}z_{\alpha}^2$ and $z_{\mathbf{1}}^{-2}z_{\alpha}^{-2}$ respectively. Thus we can arrange so that
\[  d_1 = d_2 = \mathbf{i}\frac{1}{2} \sqrt{\frac{1}{2} \left(-1+\sqrt{5}\right)}. \]
Hence, up to choice of our basis elements, we have a unique solution of all free parameters determining our category. 
Thus Proposition~\ref{prop:UniqueSD} gives that we have at most one unitary fusion category with fusion ring $R(2)$, and $\lambda_\alpha = -1$.

We know that the even part of the $2D2$ subfactor is a unitary fusion category with fusion ring $R(2)$; hence this must be the unique example.
\end{proof}
Let us write $\cC_2$ for the categorification of $R(2)$ we have classified in this section.
\begin{rem}
We wish to point out the above solutions to our free variables can be used to construct a system of dualizable endomorphisms of the Cuntz algebra $O_5\rtimes \bbZ/2\bbZ$. This gives an independent construction of the category $\cC_2$.
\end{rem}

To finish up, we connect $\cC_2$ to the even part of the $3^{\bbZ/4\bbZ}$ category.

\begin{cor}
There is a monoidal $\bbZ/2\bbZ$ action on $\cC_2$, such that equivariantisation by this action gives the $3^{\bbZ/4\bbZ}$ category of \cite{MR3827808,MR3402358}.
\end{cor}
\begin{proof}
Using the same gauge choice as in the previous Theorem, we can define an order two monoidal equivalence on $\cC_2$ by 
\begin{align*}
    \tikzmath{
\draw[mid>] (-.3,-.6) -- (0,0);
\draw[mid>] (.3,-.6) -- (0,0);
\draw[mid<] (0,.6) -- (0,0);
\filldraw (0,0) node[left]{$\scriptstyle 0$} circle (.05cm);
} \leftrightarrow \tikzmath{
\draw[mid>] (-.3,-.6) -- (0,0);
\draw[mid>] (.3,-.6) -- (0,0);
\draw[mid<] (0,.6) -- (0,0);
\filldraw (0,0) node[left]{$\scriptstyle 1$} circle (.05cm);
}, \qquad \tikzmath{
\draw[thick, red, mid>] (0,0) -- (-.6,.6);
\draw[mid>] (-.6,-.6) -- (0,0);
\draw[mid>] (.6,-.6) -- (0,0);
\draw[mid>] (0,0) -- (.6,.6);
\filldraw (0,0) node[left]{$\scriptstyle 0$} circle (.05cm);
} \leftrightarrow \tikzmath{
\draw[thick, red, mid>] (0,0) -- (-.6,.6);
\draw[mid>] (-.6,-.6) -- (0,0);
\draw[mid>] (.6,-.6) -- (0,0);
\draw[mid>] (0,0) -- (.6,.6);
\filldraw (0,0) node[left]{$\scriptstyle 1$} circle (.05cm);
},\quad \text{and} \quad \tikzmath{
\draw[thick, red, mid>] (0,-.5) -- (0,0);
\draw[thick, red, mid>] (0,.5) -- (0,0);
\draw[thick, red] (0,0) -- (.1,0);
} \leftrightarrow -\tikzmath{
\draw[thick, red, mid>] (0,-.5) -- (0,0);
\draw[thick, red, mid>] (0,.5) -- (0,0);
\draw[thick, red] (0,0) -- (.1,0);
}.
\end{align*}
By equivariantising by this order two monoidal auto-equivalence we obtain a unitary fusion category generated by the four morphisms 
\[  \tikzmath{
\draw[mid>] (-.3,-.6) -- (0,0);
\draw[mid>] (.3,-.6) -- (0,0);
\draw[mid<] (0,.6) -- (0,0);
\filldraw (0,0) node[left]{$\scriptstyle 0$} circle (.05cm);
}+ \tikzmath{
\draw[mid>] (-.3,-.6) -- (0,0);
\draw[mid>] (.3,-.6) -- (0,0);
\draw[mid<] (0,.6) -- (0,0);
\filldraw (0,0) node[left]{$\scriptstyle 1$} circle (.05cm);
},\qquad    \tikzmath{
\draw[thick, red, mid>] (0,0) -- (-.6,.6);
\draw[mid>] (-.6,-.6) -- (0,0);
\draw[mid>] (.6,-.6) -- (0,0);
\draw[mid>] (0,0) -- (.6,.6);
\filldraw (0,0) node[left]{$\scriptstyle 0$} circle (.05cm);
}+ \tikzmath{
\draw[thick, red, mid>] (0,0) -- (-.6,.6);
\draw[mid>] (-.6,-.6) -- (0,0);
\draw[mid>] (.6,-.6) -- (0,0);
\draw[mid>] (0,0) -- (.6,.6);
\filldraw (0,0) node[left]{$\scriptstyle 1$} circle (.05cm);
}, 
\qquad 
\tikzmath{
\draw[thick, red, mid>] (-.4,.2) -- (-.4,.6);
\draw[thick, red] (-.5,.2) -- (-.4,.2) arc (0:-180:.2cm);
\draw[thick, red, mid>] (-.8,.2) -- (-.8,.6);
\draw[mid>] (-.3,-.6) -- (0,0);
\draw[mid>] (.3,-.6) -- (0,0);
\draw[mid<] (0,.6) -- (0,0);
\filldraw (0,0) node[left]{$\scriptstyle 0$} circle (.05cm);
}
-
\tikzmath{
\draw[thick, red, mid>] (-.4,.2) -- (-.4,.6);
\draw[thick, red] (-.5,.2) -- (-.4,.2) arc (0:-180:.2cm);
\draw[thick, red, mid>] (-.8,.2) -- (-.8,.6);
\draw[mid>] (-.3,-.6) -- (0,0);
\draw[mid>] (.3,-.6) -- (0,0);
\draw[mid<] (0,.6) -- (0,0);
\filldraw (0,0) node[left]{$\scriptstyle 1$} circle (.05cm);
},
\quad \text{and}\quad 
\tikzmath{
\draw[thick, red, mid>] (0,0) -- (-.8,.8);
\draw[thick, red, mid>] (.2,.5) -- (.2,.8);
\draw[thick, red] (.1,.5) -- (.2,.5) arc (0:-180:.2cm);
\draw[thick, red, mid>] (-.2,.5) -- (-.2,.8);
\draw[mid>] (-.6,-.6) -- (0,0);
\draw[mid>] (.6,-.6) -- (0,0);
\draw[mid>] (0,0) -- (.8,.8);
\filldraw (0,0) node[left]{$\scriptstyle 0$} circle (.05cm);
}
-
\tikzmath{
\draw[thick, red, mid>] (0,0) -- (-.8,.8);
\draw[thick, red, mid>] (.2,.5) -- (.2,.8);
\draw[thick, red] (.1,.5) -- (.2,.5) arc (0:-180:.2cm);
\draw[thick, red, mid>] (-.2,.5) -- (-.2,.8);
\draw[mid>] (-.6,-.6) -- (0,0);
\draw[mid>] (.6,-.6) -- (0,0);
\draw[mid>] (0,0) -- (.8,.8);
\filldraw (0,0) node[left]{$\scriptstyle 1$} circle (.05cm);
}
,
\]
and the isomorphism
\[ 
\tikzmath{
\draw[red, thick, mid>] (0,0) -- (0,.4);
\draw[red, thick, mid>] (.4,0) -- (.4,.4) arc (0:180:.2cm);
\draw[red, thick] (0,.4) -- (.1,.4);
\draw[red, thick, mid>] (.6,0) -- (.6,.4);
\draw[red, thick, mid>] (1,0) -- (1,.4) arc (0:180:.2cm);
\draw[red, thick] (.6,.4) -- (.7,.4);
}
:
\alpha^{\otimes 4} \xrightarrow{\sim} \mathbf{1} .
\]
This is the presentation of the $3^{\bbZ/4\bbZ}$ category from \cite{MR3827808}.
\end{proof}

%%%%%%%%%%%%%%%%%%%%%%%%%%%%%%%%%%%%%%%%%%%%%%%%%%%%%%
%%%%%%%%%%%%%%%%%%%%%%%%%%%%%%%%%%%%%%%%%%%%%%%%%%%%%%
%%%%%%%%%%%%%%%%%%%%%%%%%%%%%%%%%%%%%%%%%%%%%%%%%%%%%%
\section{The non-self-dual case}
\label{sec:NSD}

%%%%%%%%%%%%%%%%%%%%%%%%%%%%%%%%%%%%%%%%%%%%%%%%%%%%%%%%%%%%%%%%%%%
In this section we focus on the unitary categorification of the fusion rings with four simple objects
$\mathbf{1},\alpha, \rho, \alpha\rho$
and fusion rules
\begin{equation}
\label{eq:S(m)}
\alpha \otimes \alpha\cong \mathbf{1}
\qquad\qquad
\rho\otimes \rho \cong \alpha \oplus m \rho \oplus m\alpha \rho.
\tag{$S(m)$}
\end{equation}
Let us write \ref{eq:S(m)} for such a fusion ring. 
By \cite{MR3229513} we know that \ref{eq:S(m)} has a categorification only if $m = 0,1, 2$.

Our main result of this section is as follows.
\begin{thm}\label{thm:NSD}
Let $\cD_m$ be a unitary fusion category with $K_0(\cC) \cong$ \ref{eq:S(m)}. 
Then either 
\begin{itemize}
\item 
$m=0$, in which case $\cD_0$ is equivalent to one of the four monoidally distinct categories $\Hilb(\bbZ/4\bbZ,\omega)$ where $\omega \in H^3(\bbZ/4\bbZ, \mathbb{C}^\times)$, or
\item
$m=1$, in which case $\cD_1$ is equivalent to the monoidally distinct even parts of the two complex conjugate subfactors with principal graphs $\mathcal{S}'$ from \cite{MR3306607,MR3827808}.
\end{itemize}
In particular the case $m=2$ from \cite[Thm.~1.1(6)]{MR3229513} is not categorifiable.
\end{thm}
\begin{proof}
The $m=0$ case is easily seen to be pointed, and hence the claim of the above theorem follows from \cite[Remark 4.10.4]{MR3242743}. 
Thus it suffices to restrict our attention to the cases of $m=1$ and $m=2$.

The general outline of this section follows for the most part as in the self-dual case. 
In \S\ref{NSD:NumericalData},
We begin by writing down a list of numerical data (essentially the 6j symbols of the category) which fully describe a unitary fusion category with fusion ring \ref{eq:S(m)}. 
In \S\ref{NSD:CenterAnalysis}, by studying the Drinfeld centre via the tube algebra of the category, we are able to deduce the precise values of some of this numerical data. 
To reduce the complexity of our numerical data, in \S\ref{NSD:Symmetries}, we use tetrahedral symmetries to essentially cut down the number of free variables in our numerical data by a factor of $24$. 
Finally, in \S\ref{NSD:Classification}, we solve for this numerical data by evaluating various morphisms in our categories in multiple ways to obtain equations. 

In the case $m=1$, we reduce our numerical data to two possible solutions, which shows there are at most two distinct unitary fusion categories categorifying $S(1)$. 
From the subfactor classification literature \cite{MR3306607,MR3827808} we know that two such categories exist.
We then show there are no solutions to the numerical data in the case $m=2$, and hence there are no such unitary fusion categories. 
\end{proof}

%%%%%%%%%%%%%%%%%%%%%%%%%%%%%%%%%%%%%%%%%%%%%%
\subsection{Numerical data}
\label{NSD:NumericalData}

We now produce a set of numerical data which completely describes a categorification of the ring $S(m)$. 
Let us write $\cD_m$ for such a unitary fusion category. 
We will show that the category $\cD_m$ can be described by the following data:
\begin{itemize}
    \item an 8th root of unity $\nu = e^{\pm i\pi\frac{1}{4}}$,
    \item $m$ choices of signs $\chi_i \in \{-1, 1\}$, 
    \item $m$ choices of 3rd roots of unity $\omega_i \in \{1, e^{2i\pi\frac{1}{3}},e^{2i\pi\frac{2}{3}}\}$, and
    \item $8m^4$ complex scalars $A^{i,j}_{k,\ell},B^{i,j}_{k,\ell},C^{i,j}_{k,\ell},D^{i,j}_{k,\ell},\widehat{A}^{i,j}_{k,\ell},\widehat{B}^{i,j}_{k,\ell},\widehat{C}^{i,j}_{k,\ell},\widehat{D}^{i,j}_{k,\ell} \in \mathbb{C}$ for $0 \leq i,j,k,\ell <m$. These complex scalars are the entries of the $F$-tensors $F^{\rho,\rho,\rho}_{\rho}$ and $F^{\rho,\rho,\rho}_{\alpha \rho}$
    \end{itemize}
    While the $128$ complex scalars in the $m=2$ case seems infeasible to deal with as is, we will use tetrahedral symmetries later on to reduce this $128$ to a more workable number.

To simplify notation, we define $d:=\dim(\rho)$, which is the largest solution to $d^2=1+2md$.
If $m=1$ then $d = 1+\sqrt{2}$, and if $m=2$ then $d = 2+\sqrt{5}$.
We pick orthonormal bases for the hom spaces
\[
\tikzmath{
\draw[mid>] (-.3,-.6) node[below]{$\scriptstyle \rho$} -- (0,0);
\draw[mid>] (.3,-.6) node[below]{$\scriptstyle \rho$} -- (0,0);
\draw[mid<, thick, blue] (0,.6) node[above]{$\scriptstyle \alpha$} -- (0,0);
\filldraw[blue] (0,0) circle (.05cm);
}
\in \cD_m(\rho\otimes \rho\to \alpha)
,\qquad\tikzmath{
\draw[mid>] (-.3,-.6) node[below]{$\scriptstyle \rho$} -- (0,0);
\draw[mid>] (.3,-.6) node[below]{$\scriptstyle \rho$} -- (0,0);
\draw[mid<] (0,.6) node[above]{$\scriptstyle \rho$} -- (0,0);
\filldraw (0,0) node[left]{$\scriptstyle i$} circle (.05cm);
}\in
\cD_m(\rho\otimes \rho\to \rho)
\qquad
\tikzmath{
\draw[thick, blue, mid>] (0,0) -- (-.6,.6) node[above]{$\scriptstyle \alpha$};
\draw[mid>] (-.6,-.6) node[below]{$\scriptstyle \rho$} -- (0,0);
\draw[mid>] (.6,-.6) node[below]{$\scriptstyle \rho$} -- (0,0);
\draw[mid>] (0,0) -- (.6,.6) node[above]{$\scriptstyle \rho$};
\filldraw (0,0) node[left]{$\scriptstyle i$} circle (.05cm);
}
\in\cD_m(\rho\otimes \rho\to \alpha\rho)
\qquad
0 \leq i \leq m
\]
so we have the local relations
\[
\tikzmath{
\draw[mid>] (-.3,-.6) node[below]{$\scriptstyle \rho$} -- (-.3,1.2);
\draw[mid>] (.3,-.6) node[below]{$\scriptstyle \rho$} -- (.3,1.2);
}
= 
\tikzmath{
\draw[mid>] (-.3,-.6) node[below]{$\scriptstyle \rho$} -- (0,0);
\draw[mid>] (.3,-.6) node[below]{$\scriptstyle \rho$} -- (0,0);
\draw[mid<, thick, blue] (0,.6) -- (0,0);
\filldraw[blue] (0,0) circle (.05cm);
\filldraw[blue] (0,.6) circle (.05cm);
\draw[mid<] (.3,1.2) node[above]{$\scriptstyle \rho$} -- (0,0.6);
\draw[mid<] (-.3,1.2) node[above]{$\scriptstyle \rho$} -- (0,0.6);
} 
+ 
\sum 
\tikzmath{
\draw[mid>] (-.3,-.6) node[below]{$\scriptstyle \rho$} -- (0,0);
\draw[mid>] (.3,-.6) node[below]{$\scriptstyle \rho$} -- (0,0);
\draw[mid<] (0,.6) -- (0,0);
\filldraw (0,0) node[left]{$\scriptstyle i$} circle (.05cm);
\filldraw (0,.6) node[left]{$\scriptstyle i$} circle (.05cm);
\draw[mid<] (.3,1.2) node[above]{$\scriptstyle \rho$} -- (0,0.6);
\draw[mid<] (-.3,1.2) node[above]{$\scriptstyle \rho$} -- (0,0.6);
}
+ 
\sum 
\tikzmath{
\draw[mid>] (-.3,-.6) node[below]{$\scriptstyle \rho$} -- (0,0);
\draw[mid>] (.3,-.6) node[below]{$\scriptstyle \rho$} -- (0,0);
\draw[mid>] (0,0) to[out=45,in=-45] (0,.6);
\draw[thick, blue, mid>] (0,0) to[out=135,in=-135] (0,.6);
\filldraw (0,0) node[left]{$\scriptstyle i$} circle (.05cm);
\filldraw (0,.6) node[left]{$\scriptstyle i$} circle (.05cm);
\draw[mid<] (.3,1.2) node[above]{$\scriptstyle \rho$} -- (0,0.6);
\draw[mid<] (-.3,1.2) node[above]{$\scriptstyle \rho$} -- (0,0.6);
}
\]
\begin{equation}
\label{eq:NSD-Spaghetti}
\left.
\begin{aligned}
\tikzmath{
\draw[thick, blue, mid>] (0,0) circle (.4cm);
}
&=
\tikzmath{
\draw[thick, blue, mid<] (0,0) circle (.4cm);
}
=
1
&&&
\tikzmath{
\draw[mid>] (0,0) circle (.4cm);
}
&=
\tikzmath{
\draw[mid<] (0,0) circle (.4cm);
}
=
d
\\
\tikzmath{
\draw[mid>] (-.3,-.6) node[below]{$\scriptstyle \rho$} -- (-.3,-.3) arc (180:90:.3cm);
\filldraw (0,0) node[left,yshift=.1cm]{$\scriptstyle i$} circle (.05cm);
\draw[mid>] (0,0) -- (0,.3) arc (180:0:.45cm) -- (.9,-.3) arc (0:-180:.3cm) arc (0:90:.3cm);
}
&=
\tikzmath{
\draw[mid>] (.3,-.6) node[below]{$\scriptstyle \rho$} -- (.3,-.3) arc (0:90:.3cm);
\filldraw (0,0) node[left,yshift=.12cm]{$\scriptstyle i$} circle (.05cm);
\draw[mid>] (0,0) -- (0,.3) arc (0:180:.45cm) -- (-.9,-.3) arc (-180:0:.3cm) arc (180:90:.3cm);
}
=0\,\,
&
\tikzmath{
\draw[mid>] (-.3,-.3) arc (180:90:.3cm);
\draw[mid>] (-.3,-.3) arc (-180:90:.3cm);
\draw[thick] (-.3,-.3) -- (-.2,-.3);
\filldraw (0,0) node[left,yshift=.1cm]{$\scriptstyle i$} circle (.05cm);
\draw[mid>] (0,0) -- (0,.4) node[above]{$\scriptstyle \rho$};
}
&=
0\,\,
&
\tikzmath{
\draw[mid>] (-.4,-.4) node[below]{$\scriptstyle \rho$} -- (0,0);
\draw[blue, thick, mid>] (0,0) -- (-.4,.4) node[above]{$\scriptstyle \alpha$};
\draw[mid>] (0,0) .. controls ++(45:.1cm) and ++(90:.7cm) ..(.4,0) .. controls ++(-90:.7cm) and ++(-45:.1cm) .. (0,0);
\filldraw (0,0) node[left]{$\scriptstyle i$} circle (.05cm);
}
&=
\tikzmath{
\draw[mid>] (.4,-.4) node[below]{$\scriptstyle \rho$} -- (0,0);
\draw[blue, thick] (0,0) -- (-.3,.4);
\draw[thick, blue, mid>] (-.3,.4) -- (-.6,.8) node[above]{$\scriptstyle \alpha$};
\draw[mid>] (0,0) .. controls ++(45:.7cm) and ++(90:.7cm) .. (-.4,0) .. controls ++(-90:.7cm) and ++(-45:.1cm) ..(0,0);
\filldraw (0,0) node[left]{$\scriptstyle i$} circle (.05cm);
}
=0\,\,
&
\tikzmath{
\draw[mid>] (-.4,-.4) to[out=90,in=-135] (0,0);
\draw[mid>] (-.4,-.4) arc (-180:0:.4cm) to[out=90,in=-45] (0,0);
\draw[thick] (-.4,-.4) -- (-.3,-.4);
\draw[mid>] (0,0) -- (.4,.4) node[above]{$\scriptstyle \rho$};
\draw[blue, thick, mid>] (0,0) -- (-.4,.4) node[above]{$\scriptstyle \alpha$};
\filldraw (0,0) node[left]{$\scriptstyle i$} circle (.05cm);
}
&=
0.
\end{aligned}
\,\,
\right\}
\end{equation}
We also choose unitary isomorphisms
 
\[ \tikzmath{
\draw[thick, blue, mid>] (0,-.5) node[below]{$\scriptstyle \alpha$} -- (0,0);
\draw[thick, blue, mid>] (0,.5) node[above]{$\scriptstyle \overline{\alpha}$} -- (0,0);
\draw[thick, blue] (0,0) -- (.1,0);
}\in \cD_m(\alpha \to \overline{\alpha})
\qquad \text{and} \qquad  
\tikzmath{
\draw[mid>] (-.5,-.5) node[below]{$\scriptstyle \rho$} -- (0,0);
\draw[mid>] (0,0) -- (.5,.5) node[above]{$\scriptstyle \rho$};
\draw[thick, blue, mid>] (.5,-.5) node[below]{$\scriptstyle \alpha$} -- (0,0);
\draw[thick, blue, mid<] (-.5,.5) node[above]{$\scriptstyle \alpha$} -- (0,0);
}
\in \cD_m(\rho \otimes \alpha \to \alpha \otimes \rho).  
\]
We normalise this last morphism so that
\[ 
\tikzmath{
\draw[mid>] (-.8,-.8) node[below]{$\scriptstyle \rho$} -- (0,0);
\draw[mid>] (0,0) -- (.5,.5) node[above]{$\scriptstyle \rho$};
\draw[thick, blue, mid<] (.5,-.5) -- (0,0);
\draw[thick, blue, mid>] (.8,-.8) node[below]{$\scriptstyle \alpha$} -- (.5,-.5);
\draw[thick, blue] (.5,-.5) -- (.62,-.38);
\draw[thick, blue, mid>] (-.5,.5) node[above]{$\scriptstyle \overline{\alpha}$} -- (0,0);
}
=
\tikzmath{
\draw[mid>] (-.5,-.5) node[below]{$\scriptstyle \rho$} -- (0,0);
\draw[mid>] (0,0) -- (.8,.8) node[above]{$\scriptstyle \rho$};
\draw[thick, blue, mid<] (-.5,.5) -- (0,0);
\draw[thick, blue, mid>] (-.8,.8) node[above]{$\scriptstyle \overline{\alpha}$} -- (-.5,.5);
\draw[thick, blue] (-.5,.5) -- (-.38,.62);
\draw[thick, blue, mid>] (.5,-.5) node[below]{$\scriptstyle \alpha$} -- (0,0);
}.
\]
We are still free to rescale the crossing up to sign.

Note that as pointed out in the proof of \cite[Theorem 5.8]{MR3229513}, we may assume that $\alpha$ has second Frobenius-Schur indicator $-1$, so
\[   
\tikzmath{
\draw[thick, blue, mid>] (0,-.5) node[below]{$\scriptstyle \alpha$} -- (0,0);
\draw[thick, blue, mid>] (0,.5) node[above]{$\scriptstyle \overline{\alpha}$} -- (0,0);
\draw[thick, blue] (0,0) -- (.1,0);
}
=
-
\tikzmath{
\draw[thick, blue, mid>] (0,-.5) node[below]{$\scriptstyle \alpha$} -- (0,0);
\draw[thick, blue, mid>] (0,.5) node[above]{$\scriptstyle \overline{\alpha}$} -- (0,0);
\draw[thick, blue] (0,0) -- (-.1,0);
}\,.
\]
Let $\mu$ be the scalar defined by
\begin{equation}
\label{eq:NSDExtraJellyfish}
\tikzmath{
\draw[thick, blue, mid>] (-.6,-.6) -- (-.6,.6);
\draw[mid>] (-.3,-.6) -- (0,0);
\draw[mid>] (.3,-.6) -- (0,0);
\draw[thick, blue, mid<] (0,.6) -- (0,0);
\filldraw[blue] (0,0) circle (.05cm);
}
=
\mu\,
\tikzmath{
\draw[thick, blue] (-60:.4cm) .. controls ++(-150:.2cm) and ++(90:.3cm) .. (-.7,-.8);
\draw[thick, blue, mid<] (.4,.6) -- (.4,0) arc (0:-60:.4cm);
\draw[mid>] (-.4,-.8) -- (0,0);
\draw[mid>] (.4,-.8) -- (0,0);
\draw[thick, blue, mid<] (0,.6) -- (0,0);
\filldraw[blue] (0,0) circle (.05cm);
}
.
\end{equation}
Note that from our normalisation we have that $\mu^2 = -1$.

\begin{lem}
Without loss of generality, we have the relation
\[
\tikzmath{
\draw[mid>] (.6,-.3) arc (-90:-180:.3cm) -- (.3,.8); 
\draw[mid>] (.6,-.3) arc (-90:0:.3cm) arc (180:90:.3cm); 
\draw[mid>] (1.5,-.8) -- (1.5,0) arc (0:90:.3cm);
\draw[thick, blue, mid>] (1.2,.3) arc (180:0:.3cm) .. controls ++(-90:.4cm) and ++(-90:.7cm) .. (.6,.-.3);
\filldraw[blue] (.6,-.3) circle (.05cm);
\filldraw[blue] (1.2,.3) circle (.05cm);
}
= 
\quad \frac{\nu}{d} \,
\tikzmath{
\draw[mid>] (0,-.5) -- (0,.5);
}
\qquad\qquad\text{where }
\nu=\exp(\pm \pi i/4)
\,.
\]
\end{lem}
\begin{proof}
First, by our normalizations for orthonormal bases of hom spaces, we observe that
$$
\tikzmath{
\draw[thick, blue, mid>] (-.4,.8) .. controls ++(90:.4cm) and ++(90:.7cm) .. (.4,.4) arc (0:-180:.2cm);
\draw[thick, blue, mid<] (-.4,-.8) .. controls ++(-90:.4cm) and ++(-90:.7cm) .. (.4,-.4) arc (0:180:.2cm);
\draw[mid>] (-.4,-.8) arc (-90:0:.2cm) arc (180:90:.2cm);
\draw[mid>] (0,.4) arc (-90:-180:.2cm) arc (0:90:.2cm);
\draw[mid>] (-.4,-.8) arc (-90:-180:.2cm) -- (-.6,.6) arc (180:90:.2cm);
\draw[mid>] (0,.4) arc (-90:0:.2cm) -- (.2,1.2) arc (180:0:.2cm) -- (.6,-1.2) arc (0:-180:.2cm) -- (.2,-.6) arc (0:90:.2cm);
\filldraw[blue] (-.4,-.8) circle (.05cm);
\filldraw[blue] (0,-.4) circle (.05cm);
\filldraw[blue] (0,.4) circle (.05cm);
\filldraw[blue] (-.4,.8) circle (.05cm);
}
=
\tikzmath{
\draw[thick, blue, mid>] (-.4,.8) .. controls ++(90:.4cm) and ++(90:.7cm) .. (.4,.4) -- (.4,-.4) .. controls ++(-90:.7cm) and ++(-90:.4cm) .. (-.4,-.8);
\draw[thick, blue, mid<] (0,-.4) -- (0,.4);
\draw[mid>] (-.4,-.8) arc (-90:0:.2cm) arc (180:90:.2cm);
\draw[mid>] (0,.4) arc (-90:-180:.2cm) arc (0:90:.2cm);
\draw[mid>] (-.4,-.8) arc (-90:-180:.2cm) -- (-.6,.6) arc (180:90:.2cm);
\draw[mid>] (0,.4) arc (-90:0:.2cm) arc (180:0:.2cm) -- (.6,-.6) arc (0:-180:.2cm) arc (0:90:.2cm);
\filldraw[blue] (-.4,-.8) circle (.05cm);
\filldraw[blue] (0,-.4) circle (.05cm);
\filldraw[blue] (0,.4) circle (.05cm);
\filldraw[blue] (-.4,.8) circle (.05cm);
}
=
\tikzmath{
\draw[thick, blue, mid>] (-.4,1.2) arc (0:180:.2cm) -- (-.8,-1.2) arc (-180:0:.2cm);
\draw[thick, blue, mid>] (0,-.4) -- (0,.4);
\draw[mid>] (-.4,-1.2) arc (-90:0:.2cm) -- (-.2,-.6) arc (180:90:.2cm);
\draw[mid>] (0,.4) arc (-90:-180:.2cm) -- (-.2,1) arc (0:90:.2cm);
\draw[mid>] (-.4,-1.2) arc (-90:-180:.2cm) -- (-.6,1) arc (180:90:.2cm);
\draw[mid>] (0,.4) arc (-90:0:.2cm) arc (180:0:.2cm) -- (.6,-.6) arc (0:-180:.2cm) arc (0:90:.2cm);
\filldraw[blue] (-.4,-1.2) circle (.05cm);
\filldraw[blue] (0,-.4) circle (.05cm);
\filldraw[blue] (0,.4) circle (.05cm);
\filldraw[blue] (-.4,1.2) circle (.05cm);
\draw[dotted, thick, red, rounded corners=5pt] (-.4,-.9) rectangle (.8,.9);
}
=
\frac{1}{d}\,
\tikzmath{
\draw[thick, blue, mid>] (-.4,1.2) arc (0:180:.2cm) -- (-.8,-1.2) arc (-180:0:.2cm);
\draw[mid>] (-.4,-1.2) arc (-90:0:.2cm) -- (-.2,1) arc (0:90:.2cm);
\draw[mid>] (-.4,-1.2) arc (-90:-180:.2cm) -- (-.6,1) arc (180:90:.2cm);
\filldraw[blue] (-.4,-1.2) circle (.05cm);
\filldraw[blue] (-.4,1.2) circle (.05cm);
}
=
\frac{1}{d}
\qquad
\Longrightarrow
\qquad
\tikzmath{
\draw[mid>] (.6,-.3) arc (-90:-180:.3cm) -- (.3,.8); 
\draw[mid>] (.6,-.3) arc (-90:0:.3cm) arc (180:90:.3cm); 
\draw[mid>] (1.5,-.8) -- (1.5,0) arc (0:90:.3cm);
\draw[thick, blue, mid>] (1.2,.3) arc (180:0:.3cm) .. controls ++(-90:.4cm) and ++(-90:.7cm) .. (.6,.-.3);
\filldraw[blue] (.6,-.3) circle (.05cm);
\filldraw[blue] (1.2,.3) circle (.05cm);
}
= 
\quad \frac{\nu}{d} \,
\tikzmath{
\draw[mid>] (0,-.5) -- (0,.5);
}
$$
for some unimodular scalar $\nu$.
By computing
\[ 
\frac{\nu}{d}
\tikzmath{
\draw[mid>] (-.3,-.6) -- (0,0);
\draw[mid>] (.3,-.6) -- (0,0);
\draw[mid<, thick, blue] (0,.6) -- (0,0);
\filldraw[blue] (0,0) circle (.05cm);
}
= 
\tikzmath{
\draw[mid>] (-.3,-.9) -- (-.3,-.3) arc (180:90:.3cm);
\draw[mid>] (.6,-.6) arc (-90:-180:.3cm) arc (0:90:.3cm); 
\draw[mid>] (.6,-.6) arc (-90:0:.3cm) -- (.9,0) arc (180:90:.3cm); 
\draw[mid>] (1.5,-.9) -- (1.5,0) arc (0:90:.3cm);
\draw[thick, blue, mid>] (0,0) -- (0,.8);
\draw[thick, blue, mid>] (1.2,.3) arc (180:0:.3cm) .. controls ++(-90:.4cm) and ++(-90:.7cm) .. (.6,.-.6);
\filldraw[blue] (0,0) circle (.05cm);
\filldraw[blue] (.6,-.6) circle (.05cm);
\filldraw[blue] (1.2,.3) circle (.05cm);
\draw[red, dotted, thick, rounded corners=5pt] (-.3,.1) rectangle (1.6,.5);
}
= 
\mu 
\tikzmath[xscale=-1]{
\draw[mid>] (-.3,-.9) -- (-.3,-.3) arc (180:90:.3cm);
\draw[mid>] (.6,-.6) arc (-90:-180:.3cm) arc (0:90:.3cm); 
\draw[mid>] (.6,-.6) arc (-90:0:.3cm) -- (.9,0) arc (180:90:.3cm); 
\draw[mid>] (1.5,-.9) -- (1.5,0) arc (0:90:.3cm);
\draw[thick, blue, mid>] (0,0) -- (0,.8);
\draw[thick, blue, mid>] (1.2,.3) arc (0:180:.3cm) .. controls ++(-90:.5cm) and ++(90:.4cm) .. (0,-.6) arc (-180:0:.3cm) (.6,.-.6);
\filldraw[blue] (0,0) circle (.05cm);
\filldraw[blue] (.6,-.6) circle (.05cm);
\filldraw[blue] (1.2,.3) circle (.05cm);
}
= 
\frac{\mu\overline{\nu}}{d}
\tikzmath{
\draw[mid>] (-.3,-.6) -- (0,0);
\draw[mid>] (.3,-.6) -- (0,0);
\draw[mid<, thick, blue] (0,.6) -- (0,0);
\filldraw[blue] (0,0) circle (.05cm);
}
\]
we find that $\nu^2 = \mu$. 
By re-scaling the crossing by a  sign, we may assume that 
$\nu = e^{\frac{\pm \pi i}{4}}$.
\end{proof}

In order to define natural orthornormal bases for the spaces $\cD_m(\rho\otimes \rho\to \rho)$ and $\cD_m(\rho\otimes \rho\to \alpha\rho)$ we define the operators
\[  
K^\mathbf{1}\left(\TI{i} \right) 
=
\tikzmath{
\draw[thick, blue, mid<] (0,0) circle (.4cm); 
\draw[mid>] (-.4,-.8) -- (0,0);
\draw[mid>] (.4,-.8) -- (0,0);
\draw[mid<] (0,.8) -- (0,0);
\filldraw (0,0) node[left]{$\scriptstyle i$} circle (.05cm);
}
\qquad \text{ and }\qquad   
K^\alpha\left(\Ta{i}\right) 
=
\tikzmath{
\draw[thick, blue, mid>] (0,0) .. controls ++(135:.7cm) and ++(90:.6cm) .. (.6,0) arc (0:-180:.6cm) -- (-.6,.8);
\draw[mid>] (-.8,-.8) -- (-.4,-.4);
\draw[mid>] (-.4,-.4) -- (0,0);
\draw[mid>] (.8,-.8) -- (.4,-.4);
\draw[mid>] (.4,-.4) -- (0,0);
\draw[mid>] (0,0) -- (.4,.4);
\draw[mid>] (.4,.4) -- (.8,.8);
\filldraw (0,0) node[left]{$\scriptstyle i$} circle (.05cm);
}
,
\]
and the Frobenius operators:
\begin{align*}
R^{\mathbf{1}}\left(\TI{i} \right) 
&= 
\sqrt{d}\,
\tikzmath{
\draw[mid>] (0,0) arc (-90:-180:.3cm) to[out=90,in=-90] (.6,1.3);
\draw[mid>] (0,0) arc (-90:0:.3cm) arc (180:0:.3cm);
\draw[mid>] (.9,-.4) -- (.9,.3) arc (0:90:.3cm);
\draw[mid>] (0,-.4) -- (0,0);
\draw[thick, blue, mid>] (.6,.6) to[out=90,in=-90] (-.3,1.3);
\filldraw[blue] (.6,.6) circle (.05cm);
\filldraw (0,0) node[left, yshift=-.1cm]{$\scriptstyle i$} circle (.05cm);
}
&
R^{\alpha}\left(\Ta{i}\right) 
&= 
\sqrt{d}\,
\tikzmath{
\draw[mid>] (0,0) arc (-90:-180:.3cm) -- (-.3,1);
\draw[mid>] (0,0) -- (.6,.6);
\draw[mid>] (.9,-.6) -- (.9,.3) arc (0:90:.3cm);
\draw[mid>] (.6,-.6) -- (0,0);
\draw[thick, blue, mid>] (.6,.6) arc (180:0:.3cm) .. controls ++(-90:.6cm) and ++(-135:1cm) .. (0,0);
\filldraw[blue] (.6,.6) circle (.05cm);
\filldraw (0,0) node[left]{$\scriptstyle i$} circle (.05cm);
}
\\
L^{\mathbf{1}}\left(\TI{i} \right) 
&= 
\sqrt{d}\,
\tikzmath[xscale=-1]{
\draw[mid>] (0,0) arc (-90:-180:.3cm) -- (-.3,1);
\draw[mid>] (0,0) arc (-90:0:.3cm) arc (180:0:.3cm);
\draw[mid>] (.9,-.4) -- (.9,.3) arc (0:90:.3cm);
\draw[mid>] (0,-.4) -- (0,0);
\draw[thick, blue, mid>] (.6,.6) -- (.6,1);
\filldraw[blue] (.6,.6) circle (.05cm);
\filldraw (0,0) node[left, yshift=-.1cm]{$\scriptstyle i$} circle (.05cm);
}
&
L^{\alpha}\left(\Ta{i}\right) 
&= 
\sqrt{d}\,
\tikzmath[xscale=-1]{
\draw (0,0) -- (-.3,.3);
\draw[mid>] (-.3,.3) to[out=135,in=-90] (-.4,1);
\draw[mid>] (0,0) -- (.6,.6);
\draw[mid>] (.9,-.4) -- (.9,.3) arc (0:90:.3cm);
\draw[mid>] (-.4,-.4) -- (0,0);
\draw[thick, blue, mid<] (0,0) .. controls ++(-45:.4cm) and ++(-90:.4cm) .. (-.3,0) .. controls ++(90:.4cm) and ++(90:.6cm) .. (.6,.6);
\filldraw[blue] (.6,.6) circle (.05cm);
\filldraw (0,0) node[left]{$\scriptstyle i$} circle (.05cm);
}
\end{align*}
on these spaces. 
These operators satisfy the following relations:
\begin{align*}
     K^{\mathbf{1}}\circ  K^{\mathbf{1}} 
     &= 
     1
     & 
     K^{\alpha}\circ  K^{\alpha} &= -1
     \\
     R^{ \alpha}\circ R^{\mathbf{1}} & = \nu 
     &
     R^{\mathbf{1}}\circ R^{\alpha }  
     &=
     \nu^{-1}\\
     L^{\alpha}\circ L^{\mathbf{1}} 
     & = 
     \nu^{-1} K^{\mathbf{1}}
     &
     L^{\mathbf{1}}\circ L^{\alpha } 
     &=
     \nu^{-1} K^{\alpha}
     \\
     K^{\alpha}\circ R^{ \mathbf{1}} &= \mu (R^{ \mathbf{1}} \circ K^{\mathbf{1}}) 
     &
    K^{\alpha}\circ L^{ \mathbf{1}} &= \mu (L^{\mathbf{1}} \circ K^{\mathbf{1}})\\
     (R^{\alpha }\circ L^{\mathbf{1} })^3 &= -1
     &
    (R^{\mathbf{1} }\circ L^{\alpha })^3 &= \mu.
\end{align*}
As a consequence of these relations, we can diagonalise the action of the operator $K^\mathbf{1}$, and set
\[
\Ta{i} := R^{\mathbf{1}}\left(\TI{i}\right)
\]
to obtain that there exist scalars $\chi_{i} \in \{-1,1\}$ and $\omega_{i}\in \{1, e^{2\pi i \frac{1}{3}}, e^{2\pi i \frac{2}{3}}\}$ such that
\begin{align*}
    K^\mathbf{1}\left(\TI{i}\right) & = \chi_i \TI{i}\qquad &&K^\alpha\left(\Ta{i}\right)  = \mu\chi_i \Ta{i}\\
    R^\mathbf{1}\left(\TI{i}\right) & =  \Ta{i}\qquad &&R^\alpha\left(\Ta{i}\right)  = \nu \TI{i}\\
    L^\mathbf{1}\left(\TI{i}\right) & =  -\nu\omega_i^2\Ta{i}\qquad &&L^\alpha\left(\Ta{i}\right)  = -\chi_i \omega_i^2 \TI{i}.
\end{align*}
In particular, this gives us the local relations
\begin{equation}
\label{eq:NSDalphaJellyfish}
\tikzmath{
\draw[thick, blue, mid>] (-.6,-.6) -- (-.6,.6);
\draw[mid>] (-.3,-.6) node[below]{$\scriptstyle \rho$} -- (0,0);
\draw[mid>] (.3,-.6) node[below]{$\scriptstyle \rho$} -- (0,0);
\draw[mid<] (0,.6) node[above]{$\scriptstyle \rho$} -- (0,0);
\filldraw (0,0) node[left]{$\scriptstyle i$} circle (.05cm);
}
=
\chi_{i} 
\tikzmath{
\draw[thick, blue] (-60:.4cm) .. controls ++(-150:.2cm) and ++(90:.3cm) .. (-.7,-.8);
\draw[thick, blue] (90:.4cm) arc (90:-60:.4cm);
\draw[thick, blue, mid>]  (90:.4cm) to[out=180,in=-90] (-.4,.6);
\draw[mid>] (-.4,-.8) -- (0,0);
\draw[mid>] (.4,-.8) -- (0,0);
\draw[mid<] (0,.6) -- (0,0);
\filldraw (0,0) node[left]{$\scriptstyle i$} circle (.05cm);
}
\qquad \text{ and } \qquad  
\tikzmath{
\draw[thick, blue, mid>] (-.9,-.6) -- (-.9,.6);
\draw[thick, blue, mid>] (0,0) -- (-.6,.6) node[above]{$\scriptstyle \alpha$};
\draw[mid>] (-.6,-.6) node[below]{$\scriptstyle \rho$} -- (0,0);
\draw[mid>] (.6,-.6) node[below]{$\scriptstyle \rho$} -- (0,0);
\draw[mid>] (0,0) -- (.6,.6) node[above]{$\scriptstyle \rho$};
\filldraw (0,0) node[left]{$\scriptstyle i$} circle (.05cm);
}   
= 
\mu\chi_i\,
\tikzmath{
\draw[thick, blue, mid>] (0,0) -- (-.8,.8);
\draw[thick, blue] (45:.6cm) arc (45:90:.6cm) arc (270:180:.2cm);
\draw[thick, blue, mid<] (45:.6cm) arc (45:-45:.6cm);
\draw[thick, blue] (-45:.6cm) arc (-45:-135:.6cm);
\draw[thick, blue] (225:.6) .. controls ++(135:.3cm) and ++(90:.3cm) .. (-1.1,-.8);
\draw[mid>] (-.8,-.8) -- (-.4,-.4);
\draw[mid>] (-.4,-.4) -- (0,0);
\draw[mid>] (.8,-.8) -- (.4,-.4);
\draw[mid>] (.4,-.4) -- (0,0);
\draw[mid>] (0,0) -- (.4,.4);
\draw[mid>] (.4,.4) -- (.8,.8);
\filldraw (0,0) node[left]{$\scriptstyle i$} circle (.05cm);
}\,.
\end{equation}

\begin{rem}\label{rem:RescaleNSD}
Note we are free to change our basis of $\cD_m(\rho\otimes \rho \to \rho)$ by a unitary which commutes with the operator $K^\mathbf{1}$. 
In particular, if $m=2$ and $\chi_0 = \chi_1$, then we are free to pick any other orthonormal basis of $\cD_m(\rho\otimes \rho \to \rho)$, and if $\chi_0 \neq \chi_1$ then we can only re-scale each basis vector by an element of $U(1)$.
\end{rem}

With this special choice of bases, we can determine the following local relations in $\cD_m$.

\begin{lem}
\label{lem:NSDJellyfishRelations}
There are scalars
\[
A^{i,j}_{k,\ell}, B^{i,j}_{k,\ell}, C^{i,j}_{k,\ell}, D^{i,j}_{k,\ell}, \widehat{A}^{i,j}_{k,\ell},\widehat{B}^{i,j}_{k,\ell},\widehat{C}^{i,j}_{k,\ell},\widehat{D}^{i,j}_{k,\ell}\in \mathbb{C} 
\qquad\qquad 
0 \leq i,j,k,\ell < m
\]
such that the following local relations hold in $\cD_m$:
\begin{align*}
\tikzmath{
\draw[mid>] (-.6,-.6) node[below]{$\scriptstyle \rho$} -- (-.6,.6);
\draw[mid>] (-.3,-.6) node[below]{$\scriptstyle \rho$} -- (0,0);
\draw[mid>] (.3,-.6) node[below]{$\scriptstyle \rho$} -- (0,0);
\draw[thick, blue, mid<] (0,.6) node[above]{$\scriptstyle \alpha$} -- (0,0);
\filldraw[blue] (0,0) circle (.05cm);
}
&=
\frac{\nu}{d}
\tikzmath{
\draw[mid>] (-.3,-.6) -- (0,0);
\draw[mid>] (.3,-.6) -- (0,0);
\draw[mid>] (.8,-.6) -- (0,.6);
\draw[thick, blue, mid<] (.6,.6) to[out=-135,in=90] (0,0);
\filldraw[blue] (0,0) circle (.05cm);
}
+
\frac{1}{\sqrt{d}}
\sum_{i=0}^1
\tikzmath{
\draw[thick, blue, mid<]  (.4,.6) .. controls ++(-90:.3cm) and ++(135:.5cm) .. (0,0);
\draw[mid>] (-.4,-.8) -- (0,0);
\draw[mid>] (.4,-1.2) -- (0,0);
\draw[mid>] (0,-1.2) -- (-.4,-.8);
\draw[mid>] (-.8,-1.2) -- (-.4,-.8);
\draw[mid<] (-.4,.6) .. controls ++(-90:.3cm) and ++(45:.5cm) ..  (0,0);
\filldraw (0,0) node[left]{$\scriptstyle i$} circle (.05cm);
\filldraw (-.4,-.8) node[left]{$\scriptstyle i$} circle (.05cm);
}
+
\frac{\nu}{\sqrt{d}}
\sum_{i=0}^1
\tikzmath{
\draw[thick, blue, mid<] (.4,.6) .. controls ++(-90:.2cm) and ++(90:.6cm) .. (.4,0) .. controls ++(-90:.4cm) and ++(135:.6cm) .. (-.4,-.8);
\draw[mid>] (-.4,-.8) -- (0,0);
\draw[mid>] (.4,-1.2) -- (0,0);
\draw[mid>] (0,-1.2) -- (-.4,-.8);
\draw[mid>] (-.8,-1.2) -- (-.4,-.8);
\draw[mid<] (0,.6) -- (0,0);
\filldraw (0,0) node[left]{$\scriptstyle i$} circle (.05cm);
\filldraw (-.4,-.8) node[left]{$\scriptstyle i$} circle (.05cm);
}
\\
\tikzmath{
\draw[mid>] (-.6,-.6) node[below]{$\scriptstyle \rho$} -- (-.6,.6);
\draw[mid>] (-.3,-.6) node[below]{$\scriptstyle \rho$} -- (0,0);
\draw[mid>] (.3,-.6) node[below]{$\scriptstyle \rho$} -- (0,0);
\draw[mid<] (0,.6) node[above]{$\scriptstyle \rho$} -- (0,0);
\filldraw (0,0) node[left]{$\scriptstyle \ell$} circle (.05cm);
}
&=
-
\frac{\omega_\ell}{\nu}\,
\tikzmath{
\draw[mid<] (-.6,.2) -- (-.2,-.2);
\draw[thick, blue, mid>] (-.4,-.8) -- (-.2,-.2);
\draw[mid>] (.4,-1.2) -- (-.2,-.2);
\draw[mid>] (0,-1.2) -- (-.4,-.8);
\draw[mid>] (-.8,-1.2) -- (-.4,-.8);
\draw[mid<] (.2,.2) -- (-.2,-.2);
\filldraw (-.2,-.2) node[left]{$\scriptstyle \ell$} circle (.05cm);
\filldraw[blue] (-.4,-.8) circle (.05cm);
}
-\nu \chi_\ell\omega_\ell^2
\tikzmath{
\draw[mid>] (-.2,.3) -- (-.6,.7);
\draw[mid>] (-.2,.3) -- (.2,.7);
\draw[thick, blue, mid>] (-.2,-.2) -- (-.2,.3);
\draw[mid>] (-.4,-.8) -- (-.2,-.2);
\draw[mid>] (.4,-1.2) -- (-.2,-.2);
\draw[mid>] (0,-1.2) -- (-.4,-.8);
\draw[mid>] (-.8,-1.2) -- (-.4,-.8);
\filldraw[blue] (-.2,-.2) circle (.05cm);
\filldraw[blue] (-.2,.3) circle (.05cm);
\filldraw (-.4,-.8) node[left]{$\scriptstyle \ell$} circle (.05cm);
}
+
\sum_{i,j,k}
A_{k,\ell}^{i,j}
\tikzmath{
\draw[mid>] (-.4,-.8) -- (-.2,-.2);
\draw[mid>] (.4,-1.2) -- (-.2,-.2);
\draw[mid>] (0,-1.2) -- (-.4,-.8);
\draw[mid>] (-.8,-1.2) -- (-.4,-.8);
\draw[mid<] (-.2,.4)-- (-.2,-.2);
\draw[mid>] (-.2,.4) -- (.2,.8);
\draw[mid>] (-.2,.4) -- (-.6,.8);
\filldraw (-.2,.4) node[left]{$\scriptstyle k$} circle (.05cm);
\filldraw (-.2,-.2) node[left]{$\scriptstyle j$} circle (.05cm);
\filldraw (-.4,-.8) node[left]{$\scriptstyle i$} circle (.05cm);
}
+
B_{k,\ell}^{i,j}
\tikzmath{
\draw[thick, blue, mid<] (-.2,.4) to[out=-135,in=135] (-.2,-.2);
\draw[mid>] (-.4,-.8) -- (-.2,-.2);
\draw[mid>] (.4,-1.2) -- (-.2,-.2);
\draw[mid>] (0,-1.2) -- (-.4,-.8);
\draw[mid>] (-.8,-1.2) -- (-.4,-.8);
\draw[mid<] (-.2,.4) to[out=-45,in=45] (-.2,-.2);
\draw[mid>] (-.2,.4) -- (.2,.8);
\draw[mid>] (-.2,.4) -- (-.6,.8);
\filldraw (-.2,.4) node[left]{$\scriptstyle k$} circle (.05cm);
\filldraw (-.2,-.2) node[left]{$\scriptstyle j$} circle (.05cm);
\filldraw (-.4,-.8) node[left]{$\scriptstyle i$} circle (.05cm);
}
+
C_{k,\ell}^{i,j}
\tikzmath{
\draw[thick, blue, mid<] (-.4,.4) .. controls ++(90:.3cm) and ++(90:.4cm) .. (.4,0) .. controls ++(-90:.4cm) and ++(135:.6cm) .. (-.4,-.8);
\draw[thick, blue, mid>] (0,0) to[out=135,in=-90] (-.4,.4);
\draw[thick, blue] (-.4,.4) -- (-.3,.4);
\draw[mid>] (-.4,-.8) -- (0,0);
\draw[mid>] (.4,-1.2) -- (0,0);
\draw[mid>] (0,-1.2) -- (-.4,-.8);
\draw[mid>] (-.8,-1.2) -- (-.4,-.8);
\draw[mid<] (0,.8) to[out=-90,in=45] (0,0);
\draw[mid>] (0,.8) -- (.4,1.2);
\draw[mid>] (0,.8) -- (-.4,1.2);
\filldraw (0,.8) node[left]{$\scriptstyle k$} circle (.05cm);
\filldraw (0,0) node[left]{$\scriptstyle j$} circle (.05cm);
\filldraw (-.4,-.8) node[left]{$\scriptstyle i$} circle (.05cm);
}
+
D_{k,\ell}^{i,j}
\tikzmath{
\draw[thick, blue, mid<] (0,.5) .. controls ++(-135:.6cm) and ++(90:.4cm) .. (.4,0) .. controls ++(-90:.4cm) and ++(135:.6cm) .. (-.4,-.8);
\draw[mid>] (-.4,-.8) -- (0,0);
\draw[mid>] (.4,-1.2) -- (0,0);
\draw[mid>] (0,-1.2) -- (-.4,-.8);
\draw[mid>] (-.8,-1.2) -- (-.4,-.8);
\draw[mid<] (0,.5) to[out=-45,in=45] (0,0);
\draw[mid>] (0,.5) -- (.4,.9);
\draw[mid>] (0,.5) -- (-.4,.9);
\filldraw (0,.5) node[left]{$\scriptstyle k$} circle (.05cm);
\filldraw (0,0) node[left]{$\scriptstyle j$} circle (.05cm);
\filldraw (-.4,-.8) node[left]{$\scriptstyle i$} circle (.05cm);
}
\\
\tikzmath{
\draw[mid>] (-.6,-.6) node[below]{$\scriptstyle \rho$} -- (-.6,.6);
\draw[mid>] (-.3,-.6) node[below]{$\scriptstyle \rho$} -- (0,0);
\draw[mid>] (.3,-.6) node[below]{$\scriptstyle \rho$} -- (0,0);
\draw[thick, blue, mid<] (-.3,.6) node[above]{$\scriptstyle \alpha$} -- (0,0);
\draw[mid<] (.3,.6) node[above]{$\scriptstyle \rho$} -- (0,0);
\filldraw (0,0) node[left]{$\scriptstyle \ell$} circle (.05cm);
}
&=
-\chi_\ell\omega_\ell\,
\tikzmath{
\draw[mid<] (-.6,.2) -- (-.2,-.2);
\draw[thick, blue, mid>] (-.4,-.8) .. controls ++(90:.2cm) and ++(-90:.3cm) .. (0,-.2) .. controls ++(90:.2cm) and ++(-90:.2cm) .. (-.2,.2);
\draw[mid>] (.4,-1.2) -- (-.2,-.2);
\draw[mid>] (0,-1.2) -- (-.4,-.8);
\draw[mid>] (-.8,-1.2) -- (-.4,-.8);
\draw[mid<] (.2,.2) -- (-.2,-.2);
\filldraw (-.2,-.2) node[left]{$\scriptstyle \ell$} circle (.05cm);
\filldraw[blue] (-.4,-.8) circle (.05cm);
}
-
\nu^3\omega_\ell^2\,
\tikzmath{
\draw[thick, blue, mid<] (0,.8) .. controls ++(-90:.3cm) and ++(-90:.3cm) .. (.4,.8);
\draw[thick, blue, mid>] (.4,.8) .. controls ++(90:.3cm) and ++(-90:.3cm) .. (0,1.2);
\draw[thick, blue] (.4,.8) -- (.3,.8);
\draw[thick, blue, mid>] (0,0) .. controls ++(90:.3cm) and ++(90:.4cm) .. (.4,0);
\draw[thick, blue, mid<] (.4,0) .. controls ++(-90:.4cm) and ++(135:.6cm) .. (-.4,-.8);
\draw[thick, blue] (.4,0) -- (.3,0);
\draw[mid>] (-.4,-.8) -- (0,0);
\draw[mid>] (.4,-1.2) -- (0,0);
\draw[mid>] (0,-1.2) -- (-.4,-.8);
\draw[mid>] (-.8,-1.2) -- (-.4,-.8);
\draw[mid>] (0,.8) -- (.4,1.2);
\draw[mid>] (0,.8) -- (-.4,1.2);
\filldraw[blue] (0,.8) circle (.05cm);
\filldraw[blue] (0,0) circle (.05cm);
\filldraw (-.4,-.8) node[left]{$\scriptstyle \ell$} circle (.05cm);
}
+
\sum_{i,j,k}
\widehat{A}_{k,\ell}^{i,j}
\tikzmath{
\draw[thick, blue, mid>] (-.4,-.8) .. controls ++(135:.5cm) and ++(-90:.4cm) .. (.2,0) .. controls ++(90:.3cm) and ++(-90:.3cm) .. (-.2,.8);
\draw[thick, blue, mid<] (-.2,.4) to[out=-135,in=135] (-.2,-.2);
\draw[mid>] (-.4,-.8) to[out=45,in=-135] (-.2,-.2);
\draw[mid>] (.4,-1.2) -- (-.2,-.2);
\draw[mid>] (0,-1.2) -- (-.4,-.8);
\draw[mid>] (-.8,-1.2) -- (-.4,-.8);
\draw[mid<] (-.2,.4) to[out=-45,in=45] (-.2,-.2);
\draw[mid>] (-.2,.4) -- (.2,.8);
\draw[mid>] (-.2,.4) -- (-.6,.8);
\filldraw (-.2,.4) node[left]{$\scriptstyle k$} circle (.05cm);
\filldraw (-.2,-.2) node[left]{$\scriptstyle j$} circle (.05cm);
\filldraw (-.4,-.8) node[left]{$\scriptstyle i$} circle (.05cm);
}
+
\widehat{B}_{k,\ell}^{i,j}
\tikzmath{
\draw[thick, blue, mid<] (0,.5) to[out=-135,in=90] (-.3,.2);
\draw[thick, blue, mid>] (-.3,.2) .. controls ++(-90:.2cm) and ++(-90:.4cm) .. (.4,.5) .. controls ++(90:.2cm) and ++(-90:.2cm) .. (0,.9);
\draw[thick, blue] (-.3,.2) -- (-.2,.2);
\draw[mid>] (-.4,-.8) -- (0,0);
\draw[mid>] (.4,-1.2) -- (0,0);
\draw[mid>] (0,-1.2) -- (-.4,-.8);
\draw[mid>] (-.8,-1.2) -- (-.4,-.8);
\draw[mid<] (0,.5) to[out=-45,in=90] (0,0);
\draw[mid>] (0,.5) -- (.4,.9);
\draw[mid>] (0,.5) -- (-.4,.9);
\filldraw (0,.5) node[left]{$\scriptstyle k$} circle (.05cm);
\filldraw (0,0) node[left]{$\scriptstyle j$} circle (.05cm);
\filldraw (-.4,-.8) node[left]{$\scriptstyle i$} circle (.05cm);
}
+
\widehat{C}_{k,\ell}^{i,j}
\tikzmath{
\draw[thick, blue, mid>] (-.4,-.8) .. controls ++(135:.5cm) and ++(-90:.4cm) .. (.2,0) .. controls ++(90:.3cm) and ++(-90:.3cm) .. (-.2,.8);
\draw[mid>] (-.4,-.8) to[out=45,in=-135] (-.2,-.2);
\draw[mid>] (.4,-1.2) -- (-.2,-.2);
\draw[mid>] (0,-1.2) -- (-.4,-.8);
\draw[mid>] (-.8,-1.2) -- (-.4,-.8);
\draw[mid<] (-.2,.4) -- (-.2,-.2);
\draw[mid>] (-.2,.4) -- (.2,.8);
\draw[mid>] (-.2,.4) -- (-.6,.8);
\filldraw (-.2,.4) node[left]{$\scriptstyle k$} circle (.05cm);
\filldraw (-.2,-.2) node[left]{$\scriptstyle j$} circle (.05cm);
\filldraw (-.4,-.8) node[left]{$\scriptstyle i$} circle (.05cm);
}
+
\widehat{D}_{k,\ell}^{i,j}
\tikzmath{
\draw[thick, blue, mid>] (-.2,.-.2) .. controls ++(135:.5cm) and ++(-90:.2cm) ..(.1,.4) .. controls ++(90:.3cm) and ++(-90:.3cm) .. (-.2,.8);
\draw[mid>] (-.4,-.8) -- (-.2,-.2);
\draw[mid>] (.4,-1.2) -- (-.2,-.2);
\draw[mid>] (0,-1.2) -- (-.4,-.8);
\draw[mid>] (-.8,-1.2) -- (-.4,-.8);
\draw[mid<] (-.2,.4) to[out=-90,in=45] (-.2,-.2);
\draw[mid>] (-.2,.4) -- (.2,.8);
\draw[mid>] (-.2,.4) -- (-.6,.8);
\filldraw (-.2,.4) node[left]{$\scriptstyle k$} circle (.05cm);
\filldraw (-.2,-.2) node[left]{$\scriptstyle j$} circle (.05cm);
\filldraw (-.4,-.8) node[left]{$\scriptstyle i$} circle (.05cm);
}
\end{align*}
\end{lem}
\begin{proof}
The proof is omitted as it is nearly identical to the proof of Lemma~\ref{lem:SDJellyfishRelations}.
\end{proof}

\begin{rem}
As in the self-dual case described in Remark \ref{rem:SD6j}, the above complex scalars are precisely entries of certain $F$-tensors of $\cD_m$. 
\end{rem}

With these local relations, we can show that our described numerical data fully determines the category $\cD_m$.
\begin{prop}
\label{prop:UniqueNSD}
There is at most one unitary fusion category $\cD_m$ realising each tuple of data 
\[(\nu, \chi, \omega, A, B, C, D, \widehat{A},\widehat{B},\widehat{C},\widehat{D}).\]
\end{prop}
\begin{proof}
We omit the proof which is nearly identical to the proof of Proposition \ref{prop:UniqueSD}
replacing
\eqref{eq:SD-Spaghetti} with \eqref{eq:NSD-Spaghetti},
the jellyfish relations from Lemmas \ref{lem:alphaJellyfish} and  \ref{lem:SDJellyfishRelations} with those from \eqref{eq:NSDalphaJellyfish}, \eqref{eq:NSDExtraJellyfish}, and Lemma \ref{lem:NSDJellyfishRelations},
and using absorption relations similar to Lemma \ref{lem:AbsorptionRelations}.
\end{proof}
%%%%%%%%%%%%%%%%%%%%%%%%%%%%%%%%%%%%%%%%%%%%%%%%%%%%%%%%%%%%%%%%%%%
\subsection{Centre Analysis}
\label{NSD:CenterAnalysis}

As in the self-dual case, we study the centre of $\cD_m$ in order to determine information about our free variables. We restrict our attention to the case of $m=2$, as this is the most difficult case, and we need as much information about our numerical data as possible in order to make progress on the classification. While we could repeat the analysis for $m=1$, this is unnecessary as in this case the lack of multiplicity makes it easy to solve for our numerical data.

Our main result of this section is as follows.

\begin{lem}\label{lem:centreNSD}
If $m=2$ and $\chi_0 = \chi_1$, then
\[ \sum_{i} A^{i,0}_{i,0} = \frac{(2+i)-\sqrt{5}}{2} ,\qquad \sum_i A^{i,1}_{i,1} = \frac{(2-i)-\sqrt{5}}{2}, \quad \text{and}\quad \sum_i A^{i,1}_{i,0} =\sum_i A^{i,0}_{i,1} =0.\]
\end{lem}
In the case of $\chi_1 = \chi_0$, knowing the above information about the free variables $A^{i,j}_{k,l}$ will be the key starting point in showing non-existence of the category $\cD_{2}$ later on in this paper.

To show this result we study the tube algebra of $\cD_{2}$. As in the self-dual case, we only study a small sub-algebra. 
We choose the following bases:
\begin{align*}
    A_{\mathbf{1}\to \mathbf{1}} &=\operatorname{span} \left\{\,\,  
    \tikzmath{\filldraw[very thick, fill=gray!30] (0,0) circle (.2cm);}
    \,,\,\,
    \tikzmath{\filldraw[very thick, fill=gray!30] (0,0) circle (.2cm);
    \draw[thick, blue, mid<] (0,0) circle (.4cm);}
    \,,\,\,
    \tikzmath{\filldraw[very thick, fill=gray!30] (0,0) circle (.2cm);
    \draw[mid<] (0,0) circle (.4cm);}
    \,,\,\,
    \tikzmath{\filldraw[very thick, fill=gray!30] (0,0) circle (.2cm);
    \draw[blue, thick, mid<] (0,0) circle (.65cm);
    \draw[mid<] (0,0) circle (.4cm);}\,\,  
    \right\}
    \\
    A_{\mathbf{1}\to \rho} 
    &= 
    \operatorname{span}\left\{\,\,
    \tikzmath{\filldraw[very thick, fill=gray!30] (0,-.2) circle (.2cm);
    \draw[mid<] (0,0) ellipse (.4cm and .6cm);
    \draw[mid>] (0,0) .. controls ++(90:.2cm) and ++(-45:.3cm) .. (135:.46cm);
    \filldraw (135:.46cm) node[left]{$\scriptstyle 0$} circle (.05cm);
    }
    \,,\,\,
    \tikzmath{\filldraw[very thick, fill=gray!30] (0,-.2) circle (.2cm);
    \draw[mid<] (0,0) ellipse (.4cm and .6cm);
    \draw[mid>] (0,0) .. controls ++(90:.2cm) and ++(-45:.3cm) .. (135:.46cm);
    \filldraw (135:.46cm) node[left]{$\scriptstyle 1$} circle (.05cm);
    }\,,\,\,     
    \tikzmath{\filldraw[very thick, fill=gray!30] (0,-.2) circle (.2cm);
    \draw[mid<] (0,0) ellipse (.4cm and .6cm);
    \draw[thick, blue, mid<] (0,0) ellipse (.7cm and .9cm);
    \draw[mid>] (0,0) .. controls ++(90:.2cm) and ++(-45:.3cm) .. (135:.46cm);
    \filldraw (135:.46cm) node[left]{$\scriptstyle 0$} circle (.05cm);
    }
    \,,\,\,
    \tikzmath{\filldraw[very thick, fill=gray!30] (0,-.2) circle (.2cm);
    \draw[mid<] (0,0) ellipse (.4cm and .6cm);
    \draw[thick, blue, mid<] (0,0) ellipse (.7cm and .9cm);
    \draw[mid>] (0,0) .. controls ++(90:.2cm) and ++(-45:.3cm) .. (135:.46cm);
    \filldraw (135:.46cm) node[left]{$\scriptstyle 1$} circle (.05cm);
    }\,\,   
    \right\}\\
    A_{\alpha\to \alpha} 
    &= 
    \operatorname{span}\left\{\,\,  
    \tikzmath{\filldraw[very thick, fill=gray!30] (0,0) circle (.2cm);
    \draw[blue, thick, mid>] (0,.2) -- (0,.8);}\,,\,\, 
    \tikzmath{\filldraw[very thick, fill=gray!30] (0,0) circle (.2cm);
    \draw[blue, thick, mid>] (0,.2) arc (180:0:.2cm) -- (.4,0) arc (0:-180:.4cm) .. controls ++(90:.3cm) and ++(270:.3cm) .. (0,.8);}\,,\,\,
    \tikzmath{\filldraw[very thick, fill=gray!30] (0,0) circle (.2cm);
    \draw[blue, thick, mid>] (0,.2) -- (0,.6);
    \draw[mid<] (0,0) circle (.6cm);
    \draw[thick, blue] (0,.6) -- (0,1);
    \draw[thick, blue, mid>] (0,1) -- (0,1.4);
    }\,,\,\,
    \tikzmath{\filldraw[very thick, fill=gray!30] (0,0) circle (.2cm);
    \draw[blue, thick, mid>] (0,.2) arc (180:0:.2cm) -- (.4,0) arc (0:-180:.4cm) .. controls ++(90:.3cm) and ++(270:.3cm) .. (0,.8);
    \draw[thick, blue,mid>] (0,.8) -- (0,1.2);
    \draw[mid<] (0,0) circle (.6cm);}\,\,
    \right\}
    .
\end{align*}
By direct computation we obtain that: 
\begin{enumerate}[label=(\arabic*)]
    \item The irreducible representations of $A_{\mathbf{1}\to \mathbf{1}}$ are: 
$$
    \begin{array}{c|c c c c}
         &
         \tikzmath{\filldraw[very thick, fill=gray!30] (0,0) circle (.1cm);}
         &
        \tikzmath{\filldraw[very thick, fill=gray!30] (0,0) circle (.1cm);
        \draw[thick, blue, mid<] (0,0) circle (.2cm);}
        &
        \tikzmath{\filldraw[very thick, fill=gray!30] (0,0) circle (.1cm);
        \draw[mid<] (0,0) circle (.2cm);}
        &
        \tikzmath{\filldraw[very thick, fill=gray!30] (0,0) circle (.1cm);
        \draw[blue, thick, mid<] (0,0) circle (.35cm);
        \draw[mid<] (0,0) circle (.2cm);}
         \\ \hline
         \chi_0 &1 & 1& 2 + \sqrt{5}&2 + \sqrt{5} \\
         \chi_1 &1 & 1& 2 - \sqrt{5}& 2 - \sqrt{5}\\
         \chi_2 &1 &-1 &\mathbf{i}&-\mathbf{i} \\
         \chi_3 &1 & -1&-\mathbf{i} & \mathbf{i}
    \end{array}
$$
    Hence $\cI(\mathbf{1})$ contains 4 simple objects $X_i$ with dimensions
\[ \dim(X_0) = 1 , \quad \dim(X_1) = 9 + 4\sqrt{5},\quad  \text{and} \quad \dim(X_2) = \dim(X_3) = 5 + 2\sqrt{5}.\]
\item  The irreducible representations of $A_{\alpha\to \alpha}$ are:
\[  
\begin{array}{c|c c c c}
    & 
    \tikzmath{\filldraw[very thick, fill=gray!30] (0,0) circle (.2cm);
    \draw[blue, thick, mid>] (0,.2) -- (0,.8);}
    &
    \tikzmath{\filldraw[very thick, fill=gray!30] (0,0) circle (.2cm);
    \draw[blue, thick, mid>] (0,.2) arc (180:0:.2cm) -- (.4,0) arc (0:-180:.4cm) .. controls ++(90:.3cm) and ++(270:.3cm) .. (0,.8);}
    &
    \tikzmath{\filldraw[very thick, fill=gray!30] (0,0) circle (.2cm);
    \draw[blue, thick, mid>] (0,.2) -- (0,.6);
    \draw[mid<] (0,0) circle (.6cm);
    \draw[thick, blue] (0,.6) -- (0,1);
    \draw[thick, blue, mid>] (0,1) -- (0,1.4);
    \filldraw[white] (0,-.75) circle (.01cm);
    }
    &
    \tikzmath{\filldraw[very thick, fill=gray!30] (0,0) circle (.2cm);
    \draw[blue, thick, mid>] (0,.2) arc (180:0:.2cm) -- (.4,0) arc (0:-180:.4cm) .. controls ++(90:.3cm) and ++(270:.3cm) .. (0,.8);
    \draw[thick, blue,mid>] (0,.8) -- (0,1.2);
    \draw[mid<] (0,0) circle (.6cm);}
    \\\hline
         \tau_{0} &1 &   \mathbf{i}& \frac{a(1 + \mathbf{i}\mu) +  \sqrt{2a^2(  1 + \mathbf{i}\mu) + 4\mathbf{i}\mu}}{2}& \mathbf{i}\frac{a(1 + \mathbf{i}\mu) +  \sqrt{2a^2(  1 + \mathbf{i}\mu) + 4\mathbf{i}\mu}}{2} \\
         \tau_{1} &1 &   -\mathbf{i}& \frac{a(1 - \mathbf{i}\mu) +  \sqrt{2a^2(  1 - \mathbf{i}\mu)- 4\mathbf{i}\mu}}{2}& -\mathbf{i}\frac{a(1 - \mathbf{i}\mu) +  \sqrt{2a^2(  1 - \mathbf{i}\mu)- 4\mathbf{i}\mu}}{2} \\
            \tau_{2} &1 &   \mathbf{i}& \frac{a(1 + \mathbf{i}\mu) -  \sqrt{2a^2(  1 + \mathbf{i}\mu) + 4\mathbf{i}\mu}}{2}& \mathbf{i}\frac{a(1 + \mathbf{i}\mu) -  \sqrt{2a^2(  1 + \mathbf{i}\mu) + 4\mathbf{i}\mu}}{2}
            \\  \tau_{3} &1 &   -\mathbf{i}& \frac{a(1 - \mathbf{i}\mu) -  \sqrt{2a^2(  1 - \mathbf{i}\mu)- 4\mathbf{i}\mu}}{2}& -\mathbf{i}\frac{a(1 - \mathbf{i}\mu) -  \sqrt{2a^2(  1 - \mathbf{i}\mu)- 4\mathbf{i}\mu}}{2}
\end{array}
\]
where $a :=\chi_0  +\chi_1 \in \{0,\pm2\}$. 
Hence $\cI(\alpha)$ contains 4 simple objects $Y_i$ with dimensions
\[
\frac{20 +8\sqrt{5}}{ 2 + \frac{1}{2}\left| a(1 + \mathbf{i}\mu) \pm \sqrt{2a^2(1 + \mathbf{i}\mu) + 4\mathbf{i}\mu}  \right |^2   }
\quad\text{and}\quad 
\frac{20 +8\sqrt{5}}{ 2 + \frac{1}{2}\left| a(1 - \mathbf{i}\mu) \pm \sqrt{2a^2(1 - \mathbf{i}\mu) - 4\mathbf{i}\mu}  \right |^2   }  
\]

\item 
Let $_\mathbf{1} \pi _\rho$ be the action of $A_{\mathbf{1} \to \mathbf{1}}$ on  $A_{\mathbf{1} \to \rho}$. Then 
\[
\resizebox{.9\textwidth}{!}{$   
{}_\mathbf{1}\pi_\rho
\left(\,
    \tikzmath{\filldraw[very thick, fill=gray!30] (0,0) circle (.2cm);
    \draw[thick, blue, mid<] (0,0) circle (.4cm);}
\,\right)
= 
\begin{bmatrix}
0 &0 & 1 & 0\\
0 & 0 &0& 1 \\
1 & 0 & 0 &0 \\
0 & 1 & 0 &0
\end{bmatrix}. 
\qquad  
{}_\mathbf{1}\pi_\rho
\left(\,
    \tikzmath{\filldraw[very thick, fill=gray!30] (0,0) circle (.2cm);
    \draw[mid<] (0,0) circle (.4cm);}
\,\right) 
= 
\begin{bmatrix}
\phi & \phi'\\
\phi' & \phi
\end{bmatrix}
\quad \text{and} \qquad 
{}_\mathbf{1}\pi_\rho
\left(\,
    \tikzmath{\filldraw[very thick, fill=gray!30] (0,0) circle (.2cm);
    \draw[blue, thick, mid<] (0,0) circle (.65cm);
    \draw[mid<] (0,0) circle (.4cm);}
\,\right) 
= 
\begin{bmatrix}
\phi' & \phi\\
\phi & \phi'
\end{bmatrix} 
$} 
\]
where $\phi$ and $\phi'$ are the operators on $\Hom(\rho\otimes\rho\to \rho)$ defined by
\begin{align*}
\phi\left(
    \tikzmath{
    \draw[mid>] (-.3,-.6) -- (0,0);
    \draw[mid>] (.3,-.6) -- (0,0);
    \draw[mid<] (0,.6) -- (0,0);
    \filldraw (0,0) node[left]{$\scriptstyle i$} circle (.05cm);
    }
\right)
&=
\sum_{j} 
\tikzmath{
\draw[mid>] (-.2,-.4) -- (.2,.8);
\draw[mid>] (-.3,-.8) -- (-.2,-.4);
\draw[mid>] (.2,.8) -- (.2,1.2);
\draw[mid>] (-.2,-.4) -- (.4,0);
\draw[mid>] (.4,0) -- (.2,.8);
\draw[mid>] (.6,-.8) -- (.4,0);
\filldraw (-.2,-.4) node[left]{$\scriptstyle j$} circle (.05cm);
\filldraw (.2,.8) node[left]{$\scriptstyle j$} circle (.05cm);
\filldraw (.4,0) node[left, yshift=.05cm]{$\scriptstyle i$} circle (.05cm);
}
= 
\sum_{j,k} A^{j,k}_{j,i}
    \tikzmath{
    \draw[mid>] (-.3,-.6) -- (0,0);
    \draw[mid>] (.3,-.6) -- (0,0);
    \draw[mid<] (0,.6) -- (0,0);
    \filldraw (0,0) node[left]{$\scriptstyle k$} circle (.05cm);
    }
\\
\phi'\left(
    \tikzmath{
    \draw[mid>] (-.3,-.6) -- (0,0);
    \draw[mid>] (.3,-.6) -- (0,0);
    \draw[mid<] (0,.6) -- (0,0);
    \filldraw (0,0) node[left]{$\scriptstyle i$} circle (.05cm);
    }
\right)
&=
\sum_{j}
\tikzmath{
\draw[thick, blue] (-.2,-.4) .. controls ++(-135:.6cm) and ++(270:.5cm) .. (.1,-.4) .. controls ++(90:.3cm) and ++(270:.3cm) .. (-.3,.2) .. controls ++(90:.3cm) and ++(270:.3cm) .. (.4,.8) .. controls ++(90:.5cm) and ++(135:.6cm) .. (.2,.8);
\draw[mid>] (-.2,-.4) -- (.2,.8);
\draw[far>] (-.05,-1) -- (-.2,-.4);
\draw[far>] (.2,.8) -- (.3,1.5);
\draw[far>] (-.2,-.4) -- (.4,0);
\draw[mid>] (.4,0) -- (.2,.8);
\draw[mid>] (.65,-1) -- (.4,0);
\filldraw (-.2,-.4) node[left]{$\scriptstyle j$} circle (.05cm);
\filldraw (.2,.8) node[left]{$\scriptstyle j$} circle (.05cm);
\filldraw (.4,0) node[left, yshift=.1cm]{$\scriptstyle i$} circle (.05cm);
}
=
 \chi_{\mathbf{1}, i}\sum_{j,k}D^{j,k}_{j,i}
    \tikzmath{
    \draw[mid>] (-.3,-.6) -- (0,0);
    \draw[mid>] (.3,-.6) -- (0,0);
    \draw[mid<] (0,.6) -- (0,0);
    \filldraw (0,0) node[left]{$\scriptstyle k$} circle (.05cm);
    }
\end{align*}
which we can naturally identify as operators on the two spaces:
\[
\left\{ \,\,
\tikzmath{\filldraw[very thick, fill=gray!30] (0,-.2) circle (.2cm);
    \draw[mid<] (0,0) ellipse (.4cm and .6cm);
    \draw[mid>] (0,0) .. controls ++(90:.2cm) and ++(-45:.3cm) .. (135:.46cm);
    \filldraw (135:.46cm) node[left]{$\scriptstyle 0$} circle (.05cm);
    }
    \,,\,\,
    \tikzmath{\filldraw[very thick, fill=gray!30] (0,-.2) circle (.2cm);
    \draw[mid<] (0,0) ellipse (.4cm and .6cm);
    \draw[mid>] (0,0) .. controls ++(90:.2cm) and ++(-45:.3cm) .. (135:.46cm);
    \filldraw (135:.46cm) node[left]{$\scriptstyle 1$} circle (.05cm);
    }\,\,
\right\}
\qquad\text{and}\qquad
\left\{\,\,
    \tikzmath{\filldraw[very thick, fill=gray!30] (0,-.2) circle (.2cm);
    \draw[mid<] (0,0) ellipse (.4cm and .6cm);
    \draw[thick, blue, mid<] (0,0) ellipse (.7cm and .9cm);
    \draw[mid>] (0,0) .. controls ++(90:.2cm) and ++(-45:.3cm) .. (135:.46cm);
    \filldraw (135:.46cm) node[left]{$\scriptstyle 0$} circle (.05cm);
    }
    \,,\,\,
    \tikzmath{\filldraw[very thick, fill=gray!30] (0,-.2) circle (.2cm);
    \draw[mid<] (0,0) ellipse (.4cm and .6cm);
    \draw[thick, blue, mid<] (0,0) ellipse (.7cm and .9cm);
    \draw[mid>] (0,0) .. controls ++(90:.2cm) and ++(-45:.3cm) .. (135:.46cm);
    \filldraw (135:.46cm) node[left]{$\scriptstyle 1$} circle (.05cm);
    }\,\, 
\right\}
\]
by local insertion.
That is, the elements of $A_{\mathbf{1}\leftarrow \mathbf{1}}$ which involve $\phi,\phi'$ above act on $A_{1\leftarrow \rho}$ by applying $\phi,\phi'$ locally on the trivalent vertices in our standard basis of $A_{1\leftarrow \rho}$.
\end{enumerate}
With these computations in hand, we either pin down the scalars $\chi_0$ and $\chi_1$, or determine the operator $\phi$.

\begin{proof}[Proof of Lemma~\ref{lem:centreNSD}]
Recall we have three possibilities for $a\in \{-2,0,2\}$. 
If $a=0$, then we have $\chi_0 = -\chi_1$.
Thus we can restrict our attention to the case of $a= \pm 2$.

We begin by determining the decomposition of $_\mathbf{1}\pi_\rho$ into irreducible representations of $A_{\mathbf{1}\to \mathbf{1}}$. 
As $X_0$ is the tensor unit of $\cZ(\cD_m)$, its restriction contains no copies of $\rho$, and thus $_\mathbf{1}\pi_\rho$ contains no copies of $\chi_0$. 
We also know that $\Tr(_\mathbf{1}\pi_\rho)=0$, and so from the character table of $A_{\mathbf{1}\to \mathbf{1}}$ above, we must have that 
\[  {}_\mathbf{1}\pi_\rho \cong 2\chi_1\oplus k \chi_2 \oplus (2-k)\chi_3\]
with $k \in \{0,1,2\}.$
In particular, we find that 
\[
\Tr\left({}_\mathbf{1}\pi_\rho
\left(\,
    \tikzmath{\filldraw[very thick, fill=gray!30] (0,0) circle (.2cm);
    \draw[mid<] (0,0) circle (.4cm);}
\,\right) \right) = 4 - 2\sqrt{5} + 2\mathbf{i}(k - 1)
\qquad\Longrightarrow\qquad
\Tr(\phi) = 2-\sqrt{5} + \mathbf{i}(k-1).  
\]

To determine $k$ we study the restriction of the objects $X_i$ and $Y_i$. 
By the above decomposition of ${}_\mathbf{1}\pi_\rho$ and from counting dimensions we have
\begin{align*}
    \cF(X_0) &= \mathbf{1}\\
    \cF(X_1) &= \mathbf{1} \oplus 2\rho \oplus 2 \alpha \rho\\
    \cF(X_2) &= \mathbf{1} \oplus k\rho \oplus (2-k) \alpha \rho\\
    \cF(X_3) &= \mathbf{1} \oplus (2-k)\rho \oplus k \alpha \rho.
\end{align*}
By our assumption that $a = \pm 2$, one of the objects $Y_i$ must be invertible. 
Thus we can label our $Y_i$ so that 
\begin{align*}
    \cF(Y_0) &= \alpha\\
    \cF(Y_1) &= \alpha \oplus 2\rho \oplus 2 \alpha \rho\\
    \cF(Y_2) &=\alpha \oplus  (2-k)\rho \oplus k \alpha \rho\\
    \cF(Y_3) &= \alpha \oplus k\rho \oplus  (2-k) \alpha \rho.
\end{align*}
Hence we now know the restriction of all the objects in both $\cI(\mathbf{1})$ and $\cI(\alpha)$, up to the integer $k$. 
Denote by $Z_i$ the remaining simple objects in $\cZ(\cD)$, i.e. those simple objects such that
\[ 
\cF(Z_i)  = p_i \rho \oplus q_i \alpha \rho   
\]
where $p_i,q_i$ are positive integers.
This allows us to write
\begin{align*}
    \cI(\rho) &= 2X_1 + kX_2 + (2-k)X_3 + 2Y_1 + (2-k)Y_2 + kY_3 + \sum p_i Z_i\\
    \cI(\alpha\rho) &= 2X_1 + (2-k)X_2 + kX_3 + 2Y_1 + kY_2 + (2-k)Y_3 + \sum q_i Z_i.
\end{align*}
Therefore
\begin{align*}
    20 &= \dim\Hom(\cI(\rho), \cI(\rho)) = 4 k^2-8 k+16 + \sum p_i^2 \\
     20 &= \dim\Hom(\cI(\alpha\rho), \cI(\alpha\rho)) = 4 k^2-8 k+16 + \sum q_i^2 \\
     16 &= \dim\Hom(\cI(\rho), \cI(\alpha\rho)) = -4 k^2+8 k+8 + \sum p_iq_i.
\end{align*}
If $k\in \{0,2\}$ then we get
\begin{align*}  
\sum p_i^2 &= \sum q_i^2 = 4
&\text{and} 
\sum p_iq_i&= 8,     
\end{align*}
which is impossible. 
Thus we must have $k=1$, and so 
$\Tr(\phi) = 2-\sqrt{5}$.
From 
\[{}_\mathbf{1}\pi_\rho
\left(  \frac{1}{\dim(\cD)}\left(\tikzmath{\filldraw[very thick, fill=gray!30] (0,0) circle (.2cm);} +\tikzmath{\filldraw[very thick, fill=gray!30] (0,0) circle (.2cm);
    \draw[thick, blue, mid<] (0,0) circle (.4cm);} + (2+\sqrt{5}) \tikzmath{\filldraw[very thick, fill=gray!30] (0,0) circle (.2cm);
    \draw[mid<] (0,0) circle (.4cm);} + (2+\sqrt{5})  \tikzmath{\filldraw[very thick, fill=gray!30] (0,0) circle (.2cm);
    \draw[blue, thick, mid<] (0,0) circle (.65cm);
    \draw[mid<] (0,0) circle (.4cm);}\right)  \right)=0
    \]
    and
    \[ {}_\mathbf{1}\pi_\rho
\left( \tikzmath{\filldraw[very thick, fill=gray!30] (0,0) circle (.2cm);
    \draw[mid<] (0,0) circle (.4cm);}\right)^2 = {}_\mathbf{1}\pi_\rho
\left(\tikzmath{\filldraw[very thick, fill=gray!30] (0,0) circle (.2cm);
    \draw[thick, blue, mid<] (0,0) circle (.4cm);}\right) + 2\cdot{}_\mathbf{1}\pi_\rho
\left(\tikzmath{\filldraw[very thick, fill=gray!30] (0,0) circle (.2cm);
    \draw[mid<] (0,0) circle (.4cm);}\right) +2\cdot{}_\mathbf{1}\pi_\rho
\left(\tikzmath{\filldraw[very thick, fill=gray!30] (0,0) circle (.2cm);
    \draw[blue, thick, mid<] (0,0) circle (.65cm);
    \draw[mid<] (0,0) circle (.4cm);}\right), 
    \]
    we find that $\phi$ has the two distinct eigenvalues
    \[ \frac{(2+i)-\sqrt{5}}{2}  \quad \text{and}\quad \frac{(2-i)-\sqrt{5}}{2} .  \]
    
    Finally, by Remark~\ref{rem:RescaleNSD}, we are free to unitarily change our basis of $\cD_2(\rho\otimes\rho\to \rho)$ by any element of $U(2)$. 
    In particular we can choose this basis so that $\phi$ acts diagonally. 
    This gives the statement of the lemma.
\end{proof}

%%%%%%%%%%%%%%%%%%%%%%%%%%%%%%%%%%%%%%%%%%%%%%%%%%%%%%%%%%%%%%%%%%%
\subsection{Symmetries}
\label{NSD:Symmetries}

We now use the tetrahedral symmetries to determine relations between the 128 complex scalers:
\[ A^{i,j}_{k,\ell},\quad  B^{i,j}_{k,\ell} ,\quad C^{i,j}_{k,\ell},\quad D^{i,j}_{k,\ell} ,\quad \widehat{A}^{i,j}_{k,\ell}, \quad \widehat{B}^{i,j}_{k,\ell} , \quad\widehat{C}^{i,j}_{k,\ell},\quad\text{and} \quad  \widehat{D}^{i,j}_{k,\ell}.    \]
Using the same techniques as in the self-dual case we are able to show the following.
\begin{lem}\label{lem:symNSD}
The scalars 
$B^{i,j}_{k,\ell} , D^{i,j}_{k,\ell}, \widehat{A}^{i,j}_{k,\ell}   ,  \widehat{C}^{i,j}_{k,\ell} , \widehat{D}^{i,j}_{k,\ell}$
can be expressed in terms of the scalars $A^{i,j}_{k,\ell}$ as:
\begin{align*}
    B^{i,j}_{k,\ell}  &= -\nu \chi_\ell \omega_\ell^2 A^{k,i}_{j,\ell}  
    &  
    D^{i,j}_{k,\ell} &= -\nu^{-1}\omega_\ell A^{j,k}_{i,\ell} 
    &
    &
    \\
    \widehat{A}^{i,j}_{k,\ell}  &= -\nu^3\omega_\ell^2\overline{A^{j,i}_{k,\ell} } 
    & 
    \widehat{C}^{i,j}_{k,\ell} &= - \chi_\ell \omega_\ell \overline{A^{i,k}_{j,\ell}} & \widehat{D}^{i,j}_{k,\ell} &= \overline{ A^{k,j}_{i,\ell} }.
\end{align*}
The scalars $ \widehat{B}^{i,j}_{k,\ell}$ can be expressed in terms of the scalars $C^{i,j}_{k,\ell}$ as:
\[\widehat{B}^{i,j}_{k,\ell} = \nu^3\overline{C^{k,j}_{i,\ell}}.     \]
The scalars $A^{i,j}_{k,\ell}$ satisfy $\bbZ/4\bbZ$ symmetries generated by the relations
\[  A^{i,j}_{k,\ell}= -\nu \chi_i \omega_j\omega_k^2\omega_\ell^2 \overline{ A^{k,\ell}_{j,i}} = \chi_i  \chi_k\omega_i\omega_\ell  \omega_j^2\omega_k^2A^{j,i}_{\ell,k}.  \]
The scalars $C^{i,j}_{k,\ell}$ satisfy $S_3$ symmetries generated by the order three rotation
\[  C^{i,j}_{k,\ell} = \omega_\ell C^{j,k}_{i,\ell} = \omega_\ell^2 C^{k,i}_{j,\ell} \]
and the order two flip
\[C^{i,j}_{k,\ell}  =  \chi_j\chi_k\omega_i^2\omega_k C^{k,\ell}_{i,j}.  \]

Finally we have that if $\chi_0 = -\chi_1$, then 
\[   A^{i,j}_{k,\ell} = C^{i,j}_{k,\ell}  \quad \text{ if $i+j+k+\ell = 0\pmod 2$}. \]
\end{lem}
\begin{proof}
The proof of this lemma uses the exact same techniques as in the proof of Lemma~\ref{lem:mainSym}. 
The only real difference is we have different Frobenius operators in this case.
\end{proof}
%%%%%%%%%%%%%%%%%%%%%%%%%%%%%%%%%%%%%%%%%%%%%%%%%%%%%%%%%%%%%%%%%%%
\subsection{Classification}
\label{NSD:Classification}

We now complete the proof of Theorem~\ref{thm:NSD} to complete the classification in the non self-dual case. 
To prove this theorem, we break into three cases:
(1) $m=1$, (2) $m=2$ and $\chi_0 = -\chi_1$, and (3) $m=2$ and $\chi_0 = \chi_1$.

%%%%%%%%%%%%%%%%%%%%%%%%%%%%%%%%%%%%%%%%%%
\subsubsection*{The case \texorpdfstring{$m=1$}{m=1}}

If $m=1$ then from our previous analysis we only have to determine the sign $\chi_0$, the 3rd root of unity $\omega_0$, the 8th root of unity $\nu$, and the two complex scalars $a:= A^{0,0}_{0,0}$ and $c:= C^{0,0}_{0,0}$. Further, we have that if $\omega_0 \neq 1$, then $c = 0$.

By evaluating the diagrams 
\[
\tikzmath{
\draw[mid>] (-.2,1) -- (-.2,1.5);
\draw[mid>] (-.2,-1.5) -- (-.2,-1);
\draw[mid>] (-.2,-1) -- (0,-.4);
\draw[mid>] (0,.4) -- (-.2,1);
\draw[mid>] (0,-.4) arc (240:120:.47cm);
\draw[mid>] (0,-.4) arc (-60:60:.47cm);
\draw[mid>] (-.2,-1) arc (240:120:1.15cm);
\filldraw (0,-.4) node[left]{$\scriptstyle \ell'$} circle (.05cm);
\filldraw (0,.4) node[left]{$\scriptstyle \ell$} circle (.05cm);
\filldraw (-.2,-1) node[left]{$\scriptstyle k'$} circle (.05cm);
\filldraw (-.2,1) node[left]{$\scriptstyle k$} circle (.05cm);
}
\,,
\qquad 
\tikzmath{
\draw[mid>] (-.2,1) -- (.1,1.5);
\draw[thick, blue, mid>] (-.2,1) -- (-.5,1.5);
\draw[mid>] (.1,-1.5) -- (-.2,-1);
\draw[thick, blue, mid>] (-.5,-1.5) -- (-.2,-1);
\draw[mid>] (-.2,-1) -- (0,-.4);
\draw[mid>] (0,.4) -- (-.2,1);
\draw[mid>] (0,-.4) arc (240:120:.47cm);
\draw[mid>] (0,-.4) arc (-60:60:.47cm);
\draw[mid>] (-.2,-1) arc (240:120:1.15cm);
\filldraw (0,-.4) node[left]{$\scriptstyle \ell'$} circle (.05cm);
\filldraw (0,.4) node[left]{$\scriptstyle \ell$} circle (.05cm);
\filldraw (-.2,-1) node[left]{$\scriptstyle k'$} circle (.05cm);
\filldraw (-.2,1) node[left]{$\scriptstyle k$} circle (.05cm);
}
\,,\qquad\text{and}\qquad 
\tikzmath{
\draw[mid>] (0,-.6) arc (-90:90:.3cm);
\draw[mid>] (0,-.6) arc (270:90:.3cm);
\draw[mid>] (-1,-1) to[out=90,in=-135] (-.5,.5);
\draw[mid>] (0,0) to[out=90, in=-45] (-.5,.5);
\draw[mid>] (-.5,.5) -- (0,1);
\draw[thick, blue, mid>] (-.5,.5) -- (-1,1);
\draw[thick, blue, mid>] (0,-1) -- (0,-.6);
\filldraw (0,0) node[left, yshift=.1cm]{$\scriptstyle \ell$} circle (.05cm);
\filldraw (-.5,.5) node[left]{$\scriptstyle k$} circle (.05cm);
\filldraw[blue] (0,-.6) circle (.05cm);
}
\]
in two ways (see Footnote \ref{footnote:EvaluateInTwoWays})
we obtain the equations:
\begin{align*}
    |a|^2 + |c|^2 &= 1
    &
    2|a|^2 &= 1 - \frac{1}{1+\sqrt{2}}
    &
    \frac{\nu}{1+\sqrt{2}} &=  a(\chi_0\omega_0-\nu ).
\end{align*}
With the first two of these equations we can solve to find
\[ |a|^2 = 1 - \frac{1}{\sqrt{2}} \quad \text{and} \quad |c|^2 = \frac{1}{\sqrt{2}},   \]
and thus we have $\omega_0 = 1$. 
The general solution to these equations is then given by
\[ \chi_0 = 1,\qquad a =\frac{\nu}{(1-\nu)} (-1+\sqrt{2}) \quad \text{and} \quad c = e^{i\theta}2^{\frac{-1}{4}}    \]
where $\theta$ is any phase.

\begin{lem}
There are exactly two unitary fusion categories, up to monoidal equivalence, which categorify $S(1)$.
\end{lem}
\begin{proof}
By unitarily renormalising the basis element 
\[  
\TI{0}\mapsto  z \cdot \TI{0} \qquad\qquad z\in U(1),
\] we change $c$ to $z^{-2}c$. 
We can thus renormalise so that $c =2^{\frac{-1}{4}}$. 
Hence we have two solutions for our free variables, depending on the choice of $\nu = e^{\pm i \pi \frac{1}{4}}$. 
By Proposition~\ref{prop:UniqueNSD}, there are at most 2 unitary fusion categories with these fusion rules. 
These two unitary fusion categories are realised by the even parts of the two subfactors $\mathcal{S}'$ constructed in \cite{MR3306607}, which are monoidally non-equivalent and complex conjugate to each other.
Indeed, they each admit a $\bbZ/2\bbZ$-equivariantization, which produces monoidally non-equivalent $2^{\bbZ/4\bbZ}1$ near-group fusion categories which are complex conjugate \cite[Ex.~2.2]{MR4079744}, \cite[Ex.~9.5]{MR3635673}.
\end{proof}

%%%%%%%%%%%%%%%%%%%%%%%%%%%%%%%%%%%%%%%%%%
\subsubsection*{The case \texorpdfstring{$m=2$ and $\chi_0 = -\chi_1 = 1$}{m=2 and chi0=-chi1=1}}

If $m=2$ and $\chi_0 = -\chi_1=1$, then we have to determine the 3rd roots of unity $\omega_i$, the 8th root of unity $\nu$, and the free complex variables
$A^{i,j}_{k,l}$ and $C^{i,j}_{k,l}$.
We can represent these free  complex variables in the same matrix notation as in \eqref{eq:MatrixForm} in the self-dual section. 
After applying the symmetries of Lemma~\ref{lem:symNSD}, we obtain:
\[
A = 
\begin{bmatrix}
a_0 &0 & 0 &0	\\
0 &  a_1 &-\nu \omega_0^2 \overline{a_1}&0		\\
0 &  \nu \omega_1^2 \overline{a_1}&  a_1 & 0	\\
0 &0 &0 &  a_2
\end{bmatrix}
\qquad\qquad
C   = 
\begin{bmatrix}
c_0 &0 & 0 &0	\\
0 &0&0&0		\\
0 &  0&  0 & 0	\\
0 &0 &0 & c_1
\end{bmatrix}.
\]
Due to the large number of variables which are zero, it is fairly easy to derive a contradiction in this case.
\begin{lem}
There is no unitary fusion category that categorifies $S(2)$ with $\chi_0 =-\chi_1 = 1$.
\end{lem}
\begin{proof}
By evaluating 
\[
\tikzmath{
\draw[mid>] (-.2,1) -- (.1,1.5);
\draw[thick, blue, mid>] (-.2,1) -- (-.5,1.5);
\draw[mid>] (.1,-1.5) -- (-.2,-1);
\draw[thick, blue, mid>] (-.5,-1.5) -- (-.2,-1);
\draw[mid>] (-.2,-1) -- (0,-.4);
\draw[mid>] (0,.4) -- (-.2,1);
\draw[mid>] (0,-.4) arc (240:120:.47cm);
\draw[mid>] (0,-.4) arc (-60:60:.47cm);
\draw[mid>] (-.2,-1) arc (240:120:1.15cm);
\filldraw (0,-.4) node[left]{$\scriptstyle \ell'$} circle (.05cm);
\filldraw (0,.4) node[left]{$\scriptstyle \ell$} circle (.05cm);
\filldraw (-.2,-1) node[left]{$\scriptstyle k'$} circle (.05cm);
\filldraw (-.2,1) node[left]{$\scriptstyle k$} circle (.05cm);
}
\]
in two ways (see Footnote \ref{footnote:EvaluateInTwoWays}),
we obtain the equation
\[  
\chi_\ell \chi_{\ell'} \omega_\ell^2\omega_{\ell'}\sum_{i,j}A^{k,i}_{j,\ell}\overline{A^{k',i}_{j,\ell'}} + \omega_\ell\omega_{\ell'}^2 \sum_{i,j}A^{j,k}_{i,\ell}\overline{A^{j,k'}_{i,\ell'}} 
-(2- \sqrt{5})\delta_{\ell,k}\delta_{\ell',k'}
= 
\delta_{\ell,\ell'}\delta_{k,k'}.  \]
Taking $k = k' = 0$ and $\ell = \ell' = 1$ gives
$|a_1|^2  = \frac{1}{2}$,
and taking $k = k' = \ell = l'=0$ gives
$2|a_0|^2 + 2|a_1|^2 = 3 - \sqrt{5}$.
These two equations imply $2|a_0|^2<0$, a contradiction.
\end{proof}

%%%%%%%%%%%%%%%%%%%%%%%%%%%%%%%%%%%%%%%%%%
\subsubsection*{The case \texorpdfstring{$m=2$ and $\chi_0 = \chi_1$}{m=2 and chi0=chi1}}

Finally we deal with the last case where $m=2$ and $\chi_0 = \chi_1$. Let us again represent our free variables $A^{i,j}_{k,l}$ and $C^{i,j}_{k,l}$ in matrix form as in \eqref{eq:MatrixForm}. 
After applying the symmetries of Lemma~\ref{lem:symNSD}, we obtain:
\[A = 
\begin{bmatrix}
a_0 &a_1 & \omega_{0}\omega_1^2a_1 &a_2	\\
-\nu \chi_0 \omega_1^2 \overline{a_1} & a_3 &-\nu \chi_0 \omega_0^2 \overline{a_3}&a_4		\\
-\nu \chi_0 \omega_0^2 \overline{a_1} & -\nu \chi_0 \omega_1^2 \overline{a_3}&  a_3 & \omega_0^2\omega_1a_4	\\ 
-\nu \chi_0 \omega_0\omega_1 \overline{a_2} &-\nu \chi_0 \omega_0^2 \overline{a_4} &-\nu \chi_0 \omega_1^2 \overline{a_4}  &  a_5
\end{bmatrix}
\qquad
C   = 
\begin{bmatrix}
c_0 &c_1 & \omega_0c_1 &c_2	\\
c_1&\omega_1^2c_2&\omega_0^2\omega_1^2c_2&\omega_1^2c_3	\\
\omega_0^2c_1 &  \omega_1c_2&  \omega_1^2c_2 & c_3	\\
\omega_0\omega_1^2c_2 &\omega_1c_3 &c_3 & c_4
\end{bmatrix}.
\]
Recall from Lemma~\ref{lem:centreNSD} that in this case we have
\[ \sum_{i} A^{i,0}_{i,0} = \frac{(2+i)-\sqrt{5}}{2} ,\qquad \sum_i A^{i,1}_{i,1} = \frac{(2-i)-\sqrt{5}}{2}, \quad \text{and}\quad \sum_i A^{i,1}_{i,0} =\sum_i A^{i,0}_{i,1} =0,\]
which implies
\[ 
a_0 + a_3 = \frac{(2+i)-\sqrt{5}}{2},\qquad a_3 + a_5 = \frac{(2-i)-\sqrt{5}}{2}, \quad \text{and}\quad a_1 = - a_4 = -\omega_0^2\omega_1a_4.  
\]

With these linear equations in hand, it is straightforward to show non-existence in this case.

\begin{thm}
There is no unitary fusion category that categorifies $S(2)$ with $\chi_0 =\chi_1$.
\end{thm}
\begin{proof}

Evaluating the diagram
\[
\tikzmath{
\draw[mid>] (0,-.6) arc (-90:90:.3cm);
\draw[mid>] (0,-.6) arc (270:90:.3cm);
\draw[mid>] (-1,-1) to[out=90,in=-135] (-.5,.5);
\draw[mid>] (0,0) to[out=90, in=-45] (-.5,.5);
\draw[mid>] (-.5,.5) -- (0,1);
\draw[thick, blue, mid>] (-.5,.5) -- (-1,1);
\draw[thick, blue, mid>] (0,-1) -- (0,-.6);
\filldraw (0,0) node[left, yshift=.1cm]{$\scriptstyle \ell$} circle (.05cm);
\filldraw (-.5,.5) node[left]{$\scriptstyle k$} circle (.05cm);
\filldraw[blue] (0,-.6) circle (.05cm);
}
\]
in two ways (see Footnote \ref{footnote:EvaluateInTwoWays}) gives
\[ \sum_i A^{i,k}_{i,\ell} + \nu^3 \chi_\ell \omega_\ell \sum_i A^{k,i}_{i,\ell} - \delta_{\ell,k}(2 - \sqrt{5})=0.\]
In terms of our free variables this gives
\begin{align*}
    2-\sqrt{5} = a_0 + a_3 + 2\nu^3\chi_0\omega_0a_0 + 2\overline{a_3} & = \mathbf{i} - 2\nu ^3 \chi_0  \omega _0 \left(2 \nu  \chi_0  \omega _0^2 \overline{a_3}+2 a_3+(-2-i)+\sqrt{5}\right)\\
      2-\sqrt{5}=a_5 + a_3 + 2\nu^3\chi_0\omega_1a_5 + 2\overline{a_3} &= - \mathbf{i} -2 \nu ^3 \chi_0  \omega _1 \left(2 \nu  \chi_0  \omega _1^2 \overline{a_3}+2 a_3+(-2+i)+\sqrt{5}\right).
\end{align*}
This system of equations of the complex variable $a_3$ does not hold for any values of our free variables 
$\chi_0 \in \{-1, 1\}$, $\nu \in  \{e^\frac{\mathbf{i}\pi}{4}, e^\frac{-\mathbf{i}\pi}{4}\}$, and $\omega_0,\omega_1 \in \{1,e^\frac{2\mathbf{i}\pi}{3} ,e^\frac{4\mathbf{i}\pi}{3}   \}$.
\end{proof}

%%%%%%%%%%%%%%%%%%%%%%%%%%%%%%%%%%%%%%%%%%%%%%%%%%%%%%%%%%%%%%%%%%%
%%%%%%%%%%%%%%%%%%%%%%%%%%%%%%%%%%%%%%%%%%%%%%%%%%%%%%%%%%%%%%%%%%%
%%%%%%%%%%%%%%%%%%%%%%%%%%%%%%%%%%%%%%%%%%%%%%%%%%%%%%%%%%%%%%%%%%%
\appendix
\section{A multiplicity bound for \texorpdfstring{$\bbZ/2\bbZ$}{Z/2Z}-quadratic categories}
\label{appendix:PseudounitaryBound}

\begin{center}
{\tiny Ryan Johnson, Siu-Hung Ng, David Penneys, Jolie Roat, Matthew Titsworth, and Henry Tucker}
\end{center}

In this appendix, we prove Theorem \ref{thm:PseudounitaryBound}.
That is, given a pseudounitary $\bbZ/2\bbZ$ quadratic fusion category with simple objects $1,\alpha,\rho,\alpha\rho$ with $\rho$ self-dual and fusion rules determined by
\begin{equation*}
\alpha^2 \cong 1
\qquad\text{and}\qquad
\rho^2 \cong 1\oplus m \rho\oplus n \alpha \rho,
\tag*{\ref{Q:4ObjectsSelfDual}}
\end{equation*}
$(m,n)$ must be one of $(0,0),(0,1),(1,0),(1,1),(2,2)$.

%%%%%%%%%%%%%%%%%%%%%%%%%%%%%%%%%%%%%%%%%%%%%%%%%%%%%%%%%%%%%%%%%%%
\subsection{Basic number theoretic constraints}

Given a $\bbZ/2\bbZ$-quadratic category $\cC$ with fusion rules \ref{Q:4ObjectsSelfDual}, the fusion matrices are given in the ordering $1,\rho,\alpha\rho,\alpha$ by
$$
L_\rho =
\left( 
\begin{array}{c|ccc}
0 & 1 & 0 & 0
\\\hline
1 & m & n & 0
\\
0 & n & m & 1
\\
0 & 0 & 1 & 0
\end{array}
\right),\,
L_{\alpha\rho} =
\left( 
\begin{array}{c|ccc}
0 & 0 & 1 & 0
\\\hline
0 & n & m & 1
\\
1 & m & n & 0
\\
0 & 1 & 0 & 0
\end{array}
\right),
\text{ and }
L_\alpha =
\left( 
\begin{array}{c|ccc}
0 & 0 & 0 & 1
\\\hline
0 & 0 & 1 & 0
\\
0 & 1 & 0 & 0
\\
1 & 0 & 0 & 0
\end{array}
\right).
$$
Setting $d:= \dim(\rho)$, we have $\dim(\alpha)=1$ and $d^2 = 1+(m+n)d$, so that
\begin{equation}
\label{eq:dParameter}
	d = \frac{1}{2} \left(m + n + \sqrt{4 + (m+n)^2}\right).  
\end{equation}
Then
\begin{align*}
	\dim(\cC) &= 2 + 2 d^2
	%\\&
	= 2 + \frac{(m+n)^2 + 2(m+n)\sqrt{4 + (m+n)^2} + 4 + (m+n)^2}{2}
	%\\&
	= 4 + 2(m+n)d.
\end{align*}
Since $K_0(\cC)$ is abelian of dimension 4,
its irreducible representations are all 1-dimensional.
Hence by \cite[Rem.~2.11]{1309.4822}, the element  
\begin{align*}
	R 
	:= 
	I + L_\rho^2 + L_{\alpha\rho}^2 + L_\alpha^2 
	= 
	\left(\begin{array}{rrrr}
		4 & 2 \, m & 2 \, n & 0 \\
		2 \, m & 2 \, n^{2} + 2 \, m^{2} + 4 & 4 \, mn & 2 \, n \\
		2 \, n & 4 \, mn & 2 \, n^{2} + 2 \, m^{2} + 4 & 2 \, m \\
		0 & 2 \, n & 2 \, m & 4
	\end{array}\right)
\end{align*}
is central in $K_0(\cC)$, and the roots of its characteristic polynomial are called the \emph{formal codegrees} \cite{MR2576705} of $\cC$:
\begin{align*}
	f_1 &= 4 + (m+n)^2 + (m+n)\sqrt{4 + (m+n)^2}\\
	f_2 &= 4 + (m+n)^2 - (m+n)\sqrt{4 + (m+n)^2}\\
	f_3 &= 4 + (m-n)^2 + (m-n)\sqrt{4 + (m-n)^2}\\
	f_4 &= 4 + (m-n)^2 - (m-n)\sqrt{4 + (m-n)^2}.
\end{align*}

%%%%%%%%%%%%%%%%%%%%%%%%%%%%%%%%%%%%%%%%%%
\subsection{Computing the induction and forgetful functor}

We now analyze the center $Z(\cC)$,
the forgetful functor $\cF: Z(\cC)\to \cC$,
and the induction functor $\cI: \cC \to Z(\cC)$.
Recall that $\cF(\cI(c)) \cong \bigoplus_{x\in \Irr(\cC)} x\otimes c\otimes x^*$ for all $c\in \cC$, and that $\cF$ is biadjoint to $\cI$.
We use the notation
$(a,b):= \dim(\cC(a\to b))$ and $(A,B):= \dim(Z(\cC)(A\to B))$.

\begin{lem}[{\cite[Theorem 2.13]{1309.4822}}]
\label{lem:Ostrik2.13}
There are distinct simple objects 
$ 1_{Z(\cC)}, X_2, X_3, X_4 \in \Irr(Z(\cC))$
such that
$$
\cI(1_\cC) = 1_{Z(\cC)} \oplus X_2\oplus X_3\oplus X_4
\qquad\text{and}\qquad
\dim(X_k)= \frac{f_1}{f_k}.
$$
\end{lem}

Setting
$$
r := \sqrt{\frac{4 + (m+n)^2}{4 + (m-n)^2}},
$$
it is straightforward to calculate that
\begin{align*}
\dim(X_2) &=\frac{f_1}{f_2} = \frac{f_1^2}{f_1f_2} = 1 + (m+n)d = \frac{\dim(\cC)}{2}-1
\\
\dim(X_3) &=\frac{f_1}{f_3} = \frac{f_1f_4}{f_3f_4}
=1 + \frac{1}{2}(m+n)d - \frac{r}{2}(m-n)d
\\
\dim(X_4) &=\frac{f_1}{f_4} = \frac{f_1f_3}{f_3f_4}
=1 + \frac{1}{2}(m+n)d + \frac{r}{2}(m-n)d.
\end{align*}

\begin{rem}
\label{rem:K2abBetaConditions}
Since $\dim(X_2), \dim(X_3),\dim(X_4) \in \bbZ(d)$, it must be the case that either $m=n$ or $r \in \bbQ(d)$.  
If $m=0$ or $n=0$, then $r = 1$.  
One may check that if $0 \neq m \neq n \neq 0$ and $r \in \bbQ(d)$, then $m+n \geq 11$.
We will show below in Theorem \ref{thm:TwoParameterBound} that $m+n\leq 5$.
\end{rem}

\begin{prop}
The center $Z(\cC)$ has 8 distinct simple objects $1_{Z(\cC)}, X_2,X_3,X_4,Y_1,Y_2,Y_3,Y_4$ such that
$$
\cI(1_\cC) = 1_{Z(\cC)} \oplus X_2\oplus X_3\oplus X_4
\qquad\text{and}\qquad
\cI(\alpha) = Y_1\oplus Y_2\oplus Y_3\oplus Y_4.
$$
Denote the rest of the simple objects of $Z(\cC)$ by $\{Z_s\}_{s\in S}$ where $S$ is some finite set.
The matrix $F$ of the forgetful functor $\cF:Z(\cC) \to \cC$ can then be represented as follows, where zero entries are omitted:
$$
F
=
\begin{array}{c|cccc|cccc|c}
& 1_{Z(\cC)} & X_2 & X_3 & X_4 & Y_1 & Y_2 & Y_3 & Y_4 & Z_s
\\\hline
1_\cC & 1 & 1 & 1 & 1 & &&&
\\
\alpha & &&& & 1 & 1 & 1 & 1
\\\hline
\rho & & x_2 & x_3 & x_4 & y_1 & y_2 & y_3 & y_4 & z_s
\\
\alpha\rho & & x'_2 & x'_3  & x'_4 & y'_1 & y'_2 & y'_3 & y'_4 & z'_s
\end{array}
$$
and the induction matrix is given by $F^T$.
Moreover, 
\begin{align}
\sum_{j=2}^4 x_j &= 2m \;\; \text{ and } \sum_{j=2}^4 x'_j = 2n,
\label{eqn:VerticalSums1}
\\
x_2 + x_2' &= m+n,
\label{eqn:HorizontalSumsAB1}
\\
x_3+x'_3 &= \frac{1}{2}(m+n) - \frac{r}{2}(m-n), \text{ and}
\label{eqn:HorizontalSumsAB2}
\\
x_4 + x'_4 &= \frac{1}{2}(m+n) + \frac{r}{2}(m-n).
\label{eqn:HorizontalSumsAB3}
\\
\sum_{j=1}^4 y_j &= 2n \;\;\text{ and }\;\; \sum_{j=1}^4 y_j' = 2m.
\label{eqn:VerticalSumsAB2}
\end{align}
\end{prop}
\begin{proof}
First, $1_{Z(\cC)}$ decomposes as desired by Lemma \ref{lem:Ostrik2.13}.
Next, observe that  
$$
(\cI(\alpha), \cI(\alpha))
=
4
\qquad
\text{and}
\qquad
(\cF\cI(1_\cC),\alpha)=0.
$$
Since the first Frobenius-Schur indicator $\nu_1$ satisfies 
$\Tr_{Z(\cC)}(\theta_{\cI(\alpha)}) = 0$  
\cite[Rem.~4.6]{MR2381536}
(see also \cite[Thm.~2.4]{1309.4822}),
$\cI(\alpha)$ decomposes as 4 distinct simples which are distinct from $1_{Z(\cC)},X_2,X_3,X_4$.
Equations
\eqref{eqn:VerticalSums1}
and
\eqref{eqn:VerticalSumsAB2}
follow from 
calculating $\cF\cI(1_\cC)$ and $\cF\cI(\alpha)$.
Equations
\eqref{eqn:HorizontalSumsAB1},
\eqref{eqn:HorizontalSumsAB2},
and
\eqref{eqn:HorizontalSumsAB3}
now follow from
the formulas for $\dim(X_k)$ for $k=2,3,4$.
\end{proof}

%Now using number theory and the formal codegrees, we can determine the $x_i, x_i', y_s, y_s'$ in terms of $a,b,d$ and

We now compute the dimensions of all the hom spaces amongst $\cI(\rho)$ and $\cI(\alpha\rho)$ in two ways to see
\begin{align}
(\cI(\rho),\cI(\rho)) &= 4+2m^2 + 2n^2 = \sum_{j=2}^4 x_j^2+ \sum_{j=1}^4 y_j^2 +\sum_{s} z_s^2 
\label{eqn:DimHomABXX}
\\
(\cI(\rho),\cI(\alpha\rho)) &= 4mn = \sum_{j=2}^4 x_jx'_j + \sum_{j=1}^4 y_jy'_j +\sum_{s} z_sz_s' 
\label{eqn:DimHomABXY}
\\
(\cI(g\rho),\cI(\alpha\rho)) &= 4 + 2m^2 + 2n^2  = \sum_{j=2}^4 (x'_j)^2 + \sum_{j=1}^4 (y'_j)^2 +\sum_{s} (z_s')^2.
\label{eqn:DimHomABYY}
\end{align}

\begin{lem}
The non-negative integers $x_j,x'_j,y_j,y_j',z_s,z_s'$ satisfy
\begin{align}
8+\frac{5}{2}(m +n)^2 - \frac{r^2}{2}(m-n)^2 &= \sum_{j=1}^4 (y_j+y'_j)^2 + \sum_{s} (z_s+z_s')^2
\label{eqn:PlusXYsquaresAB}
\\
8 + 4(m-n)^2 &= \sum_{j=2}^4 (x_j-x'_j)^2+\sum_{j=1}^4 (y_j-y'_j)^+ \sum_{s} (z_s-z_s')^2.
\label{eqn:MinusXYsquaresAB}
\end{align}
\end{lem}
\begin{proof}
To get the first equation, sum Equations \eqref{eqn:DimHomABXX} and \eqref{eqn:DimHomABYY} and twice Equation \eqref{eqn:DimHomABXY}.
Then use Equations \eqref{eqn:HorizontalSumsAB1}, \eqref{eqn:HorizontalSumsAB2}, and \eqref{eqn:HorizontalSumsAB3} and simplify.
The second is similar.
\end{proof}

\begin{prop}
We have the following upper bound:
\begin{equation}
\sum_{s} (z_s+z_s')^2 
\leq 8+\frac{3}{2}(m+n)^2.
\label{eqn:SquaresUpperBoundAB}
\end{equation}
\end{prop}
\begin{proof}
By Equation \eqref{eqn:PlusXYsquaresAB}, the desired inequality is equivalent to
$$
\sum_{j=1}^4 (y_j+y'_j)^2 
\geq 
\frac{1}{4}\left(\sum_{j=1}^4 y_j+y'_j\right)^2 
= 
(m+n)^2
\geq
(m+n)^2 - r^2(m-n)^2,
$$
which is true.
The second equality above holds by Equation \eqref{eqn:VerticalSumsAB2}, and the first inequality above follows from the fact that for any real numbers $w,x,y,z$, we have
\begin{align*}
4(w^2{+}x^2{+}y^2{+}z^2){-}(w{+}x{+}y{+}z)^2 
&= 
(w{-}x)^2{+}(w{-}y)^2{+}(w{-}z)^2{+}(x{-}y)^2{+}(x{-}z)^2{+}(y{-}z)^2 
\geq 0.
\end{align*}
The proof is complete.
\end{proof}

\begin{prop}\label{prop:TwoThetasAB}
Denote the twists of $Y_1,\dots, Y_4$ by $\theta_1,\dots, \theta_4$.
We have $\theta_1^2=\theta_2^2=\theta_3^2=\theta_4^2=\pm 1$.
Setting $\theta=\theta_1$, without loss of generality, we have $\theta\in \{1,i\}$.
\end{prop}
\begin{proof}
We calculate the following Frobenius-Schur indicators \cite{MR2381536} of $\alpha$:
\begin{align}
0&=\Tr_{Z(\cC)}(\theta_{\cI(\alpha)})
=
\sum_{j=1}^4 \theta_j  \dim(Y_j)
\label{eqn:FirstIndicatorZAB}
\\
\pm \dim(\cC)&=\Tr_{Z(\cC)}(\theta_{\cI(Z)}^2)
=
\sum_{j=1}^4 \theta_j^2  \dim(Y_j).
\label{eqn:SecondIndicatorZAB}
\end{align}
Since $\dim(\cI(\alpha))=\dim(\cC)=\sum_{j=1}^4 \dim(Y_k)$, Equation \eqref{eqn:SecondIndicatorZAB} implies that $\theta_1^2=\theta_2^2=\theta_3^2=\theta_4^2=\pm 1$.
By Equation \eqref{eqn:FirstIndicatorZAB}, the $\theta_j$'s split up into two nonempty groups of opposite sign, so without loss of generality, $\theta_1\in \{1,i\}$.
\end{proof}

\begin{defn}
\label{defn:Gamma}
For $j=2,\dots, 4$, we let $\epsilon_j\in \{-1,+1\}$ such that $\epsilon_j\theta=\theta_j$ where $\theta=\theta_1$ is the twist of $Y_1$.  
We also set 
$\gamma := \frac{1}{2}(m+n) + \frac{r}{2}$ and 
$\overline{\gamma} := \frac{1}{2}(m+n) - \frac{r}{2}$.
Notice that 
$\gamma + \overline{\gamma} = m+n$ 
and that 
$\gamma^2 + \overline{\gamma}^2 = \frac{1}{2}(m+n)^2 + \frac{r^2}{2}(m-n)^2$.
\end{defn}

\begin{prop}
We have the following equalities:
\begin{align}
-\sum_s \theta_s(z_s{+}z_s')^2d &= \frac{3}{2}(m{+}n)^2 d {+} \frac{r^2}{2}(m-n)^2 d{+}2(m{+}n) {+}
\theta
\sum_{j=1}^4\epsilon_j(y_j{+}y'_j)(1{+}(y_j{+}y'_j)d)
\label{eqn:UsefulIndicatorAB1}
\\
-\sum_s \theta_s^2(z_s{+}z_s')^2d &= \Delta {+}\frac{3}{2}(m{+}n)^2 d {+} \frac{r^2}{2}(m-n)^2 d{+}2(m{+}n){+}\theta^2\sum_{j=1}^4(y_j{+}y'_j)(1{+}(y_j{+}y'_j)d)
\label{eqn:UsefulIndicatorAB2}
\end{align}
where $\Delta:=\pm(\Tr_{Z(\cC)}(\theta_{\cI(\rho)}^2) + \Tr_{Z(\cC)}(\theta_{\cI(\alpha\rho)}^2))
\in \{0,\pm2\dim(\sC)=\pm  (8+4(m+n)d)\}$.
\end{prop}
\begin{proof}
To get Equation \eqref{eqn:UsefulIndicatorAB1}, we add the following two equations for the first Frobenius-Schur indicators \cite{MR2381536} of $\cI(\rho)$ and $\cI(\alpha\rho)$, and we use Equations \eqref{eqn:HorizontalSumsAB1}, \eqref{eqn:HorizontalSumsAB2}, \eqref{eqn:HorizontalSumsAB3} 
in conjunction with Definition \ref{defn:Gamma}.
\begin{align*}
0&=\Tr_{\sZ(\sC)}(\theta_{\cI(\rho)})\\
&= 
x_2(1+(m+n)d)+x_3(1+\gamma d)+x_4(1+\overline{\gamma}d)
+
\theta
\sum_{j=1}^4 \epsilon_j y_j(1+(y_j+y'_j)d)
+
\sum_s \theta_sz_s(z_s+z_s')d
\\
0&=\Tr_{\sZ(\sC)}(\theta_{\cI(g\rho)})\\
&=
x'_2(1+(m+n)d)+x'_3(1+\gamma d)+x'_4(1+\overline{\gamma}d)+
\theta
\sum_{j=1}^4 \epsilon_j y_j(1+(y_j+y'_j)d)
+
\sum_s \theta_sz_s'(z_s+z_s')d.
\end{align*}
Obtaining equation \eqref{eqn:UsefulIndicatorAB2} is similar using the second Frobenius-Schur indicators \cite{MR2381536} of $\rho$ and $\alpha\rho$:
\begin{align*}
\pm \dim(\cC)&=\Tr_{Z(\cC)}(\theta_{\cI(\rho)}^2)\\
&=x_2(1{+}(m{+}n)d){+}x_3(1{+}\gamma d){+}x_4(1{+}\overline{\gamma}d)
{+}
\theta^2\sum_{j=1}^4  y_j(1{+}(y_j{+}y'_j)d)
{+}
\sum_s \theta_s^2z_s(z_s{+}z_s')d
\\
\pm \dim(\cC)&=\Tr_{Z(\cC)}(\theta_{\cI(\alpha\rho)}^2)\\
&=x'_2(1{+}(m{+}n) d){+}x'_3(1{+}\gamma d){+}x'_4(1{+}\overline{\gamma}d)
{+}
\theta^2
\sum_{j=1}^4  y'_j(1{+}(y_j{+}y'_j) d)
{+}
\sum_s \theta_s^2z_s'(z_s{+}z_s')d.
\end{align*}
Adding the above equations,
applying 
\eqref{eqn:HorizontalSumsAB1}, \eqref{eqn:HorizontalSumsAB2}, \eqref{eqn:HorizontalSumsAB3} 
and Definition \ref{defn:Gamma},
and rearranging gives the result.
\end{proof}

\begin{thm}[{\cite[Prop.~5.6 and Thm.~5.7]{MR3229513}}]
\label{thm:LarsonRootsOfUnityAB}
Suppose $u,v\in \bbZ$ and $t\in \bbN$ is square free.
\begin{enumerate}[label=(\arabic*)]
\item It requires at least $|u|+2|v|$ roots of unity to write $u+v\sqrt{2}$ as a sum of roots of unity.
\item It requires at least $|v|\varphi(2t)$ roots of unity to write $u+v\sqrt{t}$ as a sum of roots of unity.
\end{enumerate}
\end{thm}

\begin{cor}
\label{cor:LarsonRationalExtension}
Suppose $u\in \bbQ$, $v\in \bbZ$, and $t\in \bbN$ is square free.
\begin{enumerate}[label=(\arabic*)]
\item It requires at least $|u|+2|v|$ roots of unity to write $u+v\sqrt{2}$ as a sum of roots of unity.
\item It requires at least $|v|\varphi(2t)$ roots of unity to write $u+v\sqrt{t}$ as a sum of roots of unity.
\end{enumerate}
\end{cor}
\begin{proof}
Suppose $\sum_{i=1}^N \zeta_i = u+v\sqrt{2}$.
Write $u=p/q$ in lowest terms with $q>0$
so that
$q\sum_{s=1}^N \zeta_s = p+qv\sqrt{2}$.
By Theorem \ref{thm:LarsonRootsOfUnityAB}(1), $qN \geq |p|+2q|v|$, so $N \geq |u|+2|v|$.

Now suppose $\sum_{i=1}^N \zeta_i = u+v\sqrt{t}$.
Again write $u=p/q$ in lowest terms with $q>0$
so that
$q\sum_{s=1}^N \zeta_s = p+qv\sqrt{t}$.
By Theorem \ref{thm:LarsonRootsOfUnityAB}(2), $qN \geq q|v|\varphi(2t)$, so $N\geq |v|\varphi(2t)$.
\end{proof}

\begin{lem}
\label{lem:BoundTotientAB}
For all $t\in \bbN$ with $t\neq 1,2,3,6$,
$\varphi(2t) \geq \sqrt{\frac{16t}{5}}$.
\end{lem}
\begin{proof}
By \cite{301856}, for all $t\in \bbN$, $\varphi(2t)\geq 2\left(\frac{t}{3}\right)^{2/3}$.
It is straightforward to show that for $t> 42$, 
$2\left(\frac{t}{3}\right)^{2/3} \geq \sqrt{16t/5}$.
One verifies directly that for $t=4,5$ and $7\leq t \leq 42$,
$\varphi(2t) \geq \sqrt{16t/5}$.
The result follows.
\end{proof}

\begin{thm}
If there is a pseudounitary fusion category $\cC$ with the fusion rules \ref{Q:4ObjectsSelfDual}, then $(m+n) \leq 5$.
\label{thm:TwoParameterBound}
\end{thm}
\begin{proof}
We consider the two cases for $\theta\in \{1,i\}$ afforded by Proposition \ref{prop:TwoThetasAB}.
\begin{enumerate}[label=(\arabic*)]

\item
Suppose $\theta = i$.
We add Equation \eqref{eqn:UsefulIndicatorAB1} to its complex conjugate, divide by $d$, and simplify to obtain
\begin{align}
-\sum_s (\theta_s+\overline{\theta_s})(z_s+z_s')^2 
&=
(m+n)^2 + r^2(m-n)^2 + 2(m+n)\sqrt{4 + (m+n)^2}.
\label{eqn:iSumRootsAB}
\end{align}
\begin{itemize}
\item 
{\bf \underline{Case 1:}}
Suppose $(m+n)^2+4=2v_0^2$ for some integer $v_0>0$.
Then by Corollary \ref{cor:LarsonRationalExtension}(1) with $v= 2(m+n)v_0$, it requires at least 
$$
(m+n)^2 + r^2(m-n)^2+\underbrace{4(m+n)\underbrace{\sqrt{\frac{4 + (m+n)^2}{2}}}_{=v_0}}_{=2v}
\geq 
\left(1+\frac{4}{\sqrt{2}}\right)(m+n)^2
$$
roots of unity to write the right hand side of Equation \eqref{eqn:iSumRootsAB}. 
Together with Inequality \eqref{eqn:SquaresUpperBoundAB}, we see
$$
16 + 3(m+n)^2 \geq 2 \sum_{s}(z_s + z_s')^2 \geq \left(1 + \frac{4}{\sqrt{2}}\right)(m+n)^2,
$$
which implies $m+n \leq 4$.

\item
{\bf \underline{Case 2:}}
If $4 + (m+n)^2 \neq 2v^2$, then we can write $4 + (m+n)^2=v^2t$ where $v,t$ are integers with $v>0$ and $t>2$ is square free.
Then by Corollary \ref{cor:LarsonRationalExtension}(2), it requires at least 
$2(m+n)v\varphi(2t)$ roots of unity to write the right hand side of Equation \eqref{eqn:iSumRootsAB}. 
Since $4 + (m+n)^2 \equiv \pm 1 \bmod 3$, 
we know $4+(m+n)^2 \notin \{1,2,3,6\}$.
By Lemma \ref{lem:BoundTotientAB},
$$
v^2\varphi(2t)^2 
\geq 
\frac{16v^2t}{5}
\Longleftrightarrow 
v\varphi(2t) 
\geq 
4\sqrt{\frac{4 + (m+n)^2}{5}}.
$$
Now by inequality \eqref{eqn:SquaresUpperBoundAB}, we see
$$
16 + 3(m+n)^2
\geq 
2\sum_s (z_s+z_s')^2 
\geq 
2(m+n)v\varphi(2t)
\geq 
8(m+n) \sqrt{\frac{4 + (m+n)^2}{5}},
$$
which implies $m+n \leq 4$.
\end{itemize}

\item
Suppose $\theta = 1$.
Then dividing Equation \eqref{eqn:UsefulIndicatorAB2} by $d$ and simplifying, we get
\begin{align}
-\sum_s \theta_s^2(z_s+z_s')^2
&=
\frac{\Delta+4(m+n)}{d} +\frac{3}{2}(m+n)^2 + \frac{r^2}{2}(m-n)^2+\sum_{j=1}^4(y_j+y'_j)^2
\label{eqn:1SumRootsAB1}
\end{align}
There are now 2 cases depending on the value of $\Delta$.

\begin{itemize}
\item
{\bf \underline{Case 1:}}
Suppose $\Delta=0$. 
Then Equation \eqref{eqn:1SumRootsAB1} becomes
\begin{align}
-\sum_s \theta_s^2(z_s+z_s')^2 
&= 
2(a+b)\sqrt{4 + (m+n)^2} -\frac{1}{2}(m+n)^2 + \frac{r^2}{2}(m-n)^2 +\sum_{j=1}^4(y_j+y'_j)^2.
\label{eqn:1SumRootsAB2}
\end{align}

\item
{\bf \underline{Case 2:}}
Suppose $\Delta=\pm 2\dim(\sC)=\pm(8+4(m+n)d)$. 
Then equation \eqref{eqn:1SumRootsAB1} becomes
\begin{equation}
\begin{split}
-\sum_s \theta_s^2(z_s+z_s')^2 
&= 
2(m+n\pm 2)\sqrt{4 + (m+n)^2} -2(m+n\pm 2)(a+b) 
\\
&
\qquad
\pm 4(m+n) + \frac{3}{2}(m+n)^2 
+ \frac{r^2}{2}(m-n)^2 +\sum_{j=1}^4(y_j+y'_j)^2.
\end{split}
\label{eqn:1SumRootsAB3}
\end{equation}
\end{itemize}

In either of the above cases, arguing as in (1) where $\theta=i$, we see that it takes at least
$$
\min\left\{ 
8(m+n - 2)\sqrt{\frac{4 + (m+n)^2}{5}}, 4(m+n - 2)\sqrt{\frac{4 + (m+n)^2}{2}}
\right\}
\geq 
\frac{4}{\sqrt{2}}(m+n)(m+n-2)
$$
roots of unity to write the right hand sides of Equations \eqref{eqn:1SumRootsAB2} and \eqref{eqn:1SumRootsAB3}.
Now by inequality \eqref{eqn:SquaresUpperBoundAB}, we see
$$
8 + \frac{3}{2}(m+n)^2
\geq 
\sum_s (z_s+z_s')^2 
\geq 
\frac{4}{\sqrt{2}}(m+n)(m+n-2),
$$
which implies $m+n\leq 5$.
\qedhere
\end{enumerate}
\end{proof}

\bibliographystyle{alpha}
{\footnotesize{
\bibliography{bibliography}

\newcommand{\etalchar}[1]{$^{#1}$}
\begin{thebibliography}{GMP{\etalchar{+}}18}

\bibitem[AH99]{MR1686551}
Marta Asaeda and Uffe Haagerup.
\newblock Exotic subfactors of finite depth with {J}ones indices
  {$(5+\sqrt{13})/2$} and {$(5+\sqrt{17})/2$}.
\newblock {\em Comm. Math. Phys.}, 202(1):1--63, 1999.
\newblock \mathscinet{MR1686551}, \doi{10.1007/s002200050574},
  \arXiv{math.OA/9803044}.

\bibitem[AP95]{MR1328736}
Henning~Haahr Andersen and Jan Paradowski.
\newblock Fusion categories arising from semisimple {L}ie algebras.
\newblock {\em Comm. Math. Phys.}, 169(3):563--588, 1995.
\newblock \mathscinet{MR1328736}.

\bibitem[Big10]{MR2577673}
Stephen Bigelow.
\newblock Skein theory for the {$ADE$} planar algebras.
\newblock {\em J. Pure Appl. Algebra}, 214(5):658--666, 2010.
\newblock \arXiv{math.QA/0903.0144} \mathscinet{MR2577673}
  \doi{10.1016/j.jpaa.2009.07.010}.

\bibitem[BMPS12]{MR2979509}
Stephen Bigelow, Scott Morrison, Emily Peters, and Noah Snyder.
\newblock Constructing the extended {H}aagerup planar algebra.
\newblock {\em Acta Math.}, 209(1):29--82, 2012.
\newblock \mathscinet{MR2979509}, \arXiv{0909.4099},
  \doi{10.1007/s11511-012-0081-7}.

\bibitem[BN91]{MR1193933}
Jocelyne Bion-Nadal.
\newblock An example of a subfactor of the hyperfinite {${\rm II}\sb 1$} factor
  whose principal graph invariant is the {C}oxeter graph {$E\sb 6$}.
\newblock In {\em Current topics in operator algebras ({N}ara, 1990)}, pages
  104--113. World Sci. Publ., River Edge, NJ, 1991.
\newblock \mathscinet{MR1193933}.

\bibitem[BP14]{MR3157990}
Stephen Bigelow and David Penneys.
\newblock Principal graph stability and the jellyfish algorithm.
\newblock {\em Math. Ann.}, 358(1-2):1--24, 2014.
\newblock \mathscinet{MR3157990}, \doi{10.1007/s00208-013-0941-2},
  \arxiv{1208.1564}.

\bibitem[BW96]{MR1357878}
John~W. Barrett and Bruce~W. Westbury.
\newblock Invariants of piecewise-linear {$3$}-manifolds.
\newblock {\em Trans. Amer. Math. Soc.}, 348(10):3997--4022, 1996.
\newblock \mathscinet{MR1357878} \doi{10.1090/S0002-9947-96-01660-1}
  \arxiv{hep-th/9311155}.

\bibitem[EG11]{MR2837122}
David~E. Evans and Terry Gannon.
\newblock The exoticness and realisability of twisted {H}aagerup-{I}zumi
  modular data.
\newblock {\em Comm. Math. Phys.}, 307(2):463--512, 2011.
\newblock \arxiv{1006.1326} \mathscinet{MR2837122}
  \doi{10.1007/s00220-011-1329-3}.

\bibitem[EG14]{MR3167494}
David~E. Evans and Terry Gannon.
\newblock Near-group fusion categories and their doubles.
\newblock {\em Adv. Math.}, 255:586--640, 2014.
\newblock \mathscinet{MR3167494} \arxiv{1208.1500}
  \doi{10.1016/j.aim.2013.12.014}.

\bibitem[EGNO15]{MR3242743}
Pavel Etingof, Shlomo Gelaki, Dmitri Nikshych, and Victor Ostrik.
\newblock {\em Tensor categories}, volume 205 of {\em Mathematical Surveys and
  Monographs}.
\newblock American Mathematical Society, Providence, RI, 2015.
\newblock \mathscinet{MR3242743} \doi{10.1090/surv/205}.

\bibitem[EK95]{MR1316301}
David~E. Evans and Yasuyuki Kawahigashi.
\newblock On {O}cneanu's theory of asymptotic inclusions for subfactors,
  topological quantum field theories and quantum doubles.
\newblock {\em Internat. J. Math.}, 6(2):205--228, 1995.
\newblock \mathscinet{MR1316301} \doi{10.1142/S0129167X95000468}.

\bibitem[EK98]{MR1642584}
David~E. Evans and Yasuyuki Kawahigashi.
\newblock {\em Quantum Symmetries on Operator Algebras}.
\newblock Oxford Mathematical Monographs. Oxford Science Publications. The
  Clarendon Press, Oxford University Press, New York, 1998.
\newblock xvi+829 pp. ISBN: 0-19-851175-2, \mathscinet{MR1642584}.

\bibitem[FG21]{2106.16186}
J{\:u}rgen Fuchs and Tobias Gr{\o}sfjeld.
\newblock Tetrahedral symmetry of 6j-symbols in fusion categories, 2021.
\newblock \arxiv{2106.16186}.

\bibitem[FGSV99]{MR1657800}
J.~Fuchs, A.~Ch. Ganchev, K.~Szlach{\'a}nyi, and P.~Vecserny{\'e}s.
\newblock {$S\sb 4$} symmetry of {$6j$} symbols and {F}robenius-{S}chur
  indicators in rigid monoidal {$C\sp *$} categories.
\newblock {\em J. Math. Phys.}, 40(1):408--426, 1999.
\newblock \arXiv{physics/9803038}, \mathscinet{MR1657800}.

\bibitem[GJS15]{MR3354332}
Pinhas Grossman, David Jordan, and Noah Snyder.
\newblock Cyclic extensions of fusion categories via the {B}rauer-{P}icard
  groupoid.
\newblock {\em Quantum Topol.}, 6(2):313--331, 2015.
\newblock \mathscinet{MR3354332} \doi{10.4171/QT/64} \arXiv{1211.6414}.

\bibitem[GMP{\etalchar{+}}18]{1810.06076}
Pinhas Grossman, Scott Morrison, David Penneys, Emily Peters, and Noah Snyder.
\newblock The {E}xtended {H}aagerup fusion categories, 2018.
\newblock \arxiv{1810.06076}.

\bibitem[Haa94]{MR1317352}
Uffe Haagerup.
\newblock Principal graphs of subfactors in the index range
  {$4<[M:N]<3+\sqrt2$}.
\newblock In {\em Subfactors ({K}yuzeso, 1993)}, pages 1--38. World Sci. Publ.,
  River Edge, NJ, 1994.
\newblock \mathscinet{MR1317352} available at
  \url{http://tqft.net/other-papers/subfactors/haagerup.pdf}.

\bibitem[HH09]{MR2559711}
Tobias~J. Hagge and Seung-Moon Hong.
\newblock Some non-braided fusion categories of rank three.
\newblock {\em Commun. Contemp. Math.}, 11(4):615--637, 2009.
\newblock \mathscinet{MR2559711} \doi{10.1142/S0219199709003521}
  \arxiv{0704.0208}.

\bibitem[hj]{301856}
Will~Jagy (https://math.stackexchange.com/users/10400/will jagy).
\newblock Is the euler phi function bounded below?
\newblock Mathematics Stack Exchange.
\newblock \url{https://math.stackexchange.com/q/301856 (version: 2019-12-29)}.

\bibitem[IMP16]{MR3536926}
Masaki Izumi, Scott Morrison, and David Penneys.
\newblock Quotients of {$A\sb 2\ast T\sb 2$}.
\newblock {\em Canad. J. Math.}, 68(5):999--1022, 2016.
\newblock \mathscinet{MR3536926} \doi{10.4153/CJM-2015-017-4}. This is an
  abridged version of \arxiv{1308.5723}.

\bibitem[Izu00]{MR1782145}
Masaki Izumi.
\newblock The structure of sectors associated with {L}ongo-{R}ehren inclusions.
  {I}. {G}eneral theory.
\newblock {\em Comm. Math. Phys.}, 213(1):127--179, 2000.
\newblock \mathscinet{MR1782145} \doi{10.1007/s002200000234}.

\bibitem[Izu01]{MR1832764}
Masaki Izumi.
\newblock The structure of sectors associated with {L}ongo-{R}ehren inclusions.
  {II}. {E}xamples.
\newblock {\em Rev. Math. Phys.}, 13(5):603--674, 2001.
\newblock \mathscinet{MR1832764} \doi{10.1142/S0129055X01000818}.

\bibitem[Izu17]{MR3635673}
Masaki Izumi.
\newblock A {C}untz algebra approach to the classification of near-group
  categories.
\newblock In {\em Proceedings of the 2014 {M}aui and 2015 {Q}inhuangdao
  conferences in honour of {V}aughan {F}. {R}. {J}ones' 60th birthday},
  volume~46 of {\em Proc. Centre Math. Appl. Austral. Nat. Univ.}, pages
  222--343. Austral. Nat. Univ., Canberra, 2017.
\newblock \mathscinet{MR3635673} \arxiv{1512.04288}.

\bibitem[Izu18]{MR3827808}
Masaki Izumi.
\newblock The classification of {$3^n$} subfactors and related fusion
  categories.
\newblock {\em Quantum Topol.}, 9(3):473--562, 2018.
\newblock \mathscinet{MR3827808} \doi{10.4171/QT/113}.

\bibitem[KO02]{MR1936496}
Alexander Kirillov, Jr. and Viktor Ostrik.
\newblock On a {$q$}-analogue of the {M}c{K}ay correspondence and the {ADE}
  classification of {$\mathfrak{sl}_2$} conformal field theories.
\newblock {\em Adv. Math.}, 171(2):183--227, 2002.
\newblock \mathscinet{MR1936496} \arXiv{math.QA/0101219}
  \doi{10.1006/aima.2002.2072}.

\bibitem[Lar14]{MR3229513}
Hannah~K. Larson.
\newblock Pseudo-unitary non-self-dual fusion categories of rank 4.
\newblock {\em J. Algebra}, 415:184--213, 2014.
\newblock \mathscinet{MR3229513}, \doi{10.1016/j.jalgebra.2014.05.032}.

\bibitem[Liu15]{MR3345186}
Zhengwei Liu.
\newblock Composed inclusions of {$A\sb 3$} and {$A\sb 4$} subfactors.
\newblock {\em Adv. Math.}, 279:307--371, 2015.
\newblock \mathscinet{MR3345186} \doi{10.1016/j.aim.2015.03.017}
  \arxiv{1308.5691}.

\bibitem[LLB20]{2012.14424}
Chien-Hung Lin, Michael Levin, and Fiona~J. Burnell.
\newblock Generalized string-net models: A thorough exposition, 2020.
\newblock \arxiv{2012.14424}.

\bibitem[LMP15]{MR3306607}
Zhengwei Liu, Scott Morrison, and David Penneys.
\newblock 1-{S}upertransitive {S}ubfactors with {I}ndex at {M}ost
  {$6\frac{1}{5}$}.
\newblock {\em Comm. Math. Phys.}, 334(2):889--922, 2015.
\newblock \mathscinet{MR3306607}, \arXiv{1310.8566},
  \doi{10.1007/s00220-014-2160-4}.

\bibitem[LMP20]{MR4079744}
Zhengwei Liu, Scott Morrison, and David Penneys.
\newblock Lifting shadings on symmetrically self-dual subfactor planar
  algebras.
\newblock In {\em Topological phases of matter and quantum computation}, volume
  747 of {\em Contemp. Math.}, pages 51--61. Amer. Math. Soc., [Providence],
  RI, [2020] \copyright 2020.
\newblock \mathscinet{MR4079744} \doi{10.1090/conm/747/15038}
  \arxiv{1709.05023}.

\bibitem[LPR20]{2010.10264}
Zhengwei Liu, Sebastien Palcoux, and Yunxiang Ren.
\newblock Classification of {G}rothendieck rings of complex fusion categories
  of multiplicity one up to rank six, 2020.
\newblock \arxiv{2010.10264}.

\bibitem[Lus03]{MR1974442}
G.~Lusztig.
\newblock {\em Hecke algebras with unequal parameters}, volume~18 of {\em CRM
  Monograph Series}.
\newblock American Mathematical Society, Providence, RI, 2003.
\newblock \mathscinet{MR1974442} \doi{10.1090/crmm/018}.

\bibitem[LW05]{PhysRevB.71.045110}
Michael~A. Levin and Xiao-Gang Wen.
\newblock String-net condensation: A physical mechanism for topological phases.
\newblock {\em Phys. Rev. B}, 71:045110, Jan 2005.
\newblock doi{10.1103/PhysRevB.71.045110} \arxiv{cond-mat/0404617}.

\bibitem[MP15]{MR3394622}
Scott Morrison and David Penneys.
\newblock 2-supertransitive subfactors at index {$3+\sqrt{5}$}.
\newblock {\em J. Funct. Anal.}, 269(9):2845--2870, 2015.
\newblock \mathscinet{MR3394622} \doi{10.1016/j.jfa.2015.06.023}
  \arxiv{1406.3401}.

\bibitem[M{\"u}g03]{MR1966525}
Michael M{\"u}ger.
\newblock From subfactors to categories and topology. {II}. {T}he quantum
  double of tensor categories and subfactors.
\newblock {\em J. Pure Appl. Algebra}, 180(1-2):159--219, 2003.
\newblock \mathscinet{MR1966525} \doi{10.1016/S0022-4049(02)00248-7}
  \arXiv{math.CT/0111205}.

\bibitem[MW17]{MR3611056}
Scott Morrison and Kevin Walker.
\newblock The center of the extended {H}aagerup subfactor has 22 simple
  objects.
\newblock {\em Internat. J. Math.}, 28(1):1750009, 11, 2017.
\newblock \mathscinet{MR3611056} \doi{10.1142/S0129167X17500094}
  \arxiv{1404.3955}.

\bibitem[NS07a]{MR2313527}
Siu-Hung Ng and Peter Schauenburg.
\newblock Frobenius-{S}chur indicators and exponents of spherical categories.
\newblock {\em Adv. Math.}, 211(1):34--71, 2007.
\newblock \mathscinet{MR2313527} \doi{10.1016/j.aim.2006.07.017}
  \arxiv{math/0601012}.

\bibitem[NS07b]{MR2381536}
Siu-Hung Ng and Peter Schauenburg.
\newblock Higher {F}robenius-{S}chur indicators for pivotal categories.
\newblock In {\em Hopf algebras and generalizations}, volume 441 of {\em
  Contemp. Math.}, pages 63--90. Amer. Math. Soc., Providence, RI, 2007.
\newblock \mathscinet{MR2381536}, \arXiv{math/0503167}.

\bibitem[Ocn88]{MR996454}
Adrian Ocneanu.
\newblock Quantized groups, string algebras and {G}alois theory for algebras.
\newblock In {\em Operator algebras and applications, Vol.\ 2}, volume 136 of
  {\em London Math. Soc. Lecture Note Ser.}, pages 119--172. Cambridge Univ.
  Press, Cambridge, 1988.
\newblock \mathscinet{MR996454}.

\bibitem[Ocn02]{MR1907188}
Adrian Ocneanu.
\newblock The classification of subgroups of quantum {${\rm SU}(N)$}.
\newblock In {\em Quantum symmetries in theoretical physics and mathematics
  ({B}ariloche, 2000)}, volume 294 of {\em Contemp. Math.}, pages 133--159.
  Amer. Math. Soc., Providence, RI, 2002.
\newblock \mathscinet{MR1907188}.

\bibitem[Ost03]{MR1981895}
Viktor Ostrik.
\newblock Fusion categories of rank 2.
\newblock {\em Math. Res. Lett.}, 10(2-3):177--183, 2003.
\newblock \mathscinet{MR1981895} \arXiv{math.QA/0203255}.

\bibitem[Ost09]{MR2576705}
Victor Ostrik.
\newblock On formal codegrees of fusion categories.
\newblock {\em Math. Res. Lett.}, 16(5):895--901, 2009.
\newblock \arXiv{0810.3242}, \mathscinet{MR 2576705}.

\bibitem[Ost13]{1309.4822}
Victor Ostrik.
\newblock Pivotal fusion categories of rank 3, 2013.
\newblock (with an {A}ppendix written jointly with {D}mitri {N}ikshych),
  \arxiv{1309.4822}.

\bibitem[PP15]{MR3402358}
David Penneys and Emily Peters.
\newblock Calculating two-strand jellyfish relations.
\newblock {\em Pacific J. Math.}, 277(2):463--510, 2015.
\newblock \mathscinet{MR3402358} \doi{10.2140/pjm.2015.277-2}
  \arXiv{1308.5197}.

\bibitem[TY98]{MR1659954}
Daisuke Tambara and Shigeru Yamagami.
\newblock Tensor categories with fusion rules of self-duality for finite
  abelian groups.
\newblock {\em J. Algebra}, 209(2):692--707, 1998.
\newblock \mathscinet{MR1659954}.

\bibitem[Wen88]{MR936086}
Hans Wenzl.
\newblock Hecke algebras of type {$A_n$} and subfactors.
\newblock {\em Invent. Math.}, 92(2):349--383, 1988.
\newblock \mathscinet{MR936086} \doi{10.1007/BF01404457}.

\bibitem[Wen90]{MR1090432}
Hans Wenzl.
\newblock Quantum groups and subfactors of type {$B$}, {$C$}, and {$D$}.
\newblock {\em Comm. Math. Phys.}, 133(2):383--432, 1990.
\newblock \mathscinet{MR1090432}.

\end{thebibliography}
\bigskip
\noindent
(Cain Edie-Michell)
Department of Mathematics,
University of California, San Diego (UCSD),
9500 Gilman Drive \# 0112,
La Jolla, CA  92093-0112,
USA
\\
\texttt{email: cediemichell@ucsd.edu}
\\
\\
\noindent
(Masaki Izumi)
Department of Mathematics, Graduate School of Science, Kyoto University, Kyoto 606-8502, Japan
\\
\texttt{email: izumi@math.kyoto-u.ac.jp}
\\
\\
\noindent
(David Penneys)
Department of Mathematics, 100 Math Tower, 231 West 18th Ave., Columbus, OH 43210-1174, USA
\\
\texttt{email: penneys.2@osu.edu}
}}
\end{document}